\documentclass[12pt, a4paper]{book}
\usepackage[german,english]{babel}

\usepackage{wallpaper}
\usepackage{setspace}

\usepackage[shortlabels]{enumitem}
\setlist[enumerate]{parsep=0pt,listparindent=\parindent}

\usepackage{wrapfig}

\usepackage{soul}

\usepackage[frozencache,cachedir=.]{minted}

\usepackage{caption}
\usepackage{subcaption}

\usepackage{amsmath}
\usepackage{amssymb}
\usepackage{amsthm}
\usepackage{stmaryrd}

\usepackage{graphicx} 

\usepackage{mathtools}

\usepackage{mathabx,epsfig}

\usepackage{stmaryrd}

\usepackage{wasysym}

\usepackage[left=3.5cm, right=3.5cm, top=3.5cm, bottom=3.5cm]{geometry}

\usepackage{etoolbox}

\makeatletter
\newcommand*\NoIndentAfterEnv[1]{%
  \AfterEndEnvironment{#1}{\par\@afterindentfalse\@afterheading}}
\makeatother

\NoIndentAfterEnv{itemize}
\NoIndentAfterEnv{enumerate}
\NoIndentAfterEnv{proof}
\NoIndentAfterEnv{figure}
\NoIndentAfterEnv{equation}
\NoIndentAfterEnv{minted}
\NoIndentAfterEnv{Satz}
\NoIndentAfterEnv{Def}
\NoIndentAfterEnv{Ex}
\NoIndentAfterEnv{Rem}
\NoIndentAfterEnv{Algo}
\NoIndentAfterEnv{Listing}
\NoIndentAfterEnv{Prop}
\NoIndentAfterEnv{PropDef}
\NoIndentAfterEnv{Lem}
\NoIndentAfterEnv{Kor}
\NoIndentAfterEnv{Ass}
\NoIndentAfterEnv{specialSatz}

\usepackage{color}
\definecolor{darkred}{RGB}{170,0,0}
\definecolor{darkgreen}{RGB}{0,120,0}

\usepackage{soul}
\setuldepth{noDescent}

\usepackage{tikz}
\usetikzlibrary{cd}
\usetikzlibrary{babel}
\tikzset{
	labl/.style={anchor=south, rotate=90, inner sep=.5mm}
}

\usepackage{makeidx}

\usepackage{tcolorbox}

\usepackage{hyperref}
\hypersetup{
	colorlinks,
	citecolor=blue,
	filecolor=black,
	linkcolor=black,
	urlcolor=blue
}

\usepackage{booktabs}

\usepackage{nowidow}

\usepackage{calligra}
\usepackage{mathrsfs}

\usepackage[retainorgcmds]{IEEEtrantools}

\usepackage{tocloft}

\definecolor{lightgray}{rgb}{0.4,0.4,0.4}

\let\le\leqslant
\let\ge\geqslant

\DeclareMathOperator{\add}{add}

\DeclareMathOperator{\an}{an}

\DeclareMathOperator{\characteristic}{char}

\DeclareMathOperator{\CTT}{CTT}

\DeclareMathOperator{\Gal}{Gal}

\DeclareMathOperator{\Hom}{Hom}

\DeclareMathOperator{\interior}{int}

\DeclareMathOperator{\len}{len}

\DeclareMathOperator{\Nm}{Nm}

\DeclareMathOperator{\ns}{ns}

\DeclareMathOperator{\ord}{ord}

\DeclareMathOperator{\rk}{rk}

\DeclareMathOperator{\red}{red}
\DeclareMathOperator{\ret}{ret}

\DeclareMathOperator{\slope}{slope}

\DeclareMathOperator{\Spec}{Spec}

\DeclareMathOperator{\stab}{stab}

\DeclareMathOperator{\tail}{tail}
\DeclareMathOperator{\tame}{tame}

\DeclareMathOperator{\trdeg}{trdeg}

\DeclareMathOperator{\sHom}{\mathscr{H}\text{\kern -3pt {\calligra\large om}}\,}

\mathcode`l="8000
\begingroup
\makeatletter
\lccode`\~=`\l
\DeclareMathSymbol{\lsb@l}{\mathalpha}{letters}{`l}
\lowercase{\gdef~{\ifnum\the\mathgroup=\m@ne \ell \else \lsb@l \fi}}%
\endgroup

\renewcommand{\emptyset}{\varnothing}
\renewcommand{\hat}{\widehat}
\usepackage{scalerel,stackengine}
\stackMath
\newcommand\reallywidehat[1]{%
\savestack{\tmpbox}{\stretchto{%
  \scaleto{%
    \scalerel*[\widthof{\ensuremath{#1}}]{\kern-.6pt\bigwedge\kern-.6pt}%
    {\rule[-\textheight/2]{1ex}{\textheight}}
  }{\textheight}%
}{0.5ex}}%
\stackon[1pt]{#1}{\tmpbox}%
}
\renewcommand{\phi}{\varphi}
\renewcommand{\epsilon}{\varepsilon}

\usepackage{mathtools}

\newcommand{\abs}[1]{\left\lvert #1 \right\rvert}

\newcommand{\rest}[2]{\left.{#1}\right|_{#2}}

\newcommand{\BA}{{\mathbb{A}}}

\newcommand{\BC}{{\mathbb{C}}}

\newcommand{\BF}{{\mathbb{F}}}

\newcommand{\BP}{{\mathbb{P}}}
\newcommand{\BQ}{{\mathbb{Q}}}
\newcommand{\BR}{{\mathbb{R}}}

\newcommand{\BZ}{{\mathbb{Z}}}

\newcommand{\Fm}{{\mathfrak{m}}}

\newcommand{\CA}{{\cal A}}

\newcommand{\CD}{{\cal D}}

\newcommand{\CH}{{\cal H}}

\newcommand{\CM}{{\cal M}}

\newcommand{\CO}{{\cal O}}

\newcommand{\CX}{{\cal X}}
\newcommand{\CY}{{\cal Y}}

\newcommand{\rA}{{\mathsf A}}

\newcommand{\rC}{{\mathsf C}}
\newcommand{\rD}{{\mathsf D}}

\newcommand{\rU}{{\mathsf U}}
\newcommand{\rV}{{\mathsf V}}

\newcommand{\rX}{{\mathsf X}}

\newtheorem{Prin}{Principle}[chapter]
\newtheorem{Prop}[Prin]{Proposition}
\newtheorem{PropDef}[Prin]{Proposition-Definition}
\newtheorem{Lem}[Prin]{Lemma}
\newtheorem{Satz}[Prin]{Theorem}

\newtheorem{Kor}[Prin]{Corollary}
\newtheorem{Ass}[Prin]{Assumption}

\newcommand{\thistheoremname}{}
\newtheorem*{genericthm}{\thistheoremname}
\newenvironment{specialSatz}[1]
{\renewcommand{\thistheoremname}{#1}%
	\begin{genericthm}}
	{\end{genericthm}}

\theoremstyle{definition}
\newtheorem{Def}[Prin]{Definition}
\newtheorem{Ex}[Prin]{Example}
\newtheorem{Rem}[Prin]{Remark}
\newtheorem{Algo}[Prin]{Algorithm}
\newtheorem{Listing}[Prin]{Code Listing}

\definecolor{ferrarired}{rgb}{1.0, 0.11, 0.0}
\newtcolorbox{01Def}{colback=black!0,colframe=ferrarired!50}

\usepackage[style=alphabetic,citestyle=alphabetic,backend=biber,backref=true]{biblatex}
\DeclareFieldFormat[inbook]{citetitle}{#1}
\DeclareFieldFormat[inbook]{title}{#1}
\title{Semistable reduction of covers of degree $p$}
\author{Ole Ossen}

\begin{document}

\newpage

\pagenumbering{gobble}
\thispagestyle{empty}

\begin{minipage}[b]{0.4\textwidth}
\hfill
\end{minipage}
\begin{minipage}[b]{0.6\textwidth}
\includegraphics[width=7.3cm]{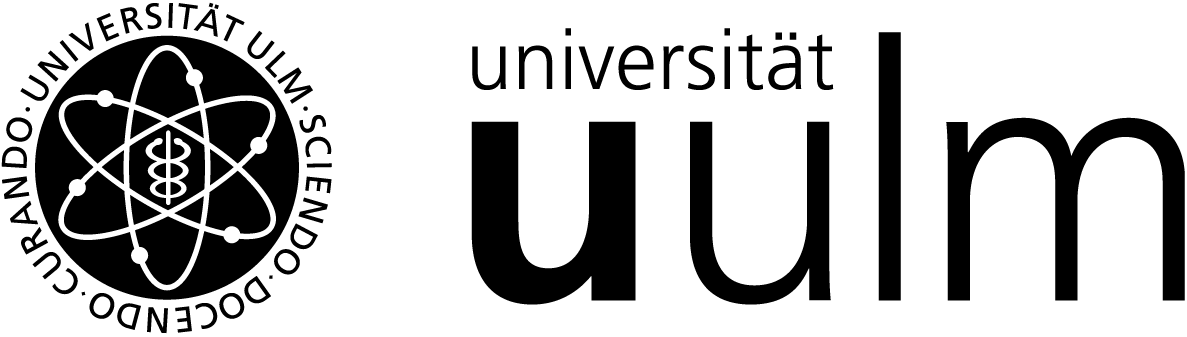}
\end{minipage}

\vspace{3.5cm}	

\begin{center}	
{\huge \bf\ Semistable reduction of covers of degree $p$}  \\
\vspace{11cm}
Dissertation zur Erlangung des Doktorgrades Dr.\,rer.\,nat.\\
der Fakultät für Mathematik und Wirtschaftswissenschaften\\
der Universität Ulm \\\vspace{1cm}

vorgelegt von\\
\textbf{Ole Ossen}\\
aus Berlin im Jahr 2024\\

\end{center}

\newpage

\vspace*{\fill}
\noindent\textbf{Amtierender Dekan}\\
Prof.\ Dr.\ Stefan Funken\\
\textbf{Gutachter}\\
Prof.\ Dr.\ Stefan Wewers\\
Prof.\ Dr.\ Jeroen Sijsling\\
\textbf{Tag der Promotion}\\
2.\ Mai 2024

\newpage

\pagenumbering{arabic}
\tableofcontents

\chapter*{Introduction}
\addcontentsline{toc}{chapter}{Introduction} 

Let $K$ be a field that is complete with respect to a non-trivial non-archimedean valuation $v_K$. We denote the valuation ring of $K$ by $\CO_K$, the residue field by $\kappa$, and we assume that $\kappa$ has positive characteristic $p\coloneqq\characteristic(\kappa)$. The goal of this thesis is to develop a method for computing semistable models of smooth projective algebraic curves $Y$ over $K$ equipped with a degree-$p$ morphism $\phi\colon Y\to\BP_K^1$.

Recall that an $\CO_K$-model of a $K$-curve $Y$ is a flat and proper $\CO_K$-scheme $\CY$ equipped with an isomorphism $\CY\otimes_{\CO_K}K\cong Y$. A model $\CY$ is called \emph{semistable} if its special fiber $\CY_s\coloneqq\CY\otimes_{\CO_K}\kappa$ is reduced and only has ordinary double points as singularities. The Semistable Reduction Theorem of Deligne and Mumford (\cite{deligne-mumford}) asserts that a semistable model $\CY$ always exists, possibly after replacing the curve $Y$ with the base change $Y_L$ to a finite separable field extension $L/K$. The proof of Deligne-Mumford suggests the following strategy for finding $L$ and $\CY$: First, find a field extension $L$ over which the Jacobian $J(Y)$ has semistable reduction by making sure that certain torsion points on $J(Y)$ are $L$-rational. Then find a minimal \emph{regular} model of $Y_L$ by resolving singularities (which is automatically semistable).

However, following this strategy currently seems unfeasible for curves of genus $g>2$. We will use a different strategy, based on the covers $\phi\colon Y\to\BP_K^1$ that the curves $Y$ we consider are always assumed to be equipped with. For the case that $K$ is a finite extension of $\BQ_3$ and $Y$ is a plane quartic curve, our strategy is fully implemented, and many examples throughout are drawn from this class of curves.

\textbf{Previous work}. Much work has been done on computing semistable reductions of curves using covers. This approach may be summarized as follows. First one chooses a finite cover $Y\to X$; often the curve $X$ is the projective line $\BP_K^1$. Then one determines a suitable semistable model $\CX$ of $X$ such that the normalization $\CY$ of $\CX$ in the function field $F_Y$ of $Y$ is a semistable model of $Y$ (again, possibly after replacing $K$ with a finite separable extension $L/K$ and $X$ and $Y$ with their base changes to $L$). While the semistable models of $\BP_K^1$ are well understood, choosing $\CX$ so that $\CY$ is also semistable is a difficult problem.

The choice of $\CX$ is comparatively simple if one can find $\phi\colon Y\to\BP_K^1$ such that the characteristic $p$ is greater than $\deg(\phi)$. In this case, a model $\CX$ of $\BP_K^1$ \emph{separating the branch locus of $\phi$} will always work. For example, this approach is leveraged in \cite{bouw-wewers} for studying superelliptic curves. It also underlies approaches to semistable reduction phrased in terms of \emph{cluster pictures}, as for example in \cite{ddmm}. We briefly review the strategy of separating branch points in \mbox{Section \ref{sec-preliminaries-admissible}}. Attention to the more difficult case $p\le\deg(\phi)$ has mostly been limited to degree-$p$ Galois covers of $\BP_K^1$, that is, to the study of \emph{superelliptic curves}. A crucial part of the present work is that it applies to \emph{all} degree-$p$ covers of $\BP_K^1$, be they Galois or not.

Much of the work on $p$-cyclic covers $Y\to\BP_K^1$ goes back to Coleman's sketch of an algorithm for computing the semistable of $Y$ presented in \cite[Section 6]{coleman}. It is written in rigid-analytic language. A complete proof of the Semistable Reduction Theorem in a similarly algorithmic and rigid-analytic spirit was eventually given in \cite{arzdorf-wewers}. Coleman also introduces the notion of \emph{$p$-approximations}, which in one way or another play a role in all the works mentioned in the following paragraph. 

The papers of Lehr and Matignon \cite{lehr}, \cite{matignon}, \cite{lehr-matignon} describe how to find a field extension $L/K$ and a semistable model $\CY$ of $Y$ in the case of $p$-cyclic $\phi\colon Y\to\BP_K^1$ conditional on the assumption that the branch locus of $\phi$ be \emph{equidistant}. This assumption was no longer necessary in the local approach of \cite{arzdorf} and \cite{arzdorf-wewers}. A modified form of this approach is used in the Sage package MCLF (\cite{mclf}), in which the computation of semistable reductions of superelliptic curves of degree $p$ is implemented. Background and information on this can be found in \cite[Section 4]{icerm}. Similar is the recent work \cite{fiore-yelton}, which is mainly concerned with the semistable reduction of hyperelliptic curves at $p=2$.

As mentioned above, we particularly consider the case of reduction of plane quartic curves at $p=3$. These admit a degree-$3$ morphism $Y\to\BP_K^1$, but this morphism is not usually Galois --- if it is, $Y$ is a \emph{Picard curve}, which case was studied in \cite{boerner-bouw-wewers}. Mention should also be made of two recent works concerning the semistable reduction of plane quartics: In \cite{super-winky-cat}, the stable reduction of plane quartics is described in terms of their \emph{Dixmier-Ohno invariants}, but conditionally on the assumption $\characteristic(\kappa)>7$. In \cite{cayley-octads}, the stable reduction of plane quartics is characterized in terms of combinatorial objects called \emph{Cayley Octads}. This approach is implemented in Magma, but as of now conjectural.

\textbf{The analytic viewpoint}. When investigating the semistable reduction of an algebraic $K$-curve $X$, it is very convenient to employ \emph{analytic geometry} over the non-archimedean field $K$. Above, mention was made of rigid-analytic geometry; we use analytic geometry in the sense of Berkovich (\cite{berkovich}). This allows local work with objects like disks and anuli that is simply impossible in the algebraic setting.

Each algebraic curve $X$ possesses an analytification $X^{\an}$. While this is a complicated topological object (an example is drawn on the cover of this thesis\footnote{Not included in the electronic version. You can take a look at \href{https://github.com/oossen/mclf/blob/plane-quartics/README.md}{my GitHub page}.}), it admits retracts to much simpler spaces, finite subgraphs of $X^{\an}$ called \emph{skeletons}. In fact, the skeletons of $X$ are the same as the intersection graphs of special fibers of semistable models of $X$. Thus the original problem of finding a semistable model of a given curve may be rephrased as the problem of finding a skeleton.

Since semistable models of $\BP_K^1$ are well-understood, so is the structure of the analytic curve $(\BP_K^1)^{\an}$. Its skeletons are simply finite trees. Given a finite morphism $\phi\colon Y\to \BP_K^1$, we stated above that our strategy consisted of finding a suitable model $\CX$ of $\BP_K^1$ such that the normalization of $\CX$ in the function field $F_Y$ is a semistable model of $Y$. In Berkovich's language, this translates to finding a suitable skeleton $\Gamma\subset(\BP_K^1)^{\an}$ such that $\phi^{-1}(\Gamma)$ is a skeleton of $Y$.

The Sage package MCLF is well-adapted to this strategy. It contains an implementation of trees $\Gamma\subset(\BP^1_K)^{\an}$, and can check if the inverse image $\phi^{-1}(\Gamma)$ of a tree is a skeleton of $Y$. The implementation of our algorithms is based on this functionality.

\textbf{Outline and results}. \mbox{Chapter \ref{cha-preliminaries}} begins by recalling (semistable) models of curves. Particularly important is the characterization of models of a curve $X$ in terms of so-called \emph{Type II} valuations on the function field $F_X$. This provides the connection to analytic geometry, which we develop to the point that we need. Special attention is given the projective line $(\BP^1_K)^{\an}$. We also introduce the notion of \emph{topological ramification}, following \cite{berkovich-etale} and \cite{ctt}. The chapter ends with a discussion of reduction of analytic curves. As mentioned above, \mbox{Section \ref{sec-preliminaries-admissible}} contains an account of the strategy of ``separating branch points''; there we also give an example illustrating that this strategy does not usually work for covers of degree $p$.

\mbox{Chapter \ref{cha-greek}} introduces \emph{admissible functions}. The most important class of these is \emph{valuative functions}, which we discuss in \mbox{Section \ref{sec-greek-valuative}}. As a set, the analytification $X^{\an}$ of a curve $X$ is obtained by adding to $X$ a point $\xi$ for each valuation $v_\xi$ on the function field $F_X$ extending the valuation $v_K$ of $K$. The valuative function associated to an element $f\in F_X^\times$ is then the function that sends a point $\xi\in X^{\an}$ to the value $v_\xi(f)$.

Next, let $\phi\colon Y\to\BP_K^1$ be a morphism of degree $p$. We discuss in detail a function $\delta$ on $Y^{\an}$, which is essentially the same as the \emph{different function} studied in \cite{ctt}. In \mbox{Section \ref{sec-greek-lambda}}, we then rephrase the main results of \cite{ctt} in terms of an admissible function $\lambda$ on $X^{\an}$. The main result is \mbox{Theorem \ref{thm-trivializing-lambda}}, which states that if $\Gamma$ is a skeleton of $\BP_K^1$ separating the branch points of $\phi$ and on whose complement $\lambda$ is locally constant, then $\phi^{-1}(\Gamma)$ is a skeleton of $Y$.

\mbox{Chapter \ref{cha-analytic}} is concerned with finding explicit formulas for the value of the function $\lambda$ at closed points $x_0\in\BP_K^1$. The main results are \mbox{Theorems \ref{thm-lambda-formula-tail-case}} \mbox{and \ref{thm-lambda-formula-general-case}}. They describe $\lambda(x_0)$ in terms of certain valuative functions on $Y^{\an}$ evaluated at a closed point $p_0\in Y$ above $x_0$. These functions are derived from the coefficients of a ``good'' equation for the curve $Y$. In \mbox{Section \ref{sec-analytic-good-equation}} we explain how one can always find such a ``good'' equation. Finally, in \mbox{Section \ref{sec-analytic-tame-locus}} we introduce the \emph{tame locus} $Y^{\tame}$ associated to a curve $Y$ equipped with a degree-$p$ morphism $Y\to\BP_K^1$. Generalizing the \emph{\'etale locus} of \cite[Section 4.2]{icerm}, it is an affinoid subdomain of $Y^{\an}$ containing much of the same information as the function $\lambda$.

\mbox{Chapter \ref{cha-quartics}} deals with the reduction of plane quartic curves at $p=3$. We begin by fixing a normal form for quartics in which a degree-$3$ morphism to $\BP_K^1$ is given by a simple projection. Then we apply the general theory of \mbox{Chapter \ref{cha-analytic}} to describe the tame locus associated to a plane quartic. We descend part of the tame locus to an affinoid subdomain of $(\BP_K^1)^{\an}$ using the ``norm trick'' of \mbox{Section \ref{sec-quartics-norm-trick}} --- a crucial step, since the structure of $(\BP_K^1)^{\an}$, including the description of affinoid subdomains and admissible functions, is much simpler than the structure of $Y^{\an}$. In \mbox{Section \ref{sec-quartics-computing-delta-directly}}, we give a method to directly compute the different function $\delta$ above any given interval in $(\BP_K^1)^{\an}$. Together with the norm trick, this allows us to concretely describe the image of the tame locus in $(\BP_K^1)^{\an}$.

In \mbox{Chapter \ref{cha-applications}}, we work towards implementing our algorithm for computing the tame locus. An important point is that throughout \mbox{Chapters \ref{cha-greek}}, \mbox{\ref{cha-analytic}, and \ref{cha-quartics}} we mostly worked over an algebraically closed field, ignoring the problem of finding a suitable field extension over which a semistable reduction exists. We describe how the computation of the tame locus may be carried out over a discretely valued subfield $K_0$ of the algebraically closed field $K$. This involves describing the structure of $(\BP^1_{K_0})^{\an}$ using so-called \emph{discoids}, which generalize closed disks. In \mbox{Section \ref{sec-applications-examples}}, we apply the resulting \mbox{Algorithm \ref{alg-interior-discoid-locus}} to compute the semistable reduction in many concrete examples of curves over $\BQ_3$. In \mbox{Section \ref{sec-modular-curve}}, we also consider the example of a plane quartic curve defined over the ramified extension $\BQ_3(\zeta_3)$ of $\BQ_3$ which arose in the study of the rational points on the modular curve $X^+_{\ns}(27)$ (\cite{rouse-sutherland-zureickbrown}).

An appendix provides a short description of our implementation for computing the semistable reduction of plane quartics in Sage. It also contains various short code fragments that were used for computations in examples throughout the thesis.

\section*{Notations and conventions}\addcontentsline{toc}{chapter}{Notations and conventions}\label{sec-notations}

For reference, we collect here various notations and conventions that are used throughout.

\textbf{Valuations}. All valuations we consider have rank $1$. Thus a valuation on a field $F$ is a map $v\colon F\to\BR\cup\{\infty\}$ satisfying
\begin{equation*}
v(x+y)\ge\min\{v(x),v(y)\},\quad v(xy)=v(x)+v(y),\quad v(x)=\infty\iff x=0
\end{equation*}
for all $x,y\in F$. The image $v(F^\times)\subseteq\BR$ is called the \emph{value group} of $v$ and denoted $\Gamma_v$.

We will sometimes use the absolute value $\abs{\cdot}_v$ associated to a valuation $v$, defined by $\abs{x}_v=e^{-v(x)}$ for $x\in K$. Equivalently, we have $v(x)=-\log\abs{x}_v$ for all $x\in F$.

\textbf{The non-archimedean ground field}. Throughout this entire thesis, $K$ will denote a field that is complete with respect to a non-trivial non-archimedean valuation $v_K$. We denote the valuation ring of $K$ by $\CO_K$, the maximal ideal of $\CO_K$ by $\Fm_K$, and the residue field of $K$ by $\kappa$. We assume that $\kappa$ is a perfect field of positive characteristic, which we denote by $\characteristic(\kappa)=p$. The field $K$ may have characteristic $0$ or characteristic $p$.

We always assume that we are in one of the following two cases:
\begin{enumerate}[(a)]
\item The field $K$ is discretely valued
\item The field $K$ is the completion of an algebraic closure of a discretely valued field $K_0$
\end{enumerate}
We will always indicate which of these holds. As a rule, $K$ is usually assumed to be algebraically closed throughout \mbox{Chapters \ref{cha-greek}}, \mbox{\ref{cha-analytic}, and \ref{cha-quartics}}.

Our most important example of a valued field satisfying (a) is given by the $p$-adic numbers $\BQ_p$, and by finite extensions of $\BQ_p$. If $K$ is a finite extension of $\BQ_p$, we will always normalize $v_K$ so that $v_K(p)=1$. The most important example of a valued field satisfying (b) is the completion of an algebraic closure of $\BQ_p$, denoted $\BC_p$.

We will denote the value group of $v_K$ by $\Gamma_K\coloneqq\Gamma_{v_K}$ and the absolute value associated to $v_K$ by $\abs{\cdot}_K\coloneqq\abs{\cdot}_{v_K}$.

For convenience, we fix a multiplicative section
\begin{equation*}
\Gamma_K\to K^\times,\qquad s\mapsto\pi^s,
\end{equation*}
of $v_K$. Thus $v_K(\pi^s)=s$ and $\pi^s\pi^t=\pi^{s+t}$ for $s,t\in\Gamma_K$.

\textbf{Residue fields}. As mentioned above, the residue field of $K$ is denoted $\kappa$. We denote the residue field of any valuation $v$ extending $v_K$ by $\kappa(v)$, reflecting the fact that this residue field is in a natural way a field extension of $\kappa$. 

The residue field of a point $P$ on a $K$-curve $X$, that is the residue field of the local ring $\CO_{X,P}$, is denoted by $k(P)$.

The completed residue field in the sense of Berkovich (see \mbox{Definition \ref{def-completed-residue-field}}) at a point $\xi$ is denoted $\CH(\xi)$. Its residue field is denoted $\kappa(\xi)$. Thus if $\xi$ is of Type II, III, or IV, then we have $\kappa(\xi)=\kappa(v_\xi)$, where $v_\xi$ denotes the valuation associated to $\xi$.

\textbf{Algebraic curves}. By a \emph{curve} over a field $F$ we mean a scheme of dimension $1$ of finite type over $F$.

If $X$ is a curve over the non-archimedean ground field $K$, then we always assume that $X$ is projective, smooth, and geometrically connected. We denote the function field of such a curve by $F_X$ and its genus by $g_X$.

\textbf{The symbols $\pm\infty$}. We will make use of the extended real line $\BR\cup\{\pm\infty\}$, endowed with the usual topology. When it comes to arithmetic involving $\pm\infty$, we will only use the obvious rules $\infty+a=\infty$ for $a\in\BR\cup\{\infty\}$, $-\infty+a=-\infty$ for $a\in\BR\cup\{-\infty\}$, $a\cdot\infty=\infty$ and $a\cdot(-\infty)=-\infty$ for $a\in\BR_{>0}$, as well as $a\cdot\infty=-\infty$ and $a\cdot(-\infty)=\infty$ for $a\in\BR_{<0}$. We will never add or multiply $\infty$ and $-\infty$.

\textbf{Fonts}. We use normal italic capital letters to denote objects over the non-archimedean ground field $K$. In particular, $X$ and $Y$ always denote curves over $K$.

Objects over the residue field $\kappa$ are often written with a bar; thus $\overline{X}$ is the common notation for a $\kappa$-curve associated to the $K$-curve $X$.

Finally, analytic curves in the sense of Berkovich are denoted by serifless letters such as $\rA,\rC,\rD,\rU,\rV,\rX$.

\section*{Acknowledgements}\addcontentsline{toc}{chapter}{Acknowledgements}

Since I am undoubtedly someone who during their PhD would have struggled greatly in an adverse and difficult environment, I am all the more grateful to have had the opposite experience.

The biggest role in this falls to my advisor Stefan Wewers who accompanied this project from originally suggesting I work on semistable reduction of plane quartics all the way to the final version of this thesis. I am very thankful not only for the continuous collaboration, but also for how well we always got along on a personal level. Likewise, much thanks is due Irene Bouw who despite a full schedule of senate sessions, teaching, and report writing always had time for encouraging and helpful discussions. I am also indebted to Jeroen Sijsling for refereeing this thesis on a harsh deadline as though this were the most natural thing in the world, and for enriching our discussions with endless tidbits about anything interesting --- be it art, language, literature, or meme culture.

I thank my fellow PhD students who also contributed greatly to the atmosphere at the institute: Jeroen, Andreas, Duc, Giovanni, Robert, Kletus, and especially the two who started at the same time as I: Robert, who single-handedly turned our trip to Marseille from what could have been a tough experience for me into an enjoyable one, and of course Tim with whom I shared my office for three and a half years. I regret that despite his great efforts in exposing me to class field theory, group cohomology, and descent on elliptic curves and despite my attempts at osmotic reception, my knowledge in these areas is still severely lacking! I also thank Bernadette Maiwald, who could never not solve any administrative problem I had.

For their continuous support during my time in Ulm, I want to thank my friends and especially my parents. Knowing that I could always count on their help, even if things somehow went south, was a great comfort.

\chapter{Models, valuations, and analytification}\label{cha-preliminaries}

In this chapter, we collect various results that we will need later. We begin by recalling the theory of \emph{models} of algebraic curves. In \mbox{Section \ref{sec-preliminaries-valuations}}, we will see how to describe models using \emph{Type II valuations}. This leads naturally to the notion of \emph{analytification} in the sense of Berkovich (\cite{berkovich}). We discuss the projective line $(\BP_K^1)^{\an}$ in detail. It is both a useful illustration of the various notions we introduce and an indispensable tool for our cover-based approach to semistable reduction.

We briefly introduce \emph{topological ramification} in \mbox{Section \ref{sec-preliminaries-topological-ramification}}. Then we turn to reduction in the context of analytic curves. Semistable models of a curve $X$ correspond to \emph{skeletons} of the analytification $X^{\an}$, and it is this skeletal viewpoint that we will take in the next two chapters. In \mbox{Section \ref{sec-preliminaries-admissible}}, we review \emph{admissible reduction}. This approach is not sufficient to compute the semistable reduction of covers of degree $p$, but it is still a necessary ingredient and a good showcase of the concepts introduced in this chapter.

\section{Models of curves}
\label{sec-preliminaries-models}

As mentioned in the section on \hyperref[sec-notations]{notations and conventions}, we work over a ground field $K$ that is complete with respect to a non-archimedean valuation $v_K$. It is always assumed to be either discretely valued or to be the completion of an algebraic closure of a discretely valued field. In this section, we allow both possibilities. 

Let us fix a smooth and geometrically irreducible projective curve $X$ \mbox{over $K$}.

\begin{Def}\label{def-model}
An \emph{$\CO_K$-model} (or simply a \emph{model}) of $X$ is a normal, flat, and proper $\CO_K$-scheme $\CX$ together with an isomorphism
\begin{equation*}
X\cong\CX_K=\CX\otimes_{\CO_K}K.
\end{equation*}
\end{Def}

Chapter 10 of \cite{liu} provides a comprehensive introduction to models of algebraic curves. A fairly self-contained account of most facts we need may be found in Chapter 3 of \cite{rueth}.

\begin{Rem}\label{rem-model-definition}
\begin{enumerate}[(a)]
\item The assumption that models be normal is non-standard, but essential for our approach. The point is that we want to describe models using certain sets of valuations (see \mbox{Proposition \ref{prop-julian-theorem}}). Normality guarantees that the local rings of a model at certain codimension-$1$ points are valuation rings.
\item By \cite[Proposition 4.3.8]{liu}, the flatness in the definition of an $\CO_K$-model $\CX$ implies that $\CX$ is integral. Conversely, a surjective morphism from an integral scheme to $\Spec\CO_K$ is automatically flat, because every torsion-free module over a valuation ring is flat (\cite[Chapter 1, §2, n\textsuperscript{o} 4,Proposition 3]{bourbaki}).
\end{enumerate}
\end{Rem}

A \emph{modification} of a model $\CX$ of $X$ is a second model $\CX'$ of $X$ together with a morphism $\phi\colon\CX'\to\CX$ such that the diagram
\begin{equation*}
\begin{tikzcd}
\CX'_K \arrow[rr, "\phi_K"] \arrow[rd, "\cong"'] &   & \CX_K \arrow[ld, "\cong"] \\
                                                 & X &                          
\end{tikzcd}
\end{equation*}
commutes.

\begin{Rem}\label{rem-partial-order-on-models}
We may regard the set of isomorphism classes of models of $X$ as a partially ordered set, where $[\CX_1]\ge[\CX_2]$ means that there exists a modification $\CX_1\to\CX_2$.
\end{Rem}
The base change of a model $\CX$ of $X$ to the residue field $\kappa$ of $K$,
\begin{equation*}
\CX_s\coloneqq\CX_\kappa=\CX\otimes_{\CO_K}\kappa,
\end{equation*}
is called the \emph{special fiber} of $\CX$. The special fiber need not be smooth over $\kappa$, nor even irreducible. But it is a projective $\kappa$-curve, that is, a one-dimensional projective $\kappa$-scheme (cf.\ \cite[Lemma 8.3.3]{liu}).

\begin{Def}
The model $\CX$ is called \emph{semistable} if the special fiber $\CX_s$ is a \emph{semistable curve}. The latter means that $\CX_s$ is reduced and has only ordinary double points as singularities.
\end{Def}

The curve $X$ is said to have \emph{semistable reduction} if there exists a semistable $\CO_K$-model of $X$. It was proved by Deligne and Mumford \cite{deligne-mumford} that this is always the case after taking a finite separable field extension $L/K$ --- of course, if $K$ is algebraically closed, this implies that every $K$-curve already has semistable reduction.

\begin{specialSatz}{Semistable Reduction Theorem}\label{thm-semistable-reduction}
There exists a finite separable field extension $L/K$ such that $X_L$ has semistable reduction.
\end{specialSatz}

Determining for certain classes of curves $X$ a field extension $L/K$ such that $X_L$ has semistable reduction and a semistable model $\CX$ of $X_L$ is the ultimate goal of this thesis.

\section{Type II valuations}\label{sec-preliminaries-valuations}

We continue the notation from the last section. Thus the ground field $K$ may still be algebraically closed or discretely valued, and $X$ denotes a smooth and geometrically irreducible projective $K$-curve.

\begin{Def}\label{def-type-ii-valuation}
A valuation 
\begin{equation*}
v\colon F_X\to\BR\cup\{\infty\}
\end{equation*}
is called a \emph{Type II valuation} if $v$ extends the valuation $v_K$ on $K$ and the transcendence degree of the residue field extension $\kappa(v)/\kappa$ is $1$.
\end{Def}

Suppose that $\CX$ is a model of $X$ and that $Z$ is an irreducible component of the special fiber $\CX_s$. Denote the generic point of $Z$ by $\xi_Z$, so that $\overline{\{\xi_Z\}}=Z$. Then since $\CX$ is normal, the local ring $\CO_{\CX,\xi_Z}$ is a valuation ring. In the case that $K$ is discretely valued and hence noetherian, this is clear, because every noetherian normal local domain of dimension $1$ is a discrete valuation ring. For the case that $K$ is algebraically closed, see \cite[Lemma 3.1]{green}.

We denote the valuation associated to the valuation ring $\CO_{\CX,\xi_Z}$ by $v_Z$; it is a Type II valuation on $F_X$ (for details see \cite[Proposition 3.4]{rueth} in the discretely valued case and again \cite[Lemma 3.1]{green} in general).

\begin{Prop}\label{prop-julian-theorem}
The map
\begin{equation*}
\CX\mapsto\{v_Z\mid\textrm{$Z\subseteq\CX_s$ irreducible component}\}
\end{equation*}
induces an isomorphism of partially ordered sets 
\begin{IEEEeqnarray*}{Cl}
&\quad\textrm{isomorphism classes of models of $X$}\\\longrightarrow&\quad\textrm{finite non-empty sets of Type II valuations on $F_X$}.
\end{IEEEeqnarray*}
\end{Prop}
\begin{proof}
For the case that $K$ is discretely valued, this is the main result of \cite[Chapter 3]{rueth}, see \cite[Corollary 3.18]{rueth}.

In general, this follows from \cite[Theorem 2]{green}. Note that the ``Local Skolem Property'' on which this theorem depends is satisfied for the complete field $K$. The ``proper, normal, integral $\CO_K$-curves'' considered in loc.\,cit.\ with function field $F_X$ are the same as our $\CO_K$-models of $X$, see \mbox{Remark \ref{rem-model-definition}}(b).
\end{proof}

The map described in \mbox{Proposition \ref{prop-julian-theorem}} is an isomorphism of partially ordered sets, meaning that it is a bijection preserving the partial orders on domain and codomain. The partial order on isomorphism classes of models of $X$ is explained in \mbox{Remark \ref{rem-partial-order-on-models}}. Sets of Type II valuations are simply ordered by inclusion.

\begin{Prop}\label{prop-julian-normalization}
Let $Y\to X$ be a cover of smooth irreducible projective curves over $K$. Let $\CX$ be an $\CO_K$-model of $X$ and let $\CY$ be the normalization of $\CX$ in the function field $F_Y$. Then $\CY$ is an $\CO_K$-model of $Y$, and the valuations corresponding to $\CY$ via the bijection of \mbox{Proposition \ref{prop-julian-theorem}} are the extensions of the valuations corresponding to $\CX$.
\end{Prop}
\begin{proof}
This is explained in \cite[Section 5.1.2]{rueth} for the case that $K$ is discretely valued. The exact same proof is valid for algebraically closed $K$ as well; all that is needed is replacing ``discrete valuation ring'' with ``valuation ring''.
\end{proof}

\mbox{Propositions \ref{prop-julian-theorem}} \mbox{and \ref{prop-julian-normalization}} are our chief tools for describing models of curves. Given a finite cover of curves $\phi\colon Y\to X$, we may describe a model of $Y$ simply using a non-empty set of Type II valuations on $X$. The corresponding model of $Y$ is obtained by normalization, as explained in \mbox{Proposition \ref{prop-julian-normalization}}.

The next lemma is an example of how valuation theory may be leveraged to gain geometric insight into the special fibers of models.

\begin{Lem}\label{lem-ramification-theory}
Let $\phi\colon Y\to X=\BP_K^1$ be a degree-$p$ morphism of $K$-curves. Suppose that $w$ is a Type II valuation on $F_Y$ with restriction $v\coloneqq\rest{w}{F_X}$. If the residue field extension $\kappa(w)/\kappa(v)$ is inseparable, then the function field $\kappa(w)$ has genus $0$.
\end{Lem}
\begin{proof}
Since the special fiber of any model of $\BP_K^1$ has arithmetic genus $0$, the function field $\kappa(v)$ must have genus $0$. 

Consider the subfield $\kappa(w)^p$ of $\kappa(w)$ consisting of $p$-th powers. The extension $\kappa(w)/\kappa(w)^p$ is inseparable of degree $p$ by \cite[Proposition II.2.11]{silverman} --- here we use the assumption that $\kappa$ is perfect. Thus $\kappa(w)^p=\kappa(v)$. Since the Frobenius maps $\kappa(w)$ isomorphically to $\kappa(w)^p$, both these function fields must have the same genus $0$.
\end{proof}

\textit{For the rest of this section, we assume that the base field $K$ is discretely valued.}

For the following definition, suppose that $L/K$ is the completion of an algebraic extension of $K$. Thus there is a unique extension $v_L$ of $v_K$ to $L$. We denote the ring of integers of $L$ by $\CO_L$. Of course, if $L/K$ is finite, then $L$ already is complete with respect to the unique extension of $v_K$ to $L$.

\begin{Def}\label{def-normlizaed-basechange}
Let $X$ be a smooth and geometrically irreducible curve over $K$ and let $\CX$ be an $\CO_K$-model of $X$. The \emph{normalized base change} $\CX_L$ of $\CX$ to $L$ is the normalization of $\CX\otimes_{\CO_K}\CO_L$.
\end{Def}

The normalized base change $\CX_L$ is an $\CO_L$-model of $X_L=X\otimes_KL$; in fact, it is the normalization of $\CX$ in the function field of $X_L$. 

Normalized base change is quite important to us. For certain curves $X$ over $K$, we will see in subsequent chapters how to compute a \emph{potentially semistable} $\CO_K$-model $\CX$ of $X$, meaning that the normalized base change $\CX_L$ to a finite extension $L/K$ is semistable. The problem of computing a semistable reduction of $X$ is then reduced to determining this field extension $L/K$. We will return to this topic in \mbox{Chapter \ref{cha-applications}}.

\begin{Def}\label{def-permanent}
An $\CO_K$-model $\CX$ of $X$ is called \emph{permanent} if $\CX_s$ is reduced.
\end{Def}

\begin{Rem}\label{rem-permanent}
The motivation for the terminology ``permanent'' is the following. Let $\CX$ be a permanent $\CO_K$-model and consider a field extension $L/K$ that is the completion of an algebraic extension of $K$. Then the special fiber of the normalized base change $\CX_L$ is simply given by $\CX_s\otimes_\kappa\lambda$, where $\lambda$ denotes the residue field of $L$ (cf.\ \cite[Section 2.2]{arzdorf-wewers}). To see this, note that $\kappa$ being perfect, \cite[Corollary 3.2.14(a)]{liu} implies that the special fiber of $\CX\otimes_{\CO_K}\CO_L$,
\begin{equation*}
(\CX\otimes_{\CO_K}\CO_L)\otimes_{\CO_L}\lambda\cong\CX\otimes_{\CO_K}\lambda\cong(\CX\otimes_{\CO_K}\kappa)\otimes_\kappa\lambda=\CX_s\otimes_\kappa\lambda,
\end{equation*}
is still reduced. From \cite[Lemma 4.1.18]{liu} it follows that $\CX\otimes_{\CO_K}\CO_L$ is normal, so that $\CX\otimes_{\CO_K}\CO_L$ is already equal to its normalized base change $\CX_L$.
\end{Rem}

We end this section by explaining some further merits of the valuation-theoretic approach to models. As mentioned above, the key point is that much information on the special fiber $\CX_s$ of a model $\CX$ is encoded in the associated set of Type II valuations. We never describe models by their global geometry; there is no need, say, to worry about how a model might be embedded in some projective space $\BP^N_{\CO_K}$. 

To elaborate further, suppose that $\CX$ is an $\CO_K$-model of $X$. Then if $Z\subseteq\CX_s$ is an irreducible component with generic point $\xi_Z$, the residue field of the valuation $v_Z$ is the same as the function field of the component $Z=\overline{\{\xi_Z\}}$. Thus knowledge of the valuations corresponding to $\CX$ means knowing the components of $\CX_s$ ``birationally''. Using the fact that the (arithmetic) genera of $X$ and $\CX_s$ coincide, it is easy to derive a criterion for $\CX$ being semistable. This is explained in detail in \cite[Section 5.3]{rueth}. 

A simple example of this is the following: If the sum of the geometric genera $g_Z$ over all components $Z\subseteq\CX_s$ equals $g_X$, then $\CX$ must be semistable. The curve $X$ is said to have \emph{abelian reduction} in this case; the arithmetic genus of $\CX_s$ comes only from the geometric genera of the components $Z\subseteq\CX_s$, not from any loops in the intersection graph of $\CX_s$ (cf.\ \mbox{Definition \ref{def-intersection-graph}}).

We can also detect non-reduced structure of $\CX_s$ using valuations. Given a component $Z\subseteq\CX_s$ with generic point $\xi_Z$, the \emph{multiplicity} of $Z$ is defined to be the length
\begin{equation*}
m(Z)\coloneqq\len(\CO_{\CX_s,\xi_Z})=\len(\CO_{\CX,\xi_Z}/\Fm_K\CO_{\CX,\xi_Z}).
\end{equation*}
(Recall that $\Fm_K$ denotes the maximal ideal of the valuation ring $\CO_K$.) The multiplicities of the components $Z\subseteq\CX_s$ measure the degree to which $\CX_s$ is non-reduced. Indeed, $\CO_{\CX_s,\xi_Z}$ is a local Artinian ring, and it is reduced if and only if its length equals $1$. But if all the local rings $\CO_{\CX_s,\xi_Z}$ are reduced, then by normality of $\CX$, the special fiber $\CX_s$ must already be reduced (\cite[\mbox{Corollary 2.2(ii)}]{arzdorf-wewers}).

\begin{Lem}\label{lem-ramification-index-and-multiplicity}
Let $\phi\colon Y\to X=\BP_K^1$ be a finite cover of $K$-curves, let $\CX$ be a permanent $\CO_K$-model of $X$ and let $\CY$ be the normalization of $\CX$ in $F_Y$. Let $Z\subseteq\CX_s$ and $W\subseteq\CY_s$ be components with $\phi(W)=Z$. Then the ramification index $e(v_W|v_Z)$ equals the multiplicity of the component $W$.
\end{Lem}
\begin{proof}
Let $\eta$ and $\xi$ be the generic points of $W$ and $Z$ respectively. Note that we have $\Fm_K\CO_{\CX,\xi}=(\pi_Z^{m(Z)})$, where $\pi_Z$ is a uniformizer for the discrete valuation ring $\CO_{\CX,\xi}$. Thus $m(Z)$ equals the ramification index of the extension of valuations $v_Z|v_K$. Similarly, $m(W)$ equals the ramification index of $v_W|v_K$. Thus the result follows from the multiplicativity of ramification indices in the tower
\begin{equation*}
\begin{tikzcd}
v_W \arrow[d, no head] \\
v_Z \arrow[d, no head]  \\
v_K     .                 
\end{tikzcd}
\end{equation*}
\end{proof}

\section{Analytic curves}\label{sec-preliminaries-analytic}

In this section, we introduce analytic curves in the sense of Berkovich (\cite{berkovich}). Berkovich defines a category of \emph{$K$-analytic spaces}, which are certain locally ringed spaces. \emph{Analytic curves} are treated in detail in Chapter 4 of \cite{berkovich}. All the analytic spaces we consider will be curves. We will denote them by serifless letters such as $\rA,\rC,\rD,\rU,\rV,\rX$. 

To every finite-type $K$-scheme $X$ is attached an \emph{analytification}, denoted $X^{\an}$, see \cite[Section 3.4]{berkovich}. It is a $K$-analytic space representing the functor that sends an analytic space $\rU$ to the set of morphisms of $K$-ringed spaces $\Hom(\rU, X)$. Thus it has attached to it a universal morphism
\begin{equation*}
\pi\colon X^{\an}\to X.
\end{equation*}
The association $X\mapsto X^{\an}$ is, of course, functorial: for every morphism $\phi\colon X\to Y$ of finite-type $K$-schemes there is an associated morphism $\phi^{\an}\colon X^{\an}\to Y^{\an}$. We will usually omit the superscript and simply write $\phi$ for $\phi^{\an}$.

Suppose for example that $\iota\colon U\to X$ is the inclusion of an open subscheme. Then $\iota^{\an}\colon U^{\an}\to X^{\an}$ is also an open immersion (\cite[Proposition 3.4.6]{berkovich}). We will often use the universal morphism $\pi\colon X^{\an}\to X$ to consider functions on $X$ as functions on $X^{\an}$, by pulling back along the map of sheaves $\CO_X\to\pi_*\CO_{X^{\an}}$. Similarly, we will consider functions on open subschemes $U\subseteq X$ as functions on $U^{\an}\subseteq X^{\an}$.

Analytic spaces are made from local building blocks, so-called \emph{affinoid spaces} (\cite[Section 2.2]{berkovich}). A subset of $X^{\an}$ is called an \emph{affinoid subdomain} if it is the image of a morphism from an affinoid space to $X^{\an}$ satisfying a certain universal property (see \cite[Section 3.1]{berkovich}). In the next section, we will get to know some examples of these (\emph{closed disks} and \emph{anuli}), and later we will study affinoid subdomains of $X^{\an}$ determined by the values of certain algebraic functions on $X^{\an}$ (see \mbox{Lemma \ref{lem-rational-domain}}). But for now, we take a different approach to defining $X^{\an}$, based entirely on valuation theory.

As a set, $X^{\an}$ may be described as follows (compare \cite[Remark 3.4.2]{berkovich}). It consists of all pairs $(P,v)$, where $P\in X$ is any point and where $v\colon k(P)\to\BR\cup\{\infty\}$ is a valuation extending $v_K$. Usually, we will denote points in $X^{\an}$ by the Greek letters $\xi$ and $\eta$ (with various subscripts). 

\begin{Def}\label{def-completed-residue-field}
For a point $\xi=(P,v)$, the completion of $k(P)$ with respect to $v$ is called the \emph{completed residue field} at $\xi$ and is denoted $\CH(\xi)$.

We denote by $\kappa(\xi)$ the residue field of the valued field $\CH(\xi)$.
\end{Def}

The map $\pi\colon X^{\an}\to X$ is, on the level of sets, simply the ``forgetful map'' $(P,v)\mapsto P$. Since there exists a unique extension of $v_K$ to every finite extension $L/K$, there is a unique point $\xi=(P,v)$ above every closed point $P\in X$.

In the case that $X$ is an irreducible curve, the only non-closed point on $X$ is the generic point. The points $\xi\in X^{\an}$ above the generic point correspond to the valuations on the function field $F_X$ extending $v_K$. We denote the valuation belonging to such a point by $v_\xi$.

We cannot associate valuations on $F_X$ to the points $\xi=(P,v)$ in $X^{\an}$ associated to closed points of $X$, but we can get close: precompose $v$ with the natural map $\CO_{X,P}\to k(P)$ to obtain a map $\tilde{v}_\xi\colon\CO_{X,P}\to\BR\cup\{\infty\}$. This map satisfies the usual rules for valuations
\begin{equation}\label{equ-valuation-rules-1}
\tilde{v}_\xi(0)=\infty,\qquad \tilde{v}_\xi(1)=0,
\end{equation}
\begin{equation}\label{equ-valuation-rules-2}
\tilde{v}_\xi(fg)=\tilde{v}_\xi(f)+\tilde{v}_\xi(g),\quad \tilde{v}_\xi(f+g)\ge\min(\tilde{v}_\xi(f),\tilde{v}_\xi(g))
\end{equation}
for all $f,g\in\CO_{X,P}$, but it also sends every element of the maximal ideal $\Fm_P$ of $\CO_{X,P}$ to $\infty$.

\begin{Def}\label{def-pseudovaluation}
A \emph{pseudovaluation} on a ring $R$ is a function $v\colon R\to\BR\cup\{\infty\}$ satisfying the rules \eqref{equ-valuation-rules-1} and \eqref{equ-valuation-rules-2}.
\end{Def}

We extend the pseudovaluation $\tilde{v}_\xi$ to the function field $F_X$ by defining
\begin{equation}\label{equ-pseudovaluation}
v_\xi\colon F_X\to\BR\cup\{\pm\infty\},\qquad f\mapsto\begin{cases}
\tilde{v}_\xi(f)&\textrm{if $f\in\CO_{X,P}$},\\
-\infty&\textrm{otherwise}.
\end{cases}
\end{equation}
Thus we have associated to each $\xi\in X^{\an}$ a function $v_\xi\colon F_X\to\BR\cup\{\pm\infty\}$. Collectively, we call all functions $v_\xi$, where $\xi\in X^{\an}$, \emph{pseudovaluations} on $F_X$, be they valuations or functions of the form \eqref{equ-pseudovaluation}.

Abhyankar's Inequality (\cite[Theorem 3.4.3]{engler-prestel}) states that
\begin{equation*}
\trdeg(\kappa(\xi)/\kappa)+\rk_\BQ(\Gamma_{v_\xi}/\Gamma_K)\le\trdeg(F_X/K)=1.
\end{equation*}
In words: The transcendence degree of $F_X/K$ is bounded by the sum of the transcendence degree of the extension of residue fields $\kappa(\xi)/\kappa$ and the rational rank of the quotient of value groups $\Gamma_{v_\xi}/\Gamma_K$. (The rational rank of an abelian group $A$ is the dimension of the $\BQ$-vector space $\BQ\otimes A$.) 

The points of $X^{\an}$ are categorized into four \emph{types} (the definitions below follow \cite[Section 1.4.4]{berkovich}; see also \cite[Section 3.3.2]{ducros} for the case of $K$ not being algebraically closed):

\begin{itemize}
\item $\xi$ is called \emph{Type I} if $\CH(\xi)$ may be embedded in the completion of an algebraic closure of $K$. In particular, the points of $X^{\an}$ associated to closed points of $X$ are Type I.
\item $\xi$ is called \emph{Type II} if $v_\xi$ is a Type II valuation (\mbox{Definition \ref{def-type-ii-valuation}}), that is, if the extension of residue fields $\kappa(\xi)/\kappa$ has transcendence degree $1$. By Abhyankar's Inequality, in this case we have $\rk_\BQ(\Gamma_{v_\xi}/\Gamma_K)=0$. Note that \cite[Theorem 3.4.3]{engler-prestel} also includes the statement that $\Gamma_{v_\xi}/\Gamma_K$ is finitely generated if $\kappa(\xi)/\kappa$ has transcendence degree $1$. Thus the index $[\Gamma_{v_\xi}:\Gamma_K]$ is finite; it is the ramification index of the extension $\CH(\xi)/K$.

\item $\xi$ is called \emph{Type III} if $\rk_\BQ(\Gamma_{v_\xi}/\Gamma_K)=1$. In this case, we necessarily have $\trdeg(\kappa(\xi)/\kappa)=0$.
\item \sloppy$\xi$ is called \emph{Type IV} if $\xi$ is not Type I, but $\rk_\BQ(\Gamma_{v_\xi}/\Gamma_K)=0$ and $\trdeg(\kappa(\xi)/\kappa)=0$.
\end{itemize}
Since for a Type II point $\xi\in X^{\an}$ the residue field $\kappa(\xi)$ has transcendence degree $1$ over $\kappa$, that is, since $\kappa(\xi)$ is an algebraic function field over $\kappa$, the following definition makes sense.

\begin{Def}\label{def-genus-of-point}
\begin{enumerate}[(a)]
\item The smooth projective $\kappa$-curve with function field $\kappa(\xi)$ is called the \emph{reduction curve} at $\xi$.
\item The \emph{genus} of a Type II point $\xi\in X^{\an}$, denoted $g(\xi)$, is the genus of the reduction curve at $\xi$.

By convention, we also define $g(\xi)=0$ for all points $\xi\in X^{\an}$ not of \mbox{Type II}.
\end{enumerate}
\end{Def}

\begin{Def}\label{def-reduction-to-residue-field}
Let $\xi\in X^{\an}$ be a Type II point and let $f\in F_X^\times$ be a function with $v_\xi(f)\in\Gamma_K$. Then we define
\begin{equation*}
[f]_\xi\coloneqq\overline{f\pi^{-v_\xi(f)}}\in\kappa(\xi)^\times.
\end{equation*}
(Recall from the section on \hyperref[sec-notations]{notations and conventions} that for $s\in \Gamma_K$, we always denote by $\pi^s$ an element of $K^\times$ with $v_K(\pi^s)=s$. The definition of $[f]_\xi$ depends on the choice of $\pi^s$, but a different choice of $\pi^s$ will only change $[f]_\xi$ by multiplication with a non-zero constant in $\kappa$.)
\end{Def}

\begin{Rem}\label{rem-divisible-groups}
Note that if $K$ is algebraically closed, we have $v_\xi(f)\in \Gamma_K$ for all $f\in F_X^\times$. Indeed, the divisible group $\Gamma_K=\BQ$ is a subgroup of $\Gamma_{v_\xi}$, which is itself a subgroup of $\BR$. Since $\xi$ is Type II, $\Gamma_{v_\xi}/\Gamma_K$ is finite, so for every $\alpha\in\Gamma_{v_\xi}$, there exists some $n\ge1$ with $n\alpha\in\BQ$. But then $\alpha\in\BQ$ as well, so $\Gamma_{v_\xi}=\Gamma_K$.
\end{Rem}

\begin{Lem}\label{lem-bracket-properties}
Let $f,g\in F_X^\times$ be non-zero functions with $v_\xi(f),v_\xi(g),v_\xi(f+g)\in \Gamma_K$. Then we have:
\begin{enumerate}[(a)]
\item $[fg]_\xi=[f]_\xi\cdot[g]_\xi$
\item $[f+g]_\xi=[f]_\xi+[g]_\xi$ if and only if $v_\xi(f)=v_\xi(g)$ and $[f]_\xi+[g]_\xi\ne0$
\end{enumerate}
\end{Lem}
\begin{proof}
We always assume that the section $s\mapsto\pi^s$ is multiplicative. Thus for (a), we can simply calculate
\begin{equation*}
[fg]_\xi=\overline{fg\pi^{-v_\xi(fg)}}=\overline{f\pi^{-v_\xi(f)}g\pi^{-v_\xi(g)}}=\overline{f\pi^{-v_\xi(f)}}\cdot\overline{g\pi^{-v_\xi(g)}}=[f]_\xi[g]_\xi.
\end{equation*}
For (b), let us first remark that under the condition $v_\xi(f)=v_\xi(g)$ we have
\begin{equation*}
[f]_\xi+[g]_\xi=0\qquad\textrm{if and only if}\qquad v_\xi(f+g)>v_\xi(f).
\end{equation*}
Indeed, we have
\begin{equation*}
[f]_\xi+[g]_\xi=\overline{f\pi^{-v_\xi(f)}}+\overline{g\pi^{-v_\xi(g)}}=\overline{(f+g)\pi^{-v_\xi(f)}},
\end{equation*}
and this is $0$ if and only if $v_\xi(f+g)>v_\xi(f)$.

Now the equivalence in (b) follows easily. Indeed, if $v_\xi(f)=v_\xi(g)$ and $[f]_\xi+[g]_\xi\ne0$, then $v_\xi(f+g)=v_\xi(f)=v_\xi(g)$, and
\begin{equation*}
[f+g]_\xi=\overline{(f+g)\pi^{-v_\xi(f+g)}}=\overline{f\pi^{-v_\xi(f+g)}}+\overline{g\pi^{-v_\xi(f+g)}}=[f]_\xi+[g]_\xi.
\end{equation*}
On the other hand, if $v_\xi(f)\ne v_\xi(g)$, say $v_\xi(f)>v_\xi(g)$, then $v_\xi(f+g)=v_\xi(g)$ and a similar calculation easily shows $[f+g]_\xi=[g]_\xi\ne[f]_\xi+[g]_\xi$. Likewise, if $[f]_\xi+[g]_\xi=0$, then clearly $[f+g]_\xi\ne[f]_\xi+[g]_\xi$. 
\end{proof}

So far, we have not discussed the topology of $X^{\an}$. We now give a definition that is well-adapted to our valuation-theoretic definition of the analytification.

Suppose that $U\subseteq X$ is an affine open. For each function $f\in\CO_X(U)$ and each point $\xi=(P,v)$ in $U^{\an}$ we may consider the value $v_\xi(f(P))$, where $f(P)$ denotes the image of $f$ in $k(P)$. Then the topology on $U^{\an}$ is the coarsest one such that the evaluation map
\begin{equation*}
U^{\an}\to\BR\cup\{\infty\},\qquad\xi\mapsto v_\xi(f(P)),
\end{equation*}
is continuous for every $f\in\CO_X(U)$ (cf.\ \cite[Section 2.1]{helminck-ren}).

\begin{Def}
A subset $I\subseteq X^{\an}$ is called an \emph{interval} if it is homeomorphic to an interval $J \subseteq\BR\cup\{\pm\infty\}$.
\end{Def}

Let $\xi_1,\xi_2\in X^{\an}$ be two points. If $I\subseteq X^{\an}$ is an interval containing $\xi_1,\xi_2$ homeomorphic to a compact interval $J\subset\BR\cup\{\pm\infty\}$, with $\xi_1,\xi_2$ mapping to the boundary points of $J$, we will use the notation
\begin{equation*}
I=[\xi_1,\xi_2]
\end{equation*}
and say that $I$ is a \emph{compact interval}. Similarly, we will write
\begin{equation*}
[\xi_1,\xi_2)\coloneqq[\xi_1,\xi_2]\setminus\{\xi_2\},\quad(\xi_1,\xi_2]\coloneqq[\xi_1,\xi_2]\setminus\{\xi_1\},\ \quad(\xi_1,\xi_2)\coloneqq[\xi_1,\xi_2]\setminus\{\xi_1,\xi_2\}.
\end{equation*}

\begin{Def}\label{def-branch}
A \emph{branch} at a point $\xi\in X^{\an}$ is an equivalence class of intervals $[\xi,\xi']$ starting at $\xi$, where $[\xi,\xi_1]$ and $[\xi,\xi_2]$ are considered equivalent if there exists an interval $[\xi,\xi_3]$ contained in both $[\xi,\xi_1]$ and $[\xi,\xi_2]$.
\end{Def}

The curve $X^{\an}$ is a \emph{special quasipolyhedron} (see \cite[Definition 4.1.1]{berkovich} for what this means and \cite[Theorem 4.3.2]{berkovich} for the result). We will not discuss all technicalities tied up in this definition here, but note that this implies in particular that the topology of $X^{\an}$ has a basis consisting of open sets $U\subseteq X^{\an}$ with the following properties:
\begin{enumerate}[(a)]
\item The boundary of $U$ is finite
\item For any distinct $\xi_1,\xi_2\in U$, there exists a unique compact interval $I\subseteq U$ with endpoints $\xi_1$ and $\xi_2$
\end{enumerate}

If these two properties are even true for the open subset $U=X^{\an}$, then $X^{\an}$ is called a \emph{simply connected special quasipolyhedron}. The most important example of a simply connected special quasipolyhedron is the analytification $(\BP_K^1)^{\an}$ of the projective line $\BP_K^1$, which we discuss in the next section.

\section{The projective line}\label{sec-preliminaries-line}

In this section, we assume that $K$ is algebraically closed. However, we will revisit the projective line in the setting of discretely valued $K$ in \mbox{Section \ref{sec-applications-discoids}}.

Our object of study is the analytification $(\BP_K^1)^{\an}$ of the projective line $\BP_K^1$. As explained in the previous section, the points of $(\BP_K^1)^{\an}$ correspond to the pseudovaluations on the rational function field $K(x)$. By definition, these are either valuations on $K(x)$ extending $v_K$, or else are obtained from evaluation at closed points in $\BP_K^1$. With the exception of the closed point $\infty\in\BP_K^1$, all points of $(\BP_K^1)^{\an}$ are induced by pseudovaluations on $K[x]$. Type I points different from $\infty$ may be identified with closed points in $\Spec K[x]$, or simply with the elements of $K$.

\begin{Def}\label{def-disk}
Let $a\in K$ be any element and let $r\in\BR$ be a real number. The \emph{closed disk} with center $a$ and radius $r$ is the subset
\begin{equation*}
\rD[a, r]\coloneqq\{\xi\in(\BP_K^1)^{\an}\mid v_\xi(x-a)\ge r\}.
\end{equation*}
\end{Def}
The closed disk $\rD[a,r]$ is an example of an affinoid subdomain of $(\BP_K^1)^{\an}$. Affinoid spaces are defined to be the \emph{spectra} of \emph{affinoid $K$-algebras} (\cite[Section 2.1]{berkovich}), which are certain commutative Banach algebras. As a set, the spectrum $\CM(\CA)$ of an affinoid $K$-algebra $\CA$ consists of all pseudovaluations $v$ on $\CA$ such that the absolute value $\abs{\cdot}_v$ associated to $v$ is bounded by the given norm on $\CA$ (\cite[Section 1.2]{berkovich}). The affinoid algebra giving rise to the closed disk $\rD[a,r]$ is the \emph{Tate algebra} of \emph{restricted power series}
\begin{equation*}
K\{\pi^{-r}t\}\coloneqq\big\{\sum_{i=0}^\infty a_it^i\in K\llbracket t\rrbracket:\lim_{i\to\infty}v_K(a_i)+ir=\infty\big\}
\end{equation*}
given as a norm the absolute value associated to the \emph{Gauss valuation} of radius $r$
\begin{equation*}
v_r\colon K\{\pi^{-r}t\}\to\BR\cup\{\infty\},\qquad\sum_{i=0}^\infty a_it^i\mapsto\min_{i\ge0}\big\{v_K(a_i)+ir\big\}.
\end{equation*}
(The Gauss valuations are so named because in proving that $v_r$, where $r\in\BQ$, is really a valuation, one can use the same strategy as in the proof of Gauss's Lemma. See for example \cite[13]{bosch-lectures}. The notation $K\{\pi^{-r}t\}$ is so chosen because this Tate algebra is a rescaled version of the Tate algebra $K\{t\}$; one can think of this rescaling as replacing the parameter $t$ with $\pi^{-r}t$, where $\pi^r$ is an element of valuation $r$.)

The \emph{spectrum} $\CM(K\{\pi^{-r}t\})$ of $K\{\pi^{-r}t\}$ is by definition the set of all pseudovaluations $v$ on $K\{\pi^{-r}t\}$ satisfying $v_r(f)\le v(f)$ for all $f\in K\{\pi^{-r}t\}$. We can identify the affinoid space $\CM(K\{\pi^{-r}t\})$ with the disk $\rD[a, r]$ as follows. Write an arbitrary element $f\in K[x]$ as a polynomial in $t=x-a$, say $f=\sum_ia_it^i$. Then we may regard $f$ as an element of $K\{\pi^{-r}t\}$ and apply any $v\in\CM(K\{\pi^{-r}t\})$ to it. In this way, $v$ defines a pseudovaluation on $K[t]$, which extends uniquely to the rational function field $K(x)$. It is easy to see that the pseudovaluations obtained in this way are exactly the elements of $\rD[a,r]$.

We denote the valuation on $K(x)$ corresponding to $v_r$ by $v_{a,r}$. It is called the \emph{Gauss valuation of radius $r$ centered at $a$}. The point corresponding to $v_{a,r}$ is the unique boundary point of $\rD[a,r]$. Given a closed disk $\rD$, we will also denote the valuation corresponding to the boundary point of $\rD$ by $v_\rD$ and call it the \emph{Gauss valuation associated to $\rD$}.

Next, consider the \emph{open disk} with center $a$ and radius $r$, the subset
\begin{equation*}
\rD(a,r)\coloneqq\{\xi\in(\BP_K^1)^{\an}\mid v_\xi(x-a)>r\}.
\end{equation*}
It is the union of all closed disks of radius $r'>r$ centered at $a$. Correspondingly, we have
\begin{equation}\label{equ-its-an-inverse-limit}
\CO_{(\BP_K^1)^{\an}}\big(\rD(a,r)\big)=\varprojlim_{r'>r}K\{\pi^{-r'}\}
\end{equation}
(\cite[Section 2.3]{berkovich}). The \emph{radius of convergence} of a power series $f=\sum_{i=0}^\infty a_it^i\in K\llbracket t\rrbracket$ is defined to be
\begin{equation}\label{equ-radius-of-convergence}
\rho(f)\coloneqq-\liminf_{i\to\infty}\frac{v_K(a_i)}{i}.
\end{equation}

\begin{Lem}\label{lem-radius-of-convergence}
\sloppy The radius $\rho(f)$ is the smallest radius for which $f$ is contained in $\CO_{(\BP_K^1)^{\an}}(\rD(a, r))$.
\end{Lem}

A similar statement can be found in \cite[Proposition 5.4.1]{gouvea}, but not phrased in the language of analytic spaces in the sense of Berkovich, and using absolute values rather than valuations.

\begin{proof}[Proof of Lemma \ref{lem-radius-of-convergence}]
That $-\rho(f)$ equals the limes inferior in \eqref{equ-radius-of-convergence} means that for every $\epsilon>0$ we have
\begin{equation}\label{equ-radius-of-convergence-liminf-1}
\frac{v_K(a_i)}{i}\ge-\rho(f)-\epsilon
\end{equation}
for almost all $i\ge0$ and that we have
\begin{equation}\label{equ-radius-of-convergence-liminf-2}
\frac{v_K(a_i)}{i}\le -\rho(f)+\epsilon
\end{equation}
for infinitely many $i\ge0$. Now if $r>\rho(f)$, we choose $\epsilon\coloneqq \frac{r-\rho(f)}{2}$. Then by \eqref{equ-radius-of-convergence-liminf-1} we have
\begin{equation*}
v_r(a_it^i)\ge -i\rho(f)-i\frac{r-\rho(f)}{2}+ir=i\frac{r-\rho(f)}{2}=i\epsilon
\end{equation*}
for almost all $i\ge0$. Since $r>\rho(f)$ was arbitrary, $f$ defines an element of \eqref{equ-its-an-inverse-limit}. Now suppose that $r<\rho(f)$. Choose $\epsilon\coloneqq\rho(f)-r$. Then by \eqref{equ-radius-of-convergence-liminf-2} we have
\begin{equation*}
v_r(a_it^i)\le -i\rho(f)+i\big(\rho(f)-r\big)+ir=0
\end{equation*}
for infinitely many $i\ge0$. Since $r<\rho(f)$ was arbitrary, $f$ is not contained in \eqref{equ-its-an-inverse-limit} for any $r<\rho(f)$.
\end{proof}

A different type of affinoid space is the \emph{closed anulus} centered at $a\in K$ with inner radius $R\in\BR$ and outer radius $r\in\BR$ satisfying $R>r$. It is the set
\begin{equation*}
\rA[a,R,r]\coloneqq\{\xi\in(\BP_K^1)^{\an}\mid R\ge v_\xi(x-a)\ge r\}.
\end{equation*}
The affinoid $K$-algebra giving rise to $\rA[a,R,r]$ is
\begin{equation*}
K\{\pi^{-r}t,\pi^Rt^{-1}\}\coloneqq\big\{\sum_{i=-\infty}^\infty a_it^i:\lim_{i\to\infty}v_K(a_i)+ir=\lim_{i\to-\infty}v_K(a_i)+iR=\infty\big\}.
\end{equation*}
(See \cite[Section 2.1]{bpr} for the definition of the norm on this Banach algebra.) Similarly, the \emph{open anulus} centered at $a\in K$ with radii $R>r$ is the subset
\begin{equation*}
\rA(a,R,r)\coloneqq\{\xi\in(\BP_K^1)^{\an}\mid R>v_\xi(x-a)> r\}.
\end{equation*}
Anuli have two boundary points, corresponding to their two radii. In fact, the two boundary points of $\rA[a,R,r]$ are the Gauss valuations of radii $r$ and $R$ centered at $a$; the anulus $\rA[a,R,r]$ equals $\rD[a,r]\setminus\rD(a,R)$, and similarly for $\rA(a,R,r)$.

\begin{Rem}\label{rem-generalized-domains}
So far, we have defined open and closed disks and open and closed anuli as subdomains of $(\BP_K^1)^{\an}$. In the sequel, we will call any affinoid space isomorphic to $\rD[a,r]$ a closed disk. Similar remarks apply to open disks and closed and open anuli.
\end{Rem}

We come now to the elegant description of the points of $(\BP_K^1)^{\an}$ using closed disks, following \cite[Section 1.4]{berkovich}.

\begin{Def}\label{def-nested-family}
Let $\CD=\{\rD_i\}_{i\in I}$ be a family of closed disks in $(\BA_K^1)^{\an}$. We call $\CD$ a \emph{nested} family of disks if $\rD_i\subseteq\rD_j$ whenever the radius of $\rD_i$ is greater than the radius of $\rD_j$.
\end{Def}

If $\CD$ is a nested family of disks, we may define a pseudovaluation on $K[x]$ associated to $\CD$ by
\begin{equation*}
v_\CD(f)=\sup_{\rD\in\CD}v_\rD(f),\qquad f\in K[x].
\end{equation*}
It is shown in \cite[Section 1.4]{berkovich} that each pseudovaluation on $K[x]$ is of the form $v_\CD$ for some family of nested disks $\CD$ and that for each family $\CD$, one of the following is true:
\begin{itemize}
\item There exist disks in $\CD$ of arbitrarily large radius. In this case, $\bigcap_i\rD_i$ consists of a single closed point, and $v_\CD$ is the pseudovaluation associated to that point.
\item $\bigcap_i\rD_i$ is equal to one of the $\rD_i\in\CD$. In this case, $v_\CD$ is the Gauss valuation asociated to this disk. If the radius of $\rD_i$ is a rational number, then $v_\CD$ is a Type II valuation. Otherwise, $v_\CD$ corresponds to a Type III point.
\item $\bigcap_i\rD_i=\emptyset$. In this case, $v_\CD$ corresponds to a Type IV point.
\end{itemize}

In the case that $K$ is discretely valued, not all points in $(\BA_K^1)^{\an}$ correspond to families of nested disks. However, there is a similar classification involving so-called \emph{discoids}. We will explain this in \mbox{Section \ref{sec-applications-discoids}}.

Disks are also helpful for describing paths in $(\BP_K^1)^{\an}$. Recall that we have asserted above that $(\BP_K^1)^{\an}$ is a simply connected special quasipolyhedron. In particular, for any pair of distinct points $\xi_1,\xi_2\in(\BP_K^1)^{\an}$, there exists a unique compact interval $[\xi_1,\xi_2]$. Let us illustrate this for the case that $\xi_1$ and $\xi_2$ are a pair of closed points, say corresponding to $a,b\in K$. Similar considerations apply for other types of points.

We write $r\coloneqq v_K(a-b)$. Then the path connecting $\xi_1$ and $\xi_2$ is given by 
\begin{equation}\label{equ-unique-path}
\{\xi\mid v_\xi=v_{a,r'}\textrm{ for some }r'\ge r\}\cup\{\xi\mid v_\xi=v_{b,r'}\textrm{ for some }r'\ge r\}.
\end{equation}
Intuitively, we consider increasingly larger disks centered at $a$ up to the point where these disks contain $b$, and similarly we consider disks centered at $b$ up to the point where they contain $a$. Then the path connecting $\xi_1$ and $\xi_2$ consists of the boundary points of all these disks.

\begin{figure}[!htb]
\centering
\begin{minipage}{.5\textwidth}
\centering
\includegraphics[scale=0.32]{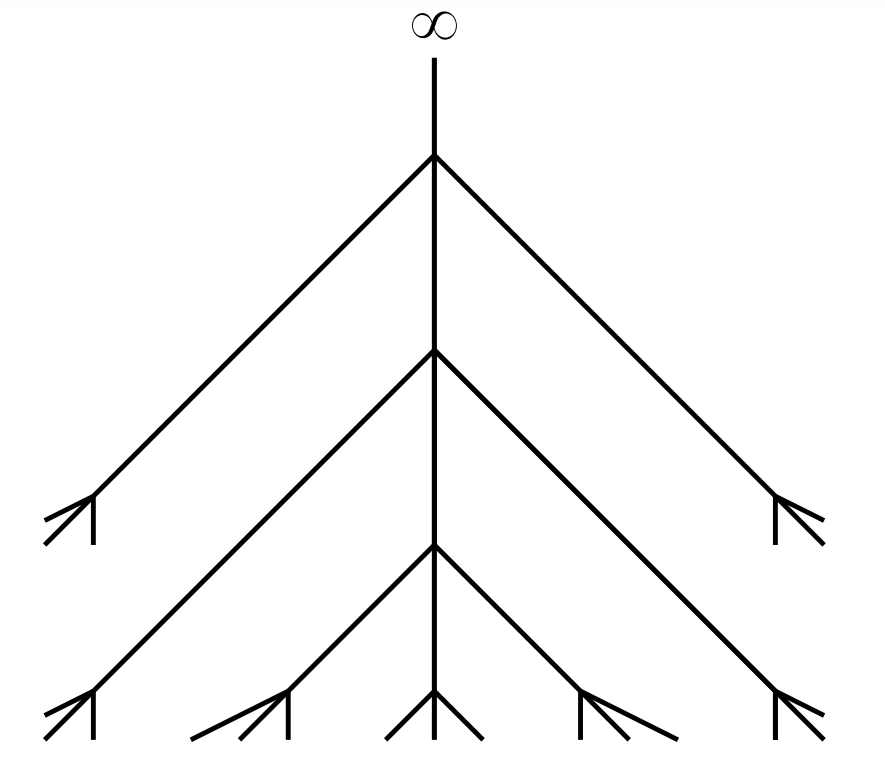}
\caption{A sketch of $(\BP_K^1)^{\an}$}
\label{fig-naked-hairy-tree}
\end{minipage}%
\begin{minipage}{.5\textwidth}\centering
\includegraphics[scale=0.32]{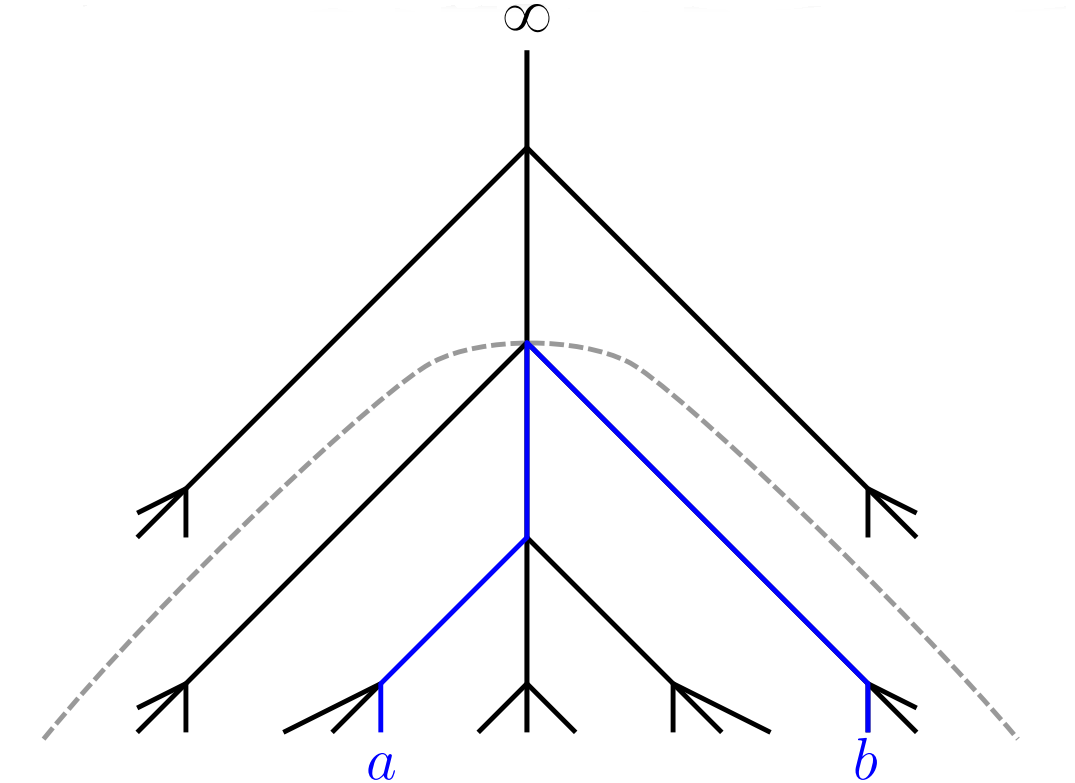}
\captionof{figure}{The same image with the path connecting two closed points highlighted}
\label{fig-hairy-tree-with-path}
\end{minipage}
\end{figure}

The ``correct'' way to visualize $(\BP_K^1)^{\an}$ is shown in the famous drawing \cite[Figure 1]{baker} by Baker (based in turn on an illustration of Silverman). A similar drawing is in \mbox{Figure \ref{fig-naked-hairy-tree}}. The curve $(\BP_K^1)^{\an}$ is an infinite tree, with branching occurring at all Type II points. In fact, let $\rD$ be a closed disk of radius $r$ with boundary point $\xi$. Then the branches at $\xi$ correspond to the open disks contained in $\rD$ of radius $r$, plus one branch pointing towards the point $\infty$. Of course, only a small part of the branching behavior is indicated in \mbox{Figure \ref{fig-naked-hairy-tree}}.

The Type I points are the endpoints at the bottom of \mbox{Figure \ref{fig-naked-hairy-tree}}, except for the point $\infty$, which is at the very top. To simplify the picture, we have left out the Type IV points. The intervals between two given Type I or Type II points consist of further points of Type II or of Type III. The Type II points are boundary points of disks of rational radius, while the Type III points are boundary points of disks of irrational radius.

In \mbox{Figure \ref{fig-hairy-tree-with-path}}, the path \eqref{equ-unique-path} connecting two Type I points is shown. The smallest closed disk containing both points is indicated by the dashed gray line.

\section{Topological ramification}\label{sec-preliminaries-topological-ramification}

In this section, $K$ is algebraically closed. Let $\phi\colon Y\to X$ be a finite separable morphism of $K$-curves and let $\xi\in X^{\an}$ and $\eta\in Y^{\an}$ be points with $\phi(\eta)=\xi$.

\begin{Def}\label{def-topological-ramification-index}
The \emph{topological ramification index} (called \emph{geometric ramification index} in \cite{berkovich-etale}; see \cite[Remark 3.2.2]{ctt} for a brief discussion of the terminology) of $\phi$ at $\eta$ is the degree
\begin{equation*}
n_\eta\coloneqq[\CH(\eta):\CH(\xi)].
\end{equation*}
\end{Def}

\begin{Rem}\label{rem-geometric-ramification-index-type-ii}
If $\xi$ and $\eta$ are of Type II, then we have $n_\eta=[\kappa(\eta):\kappa(\xi)]$. Indeed, by \cite[Theorem 6.3.1(iii)]{temkin}, $\CH(\xi)$ is \emph{stable}, meaning that
\begin{equation*}
n_\eta=[\kappa(\eta):\kappa(\xi)]\cdot[\Gamma_{v_\eta}:\Gamma_{v_\xi}].
\end{equation*}
But by \mbox{Remark \ref{rem-divisible-groups}}, both $\Gamma_{v_\eta}$ and $\Gamma_{v_\xi}$ equal the subgroup $\Gamma_K$. Thus $[\Gamma_{v_\eta}:\Gamma_{v_\xi}]=1$.
\end{Rem}

\begin{Rem}\label{rem-local-degree-formula}
The topological ramification indices are closely related to the notion of \emph{degree}. Given an open subspace $\rU\subseteq X^{\an}$ and a connected component $\rV$ of $\phi^{-1}(\rU)$, the induced morphism $\pi\colon\rV\to\rU$ is also finite. As in the case of schemes, the degree $\deg(\pi)$ of $\pi$ is defined to be the rank of the locally free $\CO_\rU$-module $\pi_*\CO_\rV$ (\cite[Remarks 6.3.1(iii)]{berkovich-etale}). For any point $\xi\in\rU$ we have the equality
\begin{equation*}
\deg(\pi)=\sum_{\eta\in\pi^{-1}(\xi)}n_\eta.
\end{equation*}
In particular, we have
\begin{equation*}
\deg(\phi)=\sum_{\eta\in\phi^{-1}(\xi)}n_\eta.
\end{equation*}
\end{Rem}

\begin{Def}
If $n_\eta>1$, then $\eta$ is called a \emph{topological ramification point} of $\phi$ and $\xi=\phi(\eta)$ is called a \emph{topological branch point}.

If furthermore $n_\eta$ is divisible by $p$, then $\eta$ is called a \emph{wild topological ramification point} of $\phi$ and $\xi$ is called a \emph{wild topological branch point}. Otherwise $\eta$ is called a \emph{tame topological ramification point} and $\xi$ a \emph{tame topological branch point}.
\end{Def}

\begin{Rem}\label{rem-topological-branch-locus}
\begin{enumerate}[(a)]
\item \sloppy If $\phi$ has degree $p$, then it follows immediately from \mbox{Remark \ref{rem-local-degree-formula}} that a wild topological branch point of $\phi$ has a unique preimage.
\item A Type II point $\eta\in Y^{\an}$ may be a wild topological ramification point even if the extension $\CH(\eta)/\CH(\xi)$ (where $\xi=\phi(\eta)$) is unramified. This choice of definition is briefly justified in \cite[Remark 3.2.4]{ctt}. We will say more on this in \mbox{Lemma \ref{lem-bijection-over-topological-branch-locus}}.
\item The function $\eta\mapsto n_\eta$ is upper semicontinuous (\cite[Lemma 3.6.10]{ctt}). In particular, if $\phi$ is a morphism of degree $p$, then the topological ramification locus and the wild topological ramification locus associated to $\phi$ are closed. Likewise, the (wild) topological branch locus, which is the image of the (wild) topological ramification locus under the finite morphism $\phi$, is closed.
\end{enumerate}
\end{Rem}

It is convenient to extend the definition of the topological ramification index to the set of branches at a given point. As explained in \cite[Section 3.6.11]{ctt}, if $v$ is a branch at a point $\eta$ represented by an interval $[\eta,\eta']$, then the topological ramification index is constant on each subinterval $(\eta,\eta'']\subset[\eta,\eta']$ with $\eta''$ chosen close enough to $\eta$. We define the \emph{topological ramification index of $v$} to be this constant value and denote it by $n_v$.

\section{Skeletons}\label{sec-preliminaries-skeletons}

In this section, $K$ is algebraically closed. Let $X$ be a $K$-curve with analytification $X^{\an}$. We begin by reviewing the central definitions, following \cite[Section 3.5.1]{ctt}.

\begin{Def}\label{def-abstract-graph}
Let $G=(V, E, s, t)$ be a finite multigraph (i.\,e.\ $V$ and $E$ are finite sets, and $s,t\colon E\to V$ are any pair of functions).

The \emph{topological realization} of $G$ is the topological space obtained as follows: First introduce a point $p_v$ for each vertex $v\in V$ and a compact interval $I_e$ for each edge $e\in E$. Then, for each edge $e$, identify the two endpoints of $I_e$ with $p_{s(e)}$ and $p_{t(e)}$ respectively.

We identify each vertex $v\in V$ with the point $p_v$ in the topological realization of $G$ and each edge $e\in E$ with the image of the interval $I_e$.
\end{Def}

\begin{Def}
Consider a pair $(\Gamma,S)$, where $\Gamma\subseteq X^{\an}$ is a closed subset and where $S\subseteq\Gamma$ is a finite subset consisting of Type I and Type II points, with at least one Type II point contained in $S$.

We say that $(\Gamma, S)$ is a \emph{graph} if there exist a finite multigraph $G=(V,E, s, t)$ and a homeomorphism $\phi$ between $\Gamma$ and the topological realization of $G$ with $\phi(S)=V$.

We then also call the elements of $S$ the \emph{vertices} of $(\Gamma, S)$ and call $\phi^{-1}(e)$, $e\in E$, the \emph{edges} of $(\Gamma,S)$.
\end{Def}

\begin{Rem}\label{rem-models-and-skeletons}
The condition that $S$ contain at least one Type II point is omitted in \cite{ctt}. Including it has the advantage that we may associate to every graph $(\Gamma, S)$ a model of $X$, namely the model corresponding by \mbox{Proposition \ref{prop-julian-theorem}} to the set of valuations
\begin{equation*}
\{v_\xi\mid\textrm{$\xi\in S$ is of Type II}\}.
\end{equation*}
If $S$ contained no Type II point, the graph would not correspond to any model, this set being empty.
\end{Rem}

\begin{Def}\label{def-skeletons}
A graph $(\Gamma, S)$ is called a \emph{skeleton} of $X^{\an}$ if the complement $X^{\an}\setminus\Gamma$ is a disjoint union of open disks and the set $S$ contains all points in $X^{\an}$ of positive genus.
\end{Def}

In the sequel, we will simply denote a graph or skeleton $(\Gamma,S)$ by $\Gamma$, and write $S(\Gamma)$ for the vertex set $S$.

\begin{Rem}\label{rem-skeletons-exist}
A skeleton $\Gamma\subset X^{\an}$ always exists. This fact is closely related to the \hyperref[thm-semistable-reduction]{Semistable Reduction Theorem}. We will see below in \mbox{Proposition \ref{prop-skeletons-and-reduction}} that the model determined by a skeleton (as explained in \mbox{Remark \ref{rem-models-and-skeletons}}) is semistable. Conversely, it can be shown that for a semistable model $\CX$ of $X$, the incidence graph of the special fiber $\CX_s$ is contained in $X^{\an}$ as a skeleton (\cite[Section 4.3]{berkovich}, \cite[Corollary 4.7]{bpr}).
\end{Rem}

\begin{Def}\label{def-leaf}
Let $\Gamma$ be a skeleton of $X^{\an}$. A point $\xi\in S(\Gamma)$ is called a \emph{leaf} of $\Gamma$ if $\xi$ maps to a leaf of $G=(V,E, s, t)$, where $G$ is a finite multigraph to whose topological realization $\Gamma$ is homeomorphic.
\end{Def}

\begin{Rem}\label{rem-unique-path}
Let $\Gamma\subset X^{\an}$ be a skeleton and let $\xi\in X^{\an}$ be any point. Since $X^{\an}\setminus\Gamma$ is a disjoint union of open disks, it follows that there is a unique interval $[\xi,\xi']$ in $X^{\an}$ with $[\xi,\xi']\cap\Gamma=\{\xi'\}$. (If $\xi\in\Gamma$, then $\xi'=\xi$ and $[\xi,\xi']=\{\xi\}$.) We will frequently make use of this property of skeletons.
\end{Rem}

\begin{Def}\label{def-retraction}
Let $\Gamma\subset X^{\an}$ be a skeleton. The \emph{canonical retraction map} associated to $\Gamma$ is the map
\begin{equation*}
\ret_\Gamma\colon X^{\an}\to\Gamma
\end{equation*}
that sends each point $\xi\in\Gamma$ to itself and sends each point $\xi\in X^{\an}\setminus\Gamma$ to the unique boundary point of the connected component of $X^{\an}\setminus\Gamma$ containing $\xi$. Alternatively, $\ret_\Gamma(\xi)$ equals $\xi'$, where $[\xi,\xi']$ is the unique interval in $X^{\an}$ with $[\xi,\xi']\cap\Gamma=\{\xi'\}$.
\end{Def}

\begin{Prop}\label{prop-tree-spanned-by-points}
Let $S\subset (\BP_K^1)^{\an}$ be a finite set of Type I and Type II points, satisfying $\abs{S}\ge3$ or containing at least one Type II point. Then there exists a skeleton of $(\BP_K^1)^{\an}$ whose underlying set is the union of all intervals $[\xi_1,\xi_2]$, where $\xi_1,\xi_2\in S$.
\end{Prop}
\begin{proof}
Let us first consider the case that $S$ contains at least one Type II point. If $\xi\in (\BP_K^1)^{\an}$ is a Type II point, then the subset $\{\xi\}$ is a skeleton of $\BP_K^1$. By induction, it thus suffices to show the following claim: If $\xi$ is a Type I or Type II point and $\Gamma'$ is a skeleton of $(\BP_K^1)^{\an}$ whose underlying set is the union of all intervals 
\begin{equation*}
[\xi_1,\xi_2],\qquad\xi_1,\xi_2\in S',
\end{equation*}
for a finite set $S'$ containing points of Type I and Type II, then one obtains again a skeleton $\Gamma$ from $\Gamma'$ by adding all intervals of the form $[\xi,\xi']$ for $\xi'\in S'$.

This follows easily from the fact that $(\BP_K^1)^{\an}$ is uniquely path-connected. Indeed, there is a unique path $[\xi,\xi']$ such that $\xi'\in\Gamma'$ and $\xi''\not\in\Gamma'$ for any $\xi''\in(\xi,\xi')$. We may construct $\Gamma$ by adding $[\xi,\xi']$ to the underlying set of $\Gamma'$ and adding $\xi$ and $\xi'$ to the vertex set of $\Gamma'$ (if the latter is not already contained in it). On the level of the underlying finite multigraph of $\Gamma'$, we simply attach a new leaf to $\Gamma'$, after potentially subdividing one of its edges.

The case where $S$ does not necessarily contain a Type II point, but contains at least three Type I points, works essentially the same. The point is that three distinct points $\xi_1,\xi_2,\xi_3$ in $(\BP_K^1)^{\an}$ determine a unique Type II point $\xi\in(\BP_K^1)^{\an}$ such that $\xi_1,\xi_2,\xi_3$ lie in distinct connected components of $(\BP_K^1)^{\an}\setminus\{\xi\}$. The skeleton $\Gamma$ is obtained from $\{\xi\}$ by the same process as the one already discussed above.
\end{proof}

\begin{figure}[htb]\centering\includegraphics[scale=0.44]{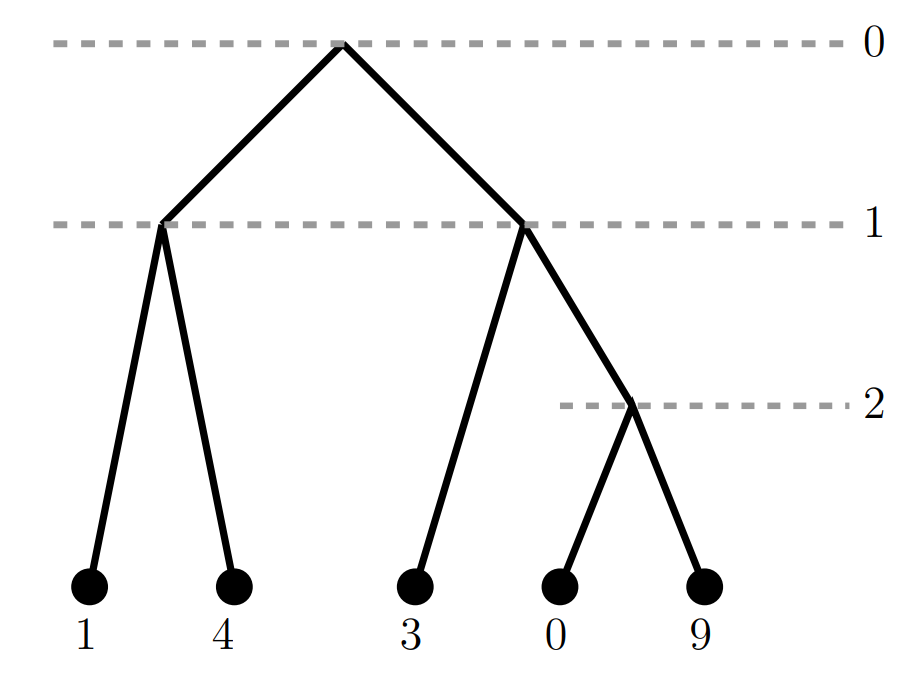}
\caption{The tree spanned by the set of Type I points discussed in \mbox{Example \ref{ex-basic-tree-example}}\label{fig-simple-tree}}
\end{figure}  

The skeleton $\Gamma$ described in \mbox{Proposition \ref{prop-tree-spanned-by-points}} is called the \emph{tree spanned by $S$}. We may visualize it by sketching the topological realization of the multigraph to which it is homeomorphic. We will do this in a similar style to the images in \mbox{Section \ref{sec-preliminaries-line}}. Thus we will put the Type I leaves different from $\infty\in(\BP_K^1)^{\an}$ at the bottom. Then the edges beginning at such a leaf $x_0$ consist of the boundary points of disks of decreasing radii centered at $x_0$. The point $\infty$ will be drawn at the top of the graph, if it is contained in $S$. 

\begin{figure}[htb]\centering\includegraphics[scale=0.35]{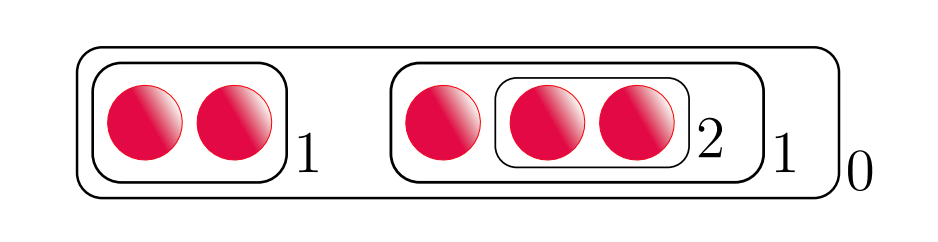}
\caption{The cluster picture visualization of the set of Type I points discussed in \mbox{Example \ref{ex-basic-tree-example}}}\label{fig-cluster-picture}
\end{figure}

\begin{Ex}\label{ex-basic-tree-example}
We consider $K=\BC_3$ and take for $S$ the set of Type I points corresponding to $0,1,3,4,9\in\BC_3$.

The tree $\Gamma$ spanned by $S$ is sketched in \mbox{Figure \ref{fig-simple-tree}}. On top of the points in $S$, the vertex set of $\Gamma$ contains three Type II points, which are the boundary points of the disk of radius $1$ centered at $1$ (or equivalently $4$), the disk of radius $1$ centered at $0$ (or $3$, or $9$), and the disk of radius $2$ centered at $0$ (or $9$). This is indicated by the dashed gray lines.

The tree spanned by $S$ can also be visualized using a \emph{cluster picture}. The connection between both viewpoints was already described and put to use by Bosch for computing semistable reductions in \cite{bosch}. Recently, cluster pictures have been used extensively in \cite{ddmm} and related works.

A \emph{cluster} is a subset of $S$ that is obtained by intersecting $S$ with a closed disk. In \mbox{Figure \ref{fig-cluster-picture}}, the five red balls correspond to the elements of $S$. The proper clusters (containing at least two points each) are drawn as rectangles, with the index indicating the radius of the smallest disk cutting out the cluster.
\end{Ex}

To end this section, we define a metric structure on $X^{\an}$, following \cite[Section 5]{bpr}. Consider first the standard open anulus of radius $r\in\BQ_{>0}$,
\begin{equation}\label{equ-open-anulus}
\rA(r)\coloneqq\rA(0,r,0)=\{\xi\in(\BP_K^1)^{\an}\mid r> v_\xi(x)>0\}.
\end{equation}
Denote by $\xi_0,\xi_r\in(\BP_K^1)^{\an}$ the two boundary points of $\rA(r)$. Then $\rA(r)$ contains the interval $(\xi_0,\xi_r)$ --- in the language of \cite{bpr}, this interval is called the \emph{skeleton} of $\rA(r)$. We define the length of this interval to be
\begin{equation*}
l\big((\xi_0,\xi_r)\big)\coloneqq r.
\end{equation*}
Next, suppose that $[\xi,\xi']\subset X^{\an}$ is any interval consisting of Type II and Type III points. By \mbox{Remark \ref{rem-skeletons-exist}}, it is possible to choose a skeleton $\Gamma$ containing $[\xi,\xi']$. Removing the vertices of $\Gamma$ yields a decomposition of $(\xi,\xi')$ into disjoint open intervals
\begin{equation*}
(\xi,\xi')\setminus S(\Gamma)=\coprod_{i=1}^k(\xi_i,\xi_{i+1}),\qquad \xi_1=\xi,\quad\xi_{k+1}=\xi'.
\end{equation*}
By \cite[Remark 3.5.2(ii)]{ctt}, there exist open anuli $\rA_i\subset X^{\an}$ such that under an identification of $\rA_i$ with a standard open anulus $\rA(r_i)$ as in \eqref{equ-open-anulus}, the interval $(\xi_i,\xi_{i+1})$ is identified with the skeleton of $\rA(r_i)$, an interval of length $r_i\in\BQ_{>0}$. Thus we define $l(\xi_i,\xi_{i+1})$ to be $r_i$, and
\begin{equation*}
l\big((\xi,\xi')\big)\coloneqq\sum_{i=1}^kl\big((\xi_i,\xi_{i+1})\big)=\sum_{i=1}^kr_i.
\end{equation*}
It can be shown (see \cite[Section 5]{bpr}) that this definition is independent of the choice of $\Gamma$. We extend the definition of the length to closed and half-open intervals consisting of Type II and Type III points by
\begin{equation*}
l\big([\xi,\xi']\big)\coloneqq l\big((\xi,\xi']\big)\coloneqq l\big([\xi,\xi')\big)\coloneqq l\big((\xi,\xi')\big).
\end{equation*}
Finally, we extend the length function to intervals $[\xi,\xi']$, $\xi\ne\xi'$, where one of the two endpoints is a Type I point, by setting $l([\xi,\xi'])\coloneqq\infty$.

The following definition coincides with the one given in \cite[\mbox{Section 3.6.2}]{ctt}.
\begin{Def}\label{def-radius-parametrization}
A homeomorphism $\phi\colon I\to J$ between an interval $I\subset X^{\an}$ and an interval $J\subseteq\BR\cup\{\pm\infty\}$ is called a \emph{radius parametrization} of $I$ if 
\begin{enumerate}[(a)]
\item the length $l(I')$ is the same as the length of $\phi(I')\subseteq\BR\cup\{\pm\infty\}$ (with respect to the usual metric on $\BR$) for every subinterval $I'\subseteq I$,
\item for some Type II point $\xi\in I$, we have $\phi(\xi)\in\BQ$.
\end{enumerate}
\end{Def}
Note that this implies that $\phi(\xi)\in\BQ$ for \emph{every} Type II point $\xi\in I$. When $I$ contains a point of Type I and at least one other point, then $J$ must have $\infty$ or $-\infty$ as boundary point. The simplest example of a radius parametrization is the identification of the skeleton of the standard open anulus $\rA(r)$ from \eqref{equ-open-anulus} with the interval $(0,r)$, which identifies $s\in(0,r)$ with the point corresponding to the Gauss valuation of radius $s$.

\begin{Def}\label{def-radius-function}
The \emph{radius function} associated to a skeleton $\Gamma$ of $X$ is the map
\begin{equation*}
r_\Gamma\colon X^{\an}\to[0,\infty]
\end{equation*}
defined by
\begin{equation*}
r_\Gamma(\xi)\coloneqq l([\xi,\xi']),
\end{equation*}
where $[\xi,\xi']$ is the unique interval with $[\xi,\xi']\cap\Gamma=\{\xi'\}$ (\mbox{Remark \ref{rem-unique-path}}).
\end{Def}
In particular, if $\xi\in\Gamma$, then $\xi'=\xi$ and so $r_\Gamma(\xi)=0$. In general, $r_\Gamma$ measures the distance of points to the skeleton $\Gamma$.

\section{Reduction maps}

The field $K$ is still assumed to be algebraically closed and $X$ denotes a $K$-curve. In this section, we discuss reduction maps associated to $\CO_K$-models \mbox{of $X$}.

Let $\CX$ be an $\CO_K$-model of $X$. Suppose that $x\colon\Spec K\to X$ is a closed point. By the valuative criterion for properness, $x$ extends to a unique morphism $\Spec\CO_K\to\CX$:
\begin{equation*}
\begin{tikzcd}
\Spec K \arrow[r] \arrow[d]                                  & \CX \arrow[d] \\
\Spec\CO_K \arrow[r, no head,equals] \arrow[ru, "\exists!", dashed] & \Spec\CO_K   
\end{tikzcd}
\end{equation*}
The image of the closed point of $\Spec\CO_K$ under this map is called the \emph{reduction} of $P$ with respect to $\CX$; it is denoted $\red_\CX(x)$. The resulting map
\begin{equation*}
\red_\CX\colon\textrm{closed points of $X$}\longrightarrow\textrm{closed points of $\CX_s$}
\end{equation*}
is surjective (\cite[Corollary 10.1.38]{liu}).

Recall from \mbox{Section \ref{sec-preliminaries-analytic}} that the underlying set of $X^{\an}$ consists of pairs $\xi=(P,v)$, where $P\in X$ is any point and where $v$ is a valuation on the residue field $k(P)$ extending $v_K$. Denote by $\CO_v$ the valuation ring of $v$. Then as before, the valuative criterion for properness gives an extension of the natural morphism $\Spec k(P)\to\CX$ to a morphism $\Spec\CO_v\to\CX$:
\begin{equation*}
\begin{tikzcd}
\Spec k(P) \arrow[r] \arrow[d]                      & \CX \arrow[d] \\
\Spec\CO_v \arrow[r] \arrow[ru, "\exists!", dashed] & \Spec\CO_K   
\end{tikzcd}
\end{equation*}
We denote the image of the closed point of $\CO_v$ under this morphism by $\red_\CX(\xi)$. Thus we have extended $\red_\CX$ to a map
\begin{equation*}
\red_\CX\colon X^{\an}\to\CX_s.
\end{equation*}
The image of a Type II point $\xi\in X^{\an}$ under $\red_\CX$ is the generic point of the component $Z\subseteq\CX_s$ giving rise to $v_\xi$ as explained in Section \ref{sec-preliminaries-valuations}.

The above reduction maps are associated to models of $X$; there is no canonical reduction map associated to $X$ alone. There is, however, a canonical reduction map associated to any affinoid subdomain of $X^{\an}$.

Recall from Section \ref{sec-preliminaries-line} that affinoid subdomains are the building blocks of analytic curves in the sense of Berkovich. They are defined to be the \emph{spectra} of certain commutative Banach algebras, the \emph{affinoid $K$-algebras}, of which we saw two examples in Section \ref{sec-preliminaries-line}. In fact, the affinoid subdomains whose reductions we consider below will be \emph{strictly affinoid spaces}, which are the spectra of \emph{strictly affinoid $K$-algebras} (\cite[Section 2.1]{berkovich}). Essentially, these are the affinoid spaces arising from affinoid spaces in the sense of rigid-analytic geometry (\cite[Chapter 3]{bosch-lectures}).

If $\CA$ is a strictly $K$-affinoid algebra, we define
\begin{equation*}
\tilde{\CA}\coloneqq\CA^\circ/\CA^{\circ\circ},\qquad\textrm{where}
\end{equation*}
\begin{equation*}
\CA^\circ=\{f\in\CA\mid\forall\xi\in\CM(\CA):v_\xi(f)\ge0\},
\end{equation*}
\begin{equation*}
\CA^{\circ\circ}=\{f\in\CA\mid\forall\xi\in\CM(\CA):v_\xi(f)>0\}
\end{equation*}
(cf.\ \cite[Section 2.4, Theorem 1.3.1]{berkovich}). The association $\CA\mapsto\tilde{\CA}$ is functorial. The \emph{canonical reduction} of the affinoid space $\CM(\CA)$ is defined to be the scheme $\Spec(\tilde{\CA})$. It is a reduced scheme of finite type over the residue field $\kappa$. There is a reduction map
\begin{equation*}
\red_{\CM(\CA)}\colon\CM(\CA)\to\Spec(\tilde{\CA})
\end{equation*}
defined as follows. A given point $\xi\in\CM(\CA)$ has an associated character $\CA\to\CH(\xi)$; applying the functor ${(\cdot)}^\sim$ yields a character $\tilde{\CA}\to\kappa(\xi)$, which in turn defines a point in $\Spec({\tilde{\CA}})$.

\begin{Ex}\label{ex-canonical-reduction}
\begin{enumerate}[(a)]
\item Consider the closed unit disk
\begin{equation*}
\rD\coloneqq\CM(\CA),\qquad\CA=K\{T\}.
\end{equation*}
Its canonical reduction is the affine line $\BA_\kappa^1$, and the reduction map sends a closed point $x\in\CO_K$ to its reduction $\overline{x}\in \kappa$.
\item Consider the closed anulus of radius $r\in\BQ_{>0}$,
\begin{equation*}
\rA[r]\coloneqq\rA[0,r,0]=\CM(\CA),\qquad\CA=K\{t,\pi^rt^{-1}\}.
\end{equation*}
The canonical reduction of $\rA[r]$ is the ``cross'' $\overline{X}\coloneqq\Spec\kappa[T,U]/(TU)$. We may identify the set of rational points
\begin{equation*}
\{x\in K\mid 0\le v_K(x)\le r\}
\end{equation*}
with a subset of $\rA[r]$. The reduction of such a point $x\in K$ is given by
\begin{itemize}
\item the reduction $\overline{x}\in \kappa$ of $x$, considered as a point in the affine line $\Spec \kappa[T]\subset\overline{X}$, if $v_K(x)=0$,
\item the nodal point of $\overline{X}$, if $0<v_K(x)<r$,
\item the reduction $\overline{x\pi^{-r}}\in \kappa$ (where $\pi^{-r}$ denotes an element of $K$ with valuation $-r$), considered as a point in the affine line $\Spec \kappa[U]\subset\overline{X}$, if $v_K(x)=r$.
\end{itemize}
\end{enumerate}
\end{Ex}

We now explain how to piece together canonical reduction maps associated to strictly affinoid subdomains $\CM(\CA)\subset X^{\an}$ to obtain reduction maps on all of $X^{\an}$, following \cite[Section 4.3]{berkovich}.

Suppose that $\{\rU_i\}_{i\in I}$ is a finite covering of $X^{\an}$ by strictly affinoid subdomains. Let us write $\rU_{ij}\coloneqq\rU_i\cap\rU_j$ and denote by $\overline{X}_i$ and $\overline{X}_{ij}$ the canonical reductions of $\rU_i$ and $\rU_{ij}$ respectively. The covering $\{\rU_i\}_{i\in I}$ is called a \emph{formal covering} of $X^{\an}$ if the induced morphism $\overline{X}_{ij}\to\overline{X}_i$ is an open embedding for all $i,j\in I$. In this case it is possible to glue the $\overline{X}_i$ along the $\overline{X}_{ij}$ to a $\kappa$-curve $\overline{X}$. Moreover, the canonical reductions $\red_{\rU_i}$ combine to a reduction map $X^{\an}\to\overline{X}$. In the following proposition, we will see how to define a formal covering associated to a skeleton of $X^{\an}$.

\begin{Prop}\label{prop-skeletons-and-reduction}
If $\Gamma\subset X^{\an}$ is a skeleton, then the model of $X$ determined by $\Gamma$ (see \mbox{Remark \ref{rem-models-and-skeletons}}) is semistable.
\end{Prop}
\begin{proof}
\sloppy A version of this folklore theorem was written down in \cite[\mbox{Theorem 4.11}]{bpr}. Below, we use a similar construction, but our version involves models in the sense of \mbox{Section \ref{sec-preliminaries-models}}, not formal models.

We assume that $\Gamma$ has at least $2$ edges (if this is not the case, the argument only needs to be slightly adjusted, as is done in \cite[\mbox{Section 4.15.2--3}]{bpr}). We then define a formal covering of $X^{\an}$ as follows. For each edge $e=[\xi_1,\xi_2]$ of $\Gamma$, put
\begin{equation*}
\rU_e\coloneqq \ret_\Gamma^{-1}(e).
\end{equation*}
The $\rU_e$ satisfy the following (\cite[Section 4.15.1]{bpr}):
\begin{itemize}
\item $\rU_e$ is a strictly affinoid subdomain of $X^{\an}$.
\item The canonical reduction of $\rU_e$ is a $\kappa$-curve $\overline{X}_e$ with one singularity, a node. If the start and end point of $e$ are not the same, then $\overline{X}_e$ has two irreducible components intersecting in this node. Otherwise, $\overline{X}_e$ is irreducible.
\item\sloppy The complement of $\{\xi_1,\xi_2\}$ in $\rU_e$ has as connected components an open anulus $\rA\coloneqq \ret_\Gamma^{-1}((\xi_1,\xi_2))$ and infinitely many open disks. The anulus $\rA$ is the inverse image of the nodal point on $\overline{X}_e$ while each open disk is the inverse image of a smooth point on $\overline{X}_e$.
\end{itemize}

Suppose now that $e$ and $f$ are two distinct adjacent edges, say with $e\cap f=\{\xi\}$. Then $\rU_e\cap\rU_f$ is the retraction $\ret_\Gamma^{-1}(\xi)$. It follows from the above bullet points that $\red_{\rU_e}(\rU_e\cap\rU_f)$ equals one of the components of $\overline{X}_e$ minus the nodal point. Thus the $\rU_e$ form a formal covering of $X^{\an}$. Moreover, it is clear that by gluing the canonical reductions $\overline{X}_e$, we obtain a semistable $\kappa$-curve $\overline{X}$.

By \cite[Proposition 2.4.4]{berkovich}, the function fields of the irreducible components of the reduction $\overline{X}$ are the same as the residue fields $\kappa(\xi)$, where $\xi\in\Gamma$ is a Type II point. On the last two pages of \cite{bosch-luetkebohmert} it is explained how to construct an honest $\CO_K$-model $\CX$ from the formal covering $\{\rU_e\}_e$, whose special fiber $\CX_s$ coincides with the semistable $\kappa$-curve $\overline{X}$. By Proposition \ref{prop-julian-theorem}, this $\CX$ must be the model determined by $\Gamma$ as explained in \mbox{Remark \ref{rem-models-and-skeletons}}, so we are done.
\end{proof}

\begin{Rem}\label{rem-reduction-interpretation}
\sloppy Let $\CX$ be the model associated to a skeleton $\Gamma$ of $X$. In \mbox{Proposition \ref{prop-skeletons-and-reduction}}, we have reconstructed the reduction map $\red_\CX\colon X^{\an}\to\CX_s$ working locally on $X^{\an}$. From the proof of \mbox{Proposition \ref{prop-skeletons-and-reduction}}, it is clear how we may describe the reduction $\red_\CX(x_0)$ of a closed point $x_0\in X$ in terms of the retraction $\ret_\Gamma(x_0)$. In summary, we have the following:
\begin{itemize}
\item For every Type II vertex $\xi$ of $\Gamma$, the special fiber $\CX_s$ has an irreducible component $\overline{X}_\xi$
\item Each edge $e$ connecting two Type II vertices $\xi_1,\xi_2$ of $\Gamma$ corresponds to a nodal point of intersection of the components $\overline{X}_{\xi_1},\overline{X}_{\xi_2}$ (here we also allow $\xi_1=\xi_2$, in which case the component $\overline{X}_{\xi_1}=\overline{X}_{\xi_2}$ has a singularity)
\item If $\ret_\Gamma(x_0)$ lies in the interior of this edge $e$, then $\red_\CX(x_0)$ is the point of intersection of $\overline{X}_{\xi_1}$ and $\overline{X}_{\xi_2}$ (or, if $\xi_1=\xi_2$, then $\red_\CX(x_0)$ is the singular point on $\overline{X}_{\xi_1}=\overline{X}_{\xi_2}$)
\item If $\ret_\Gamma(x_0)=\xi$ for a Type II vertex $\xi$ of $\Gamma$, then $\red_\CX(x_0)$ is a smooth point on $\overline{X}_\xi$
\item If $x_0$ lies on $\Gamma$, then it is a leaf of $\Gamma$, and $\red_\CX(x_0)$ is a smooth point on $\overline{X}_\xi$, where $\xi$ is the Type II vertex of $\Gamma$ adjacent to $x_0$
\item Suppose that $x_0,x_1\in X$ are two closed points with the same retraction $\xi$ to $\Gamma$. Then $\red_\CX(x_0)=\red_\CX(x_1)$ if and only if $[x_0,\xi]$ and $[x_1,\xi]$ represent the same branch at $\xi$ in the sense of \mbox{Definition \ref{def-branch}}
\end{itemize}
In light of the last point, it makes sense to talk of the \emph{reduction of a branch} at a vertex $\xi$ of $\Gamma$.
\end{Rem}

\begin{Kor}\label{cor-tree-spanned-by-points}
Let $\Gamma\subset X^{\an}$ be the tree spanned by a finite set of Type I points $S\subset X=\BP_K^1$. Let $\overline{X}$ denote the reduction of $X$ associated to $\Gamma$ as constructed in Proposition \ref{prop-skeletons-and-reduction}. Then the points $x_0\in S$ reduce to pairwise distinct smooth points on $\overline{X}$.
\end{Kor}
\begin{proof}
This is immediate from \mbox{Remark \ref{rem-reduction-interpretation}} and the construction of $\Gamma$ in \mbox{Proposition \ref{prop-tree-spanned-by-points}} --- all $x_0$ are leaves of $\Gamma$.
\end{proof}

\begin{Def}\label{def-intersection-graph}
The \emph{intersection graph} of a semistable curve $C$ is the multigraph which has one vertex for each irreducible component of $C$ and has one edge for each nodal point $x$ on $C$, the endpoints of this edge corresponding to the components on which $x$ lies.
\end{Def}

It follows from \mbox{Remark \ref{rem-reduction-interpretation}} that a skeleton $\Gamma$ of $X$ may be recovered as the topological realization of the intersection graph of the semistable curve $\overline{X}$, where $\overline{X}$ is the reduction associated to the formal covering induced by $\Gamma$.

\section{Admissible reduction}\label{sec-preliminaries-admissible}

We continue to assume that $K$ is algebraically closed, but this assumption is not essential for this section. In fact, the main goal of works like \cite{ddmm} and \cite{bouw-wewers} cited below is computing arithmetic invariants of curves, and therefore a lot of emphasis is put on computing the minimal field extension over which a semistable reduction exists. We will consider the question of computing this field extension in \mbox{Chapter \ref{cha-applications}}; for now we stick to algebraically closed $K$ for expository simplicity.

The goal of this section is to review a technique that is useful for computing the semistable reduction of many curves. The idea is to choose a finite morphism $\phi$ from the curve $Y$ under consideration to $\BP^1_K$ and then take as a model of $Y$ the normalization $\CY$ in $F_Y$ of a model $\CX$ of $\BP^1_K$ \emph{separating the branch points} of $\phi$. Our reason for discussing this technique here is twofold: First, it serves as an illustration and example for many notions introduced earlier in this chapter. And second, we will see in \mbox{Example \ref{ex-admissible-no-good}} that the method of separating branch points does not in general work without \mbox{Assumption \ref{ass-admissible}} below, thus motivating our developing a more general method in subsequent chapters. If \mbox{Assumption \ref{ass-admissible}} is satisfied, the induced map on special fibers $\CY_s\to\CX_s$ is an \emph{admissible cover} between semistable curves. (See \cite[Section 4]{harris-mumford} --- we will not use this terminology past this section.)

This approach has been used in many works treating the semistable reduction of curves. For example, \cite{bouw-wewers} deals with the case of superelliptic curves. Separating branch points is also the basis of \cite{ddmm} and related works, which treat the case of hyperelliptic curves in great detail.

Let $\phi\colon Y\to X=\BP_K^1$ be a cover of smooth and irreducible curves of degree $n\ge2$. In this section only we make the following assumption. In subsequent chapters, we will usually consider the situation that $\phi$ is of degree exactly $p$.

\begin{Ass}\label{ass-admissible}
The residue characteristic $p=\characteristic(\kappa)$ is strictly greater than $n$.
\end{Ass}

Let $\Gamma\subset X^{\an}$ be the tree spanned by the branch points of $\phi$, as constructed in \mbox{Proposition \ref{prop-tree-spanned-by-points}}. Let $\CX$ be the $\CO_K$-model associated to $\Gamma$, as explained in \mbox{Remark \ref{rem-models-and-skeletons}}. Then by \mbox{Corollary \ref{cor-tree-spanned-by-points}} the branch points of $\phi$ specialize to pairwise distinct points on $\CX_s$. Finally, let $\CY$ denote the normalization of $\CX$ in the function field $F_Y$; it is an $\CO_K$-model of $Y$.

\begin{Prop}\label{prop-paul-theorem}
The model $\CY$ is semistable.
\end{Prop}
\begin{proof}
\sloppy In the case that $\phi$ is a Galois cover, this follows from \cite[\mbox{Theorem 2.3}]{liu-lorenzini}. The general case is proved in detail in \cite[Section 3.3]{helminck}. In either case, the needed tameness assumptions follow from \mbox{Assumption \ref{ass-admissible}}.
\end{proof}

Since the most important concrete class of algebraic curves that we treat is that of plane quartic curves, many of our examples are drawn from this class, including the following \mbox{Example \ref{ex-admissible-no-good}}. Though this is not necessary for the continuity of the argument, the reader might wish to consult \mbox{Section \ref{sec-quartic-normal-form}} for a brief reminder on plane quartics and for an explanation of the normal form in which all our plane quartics are given.

\begin{Ex}\label{ex-admissible-no-good}
Consider the smooth plane quartic curve $Y$ over $K\coloneqq\BC_3$ given by
\begin{equation}\label{equ-my-favorite-example}
Y\colon\quad F(x,y)= y^3+x^3y+x^4+1=0.
\end{equation}
The discriminant of $F$ (considered as a polynomial in $y$ whose coefficients are polynomials in $x$) is 
\begin{equation*}
\Delta_F=-4x^9-27x^8-54x^4-27.
\end{equation*}
The branch locus of the degree-$3$ cover $\phi\colon Y\to\BP_K^1$ given by $(x,y)\mapsto x$ consists of the zeros of $\Delta_F$ and the point $\infty$ (cf.\ \mbox{Lemma \ref{lem-quartic-branch-locus}} below)

In fact, the branch locus of $\phi$ is \emph{equidistant} (a terminology found for example in \cite[Section 2]{lehr-matignon} or \cite[Section 4.3]{icerm}). This means that there exists a smooth $\CO_K$-model $\CX$ of $X$ such that the branch points of $\phi$ specialize to distinct points on the special fiber $\CX_s$. The set of Type II valuations associated to $\CX$ via \mbox{Proposition \ref{prop-julian-theorem}} contains only one element, which we denote $\xi$. The tree spanned by the branch points of $\phi$ is sketched in \mbox{Figure \ref{fig-equidistant-branch-locus}}; the branch point $\infty$ is at the top, the nine finite branch points at the bottom. 

\begin{figure}[htb]\centering\includegraphics[scale=0.36]{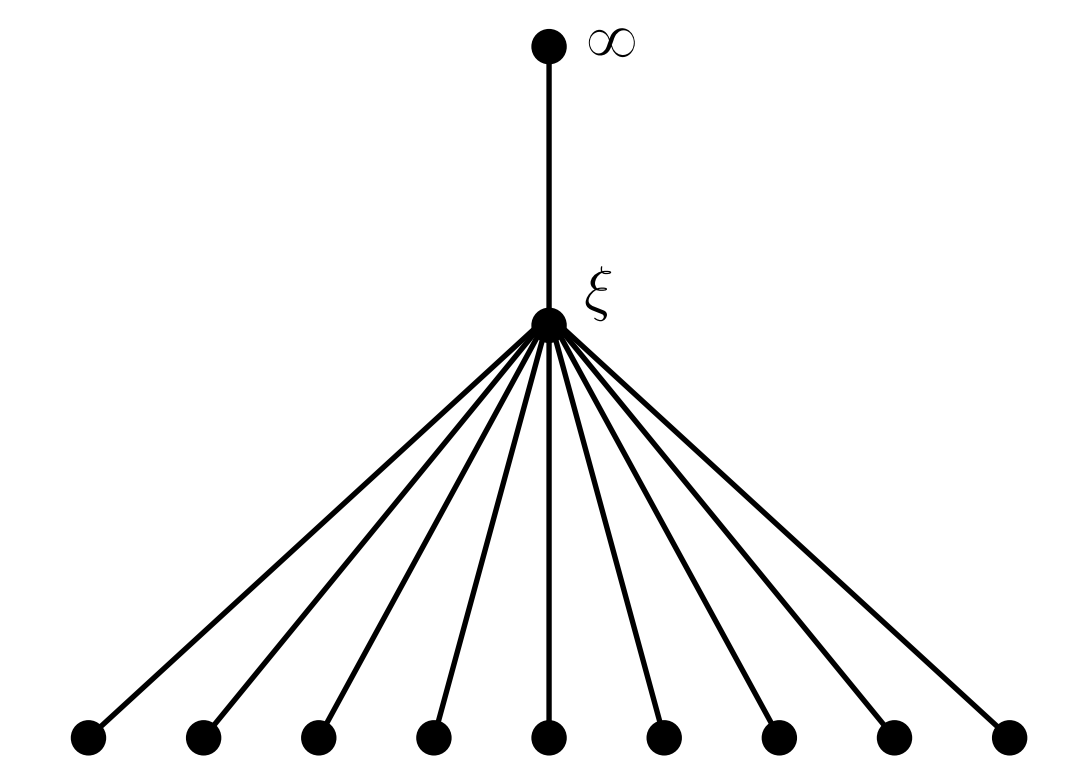}
\caption{The tree spanned by the branch locus of the cover $\phi$ discussed in \mbox{Example \ref{ex-admissible-no-good}}}\label{fig-equidistant-branch-locus}
\end{figure}  

In fact, the point $\xi$ is the boundary point of the disk of radius $1/2$ centered at one of the finite branch points, that is, at one of the zeros of the discriminant of the defining polynomial in \eqref{equ-my-favorite-example}. 

However, the normalization $\CY$ of the model $\CX$ determined by the point $\xi$ does not have semistable reduction --- thus, the conclusion of \mbox{Proposition \ref{prop-paul-theorem}} is false for this example. To show this, we claim that the valuation $v_\xi$ has a unique extension to $F_Y$ and that the associated extension of residue fields is inseparable. Then it follows from \mbox{Lemma \ref{lem-ramification-theory}} that $\CY_s$ is an irreducible curve of geometric genus $0$. In fact, the morphism $\CY_s\to\CX_s$ is a homeomorphism (\cite[Exercise 5.3.9]{liu}), so $\CY_s$ cannot have nodal singularities. Thus $\CY_s$ is not a semistable curve.

It remains to prove the claim that $v_\xi$ has a unique extension to $F_Y$ with inseparable residue field extension. Let $(x_0,y_0)$ be one of the finite ramification points of $\phi$. Since $F(x_0,y_0)=0$, the minimal polynomial of the element $z\coloneqq (y-y_0)\pi^{1/2}\in F_Y$ over the subfield $K(t)$, where $t=x-x_0$, has the form
\begin{equation}\label{equ-non-reduced-minpoly}
T^3+a_0T^2+(b_3t^3+b_2t^2+b_1t+b_0)T+c_4t^4+c_3t^3+c_2t^2+c_1t.
\end{equation}
Here the coefficients $a_0,b_0,\ldots,b_3,c_1,\ldots,c_4$ lie in a finite extension $K_0$ of $\BQ$ --- all that is needed is that the point $(x_0,y_0)$ be $K_0$-rational and that $K_0$ contain the element $\pi^{1/2}$ of valuation $1/2$. It is therefore not hard to determine the valuations of these coefficients, for example using the functionality for extending valuations along extensions of number fields in Sage --- see \mbox{Code Listing \ref{list-admissible}} in the appendix for details. The resulting valuations are
\begin{equation*}
v_K(a_0)=\frac{1}{2},\quad v_K(b_0)=1,\ v_K(b_1)=\frac{2}{3},\ v_K(b_2)=\frac{1}{3},\ v_K(b_3)=-1,
\end{equation*}
\begin{equation*}
v_K(c_1)=-\frac{1}{2},\ v_K(c_2)=-\frac{1}{6},\ v_K(c_3)=-\frac{3}{2},\ v_K(c_4)=-\frac{3}{2}.
\end{equation*}
(Note that these valuations are independent of the choice of ramification point $(x_0,y_0)$, since the nine finite branch points of $\phi$ are Galois conjugate over $\BQ_3$.) Recall that the valuation $v_\xi$ is determined by
\begin{equation*}
v_\xi\big(\sum_ia_it^i\big)=\min_i\big\{v_K(a_i)+\frac{i}{2}\big\}.
\end{equation*}
We write $\overline{t}\coloneqq[t]_{v_\xi}\in\kappa(\xi)$ and let $\eta\in Y^{\an}$ be a point with $\phi(\eta)=\xi$. Reducing the coefficients of \eqref{equ-non-reduced-minpoly} to $\kappa(v_\xi)$ shows that $\overline{z}\coloneqq[z]_{v_\eta}$ satisfies an equation of the form
\begin{equation*}
\overline{z}^3+\overline{c}_3\overline{t}^3+\overline{c}_1\overline{t}=0,
\end{equation*}
where $\overline{c}_3,\overline{c}_1\in\kappa^\times$. Clearly $\overline{c}_3\overline{t}^3+\overline{c}_1\overline{t}$ is not a third power, so $\kappa(\eta)$ must be an inseparable residue field extension of $\kappa(v_\xi)$ of \mbox{degree $3$}. The claim is proved.

We will see in \mbox{Example \ref{ex-worked-example-5}} that $Y$ has good reduction.
\end{Ex}

\chapter{Admissible functions}

Throughout this entire chapter, the ground field $K$ is algebraically closed.

We will introduce \emph{admissible functions} on an analytic curve $X^{\an}$. These are locally constant outside a skeleton of $X^{\an}$. Thus despite the complicated structure of $X^{\an}$ we can work with and understand the structure of admissible functions. In \mbox{Section \ref{sec-greek-valuative}}, we discuss an important class of admissible functions: \emph{valuative functions} attached to elements of the function field $F_X$. For the case that $X$ is the projective line, an implementation of admissible functions exists in the Sage package MCLF. This is a crucial ingredient for the implementation of our algorithms in \mbox{Chapter \ref{cha-applications}}.

Suppose that we are given a degree-$p$ morphism $\phi\colon Y\to X$. In \mbox{Sections \ref{sec-greek-mu}}, \mbox{\ref{sec-greek-delta}, and \ref{sec-greek-lambda}}, we introduce three important functions on $X^{\an}$, denoted by Greek letters. The first is the admissible function $\mu$, which measures the distance to a given skeleton of $X^{\an}$. The second is the function $\delta$. It is essentially equivalent to the \emph{different function} of \cite{ctt} and measures the topological ramification (\mbox{Section \ref{sec-preliminaries-topological-ramification}}) of $\phi\colon Y^{\an}\to X^{\an}$. While $\delta$ is not admissible, the function $\lambda$ we derive from it is. We will see in \mbox{Theorem \ref{thm-trivializing-lambda}} that $\lambda$ controls the semistable reduction of $Y^{\an}$ to a large degree and exploit this fact in \mbox{Chapter \ref{cha-analytic}}.

\label{cha-greek}

\section{Definition}\label{sec-greek-definition}

Let $I\subseteq\BR\cup\{\pm\infty\}$ be an interval. A function $f\colon I\to\BR$ is called an \emph{affine function} if it is of the form
\begin{equation*}
t\mapsto \alpha t+\beta,\qquad \alpha,\beta\in\BQ.
\end{equation*}
It is called \emph{piecewise affine} if it is continuous and for all but finitely many points $t_0\in I$ there exists an open interval $t_0\in J\subseteq I$ such that $\rest{f}{J}$ is affine. The finitely many points $t_0$ where no such interval can be found are called the \emph{kinks} of $f$.

\sloppy Let $X$ be a $K$-curve with analytification $X^{\an}$. By choosing a radius parametrization of a compact interval $I\subset X^{\an}$, say
\begin{equation*}
\phi\colon I\overset{\sim}{\longrightarrow}[a,b]\subseteq\BR\cup\{\pm\infty\},
\end{equation*}
it makes sense to talk of \emph{(piecewise) affine functions on $I$}; they are simply functions of the form $h=f\circ\phi\colon I\to\BR$, where $f$ is (piecewise) affine in the previous sense. Similarly, the \emph{kinks} of a piecewise affine function $h=f\circ\phi\colon I\to\BR\cup\{\pm\infty\}$ are given by $\{\phi^{-1}(t_0)\mid\textrm{$t_0$ is a kink of $f$}\}$. Clearly all these notions are independent of the radius parametrization chosen.

Let $U\subseteq X^{\an}$ be a subset. A function $h\colon U\to\BR\cup\{\pm\infty\}$ is called \emph{(piecewise) affine} if its restriction to every compact interval $I\subset U$ is (piecewise) affine.

\begin{Def}\label{def-admissible-functions}
An \emph{admissible function} on $X^{\an}$ is a continuous function $h\colon X^{\an}\to\BR\cup\{\pm\infty\}$ such that there exists a skeleton $\Gamma\subset X^{\an}$ with the following properties:
\begin{enumerate}[(a)]
\item The restriction of $h$ to every edge of $\Gamma$ is an affine function
\item The restriction of $h$ to the complement $X^{\an}\setminus\Gamma$ is locally constant
\end{enumerate}
Such a skeleton $\Gamma$ is said to \emph{trivialize} the admissible function $h$. A point $\xi\in X^{\an}$ with $h(\xi)=-\infty$ is called a \emph{pole} of $h$.
\end{Def}
We note that admissible functions are clearly piecewise affine, but not conversely.

\begin{Lem}\label{lem-construct-admissible-functions}
Let $h_1,h_2\colon X^{\an}\to\BR\cup\{\pm\infty\}$ be admissible functions. Then there exist unique admissible functions
\begin{equation*}
h_{+},h_{\max},h_{\min}\colon X^{\an}\to\BR\cup\{\pm\infty\}
\end{equation*}
such that for all points $\xi\in X^{\an}$ with $h_1(\xi)\ne\pm\infty$, $h_2(\xi)\ne\pm\infty$ we have
\begin{equation*}
h_+(\xi)=h_1(\xi)+h_2(\xi),\qquad h_{\max}(\xi)=\max(h_1(\xi),h_2(\xi)),
\end{equation*}
\begin{equation*}
h_{\min}(\xi)=\min(h_1(\xi),h_2(\xi)).
\end{equation*}
\end{Lem}
\begin{proof}
Let $\Gamma\subset X^{\an}$ be a skeleton trivializing both $h_1$ and $h_2$. It follows from part (a) of \mbox{Definition \ref{def-admissible-functions}} that $h_1(\xi)=\pm\infty$ or $h_2(\xi)=\pm\infty$ is only possible if $\xi$ is a Type I point on $\Gamma$. At all other points, we may simply define $h_+$ by pointwise addition, $h_{\max}$ by pointwise maximum, and $h_{\min}$ by pointwise minimum.

If $\xi$ is a Type I point contained in $\Gamma$, then using part (a) again of \mbox{Definition \ref{def-admissible-functions}} choose an interval $[\xi',\xi]\subseteq\Gamma$ on which $h_1$ and $h_2$ are affine, say of the form
\begin{equation*}
h_1(t)=\alpha_1t+\beta_1,\qquad h_2(t)=\alpha_2t+\beta_2
\end{equation*}
with respect to a radius parametrization. Then $h_+$ will be continuous if and only if we set 
\begin{equation*}
h_+(\xi)=\begin{cases}
\infty,&\alpha_1+\alpha_2>0,\\
-\infty,&\alpha_1+\alpha_2<0,\\
\beta_1+\beta_2,&\alpha_1+\alpha_2=0.
\end{cases}
\end{equation*}
Similar definitions work for $h_{\max}$ and $h_{\min}$.
\end{proof}

In the sequel, we will denote the functions $h_+$, $h_{\max}$, and $h_{\min}$ by $h_1+h_2$, $\max(h_1,h_2)$, and $\min(h_1,h_2)$ respectively. 

\begin{Rem}\label{rem-making-new-admissible-functions}
Note that given an admissible function $h\colon X^{\an}\to\BR\cup\{\pm\infty\}$ and a rational number $\alpha\in\BQ$, the function
\begin{equation*}
\alpha h\colon X^{\an}\to\BR\cup\{\pm\infty\},\qquad\xi\mapsto\alpha h(\xi),
\end{equation*}
is obviously admissible as well. Applying \mbox{Lemma \ref{lem-construct-admissible-functions}} inductively, we may thus define, given admissible functions $h_1,\ldots,h_r\colon X^{\an}\to\BR\cup\{\pm\infty\}$ and rational numbers $\alpha_1,\ldots,\alpha_r\in\BQ$, the admissible functions
\begin{equation*}
\alpha_1h_1+\ldots+\alpha_rh_r,\qquad\max(\alpha_1h_1,\ldots,\alpha_rh_r),\qquad\min(\alpha_1h_1,\ldots,\alpha_rh_r).
\end{equation*}
\end{Rem}

Let $h\colon X^{\an}\to\BR\cup\{\pm\infty\}$ be a piecewise affine function, let $\xi\in X^{\an}$ be a point, and let $v$ be a branch at $\xi$. It follows from the definition that there exists an interval $[\xi,\xi']$ representing the branch $v$ on which $h$ is affine. Choose a radius parametrization
\begin{equation*}
[\xi,\xi']\overset{\sim}{\longrightarrow}[a,b]
\end{equation*}
and suppose that $h$ corresponds to the function $t\mapsto\alpha t+\beta$ on $[a,b]$, where $\alpha,\beta\in\BQ$.

\begin{Def}\label{def-slope}
In the above situation, the rational number $\alpha$ is called the \emph{slope} of $h$ in the direction of $v$, denoted
\begin{equation*}
\slope_v(h)\coloneqq\alpha.
\end{equation*}
\end{Def}

\begin{Rem}\label{rem-slope-convention}
In \cite[Section 3.6.3]{ctt}, \emph{piecewise monomial functions} are considered. These are to our piecewise affine functions as the absolute value $\abs{\cdot}_K$ is to the valuation $v_K$. To be precise, a monomial function on an interval $I=[\xi_1,\xi_2]\subset X^{\an}$ is defined to be a function of the form
\begin{equation}\label{equ-monomial-function}
\xi\mapsto\beta\psi(\xi)^n,
\end{equation}
where $\psi\colon I\to[a,b]$ is a radius parametrization in the sense of \cite[\mbox{Section 3.6.2}]{ctt}, where $n\in\BZ$ and where $\beta$ is the absolute value of some element in $K$. A radius parametrization in the sense of \cite[\mbox{Section 3.6.2}]{ctt} is of the form
\begin{equation*}
\xi\mapsto e^{\phi(\xi)},
\end{equation*}
where $\phi$ is a radius parametrization in our sense. Composing the monomial function \eqref{equ-monomial-function} with $-\log$, we obtain the function
\begin{equation}\label{equ-associated-affine-function}
\xi\mapsto-\log(\beta\psi(\xi)^n)=-\log(\beta)-n\log(\psi(\xi))=-n\phi(\xi)-\log(\beta).
\end{equation}
Note that $\log(\beta)\in\BQ$. Thus this function is an affine function. In the same way, the piecewise monomial functions of \cite[Section 3.6.3]{ctt} correspond to our piecewise affine functions.

Let $v$ be the branch at $\xi_1$ represented by the interval $[\xi_1,\xi_2]$. In \cite[Section 3.6.5]{ctt}, the slope of the monomial function \eqref{equ-monomial-function} in the direction of $v$ is defined to be $n$. This is in contrast to the slope of the associated affine function \eqref{equ-associated-affine-function}, which according to our \mbox{Definition \ref{def-slope}} is $-n$. We will occasionally have to account for this discrepancy when appealing to results of \cite{ctt}.
\end{Rem}

\begin{Lem}\label{lem-rational-domain}
Let $h_1,\ldots,h_r\colon X^{\an}\to\BR\cup\{\pm\infty\}$ be admissible functions. Then the subset
\begin{equation*}
\rU=\bigcup_{i=1}^r\{\xi\in X^{\an}\mid h_i(\xi)\ge0\}
\end{equation*}
equals $X^{\an}$ or is an affinoid subdomain of $X^{\an}$.
\end{Lem}
\begin{proof}
Since the maximum of the admissible functions $h_1,\ldots,h_r$ is again admissible, it suffices to consider the case of only one function $h\coloneqq h_1$. Let $\Gamma$ be a skeleton trivializing $h$. Then $\rU$ may be expressed using the canonical retraction $X^{\an}\to\Gamma$ as
\begin{equation*}
\rU=r_\Gamma^{-1}(\{\xi\in\Gamma\mid h(\xi)\ge0\}.
\end{equation*}
Thus $\rU$ is obtained from $X^{\an}$ by deleting finitely many open disks and open anuli. By \cite[Lemma 4.12]{bpr}, the result is an affinoid subdomain of $X^{\an}$ (or all of $X^{\an}$ if the number of deleted disks and anuli is zero).
\end{proof}

\section{Valuative functions}\label{sec-greek-valuative}

\begin{Def}
A \emph{valuative function} is a function of the form
\begin{equation*}
\hat{f}_\alpha\colon X^{\an}\to\BR\cup\{\pm\infty\},\qquad\xi\mapsto\alpha v_\xi(f),
\end{equation*}
where $f\in F_X^\times$ is a non-zero rational function, and where $\alpha\in\BQ$ is a rational number.

If $\alpha=1$, we will denote this function simply by $\hat{f}\coloneqq\hat{f}_1$.
\end{Def}

\begin{Lem}\label{lem-valuative-is-admissible}
Every valuative function $\hat{f}_\alpha\colon X^{\an}\to\BR\cup\{\pm\infty\}$ is admissible. It is trivialized by every skeleton $\Gamma\subset X^{\an}$ containing all zeros and poles of $f$.
\end{Lem}
\begin{proof}
This follows from \cite[Theorem 5.15]{bpr}.
\end{proof}

\begin{Lem}\label{lem-generalized-valuative-functions}
Let $f_1,\ldots,f_r\in F_X^\times$ be non-zero rational functions, and let $\alpha_1,\ldots,\alpha_r\in\BQ$ be rational numbers. Then the admissible function $h\coloneqq\alpha_1\hat{f}_1+\ldots+\alpha_r\hat{f}_r$ (see \mbox{Remark \ref{rem-making-new-admissible-functions}}) is a valuative function.

For all points $\xi\in X^{\an}$ that are not poles or zeros of any of the $f_1,\ldots,f_r$, we have
\begin{equation*}
h(\xi)=\alpha_1\hat{f}_1(\xi)+\ldots+\alpha_r\hat{f}_r(\xi).
\end{equation*}
\end{Lem}
\begin{proof}
The equality $h(\xi)=\sum_i\alpha_i\hat{f}_i(\xi)$ follows from Lemma \ref{lem-construct-admissible-functions}, since the poles and zeros of $f_i$ are exactly the points where $\hat{f}_i(\xi)=\pm\infty$.

To show that $h$ is a valuative function, choose a common denominator $d$ for the $\alpha_1,\ldots,\alpha_r$ and write $\alpha_i=\nu_i/d$ for $i=1,\ldots,r$. Define the function
\begin{equation*}
f\coloneqq\prod_{i=1}^rf_i^{\nu_i}.
\end{equation*}
Then the valuative function $\hat{f}_{1/d}$
satisfies 
\begin{equation*}
\hat{f}_{1/d}(\xi)=\frac{v_\xi(f)}{d}=\frac{\big(\nu_1v_\xi(f_1)+\ldots+\nu_rv_\xi(f_r)\big)}{d}=\alpha_1\hat{f}_1(\xi)+\ldots+\alpha_r\hat{f}_r(\xi)=h(\xi)
\end{equation*}
for all $\xi$ that are not poles or zeros of any $f_i$. By the uniqueness part of \mbox{Lemma \ref{lem-construct-admissible-functions}}, we must have $h=\hat{f}_{1/d}$.
\end{proof}
We will often define a valuative function by a linear combination such as
\begin{equation*}
\alpha_1\hat{f}_1+\ldots+\alpha_r\hat{f}_r,
\end{equation*}
which is justified by \mbox{Lemma \ref{lem-generalized-valuative-functions}}.

The following example illustrates the subtlety of defining a sum of valuative functions, and how it can be overcome using the trick in the proof of \mbox{Lemma \ref{lem-generalized-valuative-functions}}.

\begin{Ex}
Consider the curve $X=\BP_K^1$ and the rational functions $f=x/(x-1)$ and $g=(x-1)^2$. The valuative function
\begin{equation*}
h\coloneqq \hat{f}+\hat{g}
\end{equation*}
may then be defined in terms of the product $fg=x(x-1)$; it is the function
\begin{equation*}
(\BP_K^1)^{\an}\to\BR\cup\{\pm\infty\},\qquad\xi\mapsto v_\xi\big(x(x-1)\big).
\end{equation*}
It is well-defined at all points $\xi\in(\BP_K^1)^{\an}$. For example, $h(0)=h(1)=\infty$. However, we have
\begin{equation*}
\hat{f}(1)=-\infty,\qquad \hat{g}(1)=\infty,
\end{equation*}
and it is not possible to derive the value $h(1)$ from this.
\end{Ex}
Let $\phi\colon Y\to X$ be a finite separable morphism of curves. We consider the norm map of the field extension of function fields $F_Y/F_X$,
\begin{equation*}
\Nm\colon F_Y\to F_X.
\end{equation*}
To end this section, we will study how to relate valuative functions on $Y^{\an}$ and valuative functions on $X^{\an}$ using the norm map.

\begin{Def}\label{def-norm-of-valuative-function}
Let $h\coloneqq\hat{f}_\alpha\colon Y^{\an}\to\BR\cup\{\pm\infty\}$ be a valuative function. The \emph{norm} of $h$, denoted $\Nm h$, is defined to be the valuative function
\begin{equation*}
(\reallywidehat{\Nm f})_\alpha\colon X^{\an}\to\BR\cup\{\pm\infty\}.
\end{equation*}
\end{Def}

\begin{Rem}
Let $f_1,\ldots,f_r\in F_Y^\times$ be rational functions and let $\alpha_1,\ldots,\alpha_r$ be rational numbers. It is an easy consequence of the multiplicativity of the map $\Nm$ that the norm $\Nm h$ of the valuative function $h\coloneqq\alpha_1\hat{f}_1+\ldots+\alpha_r\hat{f}_r$ is the valuative function
\begin{equation*}
\alpha_1(\reallywidehat{\Nm f_1})+\ldots+\alpha_r(\reallywidehat{\Nm f_r}).
\end{equation*}

\end{Rem}

Now assume that $\deg(f)=p$ (the case that we will be concerned with in subsequent chapters) and that $\xi\in X^{\an}$ is a point with $p$ distinct preimages, say $\phi^{-1}(\xi)=\{\eta_1,\ldots,\eta_p\}$. Said differently, the valuation $v_\xi$ has $p$ distinct extensions $v_{\eta_1},\ldots,v_{\eta_p}$ to $F_Y$. Let furthermore $h\colon Y^{\an}\to\BR\cup\{\pm\infty\}$ be a valuative function. 

\begin{Lem}\label{lem-norm-of-valuative-function}
We have
\begin{equation*}
(\Nm h)(\xi)=\sum_{i=1}^p h(\eta_i).
\end{equation*}
\end{Lem}
\begin{proof}
This is obvious in case $h=\hat{f}_\alpha$, where $f\in F_X^\times$. Indeed, in this case $\Nm f= f^p$ and $v_\xi(f)=v_{\eta_i}(f)$ for $i=1,\ldots,p$.

Otherwise we have $h=\hat{f}_\alpha$ for a function $f\in F_Y^\times\setminus F_X^\times$. Then $f$ is necessarily a generator of the field extension of function fields $F_Y/F_X$. By \cite[\mbox{Theorem II.8.2}]{neukirch}, the extensions of $v_\xi$ to $F_Y$ are obtained as follows. There are $p$ embeddings of $F_Y$ into an algebraic closure $\overline{(F_X)_{v_\xi}}$ of the completion $(F_X)_{v_\xi}$. The $p$ extensions $v_{\eta_1},\ldots,v_{\eta_p}$ are obtained by restricting the unique extension of $v_\xi$ to $\overline{(F_X)_{v_\xi}}$ to $F_Y$ along these $p$ embeddings.

If we denote the images of $f$ under the $p$ embeddings $F_Y\to\overline{(F_X)_{v_\xi}}$ by $f_1,\ldots,f_p$, then because $F_Y/F_X$ is separable of degree $p$ we have $\Nm f=f_1\cdots f_p$. Now everything follows by taking valuations in this equality.
\end{proof}

\section{The function $\mu$}\label{sec-greek-mu}

Let $\Gamma_0$ be a skeleton of $\BP_K^1$ containing the point $\infty$. For every closed point $x_0\in(\BP_K^1)^{\an}\setminus\Gamma_0$, the connected component of $(\BP_K^1)^{\an}\setminus\Gamma_0$ containing $x_0$ is an open disk; let us denote it by $\rD(x_0)$. We define the quantity
\begin{equation*}
\mu_{\Gamma_0}(x_0)\coloneqq \textrm{the radius of the disk $\rD(x_0)$}.
\end{equation*}
We will now derive a simple formula for $\mu_{\Gamma_0}(x_0)$ in the case that $\Gamma_0$ is the skeleton spanned by the zero set of a polynomial $f$ and the point $\infty$. That is, $\Gamma_0$ is spanned by the set of Type I points
\begin{equation*}
\{\infty\}\cup\{\xi\in(\BP_K^1)^{\an}\mid v_\xi(f)=\infty\}
\end{equation*}
for a polynomial
\begin{equation*}
f=\sum_{i=0}^da_ix^i\in K[x].
\end{equation*}
Suppose first that $x_0$ is the closed point $x=0$, corresponding to the pseudovaluation ``evaluation at $0$'',
\begin{equation*}
g\mapsto v_{x_0}(g)=v_K(g(0)).
\end{equation*}
The open disk of radius $r$ centered at $x_0$ is contained in the complement $(\BP_K^1)^{\an}\setminus\Gamma_0$ if and only if no zero of $f$ has valuation greater than $r$. It follows from the theory of the Newton polygon that
\begin{equation}\label{equ-simple-mu-np-calculation}
\mu_{\Gamma_0}(x_0)=\max_{i\ge1}\frac{v_K(a_0)-v_K(a_i)}{i}.
\end{equation}
(A brief reminder on Newton polygons may be found in Section \ref{sec-analytic-newton}. See in particular Proposition \ref{prop-newton-polygon-main-fact} for the fact about the valuations of zeros of a polynomial we used here.) Notice that \eqref{equ-simple-mu-np-calculation} makes sense because $x_0$ being contained in the complement of $\Gamma_0$ we have $v_K(a_0)<\infty$.

In the following lemma we generalize the above argument, which allows us to make sense of the formula \eqref{equ-simple-mu-np-calculation} for all points $\xi\in X^{\an}$, including points contained in the skeleton $\Gamma_0$.

\begin{Lem}\label{lem-mu-np}
For every closed point $x_0\in X^{\an}\setminus\Gamma_0$ we have
\begin{equation}\label{equ-not-so-simple-mu-np}
\mu_{\Gamma_0}(x_0)=\max_{i\ge1}\frac{\hat{f}_0(x_0)-\hat{f}_i(x_0)}{i},
\end{equation}
where the $f_i$ are the polynomials
\begin{equation*}
f_i=\frac{f^{(i)}}{i!}\in K[x],\qquad i\ge0.
\end{equation*}
In particular, $\mu_{\Gamma_0}$ may be extended to an admissible function on $X^{\an}$.
\end{Lem}
\begin{proof}
We may regard $x_0$ as an element of $K$ (which is not a zero of the polynomial $f$). The Taylor expansion of $f$ around $x_0$ reads
\begin{equation*}
f(x_0+t)=\sum_{i=0}^df_i(x_0)t^i.
\end{equation*}
Now the open disk of radius $r$ centered at $x_0$ is contained in the complement $(\BP^1_K)\setminus\Gamma_0$ if and only if the polynomial $f(x_0+t)$ has no root of valuation larger than $r$. By the theory of the Newton polygon, this is the case if and only if 
\begin{equation*}
r\ge\max_{i\ge1}\frac{v_K(f_0(x_0))-v_K(f_i(x_0))}{i}.
\end{equation*}
Thus \eqref{equ-not-so-simple-mu-np} holds. The function
\begin{equation*}
\max_{i\ge1}\frac{\hat{f}_0-\hat{f}_i}{i}\colon\quad X^{\an}\to\BR\cup\{\pm\infty\}
\end{equation*}
is admissible by \mbox{Lemma \ref{lem-construct-admissible-functions}}, so the last statement of the lemma immediately follows.
\end{proof}

\begin{Rem}\label{rem-mu-properties}
\begin{enumerate}[(a)]
\item By \mbox{Lemma \ref{lem-mu-np}} we can and will regard $\mu_{\Gamma_0}$ as a function
\begin{equation*}
\mu_{\Gamma_0}\colon X^{\an}\to\BR\cup\{\pm\infty\}.
\end{equation*}
If the skeleton $\Gamma_0$ is clear from context, we will simply write $\mu\coloneqq\mu_{\Gamma_0}$.
\item The only pole of the function $\mu_{\Gamma_0}$ is $\infty$. This is because the function $\hat{f}_d$ (where $d=\deg(f)$) is constant of finite value, so that at least one term appearing in the maximum in \eqref{equ-not-so-simple-mu-np} is $>-\infty$ for any closed point $x_0\ne\infty$.
\end{enumerate}
\end{Rem}

\section{The function $\delta$}\label{sec-greek-delta}

Let $\phi\colon Y\to X$ be a separable degree-$p$ morphism of $K$-curves. In \cite{ctt}, Cohen, Temkin, and Trushin define a \emph{different function}
\begin{equation*}
\delta_\phi\colon Y^{\an}\to[0, 1].
\end{equation*}
Unlike the authors of \cite{ctt} we prefer to work with additive valuations rather than with multiplicative absolute values. Just as valuations and absolute values are related by logarithms and exponentiation, so is the different function used by us to the one of \cite{ctt}. Denoting the different function of \cite{ctt} by $\delta_\phi^{\CTT}$ and the different function used by us by $\delta_\phi^{\add}$, we have the defining relations
\begin{equation*}
\delta_\phi^{\CTT}=\exp\big(-\frac{p-1}{p}\delta_\phi^{\add}\Big),\qquad\delta_\phi^{\add}=-\frac{p}{p-1}\log(\delta_\phi^{\CTT}).
\end{equation*}
Thus $\delta_\phi^{\add}$ is a function
\begin{equation*}
\delta_\phi^{\add}\colon Y^{\an}\to[0,\infty].
\end{equation*}
The function $\delta_\phi^{\CTT}$ is piecewise monomial, so by \mbox{Remark \ref{rem-slope-convention}} the function $\delta_\phi^{\add}$ is piecewise affine. The reason for the factor $p/(p-1)$ in the definition of $\delta_\phi^{\add}$ will become apparent in the next section; see in particular \mbox{Proposition-Definition \ref{propdef-lambda-definition}}. If $\phi$ is a Galois cover, then this normalization makes $\delta_\phi^{\add}$ equal to the \emph{depth Swan conductor} in the sense of Kato (\cite[\mbox{Definition 1.5.2}]{brewis}). See \cite[Section 2.3]{obus-wewers} and in particular \cite[\mbox{Lemma 2.3}]{obus-wewers} for a comparison of this Swan conductor and the different function of \cite{ctt}.

From now on, we will only write $\delta_\phi$ instead of $\delta_\phi^{\add}$. When $\phi$ is clear from context, we will simply denote this function by $\delta$.

\begin{Rem}\label{rem-delta-at-type-1}
$\delta_\phi$ never takes the value $\infty$ at points not of Type I. The behavior of $\delta_\phi$ at Type I points, i.\,e.\ closed points in $X$, is explained in \cite[Theorem 4.6.4]{ctt}. At Type I points that are not totally ramified with respect to $\phi$, the function $\delta_\phi$ vanishes. If $x_0$ is totally ramified, then $\delta_\phi(x_0)>0$. In fact, $\delta_\phi(x_0)=\infty$ if $\characteristic(K)=p$ and $\delta_\phi$ is constant in a neighborhood of $x_0$ otherwise.
\end{Rem}

\begin{Lem}\label{lem-bijection-over-topological-branch-locus}
Let $\eta\in Y^{\an}$ be a point of Type II with image $\xi\coloneqq\phi(\eta)$.

\begin{enumerate}[(a)]
\item If $\delta_\phi(\eta)>0$, then $\eta$ is a wild topological ramification point. In particular, we have $\phi^{-1}(\xi)=\{\eta\}$.
\item If $\delta_\phi(\eta)=0$ and $\eta$ is a wild topological ramification point, then the extension of completed residue fields $\CH(\eta)/\CH(\xi)$ is unramified of \mbox{degree $p$}.
\end{enumerate}
\end{Lem}
\begin{proof}
If $\delta_\phi(\eta)>0$, then by \cite[Lemma 4.2.2]{ctt} the extension of completed residue fields $\CH(\eta)/\CH(\xi)$ is wildly ramified. It follows from \mbox{Remark \ref{rem-geometric-ramification-index-type-ii}} that the ramification index $[\Gamma_{v_\eta}:\Gamma_{v_\xi}]$ equals $1$ and that $[\CH(\eta):\CH(\xi)]=[\kappa(\eta):\kappa(\xi)]$. Since $\CH(\eta)/\CH(\xi)$ is wildly ramified, this degree must equal $p$. Thus $\eta$ is a wild topological ramification point, proving (a).

Now suppose that $\delta_\phi(\eta)=0$ and that $\eta$ is a wild topological ramification point. Again by \cite[Lemma 4.2.2]{ctt}, it follows that $\CH(\eta)/\CH(\xi)$ is not wildly ramified. But by assumption, this extension is of degree $p$. As above, the ramification index $[\Gamma_{v_\eta}:\Gamma_{v_\xi}]$ equals $1$, so $\CH(\eta)/\CH(\xi)$ cannot be tamely ramified. Since it is not wildly ramified, it is unramified. This proves (b).
\end{proof}

Since $\phi$ has degree $p$, it is possible to regard $\delta_\phi$ as a function not only on $Y^{\an}$, but also on $X^{\an}$. This is because $\delta_\phi(\eta)=0$ for every point $\eta\in Y^{\an}$ outside the wild topological ramification locus (\cite[Section 4.1.1]{ctt}). On the other hand, if $\xi\in X^{\an}$ is contained in the wild topological branch locus, then $\phi^{-1}(\xi)$ contains only one point. Thus it makes sense to define
\begin{equation*}
\delta_\phi(\xi)\coloneqq\delta_\phi(\eta)\qquad\textrm{whenever $\phi(\eta)=\xi$}.
\end{equation*}
However, some care is needed when considering $\delta_\phi$ as a function on both $Y^{\an}$ and $X^{\an}$, because the slopes of $\delta_\phi$ at a point $\eta\in Y^{\an}$ and at $\xi\coloneqq\phi(\eta)$ will not agree. The following example illustrates this.

\begin{Ex}\label{ex-kummer-covering}
Consider the Kummer cover
\begin{equation*}
\phi\colon\quad\BP_K^1=Y\to X=\BP_K^1,\qquad x\mapsto x^p.
\end{equation*}
of curves over $K=\BC_p$. The function $\delta_\phi\colon Y^{\an}\to[0,\infty)$ is constant with value $p/(p-1)$ on the interval $I$ connecting the closed points $x=0$ and $\infty$; and if for an interval $J=[\eta,\eta']$ of length $1/(p-1)$ with $J\cap I=\{\eta'\}$ we choose a radius parametrization $J\cong[0,1/(p-1)]$, then $\rest{\delta}{J}$ corresponds to the function
\begin{equation*}
[0,1/(p-1)]\to[0,\infty),\qquad t\mapsto pt,
\end{equation*}
(cf.\ \cite[Remark 4.2.3(ii)]{ctt}). Thus $\delta$ increases with slope $p$ on $J$.

Let us fix such an interval $J=[\eta,\eta']$, say the one with $\eta'$ corresponding to the Gauss valuation centered at $x=0$ of radius $0$ and with $\eta$ corresponding to the Gauss valuation centered at $x=1$ of radius $1/(p-1)$. Then $\xi'\coloneqq\phi(\eta')$ corresponds to the Gauss valuation centered at $x=0$ of radius $0$ as well, but $\xi\coloneqq\phi(\eta)$ corresponds to the Gauss valuation centered at $x=1$ of radius $p/(p-1)$. Thus the interval $\phi(J)=[\xi',\xi]$ has length $p/(p-1)$, and the slope of $\delta$ considered as a function on $\phi(J)$ is $1$.
\end{Ex}

In general, we have the following lemma. In its statement, we consider a pair of skeletons $\Gamma\subset X^{\an}$ and $\Sigma\subset Y^{\an}$, where $\Sigma$ is defined as the inverse image of $\Gamma$ under $\phi$. By this we mean that we have an equality both of the underlying sets $\phi^{-1}(\Gamma)=\Sigma$ and of the vertex sets $\phi^{-1}(S(\Gamma))=S(\Sigma)$.

\begin{Lem}\label{lem-anulus-squeezing}
Let $\Gamma\subset X^{\an}$ and $\Sigma\coloneqq\phi^{-1}(\Gamma)\subset Y^{\an}$ be skeletons. Suppose that $\eta_1,\eta_2\in S(\Sigma)$ are adjacent vertices with images $\xi_1=\phi(\eta_1)$, $\xi_2=\phi(\eta_2)$. Then the topological ramification index is constant on $(\eta_1,\eta_2)$, say equal to $n$, and
\begin{equation*}
l((\xi_1,\xi_2))=nl((\eta_1,\eta_2)).
\end{equation*}
\end{Lem}
\begin{proof}
\sloppy This follows from \cite[Lemma 3.5.8]{ctt}, because the retractions $\ret_\Sigma^{-1}((\eta_1,\eta_2))$ and $\ret_\Gamma^{-1}((\xi_1,\xi_2))$ are open anuli (compare the proof of \mbox{Proposition \ref{prop-skeletons-and-reduction}}).
\end{proof}

An important tool for understanding the behavior of the function $\delta$ is the ``genus formula for wide open domains'' of \cite[Section 6.2.6]{ctt}. We will state this result in Proposition \ref{prop-genus-formula} below, adapted to our notation and applied to the morphism $\phi\colon Y\to X$ of algebraic curves we fixed at the beginning of this section. A \emph{wide open domain} on the analytification $X^{\an}$ of an algebraic curve $X$ is the same as the inverse image under the reduction map $\red_\CX$ associated to a model $\CX$ of $X$ of a closed point on the special fiber $\CX_s$, say
\begin{equation*}
\rD=\red_\CX^{-1}(\overline{x})
\end{equation*}
(cf.\ \cite[Remark 6.2.5(i)]{ctt}). Denote the normalization of the model $\CX$ in the function field $F_Y$ by $\CY$; it is a model of $Y$. Let us write $\{\overline{y}_1,\ldots,\overline{y}_r\}$ for the inverse image of $\overline{x}$ under the induced map on special fibers $\CY_s\to\CX_s$. Then the inverse image of $\rD=\red_\CX^{-1}(\overline{x})$ under $\phi$ is the disjoint union of wide open domains
\begin{equation*}
\phi^{-1}(\rD)=\rC_1\amalg\ldots\amalg\rC_r,\qquad \rC_i=\red_\CY^{-1}(\overline{y}_i).
\end{equation*}
Given a wide open domain $\rD=\red_\CX^{-1}(\overline{x})$, let us denote by $\xi_1,\ldots,\xi_m$ the generic points of the components of $\CX_s$ on which $\overline{x}$ lies. For example, if $\CX$ is semistable, then $m=1$ (if $\overline{x}$ is a smooth point) or $m=2$ (if $\overline{x}$ is an ordinary double point). Each of the $\xi_i$, $i=1,\ldots,m$, corresponds to a Type II point on $X^{\an}$, which we also denote $\xi_i$. Following \cite[Section 6.2.6]{ctt}, we denote by $\rD_\infty$ the set of all branches at any of the points $\xi_1,\ldots,\xi_m$ pointing inside of $\rD$.

The statement of the genus formula below also includes the \emph{genera} of wide open domains. The genus of a wide open domain $\rD$ is by definition (\cite[Section 6.2.4]{ctt}) the quantity
\begin{equation*}
g(\rD)=h^1(\rD)+\sum_{\xi\in\rD}g(\xi).
\end{equation*}
That is, $g(\rD)$ is the sum of the first Betti number of $\rD$ (essentially the number of loops in the topological space $\rD$) and the sum of the genera of all points of positive genus contained in $\rD$.

\begin{Prop}\label{prop-genus-formula}
Suppose that $\rD\subseteq X^{\an}$ is a wide open domain and that $\rC\subseteq Y^{\an}$ is a connected component of $\phi^{-1}(\rD)$. Write $n$ for the degree of the induced morphism $\rC\to\rD$ and $n_v$ for the topological ramification index of $\phi$ at a branch $v$ (as defined in \mbox{Section \ref{sec-preliminaries-topological-ramification}}). Then we have
\begin{equation*}
2g(\rC)-2-n(2g(\rD)-2)=\sum_{Q}d_Q+\sum_{v\in\rC_\infty}n_v-1-\frac{p-1}{p}\slope_v(\delta),
\end{equation*}
where the first sum is over all ramification points $Q$ of $\phi$ contained in $\rC$, and where $d_Q$ is the differential length $d_Q=\len(\Omega_{Y/X,Q})$.
\end{Prop}
\begin{proof}
This is \cite[Theorem 6.2.7]{ctt}. The quantity $R_Q$ appearing in loc.\,cit.\ is replaced with the differential length $d_Q$, as per \cite[Theorem 4.6.4]{ctt}. If $\characteristic(K)=0$, it simply equals $e_Q-1$, where $e_Q$ is the ramification index of $Q$ over $\phi(Q)$.

The slopes $\slope_v(\delta)$ appear with the opposite sign in our version, as explained in \mbox{Remark \ref{rem-slope-convention}}.
\end{proof}

\emph{From now on we assume that $X=\BP_K^1$.} The following example illustrates the kind of insight into the structure of $Y^{\an}$ the genus formula then gives us.

\begin{Ex}\label{ex-genus-formula}
Let $x_0\in X^{\an}$ be a closed point different from $\infty$ and let $\rD$ be an open disk of rational radius centered at $x_0$. It is a wide open domain. Let us assume that $\delta(\xi)>0$, where $\xi\in X^{\an}$ is the boundary point of $\rD$. Then $\xi$ is a wild topological branch point and the inverse image $\rC\coloneqq\phi^{-1}(\rD)$ is connected. The degree of the induced morphism $\rC\to\rD$ equals $p$, and so does the topological ramification index of $\phi$ at the unique branch $v\in\rC_\infty$. Thus \mbox{Proposition \ref{prop-genus-formula}} implies
\begin{equation*}
\frac{p-1}{p}\slope_v(\delta)=1-p-2g(\rC)+\sum_Qd_Q,
\end{equation*}
where the sum is over the ramification points of $\phi$ in $\rC$. We see that the slope of $\delta$ is the steeper the more branch points of $\phi$ are contained in $\rD$. And the slope is the less steep the larger the genus of $\rC$ is. Though we will not use this explicitly, we note that there is a simple formula for the genus of a wide open domain in terms of the special fiber $\CX_s$, where $\CX$ is the model of $X$ such that the wide open domain is defined using the reduction map $\red_\CX$ (\cite[Remark 6.2.5(ii)]{ctt}).

Gradually increasing the radius of $\rD$ will decrease the number of branch points in $\rD$, which changes the slope of $\delta$ in a predictable way. Thus knowledge of $\delta$ would allow us to compute $g(\rC)$ and hence determine if any points of positive genus lie above $\rD$. We will return to this approach in \mbox{Section \ref{sec-quartics-computing-delta-directly}}.
\end{Ex}

\begin{Prop}\label{prop-wild-paths-to-leaves}
Let $\Gamma_0$ be a skeleton of $\BP^1_K$ containing the branch points of $\phi$. Then for every point $\eta\in Y^{\an}\setminus\phi^{-1}(\Gamma_0)$ of positive genus (necessarily of Type II), there exists a unique path $[\eta,\eta_1]$ with $[\eta,\eta_1]\cap\phi^{-1}(\Gamma_0)=\{\eta_1\}$.

Moreover, $\delta$ is strictly positive on $(\eta,\eta_1]$.
\end{Prop}
\begin{proof}
Let $\eta$ be a point as in the statement of the proposition and let $\xi\coloneqq\phi(\eta)$ be its image in $(\BP_K^1)^{\an}$. Because $(\BP_K^1)^{\an}$ is uniquely path-connected, there is a unique path $[\xi,\xi_1]$ with $[\xi,\xi_1]\cap\Gamma_0=\{\xi_1\}$. It suffices to show that $\delta$ is strictly positive on $(\xi,\xi_1]$, because then there will also be a unique path connecting $\eta$ and $\eta_1$.

To this end, let us consider an arbitrary Type II point $\xi_2\in(\xi,\xi_1)$. Let $\rD$ be the open disk containing $\xi$ with boundary point $\xi_2$ and let $\rC$ be a connected component of $\phi^{-1}(\rD)$. The genus formula for wide open domains of \mbox{Proposition \ref{prop-genus-formula}} applied to $\rC$ reads
\begin{equation*}
2g(\rC)-2+2n=\sum_{v\in\rC_\infty}n_v-1-\frac{p-1}{p}\slope_v(\delta),
\end{equation*}
where $n$ is the degree of the induced morphism $\rC\to\rD$. Since $\sum_{v\in C_\infty}n_v=n$ (Remark \ref{rem-local-degree-formula}), it follows that
\begin{equation*}
0<2g(\rC)\le1-n-\sum_{v\in\rC_\infty}\frac{p-1}{p}\slope_v(\delta),
\end{equation*}
which is only possible if $\sum_{v\in\rC_\infty}\slope_v(\delta)<0$. This implies in particular that $\delta$ is strictly positive on the interval $(\xi_2,\xi_1)$. Varying $\xi_2$ we obtain as desired that $\delta$ is strictly positive on $(\xi,\xi_1]$.
\end{proof}

\begin{Lem}\label{lem-tame-topological-branch-locus}
If $\Gamma_0$ is a skeleton of $X=\BP_K^1$ containing the branch points of $\phi$, then every tame topological branch point $\xi\in X^{\an}$ of $\phi$ is contained in $\Gamma_0$. Moreover, if $\xi\in X^{\an}\setminus\Gamma_0$ is a wild topological branch point with $\delta(\xi)=0$, then $\slope_v(\delta)>0$, where $v$ is the branch at $\xi$ pointing towards $\Gamma_0$.
\end{Lem}
\begin{proof}
This is another application of the genus formula for wide open domains. The argument is a variation of the one given in \cite[Lemma 6.3.2]{ctt}.

Suppose that $\xi\in X^{\an}\setminus\Gamma_0$ is \emph{not} a wild topological branch point such that $\slope_v(\delta)>0$ for the branch $v$ at $\xi$ pointing towards $\Gamma_0$. We will then show that $\xi$ is not a topological branch point. We have $\slope_v(\delta)=0$, by assumption, if $\xi$ is a wild topological branch point, and because the wild topological branch locus is closed (Remark \ref{rem-topological-branch-locus}(c)) otherwise. Thus we may choose an open disk containing $\xi$ contained in the complement $X^{\an}\setminus\Gamma_0$ such that $\slope_u(\delta)=0$ for the unique branch $u\in\rD_\infty$. Let $\rC$ be a component of $\phi^{-1}(\rD)$. The genus formula for wide open domains applied to $\rC$ reads
\begin{equation*}
2g(\rC)-2+2n=\sum_{v\in\rC_\infty}n_v-1,
\end{equation*}
where $n$ is the degree of the induced morphism $\rC\to\rD$. Since $\sum_{v\in\rC_\infty}n_v=n$, it follows that
\begin{equation*}
2g(\rC)\le 1-n\le0.
\end{equation*}
This is only possible if $n=1$ and $g(\rC)=0$. Thus $\rC\to\rD$ is an isomorphism. Since the component $\rC$ was arbitrary, $\xi\in\rD$ is not a topological branch point.
\end{proof}

\section{The function $\lambda$}\label{sec-greek-lambda}

We continue the notation from the previous section. Thus $\phi\colon Y\to X=\BP_K^1$ is a separable degree-$p$ morphism of $K$-curves. Let $\Gamma_0\subset X^{\an}$ be the tree spanned by a finite set of Type I points which includes the branch locus of $\phi$ and the point $\infty$. In this section we will define an admissible function
\begin{equation*}
\lambda\colon X^{\an}\to\BR\cup\{\pm\infty\}
\end{equation*}
that will be very important in the subsequent discussion. To define $\lambda$ we will need the following construction.

Let $\Sigma\subset Y^{\an}$ be the graph obtained from $\Sigma_0\coloneqq\phi^{-1}(\Gamma_0)$ by
\begin{enumerate}[1)]
\item adding to the vertex set $S(\Sigma_0)$ all points of positive genus on $\Sigma_0$,
\item attaching to $\Sigma_0$ all paths $[\eta,\eta_1]$, where $\eta\in Y^{\an}\setminus\Sigma_0$ is a point of positive genus and $\eta_1\in\Sigma_0$.
\end{enumerate}
As in the previous section, by the inverse image of $\Gamma_0$ under $\phi$ we mean the graph with underlying set $\phi^{-1}(\Gamma_0)$ and vertex set $\phi^{-1}(S(\Gamma_0))$. We write $\Gamma\coloneqq\phi(\Sigma)$. Note that the construction of $\Sigma$ makes sense because of \mbox{Proposition \ref{prop-wild-paths-to-leaves}}: For each positive genus point $\eta\in Y^{\an}\setminus\Sigma_0$, there exists a unique path connecting $\eta$ and $\Sigma_0$, which consists of wild topological ramification points.

\begin{Prop}\label{prop-minimal-skeleton}
The graph $\Sigma$ is the unique minimal skeleton of $Y$ containing $\Sigma_0=\phi^{-1}(\Gamma_0)$.
\end{Prop}
\begin{proof}
We first prove that $\Sigma$ is a skeleton. We use \cite[Theorem 6.3.4]{ctt}, according to which it suffices to show that the ramification points of $\phi$ are among the vertices of $\Sigma$ and that $\Sigma$ ``locally trivializes'' $\delta$ in the sense of \cite[Section 6.3.3]{ctt} --- but see \mbox{Remark \ref{rem-mistake-in-ctt}} below. The first condition is true, because $\Gamma_0$ contains the branch locus of $\phi$. The latter condition means that for every point $\eta\in\Sigma$ and every branch $v$ at $\eta$ pointing outside of $\Sigma$ we have
\begin{equation}\label{equ-local-trivialization}
\frac{p-1}{p}\slope_v(\delta)=1-n_v
\end{equation}
(the sign here is opposite from the one in \cite[Section 6.3.3]{ctt}, recall \mbox{Remark \ref{rem-slope-convention}}). We now check this second condition. Let us write $\xi\coloneqq\phi(\eta)$ and denote the image of the branch $v$ by $u$. Denote by $\rD$ the open disk with boundary point $\xi$ such that $u$ points into $\rD$.

We claim that $\delta$ is monotonously increasing on every interval $I=[\xi',\xi]\subset\rD$. To show this, fix such an interval $I=[\xi',\xi]$, and consider any interval $[\xi_1,\xi_2]\subset I$ on which $\delta$ is strictly positive. Then $[\xi_1,\xi_2]$ consists of wild topological branch points, so above the branch at $\xi_2$ pointing in the direction of $\xi_1$ there is a unique branch $w$. Let $\rD_{\xi_2}$ denote the open disk with boundary point $\xi_2$ containing $\xi_1$ and let $\rC_{\xi_2}$ be a component of $\phi^{-1}(\rD_{\xi_2})$. The genus formula for wide open domains, \mbox{Proposition \ref{prop-genus-formula}}, applied to $\rC_{\xi_2}$ reads
\begin{equation*}
2g(\rC_{\xi_2})+2p-2=p-1-\frac{p-1}{p}\slope_w(\delta).
\end{equation*}
It follows immediately that $\slope_w(\delta)$ is negative, proving our claim that $\delta$ is monotonously increasing on the interval $[\xi',\xi]\subset\rD$.

Now we show that $\Sigma$ locally trivializes $\delta$, that is, we show that \eqref{equ-local-trivialization} holds. Let $\rC$ be a component of $\phi^{-1}(\rD)$. Suppose first that $\delta(\xi)=0$. Then \mbox{Lemma \ref{lem-tame-topological-branch-locus}} and our claim about monotonicity of $\delta$ show that there are no topological branch points contained in $\rD$. Thus $\frac{p-1}{p}\slope_v(\delta)=0=1-n_v$.

Next, suppose that $\delta(\xi)>0$. Then \mbox{Lemma \ref{lem-tame-topological-branch-locus}} and our claim show that the locus of topological branch points in $\rD$ is non-empty, connected, adjacent to $\xi$, and consisting of \emph{wild} topological branch points. In particular, there are no loops in $\rC$. Moreover, $\rC$ contains no points of positive genus by construction of $\Sigma$, so $g(\rC)=0$. Now \mbox{Proposition \ref{prop-genus-formula}} applied to $\rC$ shows that
\begin{equation*}
2p-2=p-1-\frac{p-1}{p}\slope_v(\delta),
\end{equation*}
so $\frac{p-1}{p}\slope_v(\delta)=1-p=1-n_v$. Thus we have shown that $\Sigma$ locally trivializes $\delta$, so is a skeleton.

Finally, it is obvious that $\Sigma$ is the unique minimal skeleton containing $\Sigma_0$. Indeed, we have only added to $\Sigma_0$ points of positive genus and the unique paths connecting them to $\Sigma_0$.
\end{proof}

\begin{Def}\label{def-phi-minimal-skeleton}
The skeleton $\Sigma$ is called the \emph{$\phi$-minimal skeleton} of $Y$ with respect to $\Gamma_0$.
\end{Def}

\begin{Rem}\label{rem-mistake-in-ctt}
Adding to the vertex set of $\Sigma$ the points of positive genus lying on $\Sigma_0=\phi^{-1}(\Gamma_0)$ is necessary to ensure that $\Sigma$ is a skeleton. This condition is missing in \cite[Theorem 6.3.4]{ctt}. Using our notation, it is asserted there that a graph containing $\Sigma_0$ and locally trivializing $\delta$ is already a skeleton. For a concrete example where this is not the case, see \mbox{Example \ref{ex-worked-example-2}}.
\end{Rem}

We now define the function $\lambda$ that is the subject of this section, in terms of the functions $\delta_\phi$ (Section \ref{sec-greek-delta}), $r_{\Gamma_0}$ (Definition \ref{def-radius-function}), and $\mu_{\Gamma_0}$ (Section \ref{sec-greek-mu}).

\begin{PropDef}\label{propdef-lambda-definition}
There exists a unique piecewise affine function $\lambda_{\phi,\Gamma_0}\colon X^{\an}\to\BR\cup\{\pm\infty\}$ satisfying the following:
\begin{enumerate}[(a)]
\item If $\xi\in X^{\an}\setminus\{\infty\}$ lies in $\Gamma_0$ or $\xi$ is a topological branch point, then
\begin{equation}\label{equ-lambda-definition}
\lambda_{\Gamma_0,\phi}(\xi)=\delta_\phi(\xi)+r_{\Gamma_0}(\xi)+\mu_{\Gamma_0}(\xi).
\end{equation}
\item It is locally constant on the locus
\begin{equation*}
    \{\xi\in X^{\an}\mid\textrm{$\xi\not\in\Gamma_0$ and $\xi$ is not a topological branch point}\}.
\end{equation*}
\end{enumerate}
\end{PropDef}
\begin{proof}
Suppose that $\xi\in X^{\an}\setminus\{\infty\}$ lies in $\Gamma_0$ or is a topological branch point. Then $\lambda_{\Gamma_0,\phi}(\xi)$ is unambiguously determined by \eqref{equ-lambda-definition}, except perhaps if several of $\delta_\phi(\xi),r_{\Gamma_0}(\xi),\mu_{\Gamma_0}(\xi)$ equal $\pm\infty$. However, the only pole of $\mu_{\Gamma_0}$ is by \mbox{Remark \ref{rem-mu-properties}} given by $\xi=\infty$, and we have $\delta_\phi(\xi),r_{\Gamma_0}(\xi)>-\infty$ for all $\xi\in X^{\an}$. Thus we may use \eqref{equ-lambda-definition} to define $\lambda_{\Gamma_0,\phi}(\xi)$ if $\xi\ne\infty$.

Moreover, there clearly exists a unique value $\lambda_{\Gamma_0,\phi}(\infty)$ compatible with $\lambda_{\Gamma_0,\phi}$ being a piecewise affine function; the argument is the same as in \mbox{Lemma \ref{lem-construct-admissible-functions}}. Note that in the case $\characteristic(K)=0$, which is most important to us, we have $\delta_\phi(\infty)<\infty$ by \cite[Theorem 4.6.4(b)]{ctt}. So in this case we can also compute $\lambda_{\Gamma_0,\phi}(\infty)$ using \eqref{equ-lambda-definition}.

In any case, imposing the condition \eqref{equ-lambda-definition} is no issue, and it remains to show that we may extend $\lambda$ to all of $X^{\an}$ in a way compatible with \mbox{condition (b)}. We begin by showing that the set 
\begin{equation*}
R\coloneqq\Gamma_0\cup\{\xi\in X^{\an}\mid\textrm{$\xi$ is a topological branch point}\}
\end{equation*}
is connected.

Choose a skeleton $\Gamma$ of $X$ containing $\Gamma_0$ and such that $\Sigma\coloneqq\phi^{-1}(\Gamma)$ is a skeleton of $Y$. It follows from \cite[Theorem 7.1.4]{ctt} that the topological ramification locus of $\phi$ is a certain ``radial set''. It is the set obtained by adding to
\begin{equation*}
S\coloneqq\{\eta\in\Sigma\mid\textrm{$\eta$ is a topological ramification point}\}
\end{equation*}
all closed intervals $[\eta,\eta']$ where $\eta\in S$,  $\eta'\not\in\Sigma$, and the length of the interval $[\eta,\eta']$ is $\delta(\eta)/p$. (In \cite{ctt}, the length of these intervals is instead $\delta(\eta)/(p-1)$. As usual, the discrepancy is because of our different normalization of $\delta$. See \mbox{Example \ref{ex-kummer-covering}} for an explicit example showing the intervals of length $\delta(\eta)/p$.)

According to \mbox{Proposition \ref{prop-minimal-skeleton}}, a skeleton $\Gamma$ as above is obtained from $\Gamma_0$ by adding to $\Gamma_0$ finitely many intervals consisting of topological branch points of $\phi$; then $\Sigma=\phi^{-1}(\Gamma)$ is the $\phi$-minimal skeleton of $Y$ with respect to $\Gamma_0$. Thus
\begin{equation*}
\phi^{-1}(\Gamma_0)\cup\{\eta\in Y^{\an}\mid\textrm{$\eta$ is a topological ramification point}\}
\end{equation*}
equals the union of $\Sigma$ and the topological ramification locus of $\phi$, and is therefore connected. The image of this set under $\phi$ is the set $R$ used in the definition of $\lambda$; we have now shown that this set $R$ is connected.

Now it follows easily from the fact that $X^{\an}$ is uniquely path-connected that $\lambda$ is well-defined. For every point $\xi\in X^{\an}\setminus R$ there exists a unique path $[\xi,\xi']$ with $\xi'\in R$, but $\xi''\not\in R$ for any $\xi''\in(\xi,\xi')$, and we define $\lambda_{\phi,\Gamma_0}(\xi)=\lambda_{\phi,\Gamma_0}(\xi')$.
\end{proof}

\begin{Rem}
Since $\phi$ and $\Gamma_0$ will usually be fixed and clear from context, we will often simply write $\lambda\coloneqq\lambda_{\phi,\Gamma_0}$ in the sequel.
\end{Rem}

\begin{Rem}\label{rem-wild-and-tame-lambda-definition}
\mbox{Lemma \ref{lem-tame-topological-branch-locus}} allows us to replace ``topological branch point'' in the definition of $\lambda$ with ``\emph{wild} topological branch point''. That is, $\lambda(\xi)$ is defined by \eqref{equ-lambda-definition} if $\xi\in\Gamma_0$ or $\xi$ is a \emph{wild} topological branch point, and $\lambda$ is locally constant on the locus consisting of points that are neither contained in $\Gamma_0$ nor \emph{wild} topological branch points.
\end{Rem}

\begin{Lem}\label{lem-lambda-admissible}
The function $\lambda\colon X^{\an}\to\BR\cup\{\pm\infty\}$ is admissible.
\end{Lem}
\begin{proof}
Let $\Gamma$ and $\Sigma$ denote the skeletons from the proof of Proposition-Definition \ref{propdef-lambda-definition}. That is, $\Sigma$ is the $\phi$-minimal skeleton of $Y$ with respect to $\Gamma_0$ and $\Gamma=\phi(\Sigma)$. We will show that $\lambda$ is locally constant on the complement $X^{\an}\setminus\Gamma$. Let $\xi\in X^{\an}\setminus\Gamma$ be any point, and let $v$ be the unique branch at $\xi$ pointing in the direction of $\Gamma$. We want to show that $\slope_v(\lambda)=0$.

This is obvious from the definition of $\lambda$ if $\xi$ is not a topological branch point. We consider the remaining case that $\xi$ is a topological branch point not contained in $\Gamma$ (necessarily wild, by Lemma \ref{lem-tame-topological-branch-locus}). Since $\xi\not\in\Gamma_0$, we have 
\begin{equation*}
\slope_v(\mu_{\Gamma_0}+r_{\Gamma_0})=\slope_v(\mu_{\Gamma_0})+\slope_v(r_{\Gamma_0})=0-1=-1.
\end{equation*}
Denote by $\eta$ the unique point in the preimage of $\xi$ and denote by $w$ the unique branch above $v$. It follows from \cite[Theorem 7.1.4]{ctt} that $\delta$ (when considered as a function on $Y^{\an}$) satisfies $\slope_w(\delta)=p$. By \mbox{Lemma \ref{lem-anulus-squeezing}} it follows (now considering $\delta$ as a function on $X^{\an}$) that $\slope_v(\delta)=1$. Adding up the slopes of $\delta$, of $r_{\Gamma_0}$, and of $\mu_{\Gamma_0}$ immediately shows $\slope_v(\lambda)=0$.

Since $\delta$, $\mu_{\Gamma_0}$, and $r_{\Gamma_0}$ are all piecewise affine, we can refine $\Gamma$ so that $\lambda$ is affine on each edge of the resulting skeleton. This proves that $\lambda$ is admissible.
\end{proof}

Next, we will study the value $\lambda$ takes on closed points not contained in the skeleton $\Gamma_0$. Let us fix such a point $x_0\in(\BP_K^1)^{\an}\setminus\Gamma_0$. We introduce some notation associated to the point $x_0$ that will be used throughout.

For $r\in\BR$, we denote by
\begin{equation}\label{equ-open-disk}
\rD(r)\coloneqq\{\xi\in(\BP_K^1)^{\an}\mid v_\xi(x-x_0)> r\}
\end{equation}
the \emph{open disk of radius $r$} centered at $x_0$ and by
\begin{equation}\label{equ-closed-disk}
\rD[r]\coloneqq\{\xi\in(\BP_K^1)^{\an}\mid v_\xi(x-x_0)\ge r\}
\end{equation}
the \emph{closed disk of radius $r$} centered at $x_0$. We will denote the boundary point of $\rD[r]$ by $\xi_r$ and the valuation $v_{\xi_r}$ simply by $v_r$. By convention, we set $\xi_\infty\coloneqq x_0$. The valuation $v_r$ is the Gauss valuation centered at $x_0$ of radius $r$, defined by
\begin{equation*}
v_r\big(\sum_ia_i(x-x_0)^i\big)=\min_i(v_K(a_i)+ir).
\end{equation*}
The dependence on the closed point $x_0$ of $\rD(r)$, $\rD[r]$, $v_r$, and $\xi_r$ is not reflected in the notation. This should not cause confusion, since we will usually keep one choice of $x_0$ fixed.

We begin with a simple formula for the value of $\lambda$ at the closed point $x_0$.

\begin{Prop}\label{prop-simple-lambda-consequence}
We have 
\begin{equation*}
\lambda(x_0)=\min\big\{r\ge\mu(x_0)\mid\delta(\xi_r)=0\big\}.
\end{equation*}
\end{Prop}
\begin{proof}
Consider the unique path connecting $x_0$ and $\Gamma_0$; it is the path 
\begin{equation*}
[x_0,\xi_{\mu(x_0)}]=\{\xi_r\mid \mu(x_0)\le r\le\infty\}.
\end{equation*}
By \mbox{Remark \ref{rem-delta-at-type-1}}, we have $\delta(x_0)=0$. If $\delta(\xi)=0$ for every point $\xi$ in $[x_0,\xi_{\mu(x_0)}]$, then by Lemma \ref{lem-tame-topological-branch-locus} there are no topological branch points in the interval $[x_0,\xi_{\mu(x_0)})$. Thus $\lambda$ is by its definition constant on the whole path, the value equal to $\lambda(\xi_{\mu(x_0)})$. Since $\xi_{\mu(x_0)}\in\Gamma_0$, this value is
\begin{equation*}
\lambda(\xi_{\mu(x_0)})=\delta(\xi_{\mu(x_0)})+r_{\Gamma_0}(\xi_{\mu(x_0)})+\mu(\xi_{\mu(x_0)})=\mu(\xi_{\mu(x_0)}).
\end{equation*}
Suppose conversely that $\delta$ is not identically zero on this path, say that the radius $r$ is chosen minimal among $r\ge\mu(x_0)$ with $\delta(\xi_r)=0$. We have $r=r_{\Gamma_0}(\xi_r)+\mu(\xi_r)$. Thus
\begin{equation*}
\lambda(x_0)=\lambda(\xi_r)=\delta(\xi_r)+r_{\Gamma_0}(\xi_r)+\mu(\xi_r)=r
\end{equation*}
as desired.
\end{proof}

In the following, we characterize the value of $\lambda$ at $x_0\in X^{\an}\setminus\Gamma_0$ as the smallest radius $r$ for which we can guarantee that the open disk of radius $r$ around $x_0$ splits into $p$ disjoint open disks along $\phi$.

\begin{Prop}\label{prop-splitting-from-radius-lambda}
Let $x_0\in X^{\an}$ be a Type I point that is not contained in $\Gamma_0$. Let $\rD\coloneqq\rD(\lambda(x_0))$ be the open disk of radius $\lambda(x_0)$ around $x_0$. Then $\phi^{-1}(\rD)$ is a disjoint union of $p$ open disks,
\begin{equation*}
\phi^{-1}(\rD)=\rC_1\amalg\ldots\amalg\rC_p,
\end{equation*}
and for each $1\le i\le p$, the induced morphism $\rest{\phi}{\rC_i}\colon\rC_i\to\rD$ is an isomorphism.
\end{Prop}
\begin{proof}
It follows from \mbox{Proposition \ref{prop-simple-lambda-consequence}} that $\delta=0$ on $\rD$. \mbox{Lemma \ref{lem-tame-topological-branch-locus}} then implies that there are no topological branch points in $\rD$. By Remark \ref{rem-local-degree-formula}, the restriction of $\phi$ to each component of $\phi^{-1}(\rD)$ is a finite morphism of \mbox{degree $1$}, that is, an isomorpism. Thus every component of $\phi^{-1}(\rD)$ must be an open disk, mapped isomorphically to $\rD$ by $\phi$, and everything follows.
\end{proof}

\begin{Kor}\label{cor-lambda-radius-interpretation}
Suppose that the skeleton $\Gamma_0$ is the tree spanned by $\infty$ and the branch points of $\phi$. Then the radius $\lambda(x_0)$ is the smallest radius $r\in\BQ$ for which $\phi^{-1}(\rD(r))$ is a disjoint union of $p$ open disks.
\end{Kor}
\begin{proof}
By \mbox{Proposition \ref{prop-splitting-from-radius-lambda}}, the inverse image of $\rD(\lambda(x_0))$, and hence of $\rD(r)$ for every $r\ge\lambda(x_0)$, is a disjoint union of $p$ open disks. By \mbox{Proposition \ref{prop-simple-lambda-consequence}}, if $r<\lambda$, then we have $r<\mu(x_0)$ or $\delta(\xi_r)>0$. In the first case, $\rD(r)$ contains a branch point of $\phi$, in the second case it contains a point $\xi$ with $\delta(\xi)>0$. In either case, $\phi^{-1}(\rD)$ cannot be a disjoint union of $p$ open disks.
\end{proof}

\begin{Satz}\label{thm-trivializing-lambda}
Suppose that $\Gamma$ is a skeleton of $X$ containing $\Gamma_0$ and such that $\lambda$ is locally constant on $X^{\an}\setminus\Gamma$.

If all points of positive genus on $\Sigma\coloneqq\phi^{-1}(\Gamma)$ are contained in the vertex set $S(\Sigma)$, then $\Sigma$ is a skeleton of $Y^{\an}$.
\end{Satz}
\begin{proof}
This is similar to the proof of \mbox{Proposition \ref{prop-minimal-skeleton}}. By \cite[Theorem 6.3.4]{ctt}, it suffices to show that the ramification points of $\phi$ are among the vertices of $\Sigma$ and that $\Sigma$ ``locally trivializes'' $\delta$ in the sense of \cite[Section 6.3.3]{ctt} (and that the vertex set of $\Sigma$ contains all points of positive genus, recall \mbox{Remark \ref{rem-mistake-in-ctt}}).

By assumption, $\Gamma_0$ contains the branch locus of $\phi$, so we need only show that for every point $\eta\in\Sigma$ and every branch $v$ at $\eta$ pointing outside of $\Sigma$
\begin{equation*}
\frac{p-1}{p}\slope_v(\delta)=1-n_v.
\end{equation*}
We write $\xi\coloneqq\phi(\eta)$ and denote the image of the branch $v$ by $u$. First suppose that $\delta(\eta)>0$. Then $\phi^{-1}(\xi)=\{\eta\}$. Since $u$ points outside of $\Gamma$, the function $r_{\Gamma_0}+\mu_{\Gamma_0}$ increases with slope $1$ in the direction of $u$. Since $\lambda$ is locally constant along $u$, it follows that $\delta$ (considered as a function on $X^{\an}$) decreases with slope $1$ in the direction of $u$. By \mbox{Lemma \ref{lem-anulus-squeezing}}, $\delta$ considered as a function on $Y^{\an}$ decreases with slope $p$ in the direction of $v$. Thus we have proved that $\frac{p-1}{p}\slope_v(\delta)=1-n_v$.

Next, suppose that $\delta(\eta)=0$. Again, $r_{\Gamma_0}+\mu_{\Gamma_0}$ increases with slope $1$ in the direction of $u$. Since $\delta\ge0$, it cannot be the case that $\lambda$ is defined using \eqref{equ-lambda-definition} along $u$. Rather, points close to $\xi$ in the direction of $u$ are not topological branch points, so we have
\begin{equation*}
\frac{p-1}{p}\slope_v(\delta)=0=1-n_v
\end{equation*}
as desired.
\end{proof}

\chapter{An analytic expression for $\lambda$}\label{cha-analytic}

Throughout this entire chapter, the ground field $K$ is algebraically closed. We denote by $X$ the projective line $\BP_K^1$ over $K$. Moreover, we fix a $K$-curve $Y$ (as usual assumed to be smooth, projective, and geometrically connected) and a finite separable morphism $\phi\colon Y\to X$ of degree $p$. Our main goal is to find a characterization of the value at a closed point $x_0\in X$ of the function $\lambda$ associated to $\phi$ that was introduced in \mbox{Section \ref{sec-greek-lambda}}. By \mbox{Theorem \ref{thm-trivializing-lambda}}, knowledge of $\lambda$ is intimately linked to the problem of finding a skeleton of $Y$.

We begin in \mbox{Section \ref{sec-analytic-preliminary}} by introducing several notions and notations that we will need in the sequel. The main thrust of the chapter will then be concerned with studying suitable generators of the field extension $F_Y/F_X$ and their minimal polynomials to obtain information on $\lambda$. In the final \mbox{Section \ref{sec-analytic-tame-locus}} of this chapter, we define the \emph{tame locus} associated to $\phi$, an affinoid subdomain of $Y^{\an}$ encoding the behavior of $\lambda$. It will be key for computing skeletons in concrete examples.

\section{Preliminaries}\label{sec-analytic-preliminary}

For the following definition and for \mbox{Lemma \ref{lem-equivalent-etale-chart}}, it is not necessary to assume that $\phi\colon Y\to X$ be of degree $p$.

\begin{Def}\label{def-etale-chart}
An \emph{\'etale affine chart} for $\phi$ (or simply a \emph{chart}) is a pair $(U,y)$, where $U\subseteq\BP^1_K\setminus\{\infty\}$ is a non-empty open subscheme, and where $y$ is a generator of the field extension of function fields $F_Y/F_X$, such that
\begin{enumerate}[(i)]
\item $y$ is a regular function on  $V\coloneqq\phi^{-1}(U)$, that is, $y\in\CO_Y(V)$,
\item the morphism $V\to\BA_U^1$ corresponding to the function $y$ on the $U$-scheme $V$ is a closed embedding, and
\item the induced morphism $\rest{\phi}{V}\colon V\to U$ is \'etale.
\end{enumerate}
\end{Def}

\begin{Ex}\label{ex-hyperelliptic-chart}
Assume that $\characteristic(K)\ne2$ and let $Y$ be the hyperelliptic curve with Weierstrass equation
\begin{equation*}
Y\colon\quad y^2=f(x),\qquad f(x)\in K[x].
\end{equation*}
Define $U\coloneqq\BP_K^1\setminus(\{\infty\}\cup B)$, where $B$ is the set of zeros of $f$. Then $(U,y)$ is an \'etale affine chart for the hyperelliptic map $\phi\colon Y\to\BP_K^1$, $(x,y)\mapsto x$.
\end{Ex}

\begin{Ex}\label{ex-plane-curve-chart}
Let $Y$ be a smooth plane curve and suppose that an affine chart on $Y$ is given by $\Spec S\subset Y$, where
\begin{equation*}
S\cong K[x,y]/(F),\qquad F\in K[x,y].
\end{equation*}
After possibly replacing $x$ by $x+cy$ for a suitable $c\in K$, we may assume that $S$ is finite over $K[x]$. (This is a special case of Noether normalization, see for example \cite[Sections 12.2.4--12.2.6]{vakil}.)

The corresponding finite morphism $\Spec S\to \BA_K^1$ extends to a cover $\phi\colon Y\to\BP_K^1$. If we define  $U\coloneqq\BP_K^1\setminus(\{\infty\}\cup B)$, where $B$ is the set of branch points of $\phi$, then $(U,y)$ is an \'etale affine chart for $\phi$.
\end{Ex}

\begin{Lem}\label{lem-equivalent-etale-chart}
Let $U\subseteq\BP_K^1\setminus\{\infty\}$ be a non-empty open subscheme and let $y$ be a generator of the field extension $F_Y/F_X$ with minimal polynomial
\begin{equation*}
F(T)=T^p+f_{p-1}T^{p-1}+\ldots+f_1T+f_0.
\end{equation*}
Then the pair $(U,y)$ is an \'etale affine chart for $\phi$ if and only if the following hold:
\begin{enumerate}[(a)]
\item $U$ is disjoint from the branch locus of $\phi$
\item The coefficients $f_{p-1},\ldots,f_1,f_0$ are regular functions on $U$, that is, are contained in $\CO_X(U)$
\item The discriminant of the minimal polynomial $F$ is an invertible function on $U$, that is, a unit in $\CO_X(U)$
\end{enumerate}
\end{Lem}
\begin{proof}
Assume first that the conditions (i), (ii), and (iii) from \mbox{Definition \ref{def-etale-chart}} hold. Of course (iii) implies (a). Write $A\coloneqq\CO_X(U)$ and $B\coloneqq\CO_Y(V)$, so that $B$ is the normalization of $A$ in the function field $F_Y$. Since $y\in B$, its minimal polynomial $F$ has coefficients in $A$, so (b) holds. Moreover, (i) and (ii) mean that we have the following commutative diagrams:
\begin{equation*}
\begin{tikzcd}
V \arrow[rd, "\rest{\phi}{V}"'] \arrow[rr, hook] &   & \BA_U^1 \arrow[ld] & B &                         & {A[T]} \arrow[ll, "y\mapsfrom T"', two heads] \\
                                                 & U &                    &   & A \arrow[lu] \arrow[ru] &                                              
\end{tikzcd}
\end{equation*}
It follows that $B\cong A[T]/(F)$. By \cite[Theorem III.2.6]{neukirch}, the places in $\Spec A$ that ramify in $\Spec A[T]/(F)$ are precisely the places at which the discriminant of $F$ vanishes. Because of (iii), this implies that the discriminant is a unit in $A$, so we have (c).

Conversely, assume that the conditions (a), (b), and (c) from \mbox{Lemma \ref{lem-equivalent-etale-chart}} hold. As before, we write $A=\CO_X(U)$ and $B=\CO_Y(B)$. Now $B$ being the normalization of $A$ in $F_Y$, it follows from (b) that $y$ is a regular function on $V=\Spec B$. Of course, (iii) holds because of (a). Thus to finish the proof it suffices to show that the subring $A[y]$ of $B$ is actually equal to $B$. Because the discriminant of the minimal polynomial $F$ is a unit, it follows that $A[y]\cong A[T]/(F)$ is \'etale over $A$. But this is only possible if $A[y]=B$, since otherwise $A[y]$ is not integrally closed and cannot be \'etale over $A$ (\cite[Tag 025P]{stacks}).
\end{proof}

\emph{From now on, we assume as in the previous chapter that $\phi$ is a degree-$p$ cover}. Suppose that $(U,y)$ is a chart for $\phi$ and let $x_0\in U$ be any closed point. The open subscheme $U$ and the point $x_0$ will remain fixed for most of this chapter, but in time we will modify the generator $y$. 

The complement of the chart $U$ consists of a finite number of closed points. Let $\Gamma_0\subset(\BP_K^1)^{\an}$ be the tree spanned by these points. By construction, $\Gamma_0$ contains the point $\infty$ and the branch locus of $\phi$. In particular, $\Gamma_0$ satisfies the requirements made of the skeleton used to define the functions $\mu$ and $\lambda$ from \mbox{Chapter \ref{cha-greek}}. Thus we have at our disposal the functions
\begin{equation*}
\mu=\mu_{\Gamma_0},\quad\lambda=\lambda_{\phi,\Gamma_0}\colon\quad X^{\an}\to\BR\cup\{\pm\infty\}
\end{equation*}
with respect to $\phi$ and $\Gamma_0$. The main results of this chapter are \mbox{Theorems \ref{thm-lambda-formula-tail-case}} \mbox{and \ref{thm-lambda-formula-general-case}}, in which we find an analytic formula for the value $\lambda(x_0)$ of $\lambda$ at the closed point $x_0$.

We begin by discussing power series expansions of rational functions on $X$. As in \mbox{Section \ref{sec-greek-lambda}}, we denote by $\rD(r)$, $r\in\BQ$, the open disk of radius $r$ centered at $x_0$ and by $\rD[r]$ the closed disk. Moreover, $v_r$ denotes the Gauss valuation of radius $r$ centered at $x_0$ and $\xi_r$ the corresponding Type II point. The function
\begin{equation*}
t\coloneqq x-x_0
\end{equation*}
is a parameter for the family of disks $\rD(r)$, $r\in\BQ$, in the sense that
\begin{equation*}
\rD(r)=\{\xi\in X^{\an}\mid v_\xi(t)>r\}.
\end{equation*}
Let $f\in\CO_X(U)$ be a function. If $r\ge\mu(x_0)$, then we have $\rD(r)\subset U^{\an}$ by definition of $\mu(x_0)$. Applying the restriction map
\begin{equation*}
\CO_X(U)\longrightarrow\CO_{X^{\an}}(U^{\an})\longrightarrow\CO_{X^{\an}}(\rD(r))
\end{equation*}
we can by \mbox{Lemma \ref{lem-radius-of-convergence}} write $f$ as a power series in the parameter $t$ converging on the open disk $\rD(r)$, say
\begin{equation}\label{equ-power-series-expansion}
f=\sum_{i=0}^\infty a_it^i,\qquad \rho(f)=-\liminf_{i\to\infty}\frac{v_K(a_i)}{i}\le r.
\end{equation}
In the following, we use the notation introduced in \mbox{Section \ref{sec-preliminaries-analytic}}, writing
\begin{equation*}
[f]_r\coloneqq[f]_{v_r}\in \kappa(\xi_r).
\end{equation*}
Note that $\kappa(\xi_r)$ is the rational function field $\kappa(\overline{t})$, where the transcendental $\overline{t}$ is the element $\overline{t}\coloneqq[t]_r$. Since an element $f\in\CO_X(U)$ has no poles in $\rD(r)$, the reduction $[f]_r$ has no pole at $\overline{t}=0$. Thus $[f]_r$ can be written as a power series in $\kappa\llbracket\overline{t}\rrbracket$. In the following lemma we investigate the relationship between reduction and passing to power series expansions.

\begin{Lem}\label{lem-power-series-valuation}
Suppose that $r\ge\mu(x_0)$. Let $f$ be an element of $\CO_X(U)$ with power series expansion \eqref{equ-power-series-expansion}. Then the following hold:
\begin{enumerate}[(a)]
\item The minimum $\min_{i\ge0}\{v_K(a_i)+ir\}$ exists and equals $v_r(f)$
\item The power series expansion of $[f]_r$ is $\sum_{i\ge0}\overline{a}_i\overline{t}^i\in\kappa\llbracket t\rrbracket$, where $\overline{a}_i$ is the reduction of $a_i\pi^{ir-v_r(f)}$ to $\kappa$
\end{enumerate}
\end{Lem}
\begin{proof}
For convenience, in this proof we will use the notation
\begin{equation*}
\tilde{v}_r\big(\sum_{i=0}^\infty a_it^i\big)\coloneqq\min_{i\ge0}\{v_K(a_i)+ir\},
\end{equation*}
whenever the right-hand side is well-defined, that is, whenever the set $\{v_K(a_i)+ir\mid i\ge0\}$ has a smallest element.

To prove (a), we assume by way of contradiction that the set $\{v_K(a_i)+ir\mid i\ge0\}$ does not have a smallest element. Then there exist arbitrarily large $N>0$ such that for all $i<N$ we have $v_K(a_N)+Nr<v_K(a_i)+ir$, or equivalently
\begin{equation*}
\frac{v_K(a_i)-v_K(a_N)}{N-i}>r.
\end{equation*}
Choose $s\in\BQ$ with $\frac{v_K(a_i)-v_K(a_N)}{N-i}>s>r$. It follows from \cite[\mbox{Corollary 7.4.11}]{gouvea} that $f$ has at least $N$ zeros on the closed disk $\rD[s]\subset U$. (Strictly speaking, Gouv\^ea only treats the case where $\characteristic(K)=0$, but the proof of \cite[Corollary 7.4.11]{gouvea} is also valid in general --- the result is an easy consequence of the Weierstrass Preparation Theorem for restricted power series, \cite[Section 2.2, Corollary 9]{bosch-lectures}.) Since $f$ is an algebraic function, this cannot be true for arbitrarily large $N$. Thus we have reached a contradiction, showing that the minimum $\min_{i\ge0}\{v_K(a_i)+ir\}$ exists.

Next, let us consider a product decomposition
\begin{equation*}
f=gh,\qquad g=\sum_{i=0}^\infty b_it^i,\quad h=\sum_{i=0}^\infty c_it^i.
\end{equation*}
Assuming that $\min_{i\ge0}\{v_K(b_i)+ir\}$ and $\min_{i\ge0}\{v_K(c_i)+ir\}$ both exist, we have the relation

\begin{equation}\label{equ-tilde-valuation-property}
\tilde{v}_r(f)=\tilde{v}_r(g)+\tilde{v}_r(h).
\end{equation}
Indeed, the inequality ``$\ge$'' is obvious. To show ``$\le$'', it suffices to show that if $\tilde{v}_r(g)=\tilde{v}_r(h)=0$, then $\tilde{v}_r(f)=0$ as well. This follows by considering the images of $g$ and $h$ under the following reduction homomorphism $\pi$, which is the composite of the maps ``evaluate $t$ at $\pi^rt$'' and ``reduce coefficients to $\kappa$'':
\begin{IEEEeqnarray}{CCCCCC}\label{equ-reduction-homomorphism}
\pi\colon&K\llbracket t\rrbracket&\longrightarrow&K\llbracket t\rrbracket&\longrightarrow&\kappa\llbracket\overline{t}\rrbracket\\&\nonumber
\sum_{i=0}^\infty a_it^i&\mapsto&\sum_{i=0}^\infty a_i\pi^{ir}t^i&\mapsto&\sum_{i=0}^\infty\overline{a_i\pi^{ir}}\overline{t}^i
\end{IEEEeqnarray}
The fact that $\tilde{v}_r(g)=\tilde{v}_r(h)=0$ implies that $\pi(g)\ne0$ and $\pi(h)\ne0$. Since $\kappa\llbracket\overline{t}\rrbracket$ is an integral domain, $\pi(f)\ne0$ as well, which means $\tilde{v}_r(f)=0$. Thus we have established \eqref{equ-tilde-valuation-property}.

To finish the proof of (a), we need to show that $v_r(f)=\tilde{v}_r(f)$ for all $f\in\CO_X(U)$. This is clear for polynomials $f\in K[t]$. And since every $f\in\CO_X(U)$ can be written as a quotient $f=g/h$ of two polynomials $g,h\in K[t]$, we deduce $v_r(f)=\tilde{v}_r(f)$ in general from \eqref{equ-tilde-valuation-property} and the corresponding property of $v_r$. This finishes the proof of (a).

For (b) we may assume that $v_r(f)=\tilde{v}_r(f)=0$. Then we need to show that $[f]_r=\pi(f)$ for all $f\in\CO_X(U)$. Again, this is clearly correct for polynomials $f\in K[t]$. To prove it in general, we write $f=g/h$ for $g,h\in K[t]$ and without loss of generality assume $v_r(g)=v_r(h)=0$. Then it follows from \mbox{Lemma \ref{lem-bracket-properties}(a)} and the fact that $\pi$ is a homomorphism that
\begin{equation*}
[f]_r=[g/h]_r=[g]_r/[h]_r=\pi(g)/\pi(h)=\pi(g/h)=\pi(f).
\end{equation*}
\end{proof}

In accordance with \mbox{Lemma \ref{lem-power-series-valuation}}, we will from now on use the notation
\begin{equation*}
v_r(f)=\min_{i\ge0}\{v_K(a_i)+ir\}
\end{equation*}
for $f\in\CO_{X^{\an}}(\rD(r))$ with power series expansion $f=\sum_ia_it^i$. Of course, this only makes sense if the set $\{v_K(a_i)+ir\mid i\ge0\}$ has a minimum. If this is the case we will say that \emph{$v_r(f)$ is well-defined}. Note that $v_r$ satisfies the axioms of a valuation on the set of elements in $\CO_{X^{\an}}(\rD(r))$ where it is well-defined --- the only hard part is the relation $v_r(fg)=v_r(f)+v_r(g)$, which we verified in the proof of \mbox{Lemma \ref{lem-power-series-valuation}}.

For the following definition, let $f\in\CO_{X^{\an}}(\rD(r))$ be a function with power series expansion $\sum_ia_it^i$ for which $v_r(f)$ is well-defined (for example, if $r\ge\mu(x_0)$, then this is the case for all $f\in\CO_X(U)$).

\begin{Def}\label{def-order-and-reduction}
\begin{enumerate}[(a)]
\item The \emph{order} of the function $f$ on the open disk $\rD(r)$ is defined to be
\begin{equation*}
\ord_r(f)\coloneqq\min\big\{i\ge0\mid v_K(a_i)+ir=v_r(f)\big\}.
\end{equation*}
\item We define
\begin{equation*}
[f]_r\coloneqq\sum_{i=0}^\infty\overline{a_i\pi^{ir-v_r(f)}}\cdot\overline{t}^i\in\kappa\llbracket\overline{t}\rrbracket.
\end{equation*}
\end{enumerate}
\end{Def}

Note that by \mbox{Lemma \ref{lem-power-series-valuation}}, \mbox{Definition \ref{def-order-and-reduction}(b)} agrees with the previous \mbox{Definition \ref{def-reduction-to-residue-field}} for $f\in\CO_X(U)$. In the following lemma we record the fact that the statements of \mbox{Lemma \ref{lem-bracket-properties}} also remain valid.

\begin{Lem}\label{lem-generalized-bracket-properties}
Let $f,g\in\CO_{X^{\an}}(\rD(r))$ be functions with power series expansions $f=\sum_ia_it^i$ and $g=\sum_ib_it^i$ respectively such that $v_r(f)$ and $v_r(g)$ are well-defined. Then we have:
\begin{enumerate}[(a)]
\item $[fg]_r=[f]_r\cdot[g]_r$
\item $[f+g]_r=[f]_r+[g]_r$ if and only if $v_r(f)=v_r(g)$ and $[f]_r+[g]_r\ne0$
\end{enumerate}
\end{Lem}
\begin{proof}
In the proof of \mbox{Lemma \ref{lem-power-series-valuation}} we have seen that $[f]_r$ is obtained from $f$ by first multiplying $f$ with $\pi^{-v_r(f)}$ and then applying the reduction homomorphism \eqref{equ-reduction-homomorphism}. With these facts in mind, the proof of the present lemma is identical to that of \mbox{Lemma \ref{lem-bracket-properties}}, with $v_r$ taking the role of $v_\xi$ and the homomorphism \eqref{equ-reduction-homomorphism} taking the role of the reduction denoted $f\mapsto\overline{f}$ in \mbox{Definition \ref{def-reduction-to-residue-field}}.
\end{proof}

\begin{Lem}\label{lem-weierstrass-preparation}For an element $f\in\CO_{X^{\an}}(\rD(r))$ with power series expansion $\sum_ia_it^i$ such that $v_r(f)$ is well-defined, the following quantities are equal:
\begin{enumerate}[(a)]
\item $\ord_r(f)$
\item the order of vanishing of $[f]_r$ at $\overline{t}=0$, i.\,e.\ the $\overline{t}$-adic valuation of $[f]_r$
\item the number of zeros of $f$ on the disk $\rD(r)$, counted with multiplicity
\end{enumerate}
\end{Lem}
\begin{proof}
The equivalence of (a) and (b) is obvious from \mbox{Lemma \ref{lem-power-series-valuation}(b)}. To show that these are equivalent to (c) we may assume without loss of generality that $r=0$. Indeed, replacing $t$ by $t\pi^r$ transforms $f$ into a power series converging on the open disk $\rD(0)$; the valuation of each monomial $a_it^i$ with respect to $v_r$ is the same as the valuation of $a_i(t\pi^r)^i$ with respect to $v_0$.

Now the desired equivalence follows from the Weierstrass Preparation Theorem \cite[Chapter 7, §3, n\textsuperscript{o} 8, Proposition 6]{bourbaki}. It asserts that there is a unique factorization
\begin{equation*}
f=\pi^{v_0(f)}\cdot u\cdot g,
\end{equation*}
where $u\in\CO_K\llbracket t\rrbracket^\times$ is a unit and where $g=\sum_{i=0}^sb_it^i\in\CO_K[t]$ is a polynomial of degree $s=\ord_0(f)$ with $v_K(b_s)=0$ and $v_K(b_i)>0$ for $0\le i<s$. (Note that this is a slightly different Weierstrass Preparation Theorem from the one mentioned in the proof of \mbox{Lemma \ref{lem-power-series-valuation}}, which was about \emph{restricted} power series.) The zeros of $f$ on $\rD(0)$ are the zeros of $g$ on $\rD(0)$, and there are $s$ of these, counted with multiplicity --- the latter is an easy consequence of the theory of the Newton polygon, specifically, of Proposition \ref{prop-newton-polygon-main-fact} below.
\end{proof}

Motivated by \mbox{Lemma \ref{lem-weierstrass-preparation}}, we extend the function $\ord_r$ to $r=\infty$ by letting
\begin{equation*}
\ord_\infty(f)\coloneqq\textrm{the order of vanishing of $f$ at $x_0$}.
\end{equation*}
Suppose now that $f\in K\llbracket t\rrbracket$ is a function on the \emph{closed} disk $\rD[r]$, that is, an element of the Tate algebra $K\{\pi^{-r}t\}$. For example, this is the case if $r>\mu(x_0)$, so that $\rD[r]\subset U^{\an}$. As explained in \mbox{Section \ref{sec-preliminaries-line}}, the elements of the Tate algebra are restricted power series
\begin{equation*}
f=\sum_{i=0}^\infty a_it^i\in K\llbracket t\rrbracket,\qquad\lim_{i\to\infty}v_K(a_i)+ir=\infty.
\end{equation*}
Clearly, the minimum $\min_{i\ge0}\{v_K(a_i)+ir\}$ exists and $v_r(f)$ is the same as the Gauss valuation on restricted power series as discussed in \mbox{Section \ref{sec-preliminaries-line}}. In particular, \mbox{Definition \ref{def-order-and-reduction}} and \mbox{Lemmas \ref{lem-generalized-bracket-properties}} \mbox{and \ref{lem-weierstrass-preparation}} all apply to $f$.

In \mbox{Lemma \ref{lem-weierstrass-preparation}} we saw that the number of zeros of $f$ on the open disk $\rD(r)$ equals the order of vanishing of the reduction $[f]_r$ at $\overline{t}=0$. The analogue for the closed disk $\rD[r]$ is the following.

\begin{Lem}\label{lem-zeros-on-closed-disk}
If $f\in K\{\pi^{-r}t\}$, we have the following:
\begin{enumerate}[(a)]
\item $[f]_r$ is a polynomial
\item The degree of $[f]_r$ equals the number of zeros of $f$ on $\rD[r]$
\end{enumerate}
\end{Lem}
\begin{proof}
That $[f]_r$ is a polynomial follows immediately from the formula in \mbox{Lemma \ref{lem-power-series-valuation}} and the fact that $f$ is represented by a \emph{restricted} power series. That the degree of this polynomial equals the number of zeros of $f$ on $\rD[r]$ follows from the Weierstrass Preparation Theorem for restricted power series \cite[Section 2.2, Corollary 9]{bosch-lectures}.
\end{proof}

\begin{Rem}\label{rem-disc-is-constant}
Let $F$ be the minimal polynomial over $F_X$ of the generator $y$ that is part of the affine \'etale chart $(U,y)$ we fixed earlier in this section. By \mbox{Lemma \ref{lem-equivalent-etale-chart}}, the discriminant $\Delta_F$ is an invertible function on $U$. In particular, it has no zeros in the disk $\rD(r)$ for any $r\ge\mu(x_0)$. Thus \mbox{Lemma \ref{lem-weierstrass-preparation}} shows that the reduction $[\Delta_F]_r$ does not vanish at $\overline{t}=0$ for any $r\ge\mu(x_0)$. And \mbox{Lemma \ref{lem-zeros-on-closed-disk}} shows that $[\Delta_F]_r$ is a constant in $\kappa$ for any $r>\mu(x_0)$.
\end{Rem}

For the following lemma, note that the function
\begin{equation*}
\BQ\to\BQ,\qquad r\mapsto v_r(f)
\end{equation*}
is by \mbox{Lemma \ref{lem-valuative-is-admissible}} induced by a piecewise affine function $\BR\to\BR$. Therefore it makes sense to talk of its right derivative $\partial_+v_r(f)$
with respect to $r$.

\begin{Lem}\label{lem-order-is-derivative}
If $f\in\CO_{X^{\an}}(\rD(r))$ and $v_r(f)$ is well-defined, then we have
\begin{equation*}
\ord_r(f)=\big(\partial_+v_r(f)\big)(r).
\end{equation*}
\end{Lem}
\begin{proof}
As in the proof of \mbox{Lemma \ref{lem-weierstrass-preparation}} we may assume without loss of generality that $r=0$ and use the Weierstrass Preparation Theorem \cite[Chapter 7, §3, n\textsuperscript{o} 8, Proposition 6]{bourbaki} to write
\begin{equation*}
f=\pi^{v_0(f)}\cdot u\cdot(t-\alpha_1)\cdots(t-\alpha_s),
\end{equation*}
where $u\in\CO_K\llbracket t\rrbracket^\times$ and where the $\alpha_1,\ldots,\alpha_s$ are the $s\coloneqq\ord_0(f)$ zeros of $f$ on $\rD(0)$ (this number of zeros is counted with multiplicity, so the $\alpha_i$ may not be pairwise distinct). We may choose $\epsilon>0$ with $v_K(\alpha_1),\ldots,v_K(\alpha_s)>\epsilon$. Then for $0\le r<\epsilon$,
\begin{equation*}
v_r(f)=v_0(f)+\sum_{i=1}^sv_r(t-\alpha_i)=v_0(f)+sr.
\end{equation*}
Thus the right derivative of $v_r(f)$ at $0$ equals $s=\ord_0(f)$.
\end{proof}

Now let $r\ge\lambda(x_0)$ be a rational number. We consider the decomposition from \mbox{Proposition \ref{prop-splitting-from-radius-lambda}},
\begin{equation*}
\phi^{-1}(\rD(r))=\rC_1\amalg\ldots\amalg\rC_p.
\end{equation*}
The function $y$ that is part of the affine \'etale chart $(U,y)$ is a regular function on the affine open subset $V=\phi^{-1}(U)$, that is, an element of $\CO_Y(V)$. Thus we can apply to it the following composite maps
\begin{equation}\label{equ-yi-definition}
\CO_Y(V)\longrightarrow\CO_{Y^{\an}}(V^{\an})\longrightarrow\CO_{Y^{\an}}(\rC_i)\overset{\sim}{\longrightarrow}\CO_{X^{\an}}(\rD(r)),\qquad i=1,\ldots,p.
\end{equation}
Here the first map is induced by the analytification morphism $Y^{\an}\to Y$, the second is a restriction map of the sheaf $\CO_{Y^{\an}}$, and the last is induced by the isomorphism $\rest{\phi}{\rC_i}\colon\rC_i\to\rD(r)$. 

In this way, we obtain $p$ images of the function $y$ in $\CO_{X^{\an}}(\rD(r))$. We denote these by
\begin{equation*}
y_1,y_2,\ldots,y_p,
\end{equation*}
and consider them as power series in the parameter $t$. Note that we have $y_i(0)=y(p_0)$, where $p_0$ is the point lying above $x_0$ contained in $\rC_i$, for $i=1,\ldots,p$. In other words, the constant coefficients of the power series $y_1,\ldots,y_p$ correspond to the $p$ preimages in $\phi^{-1}(x_0)$. In particular, the $y_1,\ldots,y_p$ are pairwise distinct.

\begin{Rem}\label{rem-valuations-preserved}
Let $f\in\CO_Y(V)$ be any function and let $\sum_ia_it^i$ be its image under one of the maps \eqref{equ-yi-definition}. Then $v_r(f)$ is well-defined (i.\,e.\ the minimum $\min_{i\ge0}\{v_K(a_i)+ir\}$ exists) by the same argument as in the proof of \mbox{Lemma \ref{lem-power-series-valuation}}; the key point is that the algebraic function $f$ only has finitely many zeros. Thus \mbox{Definition \ref{def-order-and-reduction}} and the lemmas following it are applicable to $f$.

Now let us furthermore suppose that $r=\lambda(x_0)$ and that $\xi_{\lambda(x_0)}$ is a wild topological branch point of $\phi$, so that $\phi^{-1}(\xi_{\lambda(x_0)})=\{\eta\}$ for a Type II point $\eta\in Y^{\an}$. Then for $f=\sum_ia_it^i$ as above we have
\begin{equation}\label{equ-pass-to-limit}
v_\eta(f)=v_{\lambda(x_0)}(f)=\min_{i\ge0}\{v_K(a_i)+i\lambda(x_0)\}.
\end{equation}
Indeed, if $r>\lambda(x_0)$ and $\eta_r\in\rC_i$ is chosen so that $\phi(\eta_r)=\xi_r$, we have $v_{\eta_r}(f)=\min_{i\ge0}\{v_K(a_i)+ir\}$, because $\rest{\phi}{\rC_i}$ is an isomorphism. Then \eqref{equ-pass-to-limit} follows because of continuity of evaluation at $f$.
\end{Rem}

Since the power series $y_1,\ldots,y_p$ are contained in $\CO_{X^{\an}}(\rD(\lambda(x_0)))$, they have radius of convergence at most $\lambda(x_0)$. In fact, it follows from \mbox{Corollary \ref{cor-lambda-radius-interpretation}} that in the case that $\Gamma_0$ is the tree spanned by $\infty$ and the branch points of $\phi$, the radius of convergence of the $y_1,\ldots,y_p$ is exactly $\lambda(x_0)$. We explore this in the following example.

\begin{Ex}\label{ex-lambda-as-radius-of-convergence}
Consider the plane quartic
\begin{equation*}
Y\colon\quad F(y)=y^3 + 3xy^2 - 3y - 2x^4 - x^2 - 1=0
\end{equation*}
over $K=\BC_3$. The branch points of the degree-$3$ cover $\phi\colon Y\to\BP_K^1$ given by $(x,y)\mapsto x$ consist of the point $\infty$ and the zeros of the discriminant of $F$,
\begin{equation*}
\Delta_F=-108x^8 + 216x^7 - 108x^6 + 432x^5 - 135x^4 + 270x^3 + 27x^2 + 162x + 81.
\end{equation*}
We use the \'etale affine chart as explained in \mbox{Example \ref{ex-plane-curve-chart}}, that is to say, $U$ is just the complement of the branch locus of $\phi$. Thus $\Gamma_0$ is the tree spanned by the zeros of $\Delta_F$ and $\infty$. Two of the zeros of $\Delta_F$ have valuation $1/2$ and the other six have valuation $0$. In particular, we may take $x_0=0$ as a closed point in the complement of $\Gamma_0$. 

Let us fix a point $p_0=(x_0,y_0)$ with $\phi(p_0)=x_0$. We denote by
\begin{equation*}
y_1=\sum_{i=0}^\infty b_it^i\in K\llbracket t\rrbracket
\end{equation*}
the image of $y$ corresponding to $p_0$. It is easy to calculate the coefficients $b_1,b_2,\ldots$ and their valuations (see \mbox{Code Listing \ref{list-power-series-valuations}}). We have the following:

\begin{center}
\begin{tabular}{c|cccccccccc}
$i$&$1$&$2$&$3$&$4$&$5$&$6$&$7$&$8$&$9$&$10$\\\midrule
$\frac{v_K(b_i)}{i}$&$-\frac{1}{3}$&$-\frac{2}{3}$&$-\frac{7}{8}$&$-\frac{1}{2}$&$-\frac{4}{5}$&$-\frac{8}{9}$&$-\frac{6}{7}$&$-\frac{7}{8}$&$-\frac{25}{27}$&$-\frac{9}{10}$
\end{tabular}
\end{center}
It looks as though the radius of convergence $y_1$,
\begin{equation*}
\rho(y_1)=-\liminf_{i\to\infty}\frac{v_K(b_i)}{i},
\end{equation*}
equals $1$. And indeed we will see in \mbox{Example \ref{ex-worked-example-3}} that $\lambda(x_0)=1$.
\end{Ex}

We discuss now another perspective on the power series $y_1,\ldots,y_p$ that will become important later, when we will consider not only one fixed closed point $x_0$, but vary $x_0$. Let us write $A\coloneqq\CO_X(U)$ and $B\coloneqq\CO_Y(V)$, where $V=\phi^{-1}(U)$, and denote by
\begin{equation*}
F(T)=T^p+f_{p-1}T^{p-1}+\ldots+f_1T+f_0\in A[T]
\end{equation*}
 the minimal polynomial of $y$. We now introduce a ``generic'' version of the minimal polynomial $F$, the polynomial
\begin{equation}\label{equ-implicit-equation}
\tilde{F}(T)=T^p+\tilde{f}_{p-1}T^{p-1}+\ldots+\tilde{f}_1T+\tilde{f}_0\in A\llbracket t\rrbracket[T],
\end{equation}
where the $\tilde{f}_i\in A\llbracket t\rrbracket$ are given by the Taylor expansions
\begin{equation}\label{equ-taylor-expansions}
\tilde{f}_i(t)=\sum_{l=0}^\infty a_{i,l}t^l,\qquad a_{i,l}=\frac{f_i^{(l)}}{l!},\qquad 0\le i\le p-1.
\end{equation}
The use of the tildes will always indicate that we are considering such a ``generic'' polynomial, whose coefficients are themselves not just elements of $K$, but functions. Evaluation of the $\tilde{f}_i$ at $x_0$ yields the minimal polynomial of $y$ over the subfield $K(t)$ of $F_Y$. That is, the minimal polynomial of $y$ over the subfield $K(t)$ is given by
\begin{equation}\label{equ-implicit-equation-plugged-in}
G(T)=T^p+g_{p-1}T^{p-1}+\ldots+g_1T+g_0,
\end{equation}
\begin{equation*}
\textrm{where}\qquad g_i=\sum_{l=0}^\infty a_{i,l}(x_0)t^l.
\end{equation*}
We need the following lemma, which may be viewed as an algebraic version of the Implicit Function Theorem:

\begin{Lem}\label{lem-implicit-function-theorem}
Let $R$ be a ring and let $F\in R\llbracket t\rrbracket[T]$ be a monic polynomial whose coefficients are formal power series over $R$ in the indeterminate $t$, say
\begin{equation*}
F=T^n+f_{n-1}T^{n-1}+\ldots+f_1T+f_0,\qquad f_0,\ldots,f_{n-1}\in R\llbracket t\rrbracket.
\end{equation*}
Suppose that there exists a zero $a_0\in R$ of the polynomial
\begin{equation*}
F_0\coloneqq T^n+f_{n-1}(0)T^{n-1}+\ldots+f_1(0)T+f_0(0)
\end{equation*}
and that $F_0'(a_0)\in R^\times$, where $F_0'$ is the derivative
\begin{equation*}
F_0'=nT^{n-1}+(n-1)T^{n-2}f_{n-1}(0)+\ldots+f_1(0).
\end{equation*}
Then there exists a unique lift of $a_0$ to a power series $P=a_0+a_1t+a_2t^2+\ldots\in R\llbracket t\rrbracket$ with $F(P)=0$.
\end{Lem}
\begin{proof}
We construct $P$ inductively. First let us make the ansatz $P_1=a_0+a_1t+O(t^2)$. We have
\begin{equation*}
F(a_0+a_1t)=F(a_0)+F_0'(a_0)a_1t+O(t^2).
\end{equation*}
Since $F_0'(a_0)\in R^\times$, we may and must choose $a_1$ so that $F(a_0+a_1t)\in O(t^2)$. Similarly, the ansatz $P_2=a_0+a_1t+a_2t^2$ leads to
\begin{equation*}
F(a_0+a_1t+a_2t^2)=F(a_0+a_1t)+F'_0(a_0)a_2t^2+O(t^3),
\end{equation*}
which forces us to choose $a_2$ so that $F(P_2)\in O(t^3)$. Continuing in this manner uniquely defines the power series $P$ with $F(P)=0$ as desired.
\end{proof}

We apply \mbox{Lemma \ref{lem-implicit-function-theorem}} to the monic polynomial $\tilde{F}\in A\llbracket t\rrbracket[T]\subseteq B\llbracket t\rrbracket[T]$ of \eqref{equ-implicit-equation} and the zero $y\in B$. Because $\phi$ is \'etale over $U$, we have $\tilde{F}'_0(y)\in B^\times$. It follows that there exists a formal solution in $B\llbracket t\rrbracket$ to \eqref{equ-implicit-equation}, which we denote
\begin{equation}\label{equ-general-formal-solution}
\tilde{y}=y+b_1t+b_2t^2+b_3t^3+\ldots.
\end{equation}
The power series $y_1,\ldots,y_p$ constructed above are obtained from $\tilde{y}$ by evaluating each of the coefficients of $\tilde{y}$ at the points $p_0\in\phi^{-1}(x_0)$. That is,
\begin{equation*}
y_i=y(p_0)+b_1(p_0)t+b_2(p_0)t^2+\ldots,
\end{equation*}
where $p_0\in\phi^{-1}(x_0)$ is contained in $\rC_i$. This may be seen by applying \mbox{Lemma \ref{lem-implicit-function-theorem}} again, to the polynomial in \eqref{equ-implicit-equation-plugged-in}, and to the root $y(p_0)$ of this polynomial.

\section{Newton polygons}\label{sec-analytic-newton}

In this section, we recall some basic facts about Newton polygons and introduce some terminology that we will need. We will later impose further restrictions, but for now let $(F,v_F)$ be any valued field, with ring of integers $\CO_F$. We begin by recalling the central definition.

\begin{Def}
Let $G=\sum_{i=0}^na_iT^i\in F[T]$ be a monic polynomial. The \emph{Newton polygon} of $G$ is the boundary of the lower convex hull of the set
\begin{equation}\label{equ-newton-vertex-set}
\big\{(0,v_F(a_0)),(1,v_F(a_1)),\ldots,(n,v_F(a_n))\big\}\subset\BR\times(\BR\cup\{\infty\}).
\end{equation}
We denote it by $N_G$.
\end{Def}

Geometrically, $N_G$ is obtained as follows. Begin with the ray $\{(n,t)\mid t\le0\}$ pointing downward from the rightmost point $(n,0)$ in \eqref{equ-newton-vertex-set}. Rotate this ray clockwise until it hits another point $(i,v_F(a_i))$. Then remove the line segment connecting the two points $(n,0)$ and $(i,v_F(a_i))$ from the ray and continue rotating its remainder, again until hitting another point. Continue this procedure until the point $(0,v_F(a_0))$ in \eqref{equ-newton-vertex-set} is hit. The Newton polygon $N_G$ then consists of the line segments connecting the points that were hit.

The points $(i,v_F(a_i))$ that lie on $N_G$ are called the \emph{vertices} of $N_G$.

Note that our definition allows certain coefficients $a_i$, $i<n$, to vanish. If this is the case, the corresponding vertex $(i,v_F(a_i))$ is $(i,\infty)$. If there exists an $i_0>0$ such that $a_{i_0}\ne0$, but $a_i=0$ for all $i<i_0$, then $N_G$ contains a vertical line segment, connecting the vertices $(i_0,v_F(i_0))$ and $(0,\infty)$.

\begin{Def}\label{def-theta}
We denote by $\theta_G$ the negative of the steepest slope of the Newton polygon $N_G$. Thus
\begin{equation*}
\theta_G=\max_{1\le i\le p}\frac{v_F(a_0)-v_F(a_i)}{i},
\end{equation*}
with the convention that $\theta_G=\infty$ if $a_0=0$.
\end{Def}

The following proposition is the most important fact about the Newton polygon of a polynomial. 

\begin{Prop}\label{prop-newton-polygon-main-fact}
If the vertices $(i,v_F(a_i))$ and $(j,v_F(a_j))$ of $N_G$, where $i<j$, are the endpoints of one of the line segments making up $N_G$, then in its splitting field, the polynomial $G$ has exactly $j-i$ roots of valuation
\begin{equation*}
\frac{v_F(a_i)-v_F(a_j)}{j-i}.
\end{equation*}
\end{Prop}
\begin{proof}
This is essentially the same statement as \cite[Theorem II.6.3]{neukirch}; however, in loc.\,cit.\ it is assumed that $a_0\ne0$. Our version is easily reduced to that case. Indeed, write $G=T^rH$, where $H(0)\ne0$. Then $N_G$ is obtained from $N_H$ by translating $N_H$ to the right by $r$ and adding a vertical line segment. The version of \cite[Theorem II.6.3]{neukirch} connects the non-vertical slopes of $N_G$ with the non-zero roots of $G$, and the first and vertical line segment of $N_G$ correctly predicts $r$ roots of infinite valuation.
\end{proof}

\emph{From now on, we assume that $F$ has residue field of characteristic $p>0$ and that $G$ is a monic irreducible polynomial of degree $p$.}

\begin{Def}\label{def-inseparable-separable-np}
The Newton polygon $N_G$ is called \emph{inseparable} if its only vertices are the points $(0,v_F(a_0))$ and $(p,0)$. Otherwise, it is called \emph{separable}.
\end{Def}

\begin{Rem}\label{rem-inseparable-np}
If $N_G$ is inseparable, then it consists of a single line segment. If $N_G$ is separable, then it may or may not consist of a single line segment. In particular, if the valuations of the roots of $G$ are not all equal, then it follows from \mbox{Proposition \ref{prop-newton-polygon-main-fact}} that $N_G$ is separable.
\end{Rem}

\begin{Lem}\label{lem-inseparable-np}
Suppose that $v_F$ has a unique extension $v_L$ to $L\coloneqq F(\alpha)$, where $\alpha$ is a root of $G$ with $v_L(\alpha)=0$. Then the following hold:

\begin{enumerate}[(a)]
\item $G\in\CO_F[T]$, and the polynomial $\overline{G}$ obtained from $G$ by mapping each coefficient to the residue field of $v_F$ is inseparable if and only if $N_G$ is inseparable.
\item If $\overline{G}$ is reducible, then $N_G$ (and hence $\overline{G}$) is inseparable.
\end{enumerate}
\end{Lem}
\begin{proof}
That $v_F$ has a unique extension to $L$ implies that $G$ is irreducible over the completion of $F$ with respect to $v_F$ (\cite[Theorem II.8.2]{neukirch}). It follows that $N_G$ consists of a single line segment (\cite[Theorem II.6.4]{neukirch}). Since $v_L(\alpha)=0$, it is the line segment connecting $(0,0)$ and $(p,0)$. Thus we have
\begin{equation}\label{equ-straight-line-np}
v_F(a_i)\ge v_F(a_0)=0,\qquad 0<i<p.
\end{equation}
In particular, all coefficients $a_0,\ldots,a_p$ lie in $\CO_F$, so it makes sense to define the reduction $\overline{G}$. Moreover, $N_G$ is inseparable if and only all the inequalities in \eqref{equ-straight-line-np} are strict.

The reduction $\overline{G}$ cannot have more than one distinct irreducible factor, since any factorization into distinct factors would by Hensel's Lemma lift to a factorization of $G$ over the completion of $F$.

There are two possible cases: either the irreducible factor is of degree $1$ or it is of degree $p$. In the first case, write $T-\overline{\alpha}$ for this factor; then
\begin{equation*}
\overline{G}(T)=(T-\overline{\alpha})^p=T^p-\overline{\alpha}^p
\end{equation*}
is inseparable. Clearly the inequalities \eqref{equ-straight-line-np} are all strict, so $N_G$ is inseparable as well. This shows (b), and (a) in the case that $\overline{G}$ is reducible.

In the second case, $\overline{G}$ is itself irreducible, so it is inseparable if and only if its derivative $\overline{G}'$ vanishes. We have
\begin{equation*}
\overline{G}=T^p+\overline{a}_{p-1}T^{p-1}+\ldots+\overline{a}_1T+\overline{a}_0,
\end{equation*}
\begin{equation*}
\overline{G}'=(p-1)\overline{a}_{p-1}T^{p-2}+\ldots+2\overline{a}_2T+\overline{a}_1,
\end{equation*}
where $\overline{a}_i$, $0\le i<p$, is the image of $a_i$ in the residue field of $v_F$. Clearly $\overline{G}'=0$ if and only if all inequalities in \eqref{equ-straight-line-np} are strict, that is, if and only if $N_G$ is inseparable. This shows (a) in the case that $\overline{G}$ is irreducible.
\end{proof}

Our most important application of Newton polygons is to Type II valuations on the function field of a curve over the complete ground field $K$. For the rest of this section, suppose in particular that $X=\BP_K^1$. As before, denote by $\rD(r)$ the open disk of radius $r\in\BQ$ around the fixed closed point $x_0\in X$. We may identify the function field of $X$ with the rational function field $K(t)$ in the parameter $t=x-x_0$ associated to the closed point $x_0$. Denoting as usual the boundary point of $\rD(r)$ by $\xi_r$, we obtain a family of valuations $v_r=v_{\xi_r}$ on $K(t)$. Thus given a polynomial $G=\sum_ia_iT^i\in K(t)[T]$, varying $r$ we also obtain a family of Newton polygons of $G$ with respect to the different $v_r$. We denote the Newton polygon of $G$ with respect to $v_r$ by $N_G(r)$. Similarly, varying $r$, we may regard $\theta_G$ as a function of $r$,
\begin{equation*}
\theta_G\colon\BQ\to\BQ,\qquad r\mapsto\theta_G(r)=\max_{1\le i\le r}\frac{v_r(a_0)-v_r(a_i)}{i}.
\end{equation*}
We end this section with a lower bound for the value of the function $\lambda$ at the closed point $x_0$, before tackling the more delicate question of an upper bound in the next section. Recall that we have fixed in the previous section an affine \'etale chart $(U,y)$. As before, we write
\begin{equation*}
G(T)=T^p+g_{p-1}T^{p-1}+\ldots+g_1T+g_0
\end{equation*}
for the minimal polynomial of $y$ over $K(t)$.

Note that if $\lambda(x_0)=\mu(x_0)$, the following proposition is never applicable.

\begin{Prop}\label{prop-np-inseparable}
If $r\in\BQ$ is a rational number with $\lambda(x_0)>r\ge\mu(x_0)$, then the Newton polygon $N_G(r)$ is inseparable.
\end{Prop}
\begin{proof}
By \mbox{Proposition \ref{prop-simple-lambda-consequence}} we have $\delta(\xi_r)>0$. It follows from \mbox{Lemma \ref{lem-bijection-over-topological-branch-locus}} that we have $\phi^{-1}(\xi_r)=\{\eta\}$ for a point $\eta\in Y^{\an}$. In other words, $v_\eta$ is the unique extension of the valuation $v_r$ to $F_Y$. Moreover, the extension $\kappa(\eta)/\kappa(\xi_r)$ is purely inseparable of degree $p$ by \mbox{Remark \ref{rem-geometric-ramification-index-type-ii}}.

To place ourselves in the position to apply \mbox{Lemma \ref{lem-inseparable-np}}, we replace $y$ with $z=y\pi^{-s}$, where $s=v_\eta(y)$. Then $v_\eta(z)=0$, and the minimal polynomial of $z$ is
\begin{equation*}
H(T)=T^p+g_{p-1}\pi^{-s}T^{p-1}+\ldots+g_1\pi^{(1-p)s}T+g_0\pi^{-ps}.
\end{equation*}
Denote by $\overline{H}$ the polynomial obtained from $H$ by mapping each coefficient to the residue field $\kappa(\xi_r)$. If $\overline{H}$ is reducible, then $N_H(r)$ is inseparable by \mbox{Lemma \ref{lem-inseparable-np}(b)}. Otherwise, $\overline{z}\coloneqq[z]_r$, the reduction of $z$, is a generator of the inseparable field extension $\kappa(\eta)/\kappa(\xi_r)$. Thus $N_H(r)$ is inseparable by \mbox{Lemma \ref{lem-inseparable-np}(a)}.

Now it suffices to notice that $N_G(r)$ is inseparable if and only if $N_H(r)$ is. This is an easy consequence of the definition. Indeed, let us write $h_i\coloneqq g_i\pi^{(i-p)s}$, $i=0,\ldots,p-1$, for the coefficients of $H$. Then $N_H(r)$ is inseparable if and only if
\begin{equation*}
v_r(g_i)+(i-p)s=v_r(h_i)>v_r(h_0)=0,\qquad 0<i<p.
\end{equation*}
Rewriting this yields
\begin{equation*}
\frac{v_r(g_i)}{p-i}>\frac{v_r(g_0)}{p}=s,
\end{equation*}
which is true for $i=1,\ldots,p-1$ if and only if $N_G(r)$ is inseparable.
\end{proof}

\begin{Rem}\label{rem-epp-foreshadowing}
In the proof of \mbox{Proposition \ref{prop-np-inseparable}}, we had to distinguish two cases: Either the reduction $\overline{z}$ of the generator $z$ was a generator of the extension of residue fields $\kappa(\eta)/\kappa(\xi_r)$, or it was not.

This distinction, and the fact that we cannot read off from the shape of the Newton polygon in which case we are, is significant. In the next section, we will formulate a condition (see \mbox{Assumption \ref{ass-combinatorial-condition}}) for ruling out that $\overline{z}$ is not a generator.
\end{Rem}

\section{An upper bound for $\lambda$: tail case}

Recall that in \mbox{Section \ref{sec-analytic-preliminary}} we have fixed an \'etale affine chart $(U,y)$ for $\phi$. We also retain the closed point $x_0\in U$ and the associated notations (introduced in \mbox{Section \ref{sec-greek-lambda}}) for the open and closed disks of radius $r\in\BQ$ around $x_0$, $\rD(r)$ and $\rD[r]$. Also recall the family of functions $\ord_r$ of \mbox{Definition \ref{def-order-and-reduction}}, which associate an order on $\rD(r)$ to all $f\in\CO_{X^{\an}}(\rD(r))$ for which $v_r(f)$ is well-defined. As in \eqref{equ-implicit-equation-plugged-in}, we denote the minimal polynomial of $y$ over $K(t)$ by
\begin{equation*}
G(T)=T^p+g_{p-1}T^{p-1}+\ldots+g_1T+g_0,
\end{equation*}
where $t=x-x_0$ is the parameter associated to $x_0$.

The results we prove in this section and the next only hold conditional on an assumption on $G$:

\begin{Ass}\label{ass-combinatorial-condition}
For some integer $m\ge1$, we have the inequalities
\begin{equation*}
\ord_\infty(g_0)\ge m,\qquad\ord_{\mu(x_0)}(g_0)<mp.
\end{equation*}
\end{Ass}
We note that this assumption is in terms of the parameter $t$, and hence dependent on the closed point $x_0$. This should cause no confusion, since the closed point $x_0$ is still fixed throughout this section and the next.

The condition in Assumption \eqref{ass-combinatorial-condition} will not usually hold \emph{a priori}. For now, we will prove our results assuming that it holds; then in \mbox{Section \ref{sec-analytic-good-equation}}, we will see how to modify the generator $y$ so that the condition is guaranteed to hold.

The ``tail case'' in the title of this section refers to the fact that the main result of this section, \mbox{Theorem \ref{thm-lambda-formula-tail-case}}, depends on another additional assumption, namely that $\lambda(x_0)>\mu(x_0)$. In the next section, we will remove this condition at the cost of replacing the formula for $\lambda(x_0)$ of \mbox{Theorem \ref{thm-lambda-formula-tail-case}} with a more complicated formula.

The reason for dubbing the case $\lambda(x_0)>\mu(x_0)$ the ``tail case'' will become apparent once in \mbox{Section \ref{sec-analytic-tame-locus}} we introduce the \emph{tame locus} associated to the tree $\Gamma_0$. As explained in \mbox{Lemma \ref{lem-mu-lambda-tilde-distinction}}, a closed point $x_0$ in the tame locus is contained in the \emph{tail locus} if and only if $\lambda(x_0)>\mu(x_0)$. 

Recall that in \mbox{Section \ref{sec-analytic-preliminary}} we have introduced the $p$ roots $y_1,\ldots,y_p\in\rD(\lambda(x_0))$ of the minimal polynomial $G$ of $y$. The following simple lemma contains the central combinatorial argument based on \mbox{Assumption \ref{ass-combinatorial-condition}}.

\begin{Lem}\label{lem-combinatorial-argument}
Suppose that \mbox{Assumption \ref{ass-combinatorial-condition}} holds for some $m\ge1$. Then there exist $i_0,i_1\in\{1,\ldots,p\}$ such that
\begin{equation*}
\ord_r(y_{i_0})<m\le\ord_r(y_{i_1})
\end{equation*}
for all rational numbers $r\ge\lambda(x_0)$.
\end{Lem}
\begin{proof}
Viewing the $y_1,\ldots,y_p$ and the coefficients of $G$ as elements of $K\llbracket t\rrbracket$, we have
\begin{equation*}
y_1\cdots y_p=g_0.
\end{equation*}
Thus by \mbox{Assumption \ref{ass-combinatorial-condition}}, at least one of the $y_1,\ldots,y_p$, say $y_{i_1}$, must have a zero in $x_0$. Since $x_0$ is not a branch point of $\phi$, the discriminant
\begin{equation*}
\Delta_G=\prod_{i<j}(y_i-y_j)^2
\end{equation*}
does not have a zero at $x_0$, so neither do any of the $y_i$, $i\ne i_1$. Thus the zero of $y_{i_1}$ must be of order $\ge m$.

In particular, it follows from \mbox{Lemma \ref{lem-weierstrass-preparation}} that $\ord_r(y_{i_1})\ge m$ for all $r\ge\lambda(x_0)$. But we cannot have $\ord_{\lambda(x_0)}(y_i)\ge m$ for all $i\in\{1,\ldots,p\}$, since this would imply
\begin{equation*}
\ord_{\mu(x_0)}(g_0)\ge\ord_{\lambda(x_0)}(g_0)=\sum_{i=1}^p\ord_{\lambda(x_0)}(y_i)\ge mp,
\end{equation*}
contradicting the second part of \mbox{Assumption \ref{ass-combinatorial-condition}}. Thus there exists some $i_0\ne i_1$ with 
\begin{equation*}
\ord_r(y_{i_0})\le\ord_{\lambda(x_0)}(y_{i_0})<m
\end{equation*}
for all $r\ge\lambda(x_0)$ as desired.
\end{proof}

\begin{Lem}\label{lem-lambda-formula}
Suppose that \mbox{Assumption \ref{ass-combinatorial-condition}} is satisfied and let $r>\lambda(x_0)$ be a rational number. Then the Newton polygon $N_r(G)$ has a kink.
\end{Lem}
\begin{proof}
We consider the reductions
\begin{equation*}
\overline{y}_i\coloneqq[y_i]_r\in\kappa\llbracket\overline{t}\rrbracket,\qquad i=1,\ldots,p,
\end{equation*}
(\mbox{Definition \ref{def-order-and-reduction}}). By \mbox{Lemma \ref{lem-zeros-on-closed-disk}}, the $\overline{y}_i$ are polynomials, of degree equal to the number of zeros of $y_i$ in $\rD[r]$. Because of \mbox{Lemma \ref{lem-combinatorial-argument}}, there exist $i_0,i_1\in\{1,\ldots,p\}$ with $\deg(\overline{y}_{i_0})<\deg(\overline{y}_{i_1})$.

By \mbox{Lemma \ref{lem-generalized-bracket-properties}(a)}, the reduction of the discriminant $\Delta_G$ satisfies
\begin{equation*}
[\Delta_G]_r=\prod_{i<j}[y_i-y_j]_r^2.
\end{equation*}
It is a constant by \mbox{Remark \ref{rem-disc-is-constant}}. Because of \mbox{Lemma \ref{lem-zeros-on-closed-disk}}, all the $[y_i-y_j]_r$ are polynomials, so they must be constant too. Since $\overline{y}_{i_0}$ and $\overline{y}_{i_1}$ are distinct polynomials, one of which is not constant, it follows from \mbox{Lemma \ref{lem-generalized-bracket-properties}(b)} that $v_r(y_{i_0})\ne v_r(y_{i_1})$. Since the valuations $v_r(y_i)$, $i=1,\ldots,p$, are the negatives of the slopes of $N_G(r)$, we have shown that this Newton polygon must have a kink.
\end{proof}

\begin{Kor}\label{cor-separable-np}
Suppose that \mbox{Assumption \ref{ass-combinatorial-condition}} holds. Then the following are true:
\begin{enumerate}[(a)]
\item The Newton polygon $N_G(\lambda(x_0))$ is separable.
\item If $\lambda(x_0)>\mu(x_0)$, then $N_G(\lambda(x_0))$ is a straight line.
\end{enumerate}
\end{Kor}
\begin{proof}
By \mbox{Lemma \ref{lem-lambda-formula}}, the Newton polygon $N_r(G)$ has a kink for any $r>\lambda(x_0)$. By \mbox{Proposition \ref{prop-np-inseparable}}, $N_r(G)$ is inseparable for $\lambda(x_0)>r\ge\mu(x_0)$.

Thus both statements of the corollary are reduced to the following elementary facts about deformations of Newton polygons: If a Newton polygon $N_r(G)$ is inseparable, then so is $N_G(r+\epsilon)$ for $\abs{\epsilon}$ small enough. And if $N_G(r)$ has a kink, then so does $N_G(r+\epsilon)$ for $\abs{\epsilon}$ small enough.
\end{proof}

\begin{Lem}\label{lem-tail-case-separable-extension}
\sloppy Assume that $\lambda(x_0)>\mu(x_0)$ and that \mbox{Assumption \ref{ass-combinatorial-condition}} is satisfied. Then $\xi_{\lambda(x_0)}$ is a wild topological branch point, and in particular $\phi^{-1}(\xi_{\lambda(x_0)})=\{\eta\}$ for a point $\eta\in Y^{\an}$. 

Moreover, the extension $\kappa(\eta)/\kappa(\xi_{\lambda(x_0)})$ of residue fields is separable of degree $p$ with generator $\overline{y}\coloneqq[y]_{\eta}\in \kappa(\eta)$.
\end{Lem}
\begin{proof}
Because $\lambda(x_0)>\mu(x_0)$, \mbox{Proposition \ref{prop-simple-lambda-consequence}} shows that $\delta(\xi_r)>0$ for $\lambda(x_0)>r>\mu(x_0)$, while $\delta(\xi_{\lambda(x_0)})=0$. The $\xi_r$ for $\lambda(x_0)>r>\mu(x_0)$ are wild topological branch points by \mbox{Lemma \ref{lem-bijection-over-topological-branch-locus}}(a), and so is $\xi_{\lambda(x_0)}$ because the wild topological branch locus is closed (\mbox{Remark \ref{rem-topological-branch-locus}(c)}). As in the statement of the lemma, we denote the unique preimage of $\xi_{\lambda(x_0)}$ by $\eta$. Because $\delta(\xi_{\lambda(x_0)})=0$, \mbox{Lemma \ref{lem-bijection-over-topological-branch-locus}(b)} shows that $\CH(\eta)/\CH(\xi_{\lambda(x_0)})$ is an unramified extension of degree $p$, so that the extension $\kappa(\eta)/\kappa(\xi_{\lambda(x_0)})$ is indeed separable of degree $p$.

We are left to verify that $\overline{y}$ generates the field extension $\kappa(\eta)/\kappa(\xi_{\lambda(x_0)})$. Write $V\coloneqq\phi^{-1}(U)$ for the inverse image of the chart $U$. Recall from \mbox{Remark \ref{rem-valuations-preserved}} that we have described maps $\CO_Y(V)\to K\llbracket t\rrbracket$ that send $y$ to the power series $y_i$, where $i=1,\ldots,p$. By passing to quotient fields, these extend to embeddings of fields $F_Y\to K(\!( t)\!)$. Valuations are preserved in the sense that we have
\begin{equation*}
v_\eta(f)=\min_{i\ge0}\{v_K(a_i)+i\lambda(x_0)\}
\end{equation*}
for $f\in\CO_Y(V)$ with image $\sum_ia_it^i$ in $K\llbracket t\rrbracket$. Thus the map $F_Y\to K(\!( t)\!)$ descends to a map of residue fields $\kappa(\eta)\to \kappa(\!(\overline{t})\!)$ with $\overline{y}\mapsto\overline{y}_i\coloneqq[y_i]_{\lambda(x_0)}$.

Each such embedding defines an extension of the $\overline{t}$-adic valuation on $\kappa(\overline{t})$ to $\kappa(\eta)$. But \mbox{Assumption \ref{ass-combinatorial-condition}} guarantees that not all $\overline{y}_i$, $i=1,\ldots,p$, have the same $\overline{t}$-adic valuation (combine Lemmas \ref{lem-combinatorial-argument} and \ref{lem-weierstrass-preparation}). In particular, $\overline{y}$ cannot lie in $\kappa(\xi_{\lambda(x_0)})$. Since $\kappa(\eta)/\kappa(\xi_{\lambda(x_0)})$ is of degree $p$, we conclude that $\overline{y}$ generates this field extension.
\end{proof}

\begin{Rem}\label{rem-tail-case-delta-constant}
\begin{enumerate}[(a)]
\item Let us write down the minimal polynomial of the generator $\overline{y}=[y]_{\lambda(x_0)}$ considered in \mbox{Lemma \ref{lem-tail-case-separable-extension}}. It is
\begin{equation}\label{equ-tail-case-separable-extension-reduced-equation}
\overline{G}(T)=T^p+\overline{g}_{p-1}T^{p-1}+\ldots+\overline{g}_1T+\overline{g}_0=0,\qquad\overline{g}_i=\overline{g_i\pi^{-(p-i)s}},
\end{equation}
where $s=v_\eta(y)=v_{\lambda(x_0)}(g_0)/p$. To see this, note that $N_G(\lambda(x_0))$ being a straight line implies $v_{\lambda(x_0)}(g_i)\ge(p-i)s$ for $i=1,\ldots,p-1$. Thus it makes sense to define the reductions $\overline{g_i\pi^{-(p-i)s}}$. Moreover, it is clear that \eqref{equ-tail-case-separable-extension-reduced-equation} is a monic polynomial satisfied by $\overline{y}=\overline{y\pi^{-s}}$.

\item The discriminant $\Delta_{\overline{G}}$ of the minimal polynomial $\overline{G}$ of the element $\overline{y}$ is a constant. Indeed, we have
\begin{equation*}
\Delta_{\overline{G}}=\prod_{i<j}(\overline{y}_i-\overline{y}_j)^2,
\end{equation*}
where the $\overline{y}_i=[y_i]_{\lambda(x_0)}$, $i=1,\ldots,p$, are the reductions of the power series $y_i$, as considered in \mbox{Lemma \ref{lem-lambda-formula}}. Since $\kappa(\eta)/\kappa(\xi_{\lambda(x_0)})$ is separable, the $\overline{y}_i$ are pairwise distinct, so by \mbox{Lemma \ref{lem-generalized-bracket-properties}(b)} we have $\overline{y}_i-\overline{y}_j=[y_i-y_j]_{\lambda(x_0)}$ for all $i<j$. It follows from \mbox{Lemma \ref{lem-generalized-bracket-properties}(a)} that 
\begin{equation*}
\Delta_{\overline{G}}=\prod_{i<j}[y_i-y_j]^2_{\lambda(x_0)}=\big[\prod_{i<j}(y_i-y_j)^2\big]_{\lambda(x_0)}=[\Delta_G]_{\lambda(x_0)}
\end{equation*}
is a constant.
\end{enumerate}
\end{Rem}

\begin{Lem}\label{lem-lambda-formula-tail-case}Assume that $\lambda(x_0)>\mu(x_0)$ and that \mbox{Assumption \ref{ass-combinatorial-condition}} is satisfied. Then we have
\begin{equation*}
\ord_{\lambda(x_0)}(g_1)=0,
\end{equation*}
and the Newton polygon $N_G(\lambda(x_0))$ contains the point
\begin{equation*}
(1,v_{\lambda(x_0)}(g_1)).
\end{equation*}
\end{Lem}
\begin{proof}
We again write
\begin{equation*}
\overline{y}_i\coloneqq[y_i]_{\lambda(x_0)}\in\kappa\llbracket\overline{t}\rrbracket,\qquad i=1,\ldots,p.
\end{equation*}
It follows from \mbox{Corollary \ref{cor-separable-np}} that $N_G(\lambda(x_0))$ is a straight line, so we have
\begin{equation*}
v_{\lambda(x_0)}(y_1)=\ldots=v_{\lambda(x_0)}(y_p).
\end{equation*}
As in the proof of \mbox{Lemma \ref{lem-lambda-formula}}, the reduction of the discriminant $\Delta_G$,
\begin{equation*}
[\Delta_G]_{\lambda(x_0)}=\prod_{i<j}[y_i-y_j]^2_{\lambda(x_0)}
\end{equation*}
and the factors $[y_i-y_j]_{\lambda(x_0)}$ are constant. (In \mbox{Lemma \ref{lem-lambda-formula}}, we considered reduction at a radius $r>\lambda(x_0)$, but the only thing that matters for this argument is that the radius be greater than $\mu(x_0)$, which is true here as it was there.)

By \mbox{Lemma \ref{lem-combinatorial-argument}} and \mbox{Lemma \ref{lem-weierstrass-preparation}}, there exists an $i_1\in\{1,\ldots,p\}$ for which $\overline{y}_{i_1}$ vanishes at $\overline{t}=0$. We saw in Remark \ref{rem-tail-case-delta-constant}(b) that the $\overline{y}_i$ are pairwise distinct, so it follows from \mbox{Lemma \ref{lem-generalized-bracket-properties}(b)} that the $\overline{y}_i$, $i\ne i_1$, do not vanish at $\overline{t}=0$. Now \mbox{Lemma \ref{lem-weierstrass-preparation}} and \mbox{Lemma \ref{lem-order-is-derivative}} imply that the right derivatives of the functions $r\mapsto v_r(y_i)$ are zero at $r=\lambda(x_0)$ for $i\ne i_1$, while the right derivative of $r\mapsto v_r(y_{i_1})$ is positive at $r=\lambda(x_0)$. This readily implies the statements of the lemma.
\end{proof}

We now introduce notation for the coefficients of the power series $g_0,g_1,\ldots,g_{p-1}$, writing
\begin{equation*}
g_i=\sum_{k=0}^\infty a_{i,k}t^k,\qquad i=0,\ldots,p-1.
\end{equation*}
Our main theorems of this chapter, \mbox{Theorem \ref{thm-lambda-formula-tail-case}} below and \mbox{Theorem \ref{thm-lambda-formula-general-case}} in the following section, give formulas for $\lambda(x_0)$ only in terms of valuations of these coefficients, that is, in terms of certain elements of $K$ depending on $x_0$.

\begin{Satz}\label{thm-lambda-formula-tail-case}
If $\lambda(x_0)>\mu(x_0)$ and \mbox{Assumption \ref{ass-combinatorial-condition}} is satisfied for $m\ge1$, then we have
\begin{equation*}
\lambda(x_0)=\max_{m\le k<mp}\frac{pv_K(a_{1,0})-(p-1)v_K(a_{0,k})}{(p-1)k}.
\end{equation*}
\end{Satz}
\begin{proof}
By \mbox{Lemma \ref{lem-lambda-formula-tail-case}}, the Newton polygon $N_G(\lambda(x_0))$ is a straight line on which lies the vertex 
\begin{equation*}
(1,v_{\lambda(x_0)}(g_1))=(1,v_K(a_{1,0})).
\end{equation*}
Thus for every $k\ge0$ we have
\begin{equation*}
\frac{v_K(a_{1,0})}{p-1}=\frac{v_{\lambda(x_0)}(g_0)}{p}\le\frac{v_K(a_{0,k})+k\lambda(x_0)}{p}.
\end{equation*}
with equality for $k$ with $v_{\lambda(x_0)}(f_0)=v_K(a_{0,k})+k\lambda(x_0)$. This readily implies
\begin{equation*}
\lambda(x_0)\ge\frac{pv_K(a_{1,0})-(p-1)v_K(a_{0,k})}{(p-1)k}
\end{equation*}
for all $k\ge0$, again with equality for $k$ with $v_{\lambda(x_0)}(f_0)=v_K(a_{0,k})+k\lambda(x_0)$. To finish the proof, simply notice that because of \mbox{Assumption \ref{ass-combinatorial-condition}} we have
\begin{equation*}
m\le\ord_\infty(g_0)\le\ord_{\lambda(x_0)}(g_0)\le\ord_{\mu(x_0)}(g_0)<mp.
\end{equation*}
Thus
$v_{\lambda(x_0)}(g_0)=v_K(a_{0,k})+k\lambda$ for certain $k$ with $m\le k<mp$.
\end{proof}

\begin{Rem}\label{rem-lambda-formula-tail-case}
Let us record the following facts implicit in the proof of \mbox{Theorem \ref{thm-lambda-formula-tail-case}}:
\begin{enumerate}[(a)]
\item The maximum in the statement of \mbox{Theorem \ref{thm-lambda-formula-tail-case}} is achieved for precisely those $k\in\{m,\ldots,mp-1\}$ that satisfy $v_{\lambda(x_0)}(g_0)=v_K(a_{0,k})+k\lambda(x_0)$.
\item We have $v_{\lambda(x_0)}(g_1)=v_K(a_{1,0})$. This in contrast with the statement of \mbox{Theorem \ref{thm-lambda-formula-general-case}} below, where we may need to consider the minimum over various $v_K(a_{1,l})$, $l\ge0$.
\end{enumerate}
\end{Rem}

\section{An upper bound for $\lambda$: general case}

As announced in the previous section, the main result of this section, \mbox{Theorem \ref{thm-lambda-formula-general-case}}, is not dependent on the assumption $\lambda(x_0)>\mu(x_0)$. We retain the notation from the previous section; in particular $G$ still denotes the minimal polynomial of the generator $y$, the power series $g_0,\ldots,g_{p-1}$ are the coefficients of $G$, and the coefficients of the $g_i$ are denoted $a_{i,k}$, $k\ge0$.

\begin{Lem}\label{lem-theta-derivative}
Suppose that \mbox{Assumption \ref{ass-combinatorial-condition}} holds for some $m\ge1$. Then for any $r\ge\lambda(x_0)$, we have the following:
\begin{enumerate}[(a)]
\item The right derivative of the function $\theta_G$ (\mbox{Definition \ref{def-theta}}) satisfies
\begin{equation*}
\big(\partial_+\theta_G\big)(r)\ge m.
\end{equation*}
\item The function $\beta\colon[\lambda(x_0),\infty)$ given by $r\mapsto\theta_G(r)-\frac{v_r(g_0)}{p}$ is monotonically increasing, and has monotonically decreasing right derivative.
\end{enumerate}
\end{Lem}
\begin{proof}
By the same argument as in the proof of \mbox{Lemma \ref{lem-combinatorial-argument}} we see that one of the $y_1,\ldots,y_p$, say $y_{i_1}$, has a zero of order $\ge m$ in $x_0$, while the $y_i$, $i\ne i_1$, do not have a zero in $x_0$. By \mbox{Lemma \ref{lem-order-is-derivative}}, the right derivative at $r\ge\lambda(x_0)$ of $v_r(y_i)$ is the same as the order $\ord_r(y_i)$, for $i=1,\ldots,p$. Thus we have $v_r(y_{i_1})>v_r(y_i)$ and
\begin{equation*}
\big(\partial_+\theta_G\big)(r)=\big(\partial_+v_r(y_{i_1})\big)(r)\ge m
\end{equation*} 
for $i\ne i_1$ and $r\gg\lambda(x_0)$. To prove (a), we show that $v_r(y_{i_1})\ge v_r(y_i)$, where $i\ne i_1$, holds for \emph{all} $r\ge\lambda(x_0)$. 

If it were not true that $v_r(y_{i_1})\ge v_r(y_i)$ for all $r\ge\lambda(x_0)$, we could find $r>\lambda(x_0)$ and $i_2\ne i_1$ with
\begin{equation*}
v_r(y_{i_1})=v_r(y_{i_2}),\qquad(\partial_+v_r(y_{i_2}))(r)=\ord_r(y_{i_2})<(\partial_+v_r(y_{i_1}))(r)=\ord_r(y_{i_1}).
\end{equation*}
The reductions $\overline{y}_{i_1}\coloneqq[y_{i_1}]_r$ and $\overline{y}_{i_2}\coloneqq[y_{i_2}]_r$ are polynomials by \mbox{Lemma \ref{lem-zeros-on-closed-disk}}. By \mbox{Lemma \ref{lem-weierstrass-preparation}}, $\overline{y}_{i_1}$ has larger order of vanishing at $\overline{t}=0$ than $\overline{y}_{i_2}$. As in the proof of \mbox{Lemma \ref{lem-lambda-formula}} we conclude from the fact that
\begin{equation*}
[\Delta_G]_r=\prod_{i<j}[y_i-y_j]_r^2
\end{equation*}
is a unit that $v_r(y_{i_1})\ne v_r(y_{i_2})$. We have arrived at a contradiction, which finishes the proof of (a).

To prove (b), note that we have
\begin{equation*}
\beta(r)=\theta_G(r)-\frac{v_r(g_0)}{p}=\frac{(p-1)v_r(y_{i_1})-\sum_{i\ne i_1}v_r(y_i)}{p}.
\end{equation*}
We just saw that $v_r(y_{i_1})>v_r(y_i)$ for $i\ne i_1$ and $r>\lambda(x_0)$. Hence 
\begin{equation*}
[y_i]_r=[y_{i_1}-y_i]_r
\end{equation*}
must be constant for $r>\lambda(x_0)$ and $i\ne i_1$. Thus $v_r(y_i)$, where $i\ne i_1$, is constant for $r>\lambda(x_0)$, so that
\begin{equation*}
\big(\partial_+\beta\big)(r)=\frac{p-1}{p}\ord_r(y_{i_1})
\end{equation*}
is positive, but monotonically decreasing as a function in $r$. This proves (b).
\end{proof}

\begin{Satz}\label{thm-lambda-formula-general-case}
Assume that \mbox{Assumption \ref{ass-combinatorial-condition}} is satisfied for $m\ge1$. Then we have
\begin{equation*}
\lambda(x_0)=\max\{\mu(x_0),\tilde{\lambda}\},
\end{equation*}
where
\begin{equation*}
\tilde{\lambda}=\max_{m\le k<mp}\min_{0<i<p}\min_{0\le pl<(p-i)k}
\frac{pv_K(a_{i,l})-(p-i)v_K(a_{0,k})}{(p-i)k-pl}.
\end{equation*}
\end{Satz}
\begin{proof}
Suppose that $k\in\{m,\ldots,mp-1\}$. In this proof, we will use the terminology of calling a pair $(i,l)\in\BZ^2$ a \emph{$k$-admissible} pair if $i\in\{1,\ldots,p-1\}$, $0\le pl<(p-i)k$, and $a_{i,l}\ne0$. The significance of this is that for every $k$-admissible pair $(i,l)$, the two affine functions
\begin{equation*}
r\mapsto\frac{v_K(a_{i,l})+lr}{p-i},\qquad r\mapsto\frac{v_K(a_{0,k})+kr}{p}
\end{equation*}
coincide for a unique $r(k,i,l)\in\BQ$, namely
\begin{equation*}
r(k,i,l)=\frac{pv_K(a_{i,l})-(p-i)v_K(a_{0,k})}{(p-i)k-pl}.
\end{equation*}
(If $a_{0,k}=0$, then $r(k,i,l)=-\infty$.) For $r<r(k,i,l)$, the first function has larger value than the second, and for values $r>r(k,i,l)$, it has smaller value. Said differently, for any $r\ge\lambda(x_0)$ we have
\begin{equation}\label{equ-magical-equivalence-1}
\frac{v_K(a_{i,l})+lr}{p-i}\le\frac{v_K(a_{0,k})+kr}{p}
\end{equation}
if and only if
\begin{equation}\label{equ-magical-equivalence-2}
r(k,i,l)=\frac{pv_K(a_{i,l})-(p-i)v_K(a_{0,k})}{(p-i)k-pl}\le r,
\end{equation}
and similarly if we replace the two ``$\le$'' with ``$\ge$''.

To prove the theorem it is enough to show
\begin{equation}\label{equ-k-by-k}
\min\{r_{k,i,l}\mid\textrm{$(i,l)$ is $k$-admissible}\}\le\lambda(x_0)
\end{equation}
for all $k\in\{m,\ldots,mp-1\}$. Indeed, this implies $\tilde{\lambda}\le\lambda$ from which follows the statement of the theorem in case $\lambda(x_0)=\mu(x_0)$. And in case $\lambda(x_0)>\mu(x_0)$, it follows from Theorem \ref{thm-lambda-formula-tail-case} that the Newton polygon $N_G(\lambda(x_0))$ is a straight line. The equivalence of \eqref{equ-magical-equivalence-1} and \eqref{equ-magical-equivalence-2} (with ``$\ge$'' instead of ``$\le$'') then shows $\tilde{\lambda}\ge\lambda(x_0)$. Together with the inequality $\tilde{\lambda}\le\lambda(x_0)$, this shows $\tilde{\lambda}=\lambda(x_0)=\max\{\lambda(x_0),\mu(x_0)\}$, as desired.

Before proving that \eqref{equ-k-by-k} holds for all $k\in\{m,\ldots,mp-1\}$, we need another piece of preparation. Let $s\ge\lambda(x_0)$ be a real number. We will attach to $s$ an element $k_s\in\{m,\ldots,mp-1\}$ and a $k_s$-admissible pair $(i_s,l_s)$ as follows. For $i\in\{1,\ldots,p-1\}$, the point $(i,v_{s}(g_{i}))$ lies on the steepest slope of the Newton polygon $N_G(s)$ if and only if 
\begin{equation}\label{equ-lambda-theta}
\theta_G(s)=\max_{1\le j\le p}\frac{v_{s}(g_0)-v_{s}(g_{j})}{j}=\frac{v_{s}(g_0)-v_{s}(g_{i})}{i}.
\end{equation}
There exists at least one $i\in\{1,\ldots,p-1\}$ with this property because of \mbox{Lemma \ref{lem-lambda-formula}}. 

Among those $i\in\{1,\ldots,p-1\}$ such that $(i,v_{s}(g_{i}))$ lies on the steepest slope of the Newton polygon $N_G(s)$, choose one for which the quotient $\ord_{s}(g_{i})/i$ is minimal, and denote it by $i_s$. \mbox{Lemma \ref{lem-order-is-derivative}} shows that $\ord_{s}(g_{i_s})$ is the right derivative of $v_s(g_{i_s})$ at $s$. Thus for $r=s+\epsilon$ with $\epsilon>0$ small enough, we still have
\begin{equation}\label{equ-just-differentiate-bro}
\theta_G(r)=\max_{1\le j\le p}\frac{v_r(g_0)-v_r(g_{j})}{j}=\frac{v_r(g_0)-v_r(g_{i_s})}{i_s},
\end{equation}
that is, $(i_s,v_r(g_{i_s}))$ lies on the steepest slope of $N_r(G)$. Write $k_s\coloneqq\ord_{s}(g_0)$ and $l_s\coloneqq\ord_{s}(g_{i_s})$. Then we have
\begin{equation*}
m\le\ord_\infty(g_0)\le k_s\le\ord_{\mu(x_0)}(g_0)<mp.
\end{equation*}
It follows from \mbox{Lemma \ref{lem-theta-derivative}}(a) and from differentiating \eqref{equ-just-differentiate-bro} that
\begin{equation*}
\frac{k_s-l_s}{i_s}=\big(\partial_+\theta_G\big)(s)\ge m>\frac{k_s}{p}.
\end{equation*}
Rearranging this gives
\begin{equation*}
(p-i_s)k_s>pl_s,
\end{equation*}
that is, $(i_s,l_s)$ is $k_s$-admissible.
Next, let us denote by $r_0,r_1,r_2,\ldots$ the kinks of the piecewise affine function
\begin{equation*}
[\lambda(x_0),\infty)\to\BR,\qquad r\mapsto v_r(g_0),
\end{equation*}
where we consider $r_0=\lambda(x_0)$ as the first kink. Suppose that $s$ is one of these kinks. We now show \eqref{equ-k-by-k} for $k=k_s$. In the following estimation, the second ``$\le$'' is true because of Lemma \ref{lem-theta-derivative}(b), and the first ``$\le$'' is true because the negative of the steepest slope $\theta_G(\lambda(x_0))$ of $N_G(\lambda(x_0))$ must be greater than $v_{\lambda(x_0)}(g_0)/p$. We have:
\begingroup
\allowdisplaybreaks[1]
\begin{align*}
0&\le\frac{i_s}{p-i_s}\Big(\theta_G(\lambda(x_0))-\frac{v_{\lambda(x_0)}(g_0)}{p}\Big)\\
&=\frac{i_s}{p-i_s}\Big(\theta_G(s)-\frac{v_{s}(g_0)}{p}-\int_{\lambda(x_0)}^s(\partial_+\theta_G)(r)-\frac{\ord_r(g_0)}{p}\,dr\Big)\\
&\le\frac{i_s}{p-i_s}\Big(\theta_G(s)-\frac{v_{s}(g_0)}{p}-\int_{\lambda(x_0)}^s(\partial_+\theta_G)(s)-\frac{\ord_s(g_0)}{p}\,dr\Big)\\
&=\frac{i_s}{p-i_s}\Big(\frac{v_s(g_0)-v_s(g_{i_s})}{i_s}-\frac{v_s(g_0)}{p}-(s-\lambda(x_0))\big(\frac{k_s-l_s}{i_s}-\frac{k_s}{p}\big)\Big)\\
&=\frac{v_s(g_0)}{p}-\frac{v_s(g_{i_s})}{p-i_s}-(s-\lambda(x_0))\big(\frac{k_s}{p}-\frac{l_s}{p-i_s}\big)\\
&=\frac{v_K(a_{0,{k_s}})+k_ss}{p}-\frac{v_K(a_{i_s,l_s})+l_ss}{p-i_s}-(s-\lambda(x_0))\big(\frac{k_s}{p}-\frac{l_s}{p-i_s}\big)\\
&=\frac{v_K(a_{0,k_s})+k_s\lambda(x_0)}{p}-\frac{v_K(a_{i_s,l_s})+l_s\lambda(x_0)}{p-i_s}
\end{align*}%
\endgroup
Hence from the equivalence of \eqref{equ-magical-equivalence-1} and \eqref{equ-magical-equivalence-2} we get $\lambda(x_0)\ge r(k_s,i_s,l_s)$, which proves \eqref{equ-k-by-k} for $k_s$.

Finally, suppose that $k\in\{m,\ldots,mp-1\}$ is not of the form $k_s$, where $s$ is one of the kinks $r_0,r_1,r_2,\ldots$. Then there is a kink $s$ for which $k\ge k_s$, and so $(i_s,l_s)$ being $k_s$-admissible is also $k$-admissible. We have
\begin{equation*}
v_K(a_{0,k})+kr\ge v_r(g_0)\ge v_K(a_{0,k_s})+k_sr
\end{equation*}
for all $r\ge\lambda(x_0)$, so that
$r(k,i_s,l_s)\le r(k_s,i_s,l_s)$. It follows that \eqref{equ-k-by-k} holds for $k$, finishing the proof of the theorem.
\end{proof}

We end this section with an example illustrating that \mbox{Theorem \ref{thm-lambda-formula-tail-case}} is wrong without the assumption $\lambda(x_0)>\mu(x_0)$.

\begin{Ex}
Consider the plane quartic curve
\begin{equation*}
Y\colon\quad y^3 + y^2 + (3x^3 + x^2 + 3)y + x^4 - 3x^2 = 0
\end{equation*}
over $K=\BC_3$ and the degree-$3$ morphism $\phi\colon Y\to X$ given birationally by $(x,y)\mapsto x$. As in \mbox{Example \ref{ex-lambda-as-radius-of-convergence}}, we take the complement of the branch locus of $\phi$ for the chart $U$. Thus $\Gamma_0$ is the tree spanned by the branch points of $\phi$, which include $\infty$.

We fix the pair of closed points $p_0=(0,0)$, $x_0\coloneqq\phi(p_0)=0$. Observe that \mbox{Assumption \ref{ass-combinatorial-condition}} is satisfied for $m=2$. The non-zero coefficients used to define the quantity $\tilde{\lambda}$ from Theorem \ref{thm-lambda-formula-general-case} are
\begin{equation*}
a_{0,2}=-3,\quad a_{0,4}=1,\qquad a_{1,0}=3,\quad a_{1,2}=1,\quad a_{1,3}=3,\qquad a_{2,0}=1.
\end{equation*}
We first consider $k=2$. The $2$-admissible pairs (see the proof of \mbox{Theorem \ref{thm-lambda-formula-general-case}} for this terminology) are $(i,l)=(1,0),(2,0)$. We have
\begin{equation*}
\min_{0<i<3}\min_{0\le 3l<2(3-i)}\Big\{\frac{3v_K(a_{i,l})-(3-i)v_K(a_{0,2})}{(3-i)\cdot2-3l}\Big\}=\min\Big\{\frac{1}{4},-\frac{1}{2}\Big\}=-\frac{1}{2}.
\end{equation*}
Next, we consider $k=4$. The $4$-admissible pairs are $(i,l)=(1,0),(1,2),(2,0)$. We have
\begin{equation*}
\min_{0<i<3}\min_{0\le 3l<4(3-i)}\Big\{\frac{3v_K(a_{i,l})-(3-i)v_K(a_{0,4})}{(3-i)\cdot4-3l}\Big\}=\min\Big\{\frac{3}{8},0,0\Big\}=0.
\end{equation*}
It follows that $\tilde{\lambda}=\max\{-1/2,0\}=0$. The Newton polygon of the discriminant of the polynomial $y^3 + y^2 + (3x^3 + x^2 + 3)y + x^4 - 3x^2$ defining $Y$ shows that $\phi$ has a branch point of valuation $1/3$. Thus $\mu(x_0)=1/3$ and it follows from Theorem \ref{thm-lambda-formula-general-case} that
\begin{equation*}
\lambda(x_0)=\max\{\mu(x_0),\tilde{\lambda}\}=\mu(x_0)=\frac{1}{3}.
\end{equation*}
The assumption $\lambda(x_0)>\mu(x_0)$ of Theorem \ref{thm-lambda-formula-tail-case} is not satisfied, and indeed we have
\begin{equation*}
\max_{k\in\{2,4\}}\Big\{\frac{3v_K(a_{1,0})-2v_K(a_{0,k})}{2k}\Big\}=\max\Big\{-\frac{1}{2},\frac{3}{8}\Big\}=\frac{3}{8},
\end{equation*}
which is strictly larger than $\lambda(x_0)$.


\end{Ex}

\section{Guaranteeing that Assumption \ref{ass-combinatorial-condition} is satisfied}\label{sec-analytic-good-equation}

The main results of this chapter, \mbox{Theorem \ref{thm-lambda-formula-tail-case}} and \mbox{Theorem \ref{thm-lambda-formula-general-case}} only hold under the additional \mbox{Assumption \ref{ass-combinatorial-condition}} on the generator $y$ that is part of the \'etale chart $(U,y)$ fixed in \mbox{Section \ref{sec-analytic-preliminary}}. In this section we explain how to modify $y$ so that \mbox{Assumption \ref{ass-combinatorial-condition}} is guaranteed to be satisfied.

Write $A\coloneqq\CO_X(U)$, where $U$ is the open subscheme that is part of the chart $(U,y)$. A non-zero function $f\in A$ may be regarded as a rational function on $X$. Let us denote by $v_P$ the valuation on $K(x)$ attached to a closed point $P\in X$. If $v_P(f)=m<0$, then $f$ is said to have a \emph{pole} of order $m$ at $P$.

The valuations $v_P$ should not be confused with the pseudovaluations attached to closed points that we introduced in \mbox{Section \ref{sec-preliminaries-analytic}}. The latter are not actually valuations (taking the values $\pm\infty$ on various functions) and extend the valuation $v_K$, while the $v_P$ are honest valuations that are trivial on $K$.

\begin{Def}
Suppose that $f$ is a non-zero regular function on an affine open $W\subset X$ containing $U$, that is, $f\in\CO_X(W)$. The \emph{degree} of $f$ on $W$ is defined to be the product
\begin{equation*}
\deg_W(f)=\#(X\setminus W)\cdot\big( \textrm{highest order of a pole of $f$}\big).
\end{equation*} 
\end{Def}

\begin{Ex}
Consider the elliptic curve $E$ over $K=\BC_2$ with Weierstrass equation
\begin{equation*}
y^2=x^3-x.
\end{equation*}
Then $(U_1,y)$ is an affine \'etale chart for $Y$, where $U_1=\BP_K^1\setminus\{x=\infty,0,1,-1\}$ (see \mbox{Example \ref{ex-hyperelliptic-chart}}). The only pole of the function $f=x^3-x$ is $\infty$. It is a pole of order $3$, so we have $\deg_{W_1}(f)=3$, where $W_1=\BP^1_K\setminus\{x=\infty\}$. The curve $E$ may also be described by the equation
\begin{equation*}
z^2=\frac{u^3}{1-u^2},\qquad\textrm{where}\qquad z=1/y,\quad u=1/x.
\end{equation*}
Thus $(U_2,z)$ is an affine \'etale chart for $E$, where $U_2=\BP_K^1\setminus\{u=\infty,0,1,-1\}$. The function $g=u^3/(1-u^2)$ has three poles, $1$, $-1$, and $\infty$, each of order $1$. Thus we find $\deg_{W_2}(g)=3$ as well, where $W_2=\BP_K^1\setminus\{u=1,-1,\infty\}$.
\end{Ex}

\begin{Rem}\label{rem-polynomial-degree}
In \mbox{Chapter \ref{cha-quartics}}, all functions $f\in A\subset K(x)$ whose degree we consider will actually be polynomials in $x$. Thus they are regular functions on $W\coloneqq\BA_K^1=\BP_K^1\setminus\{\infty\}$, with at most a pole at $\infty$. In this case, the function $\deg_W$ coincides with the usual degree function on the polynomial ring $K[x]$.

It is also in this light that the properties proved in the following lemma should be viewed.
\end{Rem}

\begin{Lem}\label{lem-degree-properties}
For $f,g\in \CO_X(W)\setminus\{0\}$ we have:
\begin{enumerate}[(a)]
\item $\deg_W(fg)\le\deg_W(f)+\deg_W(g)$
\item If $f+g\ne0$, then $\deg_W(f+g)\le\max\big(\deg_W(f),\deg_W(g)\big)$
\item If $\rD\subset W^{\an}$ is an open disk, then $\ord_\rD(f)\le\deg_W(f)$
\end{enumerate}
\end{Lem}
\begin{proof}
\begin{enumerate}[(a)]
\item Consider a closed point $P\in X\setminus W$. We have $v_P(fg)=v_P(f)+v_P(g)$, so the highest possible order of a pole of $fg$ at $P$ is the order of a pole of $f$ at $P$ plus the order of a pole of $g$ at $P$. This immediately implies (a).
\item The inequality $v_P(f+g)\ge\min(v_P(f),v_P(g))$ shows that the order of a pole of $f+g$ at $P$ can not exceed the order of a pole of $f$ at $P$ nor the order of a pole of $g$ at $P$. This immediately implies (b).
\item The function $f$ can have no more than $\deg_W(f)$ zeros on $W$, counted with multiplicity (since the degree of a principal divisor is $0$, \cite[\mbox{Corollary II.6.10}]{hartshorne}). Thus (c) results from \mbox{Lemma \ref{lem-weierstrass-preparation}}.
\end{enumerate}
\end{proof}

As in \mbox{Section \ref{sec-analytic-preliminary}}, write
\begin{equation*}
G(T)=T^p+g_{p-1}T^{p-1}+\ldots+g_1T+g_0
\end{equation*}
for the minimal polynomial of $y$ over the subfield $K(t)$, where $t=x-x_0$ is the parameter associated to the closed point $x_0\in U$. We have studied the formal solutions $y_1,\ldots,y_p$ to the equation $G(T)=0$; they are elements of $K\llbracket t\rrbracket$. Let us now choose and fix one of these, say $y_1$, and write
\begin{equation*}
y_1=\sum_{i=0}^\infty b_it^i\in K\llbracket t\rrbracket.
\end{equation*}
Let $m\ge1$ be an integer. We denote by $u_m$ the truncation of $y_1$ consisting of the first $m$ monomials,
\begin{equation*}
u_m=b_0+b_1t+\ldots+b_{m-1}t^{m-1},
\end{equation*}
and consider the generator $w\coloneqq y-u_m$ of the field extension $F_Y/F_X$. Let us write
\begin{equation}\label{equ-minimal-polynomial-of-good-equation}
H(T)=\sum_{i=0}^ph_iT^i=G(T+u_m)
\end{equation}
for the minimal polynomial of $w$. We call the integer $m$ the \emph{order of approximation} of $w$. The following lemma will help us decide what order of approximation to choose.

\begin{Lem}\label{lem-sufficient-approximation}
Let $W\subseteq\BA_K^1$ be an affine open subscheme containing $U$ such that $g_0,g_1,\ldots,g_{p-1}\in\CO_X(W)$. If the integer $m$ satisfies
\begin{equation}\label{equ-condition-for-m}
m>\frac{\deg_W(g_i)-i}{p-i}
\end{equation}
for all $0\le i\le p-1$, then
\begin{equation*}
\ord_{\rD(\mu(x_0))}(h_0)<mp.
\end{equation*}
\end{Lem}
\begin{proof}
We have
\begin{equation*}
h_0=G(u_m)=u_m^p+g_{p-1}u_m^{p-1}+\ldots+g_1u_m+g_0.
\end{equation*}
It is convenient, to also set $g_p\coloneqq1$. Since the $g_0,\ldots,g_p$ and $u_m$ lie in $\CO_X(W)$ (the latter even in $K[t]=\CO_X(\BA_K^1)$), we obtain the desired estimate as follows:
\begin{IEEEeqnarray*}{rClr}\ord_{\rD(\mu(x_0))}(h_0)&\le&\deg_W(h_0)&\qquad\qquad\qquad\textrm{(by Lemma \ref{lem-degree-properties}(c))}\\
&\le&\max_{0\le i\le p}\{\deg_W(u_m^ig_i)\}&\textrm{(by Lemma \ref{lem-degree-properties}(a))}\\
&\le&\max_{0\le i\le p}\{i(m-1)+\deg_W(g_i)\}&\textrm{(by Lemma \ref{lem-degree-properties}(b)}\\
&<&mp&\textrm{(by \eqref{equ-condition-for-m})}
\end{IEEEeqnarray*}
\end{proof}

\begin{Kor}\label{cor-good-equation}
Suppose that there exist an affine open subscheme $W\subseteq\BA_K^1$ containing $U$ with $g_0,g_1,\ldots,g_{p-1}\in\CO_X(W)$ and an integer $m\ge1$ satisfying \eqref{equ-condition-for-m}. Then $(U,w)$ is a chart satisfying \mbox{Assumption \ref{ass-combinatorial-condition}} for this integer $m$.
\end{Kor}
\begin{proof}
Since $(U,y)$ is an affine \'etale chart, and $w=y-u_m$ for a function $u_m\in\CO_X(U)$, it follows that $(U,w)$ is an affine \'etale chart as well.

Since $u_m$ is the truncation of a formal solution of the minimal polynomial $G$, it follows that $h_0=G(u_m)$ has a zero of order at least $m$ at $t=0$. Thus the first inequality in \mbox{Assumption \ref{ass-combinatorial-condition}} holds. The second inequality holds by \mbox{Lemma \ref{lem-sufficient-approximation}}.
\end{proof}

Recall the discussion surrounding \mbox{Lemma \ref{lem-implicit-function-theorem}} in \mbox{Section \ref{sec-analytic-preliminary}}. We considered a ``generic'' version of the minimal polynomial of $y$,
\begin{equation*}
\tilde{F}(T)=T^p+\Tilde{f}_{p-1}T^{p-1}+\ldots+\Tilde{f}_1T+\Tilde{f}_0\in A\llbracket t\rrbracket[T],
\end{equation*}
from which the minimal polynomial $G$ of $y$ over the subfield $K(t)$ (where $t=x-x_0)$ is obtained by evaluating the coefficients at $x_0$. Let us also now denote by
\begin{equation*}
\Tilde{y}=y+b_1t+b_2t^2+b_3t^3+\ldots\in B\llbracket t\rrbracket
\end{equation*}
the ``generic'' solution to $\tilde{F}$, from which $y_1$ is obtained by evaluation at the corresponding point $p_0$ lying above $x_0$. Using these ingredients, we may also define a ``generic'' version of the minimal polynomial $H(T)$ of $w$. Namely, first take the truncation of $\tilde{y}$ of level $m$,
\begin{equation*}
\tilde{u}_m\coloneqq y+b_1t+\ldots+b_{m-1}t^{m-1}\in B[t],
\end{equation*}
and then define
\begin{equation*}
\tilde{H}(T)=\sum_{i=0}^p\tilde{h}_iT^i\coloneqq\tilde{F}(\tilde{u}_m+T).
\end{equation*}
Now $H(T)$ is obtained from $\tilde{H}(T)$ by evaluating at $x_0$. That is, if we write
\begin{equation*}
\tilde{h}_i=\sum_{l=0}^\infty c_{i,l}t^l,\qquad c_{i,l}\in B,
\end{equation*}
then we have
\begin{equation*}
h_i=\sum_{l=0}^\infty c_{i,l}(p_0)t^l.
\end{equation*}
As usual, we denote by $\hat{c}_{i,l}$ the associated valuative function $Y^{\an}\to\BR\cup\{\pm\infty\}$. We may then rephrase \mbox{Theorems \ref{thm-lambda-formula-tail-case}} \mbox{and \ref{thm-lambda-formula-general-case}}, as follows:

\begin{Kor}\label{cor-lambda-formula-tail-case}
If $\lambda(x_0)>\mu(x_0)$, then we have
\begin{equation*}
\lambda(x_0)=\max_{m\le k<mp}\frac{p\hat{c}_{1,0}(p_0)-(p-1)\hat{c}_{0,k}(p_0)}{(p-1)k}.
\end{equation*}
\end{Kor}

\begin{Kor}\label{cor-lambda-formula-general-case}
We have $\lambda(x_0)=\max\{\mu(x_0),\tilde{\lambda}(p_0)\}$, where
\begin{equation*}
\tilde{\lambda}(p_0)=\max_{m\le k<mp}\min_{0<i<p}\min_{0\le pl<(p-i)k}\frac{p\hat{c}_{i,l}(p_0)-(p-i)\hat{c}_{0,k}(p_0)}{(p-i)k-pl}.
\end{equation*}
\end{Kor}
Here $p_0$ is any closed point in $\phi^{-1}(x_0)$.

It should be stressed that we just took an important conceptual step. Where \mbox{Theorems \ref{thm-lambda-formula-tail-case}} \mbox{and \ref{thm-lambda-formula-general-case}} were concerned with describing the value of $\lambda$ at a fixed point $x_0$ in terms of certain constants derived from $x_0$, the preceding two corollaries describe the value of $\lambda$ at \emph{any} point $x_0$ in terms of valuative functions independent of $x_0$. Crucially, the functions $c_{i,l}$ yield an equation satisfying Assumption \ref{ass-combinatorial-condition} with respect to \emph{any} closed point $x_0\in U$, by evaluation at a closed point $p_0$ above $x_0$.

\begin{Ex}\label{ex-good-equation-superelliptic}
Consider a superelliptic curve of degree $p$, say
\begin{equation*}
Y\colon\quad F(y)=y^p-f=0,\qquad f=\sum_{i=0}^da_ix^i\in K[x].
\end{equation*}
The corresponding cover $Y\to X=\BP_K^1$ is Galois, since $K$ is algebraically closed. We take the chart $U=\Spec A$, where $A=K[x]_f$ (cf.\ \mbox{Example \ref{ex-hyperelliptic-chart}}). The ``generic'' minimal polynomial of $y$ is
\begin{equation*}
\tilde{F}(T)=T^p-\sum_{i=0}^d\frac{f^{(i)}}{i!}t^i\in A\llbracket t\rrbracket[T].
\end{equation*}
It will be useful to divide the constant coefficient of $\tilde{F}$ by $f\in A^\times$, yielding
\begin{equation*}
\tilde{G}(T)=T^p-\sum_{i=0}^d\frac{f^{(i)}}{i!f}t^i\in A\llbracket t\rrbracket[T].
\end{equation*}
Indeed, this has the effect that there is a formal solution
\begin{equation*}
\tilde{y}=1+b_1t+b_2t^2+b_3t^3+\ldots\in A\llbracket t\rrbracket
\end{equation*}
to $\tilde{G}$ provided by \mbox{Lemma \ref{lem-implicit-function-theorem}} beginning with the constant $1$. It is crucial that $\tilde{y}\in A\llbracket t\rrbracket$ --- compare the general case of \eqref{equ-general-formal-solution}, where $\tilde{y}\in B\llbracket t\rrbracket$. Following \mbox{Lemma \ref{lem-sufficient-approximation}}, we define
\begin{equation*}
m\coloneqq\Big\lfloor\frac{\deg(f)}{p}\Big\rfloor+1
\end{equation*}
and write
\begin{equation*}
\tilde{u}_m=1+b_1t+\ldots+b_{m-1}t^{m-1}\in A[t]
\end{equation*}
for the truncation of $\tilde{y}$ of order $m$. Moreover, as usual we write
\begin{equation*}
\tilde{H}(T)\coloneqq\tilde{G}(\tilde{u}_m+T)=T^p+p\tilde{u}_mT^{p-1}+\ldots+p\tilde{u}_m^{p-1}T+\tilde{u}_m^p-\sum_{i=0}^d\frac{f^{(i)}}{i!f}t^i,
\end{equation*}
and for a fixed closed point $x_0$ consider the polynomial obtained from $\tilde{H}$ by evaluating coefficients at $x_0$,
\begin{equation*}
H(T)=T^p+pu_mT^{p-1}+\ldots+pu_m^{p-1}T+u_m^p-\sum_{i=0}^d\frac{f^{(i)}(x_0)}{i!f(x_0)}t^i,
\end{equation*}
\begin{equation*}
\textrm{where}\quad u_m=1+b_1(x_0)t+\ldots+b_{m-1}(x_0)t^{m-1}.
\end{equation*}
The generator $w$ that is a solution to $H(T)$ is related to the original generator $y$ by 
\begin{equation*}
w=\frac{y}{f(x_0)^{1/p}}-u_m.
\end{equation*}
It follows from \mbox{Corollary \ref{cor-good-equation}} that the chart $(U,w)$ satisfies \mbox{Assumption \ref{ass-combinatorial-condition}}.

\end{Ex}

In a certain sense we have now reached our original goal of determining the semistable reduction of curves $Y$ equipped with a degree-$p$ morphism $\phi\colon Y\to X$. \mbox{Corollary \ref{cor-lambda-formula-general-case}} completely determines the function $\lambda=\lambda_{\phi,\Gamma_0}$. And according to \mbox{Theorem \ref{thm-trivializing-lambda}}, if we can enlarge the skeleton $\Gamma_0$ of $X$ to a skeleton $\Gamma$ trivializing $\lambda$, then $\Sigma\coloneqq\phi^{-1}(\Gamma)$ is a skeleton of $Y$. Or, using the language of models, if $\CX$ is the semistable $\CO_K$-model of $X$ obtained from $\Gamma$ as in \mbox{Remark \ref{rem-models-and-skeletons}}, then the normalization $\CY$ of $\CX$ in the function field $F_Y$ is a semistable model of $\CY$.

However, this solution is not quite satisfactory for two reasons:
\begin{itemize}
\item The general formula for $\lambda$ of \mbox{Corollary \ref{cor-lambda-formula-general-case}} is quite unwieldy, involving many terms.
\item Our formulas for $\lambda$ are in terms of admissible functions on $Y^{\an}$. For working with concrete examples and for implementation it is much more advantageous to work with admissible functions on $X^{\an}$.
\end{itemize}

In the next section, we explain an approach bypassing these issues. When studying plane quartics in the next chapter, we will only appeal to the simpler formula for $\lambda$ of \mbox{Corollary \ref{cor-lambda-formula-tail-case}}. And we will see how to describe $\lambda$ using admissible functions on $X^{\an}$.

\section{Tame locus}\label{sec-analytic-tame-locus}

We denote by $\Sigma$ the $\phi$-minimal skeleton of $Y$ with respect to $\Gamma_0$ (\mbox{Definition \ref{def-phi-minimal-skeleton}}). It is minimal among all skeletons of $Y$ containing the inverse image $\phi^{-1}(\Gamma_0)$. Attached to $\Sigma$ is the canonical retraction map of \mbox{Definition \ref{def-retraction}},
\begin{equation*}
\ret_\Sigma\colon Y^{\an}\to\Sigma.
\end{equation*}

\begin{Def}\label{def-tame-locus}
The \emph{tame locus} associated to $\phi$ and $\Gamma_0$ is the subset
\begin{equation*}
Y^{\tame}\coloneqq Y^{\tame}_{\phi,\Gamma_0}\coloneqq \ret_\Sigma^{-1}\big(\{\eta\in \Sigma\mid\delta(\eta)=0\}\big).
\end{equation*}
If $\eta\in Y^{\an}$ is a Type II point that is a leaf of $\Sigma$ (\mbox{Definition \ref{def-leaf}}) with $\delta(\eta)=0$, we call $\ret_\Sigma^{-1}(\eta)$ a \emph{tail component} of $Y^{\tame}$. The union of all the tail components is called the \emph{tail locus}, denoted $Y^{\tail}$. The complement
\begin{equation*}
Y^{\interior}\coloneqq Y^{\tame}\setminus Y^{\tail}
\end{equation*}
is called the \emph{interior tame locus}.

\end{Def}

Note that the tame locus is an affinoid subdomain of $Y^{\an}$ (or equals $Y^{\an}$) by the same argument as the one used in the proof of \mbox{Lemma \ref{lem-rational-domain}}: Its complement consists of finitely many open disks and anuli, so it is an affinoid subdomain by \cite[Lemma 4.12]{bpr}. Likewise, the interior tame locus, the tail locus, and each tail component are obtained from $Y^{\tame}$ by removing certain connected components. Thus they are themselves affinoid subdomains (or equal to $Y^{\an}$).

\begin{Lem}\label{lem-mu-lambda-tilde-distinction}Let $p_0\in Y\setminus\Sigma$ be a closed point with image $x_0=\phi(p_0)$ in $X$. Then we have the following:
\begin{enumerate}[(a)]
\item If $p_0\in Y^{\interior}$, then $\lambda(x_0)=\mu(x_0)$.
\item If $p_0\not\in Y^{\interior}$, then $\lambda(x_0)>\mu(x_0)$.
\end{enumerate}
\end{Lem}
\begin{proof}
Since $\Sigma$ is a skeleton, there exists a unique path $[p_0,\eta]$ with $[p_0,\eta]\cap\Sigma=\{\eta\}$. The image $\phi([p_0,\eta])$ is the unique path connecting $x_0$ and $\xi\coloneqq\phi(\eta)$. The point $p_0$ lies in $Y^{\tame}$ if and only if $\delta(\eta)=0$ and lies in $Y^{\tail}$ if and only if $\delta(\eta)=0$ and $\eta$ is a leaf of $\Sigma$.

According to the construction of $\Sigma$ (see \mbox{Proposition \ref{prop-minimal-skeleton}}), if $\delta(\eta)=0$ and $\eta$ is a leaf of $\Sigma$, then $\eta$ does not lie on the graph $\phi^{-1}(\Gamma_0)$. Moreover, $\delta$ is strictly positive on the interval connecting $\xi$ and $\Gamma$. Thus it follows from \mbox{Proposition \ref{prop-simple-lambda-consequence}} that $\lambda(x_0)>\mu(x_0)$.

Similarly, if $\delta(\eta)=0$, but $\eta$ is not a leaf of $\Sigma$, then $\eta$ lies on the graph $\phi^{-1}(\Gamma_0)$. Thus $\xi\in \Gamma_0$ and so $\lambda(x_0)=\mu(x_0)$.

Finally, in the case that $\delta(\eta)=\delta(\xi)>0$, it follows immediately from \mbox{Proposition \ref{prop-simple-lambda-consequence}} that $\lambda(x_0)>\mu(x_0)$.
\end{proof}

The following proposition shows that the tame locus, which was defined using a $\phi$-minimal skeleton $\Sigma$ of $Y$, retains the information contained in $\Sigma$ and thus controls the semistable reduction of $Y$.

\begin{Prop}\label{prop-tame-locus-knows-all}
Let $\Sigma_1\subset Y^{\an}$ be a graph containing $\phi^{-1}(\Gamma_0)$ whose vertex set $S(\Sigma_1)$ contains all boundary points of $Y^{\tame}$. Then $\Sigma_1$ is a skeleton of $Y$.
\end{Prop}
\begin{proof}
We check that $\Sigma_1$ contains all vertices of the $\phi$-minimal model $\Sigma$. Recall its construction in \mbox{Section \ref{sec-greek-lambda}}. It is obtained from the graph $\Sigma_0=\phi^{-1}(\Gamma_0)$ by adding to the vertex set $S(\Sigma_0)$ the points of positive genus on $\Sigma_0$ and attaching to $\Sigma_0$ all paths to points of positive genus in the complement of $\Sigma_0$.

Consider first a point $\eta$ of positive genus on $\Sigma_0$. Clearly $\eta$ is a topological ramification point with $\delta(\eta)=0$. If $\eta$ is not among the vertices of $\phi^{-1}(\Gamma_0)$, then there are only two branches at $\eta$ contained in $\Sigma_0$, which we denote $v_1,v_2$. The Riemann-Hurwitz Formula \cite[Theorem 4.5.4]{ctt} applied to $\eta$ reads
\begin{equation*}
2g(\eta)-2=-2n_\eta+n_{v_1}+\slope_{v_1}(\delta)+n_{v_2}+\slope_{v_2}(\delta)-2.
\end{equation*}
If $\eta$ is a tame topological ramification point, then $\delta=0$ in a neighborhood of $\eta$, so we get
\begin{equation*}
2g(\eta)=-2n_\eta+n_{v_1}+n_{v_2}\le0,
\end{equation*}
contradicting the assumption that $\eta$ has positive genus. If $\eta$ is a wild topological ramification point, a similar argument shows that we must have $\slope_{v_1}(\delta)>0$ or $\slope_{v_2}(\delta)>0$. In this case, $\eta$ is a boundary point of the interior tame locus.

Finally consider a point $\eta$ of positive genus not on $\Sigma_0$. Again we have $\delta(\eta)=0$, while the path connecting $\eta$ and $\Sigma_0$ consists by \mbox{Proposition \ref{prop-wild-paths-to-leaves}} of wild topological ramification points. Thus $\eta$ is a boundary point of the tail locus $Y^{\tail}$.
\end{proof}

Our eventual goal is to compute for plane quartic curves $Y$ a skeleton from knowledge of the tame locus $Y^{\an}$. For illustrative purposes, we assume conversely that the $\phi$-minimal skeleton $\Sigma\subset Y^{\an}$ is known in advance in the following examples. We will explain how to compute $\Sigma$ (and the restriction of $\delta$ to $\Sigma$) in \mbox{Chapter \ref{cha-quartics}}. In particular, we will revisit the following examples in \mbox{Section \ref{sec-quartics-examples}}.

\begin{figure}[htb]\centering\includegraphics[scale=0.33]{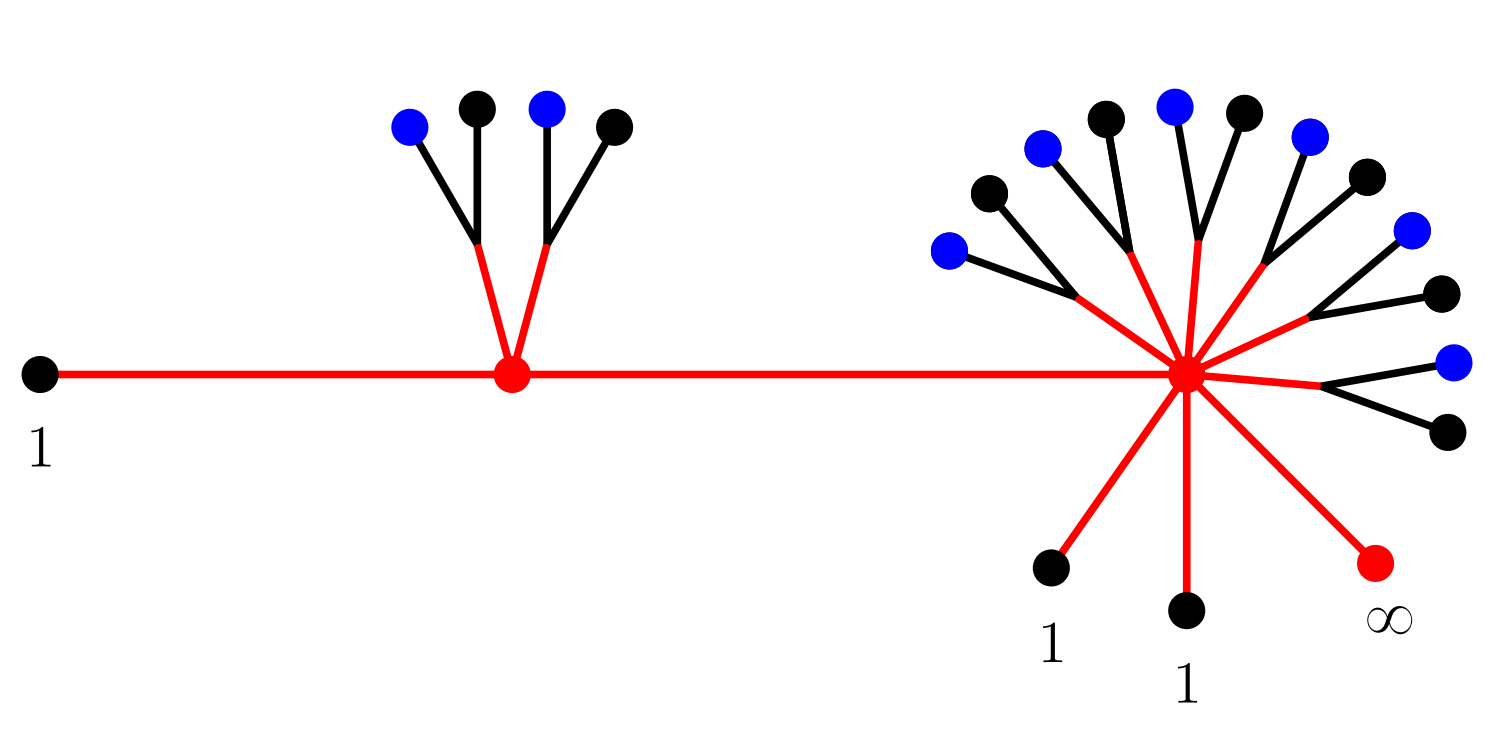}
\caption{The $\phi$-minimal skeleton of the curve considered in \mbox{Example \ref{ex-tame-locus-1}}, with the locus where $\delta>0$ in red}\label{fig-tame-locus-example-1}
\end{figure}

\begin{Ex}\label{ex-tame-locus-1}
We consider the plane quartic curve
\begin{equation*}
Y:\quad y^3 + 3xy^2 - 3y - 2x^4 - x^2 - 1=0
\end{equation*}
over $K=\BC_3$, previously studied in \mbox{Example \ref{ex-lambda-as-radius-of-convergence}}. As usual, $\phi$ denotes the associated cover $Y\to\BP^1_K$. We take for $\Gamma_0$ the skeleton spanned by the branch points of $\phi$, which include $\infty$ (as we saw in \mbox{Example \ref{ex-lambda-as-radius-of-convergence}}). We will see in \mbox{Example \ref{ex-tame-locus-1-revisited}} that the $\phi$-minimal skeleton $\Sigma$ of $Y$ with respect to $\Gamma_0$ looks as depicted in \mbox{Figure \ref{fig-tame-locus-example-1}}. 

The leaves labeled ``$1$'' are Type II points of genus $1$. The other leaves lie above the branch points of $\phi$. It follows from \mbox{Lemma \ref{lem-quartic-branch-locus}} below that $\phi$ is totally ramified above $\infty$, while the remaining eight branch points are ordinary. Above each of the ordinary branch points lies a ramification point and an unramified point. In \mbox{Figure \ref{fig-tame-locus-example-1}}, these ordinary ramification points are colored blue, while the locus where $\delta$ is positive is colored red. The latter includes the ramification point above $\infty$, labeled ``$\infty$''.

Thus the tame locus associated to $\phi$ is the disjoint union of eleven components. Eight of them are closed disks, containing one ramification point each. The other three are tail components whose unique boundary point has \mbox{genus $1$}.
\end{Ex}

\begin{figure}[htb]\centering\includegraphics[scale=0.35]{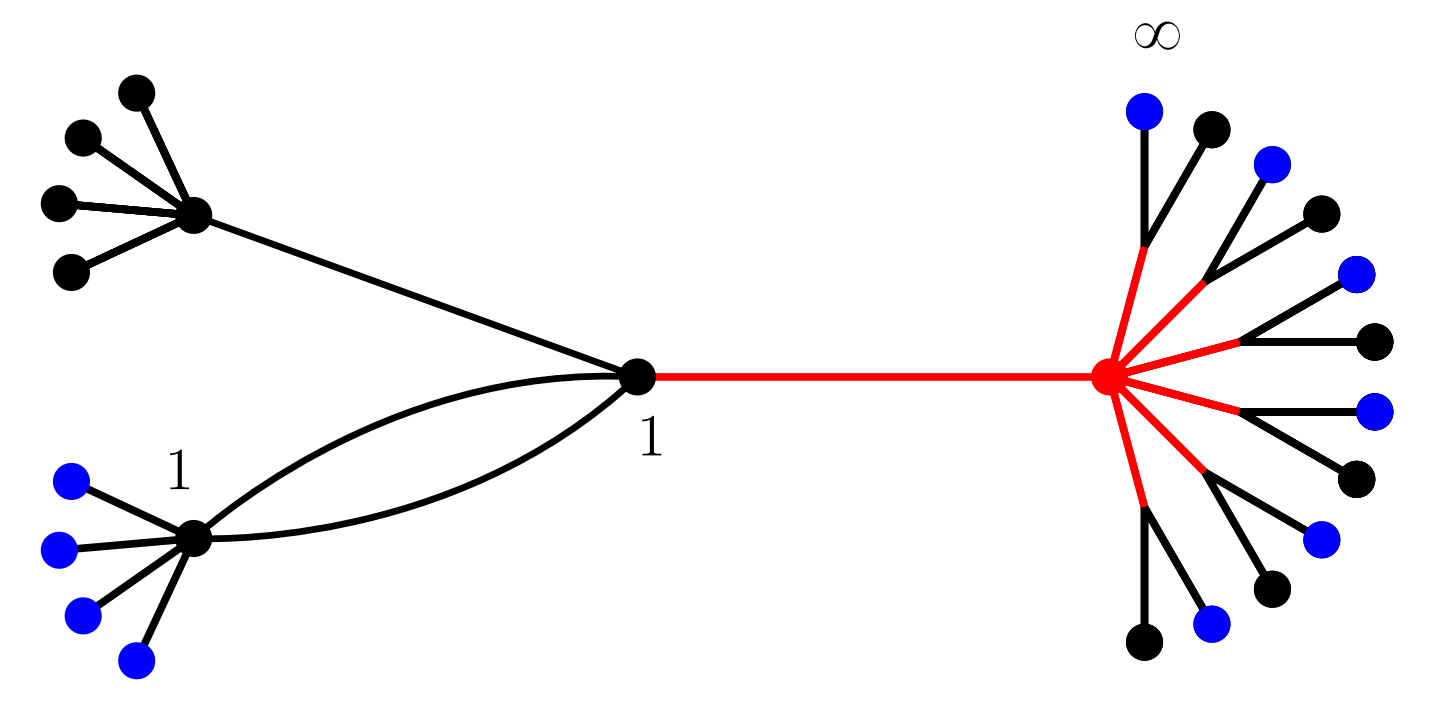}
\caption{The $\phi$-minimal skeleton of the curve considered in \mbox{Example \ref{ex-tame-locus-2}}, with the locus where $\delta>0$ in red} \label{fig-tame-locus-example-2}
\end{figure}

\begin{Ex}\label{ex-tame-locus-2}
Consider the plane quartic curve
\begin{equation*}
Y:\quad y^3-y^2+(3x^3+1)y+3x^4=0
\end{equation*}
over $K=\BC_3$. Again we take for $\Gamma_0$ the skeleton spanned by the branch points of $\phi\colon Y\to\BP_K^1$, which again include $\infty$. We will see in \mbox{Example \ref{ex-tame-locus-2-revisited}} that the $\phi$-minimal skeleton $\Sigma$ of $Y$ with respect to $\Gamma_0$ looks as depicted in \mbox{Figure \ref{fig-tame-locus-example-2}}.

Again, the vertices labeled ``$1$'' are Type II points of genus $1$. This time, there is only two of these, but the skeleton has a loop as well. All leaves lie above branch points of $\phi$, so $Y^{\tail}=\emptyset$. There are ten branch points, all ordinary, among them $\infty$. Thus above each of them lie two points, one a ramification point, the other not. The ramification points are colored blue as in \mbox{Example \ref{ex-tame-locus-1}}, while the locus where $\delta>0$ is colored red.

The interior locus has seven components of which only the one on the left is important, containing the loop and the two points of positive genus. We note that the left point of positive genus has topological ramification index $2$ (the point right above it has the same image), while the right point of positive genus is a wild topological branch point. It is part of the wild topological branch locus, but still we have $\delta=0$ there (the case of Lemma \ref{lem-bijection-over-topological-branch-locus}(b)).

An illustration of $Y^{\an}$ is drawn on the cover of this thesis\footnote{Not included in the electronic version. You can take a look at \href{https://github.com/oossen/mclf/blob/plane-quartics/README.md}{my GitHub page}.}, also depicting the locus where $\delta>0$ in red. For artistic reasons, and to visualize the myriad options to add Type II points to a skeleton of $Y^{\an}$, it is drawn with much branching and squiggliness.
\end{Ex}

\begin{Ex}\label{ex-superelliptic-tame-locus}
We consider again the superelliptic curve $Y$ given by $y^p=f$ of \mbox{Example \ref{ex-good-equation-superelliptic}}, keeping all notation introduced there. We will give a criterion for a closed point $x_0\in U$ to lie in $Y^{\tame}$ dependent on the coefficients of the minimal polynomial
\begin{equation*}
H(T)=T^p+pu_mT^{p-1}+\ldots+pu_m^{p-1}T+u_m^p-\sum_{i=0}^d\frac{f^{(i)}(x_0)}{i!f(x_0)}t^i
\end{equation*}
belong to the chart $(U,w)$ associated to $x_0$. We write
\begin{equation*}
g\coloneqq u_m^p-\sum_{i=0}^d\frac{f^{(i)}(x_0)}{i!f(x_0)}t^i=\sum_{i=m}^{mp-1}c_it^i.
\end{equation*}
Suppose first that $\lambda(x_0)>\mu(x_0)$. Recall from \mbox{Lemma \ref{lem-tail-case-separable-extension}} that there is a unique point $\eta\in Y^{\an}$ above $\xi_{\lambda(x_0)}$, whose reduction curve we now investigate. According to \mbox{Lemma \ref{lem-lambda-formula-tail-case}}, $(1,v_{\lambda(x_0)}(pu_m^{p-1}))$ is a vertex of the Newton polygon of $H(T)$ with respect to $v_{\lambda(x_0)}$, which is a straight line. It follows that
\begin{equation*}
\frac{v_{\lambda(x_0)}(pu_m^{p-i})}{p-i}=\frac{1}{p-i}+v_{\lambda(x_0)}(u_m)>\frac{1}{p-1}+v_{\lambda(x_0)}(u_m)=\frac{v_{\lambda(x_0)}(pu_m^{p-1})}{p-1}
\end{equation*}
for $i=2,\ldots p-1$, so that
$(i,v_{\lambda(x_0)}(pu_m^{p-i}))$ does not lie on $N_H(\lambda(x_0))$ for $i=2,\ldots,p-1$. By \mbox{Remark \ref{rem-tail-case-delta-constant}(a)}, the minimal polynomial of $\overline{w}\coloneqq[w]_\eta$ is then given by 
\begin{equation*}
\overline{H}(T)=T^p+[pu_m^{p-1}]_{\lambda(x_0)}T+[g]_{\lambda(x_0)}.
\end{equation*}
The discriminant of this polynomial is
\begin{equation*}
\Delta_{\overline{H}}=\det\begin{pmatrix}[pu_m^{p-1}]_{\lambda(x_0)}&&&\\
&[pu_m^{p-1}]_{\lambda(x_0)}&&\\
&&\ddots&\\
&&&[pu_m^{p-1}]_{\lambda(x_0)}
\end{pmatrix}=[pu_m^{p-1}]_{\lambda(x_0)}^p
\end{equation*}
(the matrix is the Sylvester matrix associated to $\overline{H}$ and $\overline{H}'$ as defined in \cite[Tag 00UA]{stacks}; it is a $p$-by-$p$ diagonal matrix). By \mbox{Remark \ref{rem-tail-case-delta-constant}(b)}, this discriminant $\Delta_{\overline{H}}$ is a unit. Thus the reduction curve at $\eta$ is given by an Artin-Schreier equation
\begin{equation*}
T^p+\overline{c}T+\overline{g},\qquad\overline{c}\in\kappa,\quad\overline{g}=[g]_{\lambda(x_0)}.
\end{equation*}
This equation defines a curve of positive genus if and only if $\overline{g}$ has a non-zero coefficient of index different from the unique power $q$ of $p$ contained in the set $\{m,\ldots,mp-1\}$. (See \mbox{Section \ref{sec-quartics-characteristic-p}} below and in particular \mbox{Lemma \ref{lem-artin-schreier-genus}} for this.) By \mbox{Remark \ref{rem-lambda-formula-tail-case}(a)}, this is the case if and only if there exists a $k\in\{m,\ldots,mp-1\}\setminus\{q\}$ with
\begin{equation}\label{equ-etale-locus-1}
\frac{p-(p-1)v_K(c_k)}{(p-1)k}\ge\frac{
p-(p-1)v_K(c_q)}{(p-1)q}.
\end{equation}
Now suppose that $\lambda(x_0)=\mu(x_0)$. By \mbox{Lemma \ref{lem-mu-lambda-tilde-distinction}}, we then have $x_0\in Y^{\interior}$. Unless $x_0$ already satisfies \eqref{equ-etale-locus-1}, we have
\begin{equation}\label{equ-etale-locus-2}
\mu(x_0)\ge\frac{p-(p-1)v_K(c_q)}{(p-1)q}.
\end{equation}
In summary, we have $x_0\in Y^{\tame}$ if and only if one of the inequalities \eqref{equ-etale-locus-1} or \eqref{equ-etale-locus-2} holds. These inequalities are the same as the inequalities used to compute the \emph{\'etale locus} associated to $Y$ (cf.\ \cite[Section 4.2]{icerm}) in the MCLF class \texttt{SuperpModel}. Thus the \'etale locus associated to a superelliptic curve of degree $p$ coincides with its tame locus defined in this section.
\end{Ex}

\chapter{Smooth plane quartics}\label{cha-quartics}

In this chapter, we apply the theory developed in \mbox{Chapter \ref{cha-analytic}} to smooth plane quartic curves $Y$, or just plane quartics for short. In \mbox{Section \ref{sec-quartic-normal-form}}, we introduce a normal form for plane quartics that is available whenever the curve $Y$ in question has a rational point.

Our main goal is to gain a better understanding of the tame locus associated to $Y$. In \mbox{Sections \ref{sec-quartics-norm-trick}} \mbox{and \ref{sec-quartics-computing-delta-directly}}, we study the image of the tame locus in the projective line $\BP_K^1$. The advantage of this approach is that we have a concrete understanding of the analytification $(\BP_K^1)^{\an}$, which allows us to compute this image, and hence the semistable reduction of $Y$, for concrete examples.

Every curve $Y$ of genus $3$ over $K$ is either hyperelliptic or else canonically embedded as a quartic plane curve in $\BP_K^2$ (\cite[Section 19.7]{vakil}). It is known how to compute the semistable reduction of any hyperelliptic $K$-curve (in the case that $p\ne2$, using the theory of admissible reduction as explained in \mbox{Section \ref{sec-preliminaries-admissible}}; in the case $p=2$ as explained in \cite[Section 4]{icerm}; see also \mbox{Example \ref{ex-superelliptic-tame-locus}}). As we will recall below, plane quartics are \emph{trigonal} curves, so if $p>3$, then the method explained in \mbox{Section \ref{sec-preliminaries-admissible}} will work for them. In this chapter, we explain one of the remaining cases, that of reduction at $p=3$.

\section{Normal form}\label{sec-quartic-normal-form}

In this section, we allow $K$ to be either algebraically closed or discretely valued.

A plane quartic is, by definition, a closed subscheme $Y\subset\BP^2_K$ cut out by a homogeneous polynomial
\begin{equation*}
\sum_{i+j+k=4}a_{ijk}x^iy^jz^k,\qquad a_{ijk}\in K.
\end{equation*}
We shall assume that $Y$ has a $K$-rational point. Without loss of generality we may take this rational point to be the point $[0:1:0]$. Then we must have $a_{0,4,0}=0$. We may also assume that $a_{1,3,0}=0$. Indeed, if $a_{0,3,1}\ne0$ we may replace $z$ with $z-\frac{a_{1,3,0}}{a_{0,3,1}}x$ to eliminate the $xy^3$-term. And if $a_{0,3,1}=0$, but $a_{1,3,0}\ne0$ we can simply swap $x$ and $z$.

Dehomogenizing by putting $z=1$, an affine chart of $Y$ is cut out of $\BA_K^2$ by the polynomial
\begin{equation*}\label{equ-quartic-normal-form}
F=y^3 + Ay^2 + By + C,
\end{equation*}
where $A,B,C\in K[x]$ are polynomials of degree at most $2$, at most $3$, and at most $4$ respectively. In fact, it is possible to eliminate the term $Ay^2$, which yields the normal form of \cite[Section 1]{shioda}, but we do not need this. In what follows, the discriminant of $F$,
\begin{equation*}\label{equ-disc-equation}
\Delta_F=A^2B^2-4B^3-4A^3C-27C^2+18ABC,
\end{equation*}
will play an important role.

From now on we assume that $Y$ is smooth and geometrically connected. Then $F$ is irreducible and the function field of $Y$ is $F_Y=K(x,y)$, a cubic extension of the rational function field $K(x)$, whose generator $y$ has minimal polynomial $F$. The corresponding morphism 
\begin{equation*}
\phi\colon Y\to\BP_K^1
\end{equation*}
is a degree-$3$ cover of the projective line. Geometrically, it is given by ``projection to the $x$-axis''. That is, points $(x,y)$ in $Y\cap\BA_K^2$ are mapped to their $x$-coordinate under $\phi$.

The following lemma describes the branch locus of $\phi$. Recall from \mbox{Section \ref{sec-analytic-good-equation}} that to every closed point $P$ on $X=\BP_K^1$, there is an associated valuation $v_P$ on the function field $F_X=K(x)$, trivial on $K$. It is the valuation attached to the local ring $\CO_{X,P}$, which is a discrete valuation ring. Similarly, there is a valuation $v_Q$ on the function field $F_Y$ for every closed point $Q$ on $Y$.

A closed point $P\in\BP_K^1$ is called \emph{unramified} with respect to $\phi$ if $e_Q=1$ for every point $Q$ above $P$. It is called \emph{totally ramified} with respect to $\phi$ if there is a unique point $Q$ above $P$ with $e_Q=3$. It is called an \emph{ordinary branch point} with respect to $\phi$ if there exists a point $Q$ above $P$ with $e_Q=2$. In the last case, there necessarily is another point $Q'$ above $P$, with $e_{Q'}=1$.

\begin{Lem}\label{lem-quartic-branch-locus} Suppose that $\characteristic(K)=0$. Let $P\in\Spec K[x]$ be a closed point and let $\infty$ denote the point at infinity in $\BP^1_K$.
\begin{enumerate}[(a)]
\item $P$ is unramified with respect to $\phi$ if and only if $v_P(\Delta_F)=0$, is an ordinary branch point if and only if $v_P(\Delta_F)=1$, and is totally ramified if and only if $v_P(\Delta_F)=2$.
\item $\infty$ is unramified with respect to $\phi$ if and only if $\deg(\Delta_F)=10$, is an ordinary branch point if and only if $\deg(\Delta_F)=9$, and is totally ramified if and only if $\deg(\Delta_F)=8$.
\end{enumerate}
\end{Lem}
\begin{proof}
The inverse image of $\BA^1_K=\Spec K[x]$ in $Y$ is the affine chart $V=\Spec R$, where $R=K[x,y]/(F)$. The different of the ring extension $R/K[x]$ is the principal ideal $(F'(y))$. Let $Q\in V$ be a closed point. Then according to \cite[Theorem III.2.6]{neukirch}, we have $e_Q-1=v_Q(F'(y))$.

Part (a) now quickly follows from the fact that $(\Delta_F)=(\Nm(F'(y))$ (\cite[Theorem III.2.9]{neukirch}). Indeed, if $P$ is an ordinary branch point, there are two points $Q_1,Q_2$ above $P$, say with $v_{Q_1}(F'(y))=1$ and $v_{Q_2}(F'(y))=0$. It follows that
\begin{equation*}
v_P(\Delta_F)=v_P(\Nm(F'(y))=v_{Q_1}(F'(y))+v_{Q_2}(F'(y))=1.
\end{equation*}
The rest of (a) is proved in the same way.

For (b), we use the Riemann-Hurwitz formula. In the present situation it reads
\begin{equation*}
4=2g_Y-2=\deg(\phi)(2g_X-2)+\sum_Q(e_Q-1)=-6+\sum_Q(e_Q-1),
\end{equation*}
where the sums are taken over the closed points of $Y$. It follows from (a) that
\begin{equation}\label{equ-disc-riemann-hurwitz}
\sum_{Q\in Y\setminus V}(e_Q-1)=10-\sum_{P\in\Spec K[x]}v_P(\Delta_F)=10-\deg(\Delta_F).
\end{equation}
Depending on whether $\infty$ is unramified, an ordinary branch point, or totally ramified, the sum on the left-hand side of \eqref{equ-disc-riemann-hurwitz} is equal to $0$, $1$, or $2$ respectively. Thus we get (b).
\end{proof}

\section{Degree-$p$ covers in characteristic $p$}\label{sec-quartics-characteristic-p}

Before continuing, we collect some facts about covers of curves of degree $p$ over a field $k$ of characteristic $p$. The crucial fact here is that there exist non-trivial \'etale covers of the affine line $\BA_k^1$, in contrast to the situation for fields of characteristic $0$.

The most important class of such covers is \emph{Artin-Schreier covers}, which we discuss first.

Let $k$ be an algebraically closed field of characteristic $p$, let $f=\sum_ia_ix^i\in k[x]$ be a polynomial of positive degree, and let $c\in k^\times$ be any element. Let us further suppose that $f$ is not of the form $f=g^p-cg$ for some $g\in k[x]$. Then the polynomial
\begin{equation*}
y^p-cy-f
\end{equation*}
is an irreducible element of the polynomial ring $k(x)[y]$ over the rational function field, and therefore defines a field extension of degree $p$. Let us denote by $\phi\colon Y\to\BP^1_k$ the corresponding morphism of curves.

\begin{Lem}
\label{lem-artin-schreier-genus}
If $\deg(f)$ is not a multiple of $p$, then
\begin{equation*}
g=\frac{(p-1)(\deg(f)-1)}{2}.
\end{equation*}
\end{Lem}
\begin{proof}
We may assume that $c=1$, by replacing $y$ with $yc^{1/(p-1)}$ if necessary. Now this follows immediately from \cite[Theorem 3.7.8(d)]{stichtenoth}. 
\end{proof}
If $d\coloneqq\deg(f)$ is a multiple of $p$, Lemma \ref{lem-artin-schreier-genus} is not directly applicable. But we may proceed as follows to compute the genus of $Y$. Write $d=pl$ and choose a $p$-th root $a_d^{1/p}$ of the leading coefficient $a_d$ of $f$. Then substituting $y+a_d^{1/p}x^l$ for $y$ replaces $f$ with a new polynomial of lower degree. Continuing this procedure until the degree of $f$ is no longer divisible by $p$, we may then apply \mbox{Lemma \ref{lem-artin-schreier-genus}}.

It might happen that we end up with a polynomial $f$ of degree $1$; in this case, \mbox{Lemma \ref{lem-artin-schreier-genus}} tells us that $Y$ has genus $0$.

This is essentially the procedure of \cite[Lemma 3.7.7]{stichtenoth}. Compare also \cite[Proposition 2.1.1]{farnell}, where the above change of variables is used to derive a normal form $y^p-y=f$ for Artin-Schreier covers such that $\deg(f)$ is coprime to $p$.

Now we consider the case that $p=3$. Let $k$ be an algebraically closed field of characteristic $3$ and let $A,B,C\in k[x]$ be polynomials of degree $\le2$, $\le3$, and $\le4$ respectively. Let us further suppose that the polynomial
\begin{equation*}
f=y^3+Ay^2+By+C
\end{equation*}
is an irreducible element of the polynomial ring $k(x)[y]$, so defines a field extension of degree $3$. Let us denote again by $\phi\colon Y\to\BP_k^1$ the corresponding morphism of curves.

\begin{Lem}\label{lem-quartic-is-galois}
Assume that the discriminant $\Delta_f$ of the polynomial $f$ is a unit (i.\,e. an element of $k^\times$). Then $A=0$ and $B$ is a constant in $k^\times$.
\end{Lem}

\begin{proof}
Note that we have
\begin{equation*}
\Delta_f=A^2B^2-B^3-A^3C.
\end{equation*}
Let us first suppose that $A\ne0$. Then we can define
\begin{equation*}
z\coloneqq\frac{\sqrt{\Delta_f}}{A^2y-AB},
\end{equation*}
where $\sqrt{\Delta_f}$ is a choice of square root of the constant $\Delta_f$. An easy calculation shows that this element satisfies the Artin-Schreier equation
\begin{equation*}
z^3-z-\frac{\sqrt{\Delta_f}}{A^3}.
\end{equation*}
If $A$ is a constant, then we see that $z$ and hence $y$ are contained in $k(x)$, which is impossible. If $A$ has degree $\ge1$, then $\phi$ is ramified at the roots of $A$ by \cite[Proposition 3.7.8(b)]{stichtenoth}, contradicting that $\Delta_f$ is a unit.

Thus we find that $A=0$. In this case $\Delta_f$ can only be a constant if $B$ is a constant, and we are done. 
\end{proof}

In the situation of \mbox{Lemma 
\ref{lem-quartic-is-galois}}, the discriminant $\Delta_f$ is in particular a square. Thus the Galois group of the polynomial $f$ is contained in $A_3$. Said differently, the cover $Y\to\BP_k^1$ is Galois.

\begin{Rem}\label{rem-funny-tail}
\begin{enumerate}[(a)]
\item In \mbox{Lemma \ref{lem-quartic-is-galois}}, it is not enough to assume that $\phi\colon Y\to\BP_k^1$ is unramified outside of $\infty$. For example, the cover $Y\to\BP_k^1$, where $Y$ over $k=\overline{\BF}_3$ is given by
\begin{equation*}
f(y)=y^3 + y^2 + xy -x^4=0,
\end{equation*}
is only ramified at $\infty$. But $\Delta_f=x^2(x+1)^2$ is not a unit. The problem is that $f$ does not describe $Y$ as a smooth plane curve.
\item It seems likely that \mbox{Lemma \ref{lem-quartic-is-galois}} has the following generalization to algebraically closed fields $k$ of arbitrary positive characteristic $p$: Suppose that
\begin{equation*}
f(y)=y^p+f_{p-1}y^{p-1}+\ldots+f_1y+f_0\in k[x][y]
\end{equation*}
is an irreducible polynomial with $\Delta_f\in k^\times$. Then $f$ is an Artin-Schreier polynomial, that is, $f_{p-1}=\ldots=f_2=0$ and $f_1\in k^\times$. Said differently, a degree-$p$ \'etale cover of $\BA_k^1$ by a smooth plane curve is necessarily Artin-Schreier.

If this statement were true, it should allow us to generalize much of the content of this chapter to arbitrary residue characteristic. But as of now we could not prove it.
\end{enumerate}
\end{Rem}

\section{Setup}\label{sec-quartics-setup}

For the rest of this chapter, $K$ is assumed to be algebraically closed of residue characteristic $3$.

We have seen in \mbox{Example \ref{ex-plane-curve-chart}} how to define an \'etale chart for a smooth plane curve. We now apply this to a smooth plane quartic as considered in \mbox{Section \ref{sec-quartic-normal-form}}. In particular we have at our disposal the equation for $Y$,
\begin{equation*}
0=F(T)=T^3+AT^2+BT+C,
\end{equation*}
and the discriminant $\Delta_F$ of $F$. The affine \'etale chart we use is then $(U,y)$, where
\begin{equation*}
U=D(\Delta_F)\subseteq\BA_K^1,
\end{equation*}
and where $y\in F_Y$ is a solution to $F$. As in \mbox{Section \ref{sec-analytic-preliminary}}, we define $\Gamma_0$ to be the tree spanned by the complement $\BP_K^1\setminus U$, which is a finite set of closed points. Note that $\Gamma_0$ always includes $\infty$, which may or may not be a branch point of $\phi$. We write $R\coloneqq \CO_X(U)$ and $S\coloneqq \CO_Y(V)$, where $V=\phi^{-1}(U)$.

To begin, we modify the generator $y$ so that \mbox{Assumption \ref{ass-combinatorial-condition}} is satisfied, following the strategy outlined in \mbox{Section \ref{sec-analytic-good-equation}}. Note that the functions $A,B,C$ are polynomials, so according to \mbox{Remark \ref{rem-polynomial-degree}} we have
\begin{equation*}
\deg_W(A)\le2,\qquad\deg_W(B)\le3,\qquad\deg_W(C)\le4,
\end{equation*}
where $W=\BA_K^1$. It follows from \mbox{Corollary \ref{cor-good-equation}} that order of approximation
\begin{equation*}
m>\frac{4}{3},
\end{equation*}
will suffice. Thus we fix order of approximation $m=2$.

Choose a closed point $x_0\in U$ and write $t=x-x_0$ for the associated parameter. Since $A,B,C\in K[x]$ are polynomials, the Taylor expansions introduced in equation \eqref{equ-taylor-expansions} are just polynomials as well; we write
\begin{equation}\label{equ-f-tilde}
\tilde{F}(T)=T^3+\tilde{A}T^2+\tilde{B}T+\tilde{C},
\end{equation}
where
\begin{equation*}
\tilde{A}=\sum_{l=0}^2a_lt^l,\qquad a_l=\frac{A^{(l)}}{l!},
\end{equation*}
\begin{equation*}
\tilde{B}=\sum_{l=0}^3b_lt^l,\qquad b_l=\frac{B^{(l)}}{l!},
\end{equation*}
\begin{equation*}
\tilde{C}=\sum_{l=0}^4c_lt^l,\qquad c_l=\frac{C^{(l)}}{l!}.
\end{equation*}
Note that we have $a_0=A$, $b_0=B$, and $c_0=C$.
In particular, if $p_0\in Y$ is a closed point with $\phi(p_0)=x_0$, then $y(p_0)$ satisfies
\begin{equation*}
y(p_0)^3+a_0(x_0)y(p_0)^2+b_0(x_0)y(p_0)+c_0(x_0)=0.
\end{equation*}
The minimal polynomial of $y$ over the subfield $K(t)$ of $F_Y$ is obtained from $\tilde{F}$ by evaluation at $x_0$. That is, the minimal polynomial of $y$ over $K(t)$ is
\begin{equation*}
G(T)=T^3+\Big(\sum_{l=0}^2a_l(x_0)t^l\Big)T^2+\Big(\sum_{l=0}^3b_l(x_0)t^l\Big)T+\sum_{l=0}^4c_l(x_0)t^l.
\end{equation*}
Next, we bring into play the formal solution
\begin{equation*}
\tilde{y}=y+v_1t+v_2t^2+v_3t^3+\ldots\in S\llbracket t\rrbracket.
\end{equation*}
To construct it, we follow the recipe explained in the proof of \mbox{Proposition \ref{lem-implicit-function-theorem}}. As mentioned above, we have fixed order of approximation $m=2$, so we will only need the truncation $\tilde{u}_2=y+v_1t$.

Plugging the ansatz $\tilde{u}_2=y+v_1t$, $v_1\in S$, into $\tilde{F}$ and collecting all terms of degree $1$ in $t$ shows that $v_1$ should satisfy
\begin{equation*}
3y^2v_1+y^2a_1+2yv_1a_0+yb_1+v_1b_0+c_1=0.
\end{equation*}
Thus our approximation is
\begin{equation}\label{equ-u-tilde}
\tilde{u}_2=y-\frac{y^2a_1+yb_1+c_1}{3y^2+2ya_0+b_0}t\in S[t].
\end{equation}
Given a closed point $p_0\in\phi^{-1}(x_0)$, a generator $w$ of $F_Y$ for which $(U,w)$ satisfies \mbox{Assumption \ref{ass-combinatorial-condition}} is given by
\begin{equation*}
w=y-\tilde{u}_2(p_0)=y-y(p_0)+\frac{y(p_0)^2a_1(x_0)+y(p_0)b_1(x_0)+c_1(x_0)}{3y(p_0)^2+2y(p_0)a_0(x_0)+b_0(x_0)}t.
\end{equation*}
The coefficients of the polynomial $\tilde{H}\coloneqq \tilde{F}(T+\tilde{u}_2)$ clearly lie in $S[t]$ again. It will be convenient to ``forget'' about the polynomial $\tilde{F}$ and change notation, now denoting the coefficients of $\tilde{H}$ by $\tilde{A},\tilde{B},\tilde{C}$. Thus we have
\begin{equation*}\label{equ-h-tilde}
\tilde{H}(T)=T^3+\tilde{A}T^2+\tilde{B}T+\tilde{C}.
\end{equation*}
Correspondingly we write
\begin{equation*}\label{equ-changed-notation-ABC}
\tilde{A}=\sum_{i=0}^2a_it^i,\qquad \tilde{B}=\sum_{i=0}^3b_it^i,\qquad \tilde{C}=\sum_{i=2}^4c_it^i.
\end{equation*}
Similarly, we denote the coefficients of the minimal polynomial $H(T)$ of $w$ over $K(t)$ by $A,B,C$. It is obtained from $\tilde{H}$ by evaluating at $p_0$. Thus we have
\begin{IEEEeqnarray*}{rCl}
H(T)&=&T^3+AT^2+BT+C\\&=&T^3+\Big(\sum_{l=0}^2a_l(p_0)t^l\Big)T^2+\Big(\sum_{l=0}^3b_l(p_0)t^l\Big)T+\sum_{l=0}^4c_l(p_0)t^l.
\end{IEEEeqnarray*}

\begin{Rem}\label{rem-minpoly-reduction-quartics}
Assume that $\lambda(x_0)>\mu(x_0)$. Then according to \mbox{Lemma \ref{lem-tail-case-separable-extension}}, there is a unique point $\eta\in Y^{\an}$ above $\xi_{\lambda(x_0)}$, and the extension of residue fields $\kappa(\eta)/\kappa(\xi_{\lambda(x_0)})$ is generated by $\overline{w}\coloneqq[w]_\eta\in\kappa(\eta)$. Let us repeat how to compute the minimal polynomial of $\overline{w}$, following \mbox{Remark \ref{rem-tail-case-delta-constant}}.

Set $s\coloneqq v_\eta(w)=v_{\lambda(x_0)}(C)/3$. Then the minimal polynomial of $\overline{w}$ is
\begin{equation*}
T^3+\overline{A}T^2+\overline{B}T+\overline{C},\qquad\textrm{where}
\end{equation*}
\begin{equation*}
\overline{A}=\overline{A\pi^{-s}},\qquad\overline{B}=\overline{B\pi^{-2s}},\qquad\overline{C}=\overline{C\pi^{-3s}}.
\end{equation*}
Here the bar over $A\pi^{-s}$ and the other functions denotes reduction to the residue field $\kappa(\xi_{\lambda(x_0)})$. Note that as in \mbox{Remark \ref{rem-tail-case-delta-constant}} this is well-defined because the Newton polygon of $H$ with respect to $v_{\lambda(x_0)}$ is a straight line.
\end{Rem}

\section{Tame locus}\label{sec-quartics-tame-locus}

Let $\Gamma_0$ be the tree defined in \mbox{Section \ref{sec-quartics-setup}}. As in \mbox{Section \ref{sec-analytic-tame-locus}} we denote by $\Sigma$ the $\phi$-minimal skeleton of $Y$ with respect to $\Gamma_0$. In Definition \ref{def-tame-locus} we introduced the tame locus $Y^{\tame}$ associated to $\phi$ and $\Gamma_0$, as well as the tail locus $Y^{\tail}$ and the interior tame locus $Y^{\interior}$. The goal of this section is to describe the tame locus using certain valuative functions on $Y^{\an}$.

Let $p_0\in Y$ be a closed point in the complement of $\Sigma$ and denote its image $\phi(p_0)$ in $\BP_K^1$ by $x_0$. As before, we consider the families of open disks
\begin{equation*}
\rD(r)\coloneqq\{\xi\in(\BP_K^1)^{\an}\mid v_\xi(x-x_0)> r\}
\end{equation*}
and closed disks
\begin{equation*}
\rD[r]\coloneqq\{\xi\in(\BP_K^1)^{\an}\mid v_\xi(x-x_0)\ge r\}
\end{equation*}
around $x_0$, denote the boundary point of $\rD[r]$ by $\xi_r$, and the corresponding valuation by $v_r$. Below, we will use the affine \'etale chart $(U,w)$ constructed in the previous section as well as the various coefficients of the polynomials $\tilde{A},\tilde{B},\tilde{C}$.

We define a subset $\rV\subseteq Y^{\an}$ by
\begin{equation*}
\rV=\rV_2\cup\rV_4,\qquad\textrm{where}\\
\end{equation*}
\begin{equation*}
\rV_2=\{\eta\in Y^{\an}\mid 3\hat{b}_0(\eta)+4\hat{c}_3(\eta)\ge6\hat{c}_2(\eta)\},
\end{equation*}
\begin{equation*}
\rV_4=\{\eta\in Y^{\an}\mid 8\hat{c}_3(\eta)\ge3\hat{b}_0(\eta)+6\hat{c}_4(\eta)\}.
\end{equation*}
It is an affinoid subdomain --- or equals $Y^{\an}$ --- by \mbox{Lemma \ref{lem-rational-domain}}.

\begin{Satz}\label{thm-potential-tame-locus}
We have the following inclusions:
\begin{equation*}
Y^{\tail}\subseteq\rV\subseteq Y^{\tame}
\end{equation*}
\end{Satz}
\begin{proof}
It suffices to show these inclusions for the underlying sets of closed points, and indeed for closed points $p_0$ not contained in $\Sigma$ (the set of these points being dense in $Y^{\an}$). Thus let $p_0\in Y^{\an}\setminus\Sigma$ be an arbitrary closed point with image $x_0\coloneqq\phi(p_0)$ in $X^{\an}\setminus\Gamma_0$. 

Let us assume that $p_0\not\in Y^{\interior}$. Because $Y^{\interior}=Y^{\tame}\setminus Y^{\tail}$, it suffices to show that $p_0\in Y^{\tail}$ if and only if $p_0\in\rV$. We know from \mbox{Lemma \ref{lem-mu-lambda-tilde-distinction}} that $\lambda(x_0)>\mu(x_0)$. Thus it follows from \mbox{Theorem \ref{thm-lambda-formula-tail-case}} that
\begin{equation}\label{equ-lambda-radius-quartic-case}
\lambda(x_0)=\max_{k\in\{2,3,4\}}\frac{3v_K(b_0(p_0))-2v_K(c_k(p_0))}{2k}.
\end{equation}
Following Remark \ref{rem-minpoly-reduction-quartics} we write $s\coloneqq v_\eta(w)=v_{\lambda(x_0)}(C)/3$, where $\eta$ is the unique point in the preimage of $\xi_{\lambda(x_0)}$; then we have
\begin{equation*}
\overline{w}\coloneqq[w]_\eta=\overline{w\pi^{-s}}\in \kappa(\eta).
\end{equation*}
The element $\overline{w}$ generates the field extension $\kappa(\eta)/\kappa(\xi_{\lambda(x_0)})$ and satisfies the equation
\begin{equation*}
\label{equ-potential-tame-locus-reduced-equation}
\overline{H}(\overline{w})=\overline{w}^3+\overline{A}\overline{w}^2+\overline{B}\overline{w}+\overline{C}=0,
\end{equation*}
where $\overline{A}=\overline{A\pi^{-s}}$, $\overline{B}=\overline{B}\pi^{-2s}$, $\overline{C}=\overline{C\pi^{-3s}}$. We have $p_0\in\rV$ if and only if one of the conditions
\begin{equation*}
3v_K(b_0(p_0))+4v_K(c_3(p_0))\ge6v_K(c_2(p_0)),
\end{equation*}
\begin{equation*}
8v_K(c_3(p_0))\ge3v_K(b_0(p_0))+6v_K(c_4(p_0))
\end{equation*}
holds. An easy calculation shows that the first condition holds (meaning $p_0\in\rV_2$) if and only if in the maximum in \eqref{equ-lambda-radius-quartic-case} the term involving $c_3$ is not greater than the term involving $c_2$, that is,
\begin{equation*}
\frac{3v_K(b_0(p_0))-2v_K(c_3(p_0))}{6}\le\frac{3v_K(b_0(p_0))-2v_K(c_2(p_0))}{4}.
\end{equation*}
Similarly the second condition holds (meaning $p_0\in\rV_4$) if and only if the term involving $c_3$ in \eqref{equ-lambda-radius-quartic-case} is not greater than the term involving $c_4$, that is,
\begin{equation*}
\frac{3v_K(b_0(p_0))-2v_K(c_3(p_0))}{6}\le\frac{3v_K(b_0(p_0))-2v_K(c_4(p_0))}{8}.
\end{equation*}
Thus we have $p_0\in\rV$ if and only if the maximum in \eqref{equ-lambda-radius-quartic-case} is assumed for the term involving $c_2$ or the term involving $c_4$. Using \mbox{Remark \ref{rem-lambda-formula-tail-case}} we find that this means
\begin{equation*}
\frac{v_{\lambda(x_0)}(C)}{3}=\frac{v_K(c_k(p_0))+k\lambda(x_0)}{3}=\frac{v_K(b_0(p_0))}{2}=\frac{v_{\lambda(x_0)}(B)}{2}
\end{equation*}
for $k=2$ or $k=4$. It follows that $p_0\in\rV$ if and only if in $\overline{C}=\overline{c}_4\overline{t}^4+\overline{c}_3\overline{t}^3+\overline{c}_2\overline{t}^2$ not both of $\overline{c}_2,\overline{c}_4$ vanish.

Moreover, by \mbox{Remark \ref{rem-tail-case-delta-constant}(b)}, the discriminant $\Delta_{\overline{H}}$ of the polynomial $\overline{H}$ is constant. It follows from \mbox{Lemma \ref{lem-quartic-is-galois}} that $\overline{A}=0$ and that $\overline{B}=\overline{b}_0$ is constant, that is, $\overline{H}$ is Artin-Schreier. By \mbox{Lemma \ref{lem-artin-schreier-genus}}, the genus of the reduction curve at $\eta$ is positive if and only if $p_0\in\rV$.

Said differently, if $p_0\in\rV$, then the Type II point $\eta$ is contained in the $\phi$-minimal skeleton $\Sigma$; it is a leaf, since $\lambda(x_0)>\mu(x_0)$. Thus $p_0$ is contained in a tail component with boundary point $\eta$.

If on the other hand $p_0\not\in\rV$, then the reduction curve at $\eta$ has genus $0$, so $\eta$ is not contained in the $\phi$-minimal skeleton $\Sigma$. Thus $p_0$ is not contained in the tame locus.

This concludes the proof that the point $p_0$ (with $p_0\not\in Y^{\interior}$) is contained in $\rV$ if and only if it is contained in $Y^{\tail}$ and thus the proof of the theorem.
\end{proof}

\section{Norm trick}

\label{sec-quartics-norm-trick}

In this section, we study the image $\phi(Y^{\tame})$ of the tame locus. Note that we have
\begin{equation*}
\phi(Y^{\tame})=\ret_\Gamma^{-1}\big(\{\xi\in\Gamma\mid\delta(\xi)=0\}\big),
\end{equation*}
where $\Gamma=\phi(\Sigma)$. It is the skeleton of $\BP_K^1$ obtained from $\Gamma_0$ by adding paths to all points $\xi=\phi(\eta)$, where $\eta\in Y^{\an}$ has positive genus (cf.\  \mbox{Proposition \ref{prop-wild-paths-to-leaves}}). In particular, $\phi(Y^{\tame})$ is an affinoid subdomain of $X^{\an}$, or equals $X^{\an}$. 

The structure of affinoid subdomains of $(\BP_K^1)^{\an}$ is well-understood. They are disjoint unions of so-called \emph{standard affinoid subdomains} (\cite[78]{berkovich}), which are obtained from closed disks by removing a finite number of open disks. In particular, closed disks and anuli are standard affinoid subdomains. The images of the tail components are closed disks, which we call \emph{tail disks}. Indeed, they are given by retractions
\begin{equation*}
\ret_\Gamma^{-1}(\xi),
\end{equation*}
where $\xi=\phi(\eta)$ for a point $\eta\in Y^{\an}$ with positive genus.

To study $\phi(Y^{\tame})$, we will see how to describe the image of the tail locus using admissible functions on $X^{\an}$. Combining this with the results of the next section, wherein we explain how to compute the image of the interior locus, we will be able to describe the entire tame locus $Y^{\tame}$ using admissible functions on $X^{\an}$. In this way, determining $Y^{\tame}$, and hence the semistable reduction of $Y$, becomes computationally accessible in concrete examples.

The ``norm trick'' is named for the norm map of the extension of function fields $F_Y/F_X$,
\begin{equation*}
\Nm\colon F_Y\to F_X.
\end{equation*}
Simply put, \mbox{Corollary \ref{cor-potential-tail-locus-image}} below states that for computing the image of the tail locus, we may replace all rational functions appearing in the definition of the affinoid subdomain $\rV$ from the previous section with their norms. Thus we are led to $\rU\subseteq X^{\an}$ defined as
\begin{equation*}
\rU=\rU_2\cup\rU_4,\qquad\textrm{where}
\end{equation*}
\begin{equation*}
\rU_2=\{\xi\in X^{\an}\mid 3\reallywidehat{\Nm b_0}(\xi)+4\reallywidehat{\Nm c_3}(\xi)\ge6\reallywidehat{\Nm c_2}(\xi)\} \qquad\textrm{and}
\end{equation*}
\begin{equation*}
\rU_4=\{\xi\in X^{\an}\mid 8\reallywidehat{\Nm c_3}(\xi)\ge3\reallywidehat{\Nm b_0}(\xi)+6\reallywidehat{\Nm c_4}(\xi)\}.
\end{equation*}
Recall that given a closed point $x_0\in X^{\an}\setminus\Gamma_0$ with $\lambda(x_0)>\mu(x_0)$ we have by \mbox{Theorem \ref{thm-lambda-formula-tail-case}}
\begin{equation}\label{equ-lambda-restriction-value-lambda-equation}
\lambda(x_0)=\max_{k\in\{2,3,4\}}\frac{3\hat{b}_0(p_0)-2\hat{c}_k(p_0)}{2k},
\end{equation}
where $p_0\in Y$ is any point in the preimage of $x_0$ under $\phi$. Let us denote the valuative functions appearing in this equation by
\begin{equation*}
h_k\coloneqq\frac{3\hat{b}_0-2\hat{c}_k}{2k},\qquad k\in\{2,3,4\}.
\end{equation*}
Thus we have $\lambda(x_0)=\max_{k\in\{2,3,4\}}h_k(p_0)$. The norms of these valuative functions (see \mbox{Definition \ref{def-norm-of-valuative-function}}), the functions
\begin{equation*}
\Nm h_k=\frac{3\reallywidehat{\Nm b_0}-2\reallywidehat{\Nm c_k}}{2k},\qquad k\in\{2,3,4\},
\end{equation*}
will be important below.

The following lemma states in particular that given a closed point $x_0$ with $\lambda(x_0)>\mu(x_0)$, the function $\lambda$ agrees with a valuative function on the closed disk $\rD[\lambda(x_0)]$: The restriction $\rest{\lambda}{\rD[\lambda(x_0)]}$ is equal to one of the functions $\rest{(\Nm h_2)}{\rD[\lambda(x_0)]}$, $\rest{(\Nm h_3)}{\rD[\lambda(x_0)]}$, or $\rest{(\Nm h_4)}{\rD[\lambda(x_0)]}$.

\begin{Lem}
\label{lem-lambda-restriction-valuative}
Suppose that $\lambda(x_0)>\mu(x_0)$. Then there exists an index $k_0\in\{2,3,4\}$ such that for all $\xi\in\rD[\lambda(x_0)]$ and $\eta\in\phi^{-1}(\xi)$ the following statements hold:
\begin{enumerate}[(a)]
\item $h_{k_0}(\eta)=\max_{k\in\{2,3,4\}}h_k(\eta)$
\item $(\Nm h_{k_0})(\xi)=\max_{k\in\{2,3,4\}}(\Nm h_k)(\xi)$
\item The functions $h_{k_0}$ and $\Nm h_{k_0}$ are constant on $\phi^{-1}(\rD[\lambda(x_0)])$ and $\rD[\lambda(x_0)]$ respectively, with value $\lambda(x_0)$ and $3\lambda(x_0)$ respectively
\item If (a)--(c) are only true for $k_0=3$, then the inequalities 
\begin{equation*}
h_3(\eta)\ge h_k(\eta),\qquad(\Nm h_3)(\xi)\ge(\Nm h_k)(\xi),\qquad k=2,4,
\end{equation*}
are strict
\end{enumerate}
\end{Lem}
\begin{proof}
Suppose first that the maximum in \eqref{equ-lambda-restriction-value-lambda-equation} is assumed for $k_0=4$. In the proof of \mbox{Theorem \ref{thm-potential-tame-locus}} we have verified that in this case the reduction curve at $\eta_{\lambda(x_0)}\in Y^{\an}$ (where $\eta_{\lambda(x_0)}$ is the unique preimage of $\xi_{\lambda(x_0)}$ under $\phi$) is given by an Artin-Schreier equation
\begin{equation*}
y^3+\overline{b}y+\overline{C}=0,
\end{equation*}
for a polynomial $\overline{C}$ of degree $4$ and $\overline{b}\in\kappa^\times$. By \mbox{Lemma \ref{lem-artin-schreier-genus}}, the genus of this curve is $3$. Now let $x_1\in\rD[\lambda(x_0)]$ be any closed point; choose a preimage $p_1\in\phi^{-1}(x_1)$. Then \eqref{equ-lambda-restriction-value-lambda-equation} holds with $x_1$ and $p_1$ in place of $x_0$ and $p_0$, that is,
\begin{equation}\label{equ-lambda-restriction-value-lambda-equation-2}
\lambda(x_1)=\max_{k\in\{2,3,4\}}\frac{3\hat{b}_0(p_1)-2\hat{c}_k(p_1)}{2k}=\max_{k\in\{2,3,4\}}h_k(p_1).
\end{equation}
Repeating the above argument shows that here too the maximum must be assumed for $k_0=4$. Indeed, otherwise the reduction curve at the point $\eta_{\lambda(x_0)}$ above $\xi_{\lambda(x_0)}=\xi_{\lambda(x_1)}$ could not be of genus $3$.

Thus we have shown that $\lambda(x_1)=h_4(p_1)$ for every closed point $x_1\in\rD[\lambda(x_0)]$ and every choice of preimage $p_1\in\phi^{-1}(x_1)$.

We proceed similarly if the maximum in \eqref{equ-lambda-restriction-value-lambda-equation} is not assumed for $k_0=4$, but is assumed for $k_0=2$. Then the reduction curve at $\eta_{\lambda(x_0)}$ is of genus $1$. The maximum in \eqref{equ-lambda-restriction-value-lambda-equation-2} cannot be assumed for $k_0=4$ (or the reduction curve at $\eta_{\lambda(x_0)}$ would have genus $3$), but has to be assumed for $k_0=2$ (or the reduction curve at $\eta_{\lambda(x_0)}$ would have genus $0$). Thus we deduce that $\lambda(x_1)=h_2(p_1)$ for every closed point $x_1\in\rD[\lambda(x_0)]$ and every choice of preimage $p_1\in\phi^{-1}(x_1)$.

Finally, if the maximum in \eqref{equ-lambda-restriction-value-lambda-equation} is only assumed for $k_0=3$, the reduction curve at $\eta_{\lambda(x_0)}$ is of genus $0$. The maximum in \eqref{equ-lambda-restriction-value-lambda-equation-2} can also only be assumed for $k_0=3$ (or else the reduction curve at $\eta_{\lambda(x_0)}$ would have positive genus). Thus $\lambda(x_1)=h_3(p_1)$ for every closed point $x_1\in\rD[\lambda(x_0)]$ and choice of preimage $p_1\in\phi^{-1}(x_1)$.

We have seen that there exists a $k_0\in\{2,3,4\}$ such that $h_{k_0}(p_1)$ agrees with the maximum $\lambda(x_1)$ over $k\in\{2,3,4\}$ in \eqref{equ-lambda-restriction-value-lambda-equation-2} for every closed point in $\phi^{-1}(\rD[(\lambda(x_0)])$, and hence for every point $\eta\in\phi^{-1}(\rD[(\lambda(x_0)])$ outright. This means that part (a) holds. Moreover, we have seen that if the maximum in \eqref{equ-lambda-restriction-value-lambda-equation} is only assumed for $k_0=3$, then the same is true for the maximum in \eqref{equ-lambda-restriction-value-lambda-equation-2}. It follows that the first inequality in (d) is indeed strict.

Next, note that by \mbox{Corollary \ref{cor-lambda-radius-interpretation}} the function $\lambda$ is constant on $\rD[\lambda(x_0)]$ with value $\lambda(x_0)=\lambda(\xi_{\lambda(x_0)})$. This also shows that $h_{k_0}$ is constant on $\phi^{-1}(\rD[\lambda(x_0)])$, which is the first part of (c). The second statement in (c) follows from \mbox{Lemma \ref{lem-norm-of-valuative-function}}.

To prove part (b), fix a point $\xi\in\rD[\lambda(x_0)]$, which we may take to be different from $\xi_{\lambda(x_0)}$. As a consequence of \mbox{Proposition \ref{prop-splitting-from-radius-lambda}} we have $\phi^{-1}(\xi)=\{\eta_1,\eta_2,\eta_3\}$, and these three points correspond to the three extensions of the valuation $v_\xi$ to the function field $F_Y$. As we have remarked above, we have
\begin{equation}\label{equ-lambda-restriction-value-lambda-equation-3}
h_k(\eta_j)\le h_{k_0}(\eta_j)=\frac{(\Nm h_{k_0})(\xi)}{3},\qquad k\in\{2,3,4\},\quad j\in\{1,2,3\}.
\end{equation}
Another application of \mbox{Lemma \ref{lem-norm-of-valuative-function}} shows
\begin{equation}\label{equ-lambda-restriction-value-lambda-equation-4}
\frac{(\Nm h_k)(\xi)}{3}\le\frac{(\Nm h_{k_0})(\xi)}{3},\qquad k\in\{2,3,4\}.
\end{equation}
Thus we have proven (b).

Finally note that in the case that (a) only holds for $k_0=3$, the inequalities in \eqref{equ-lambda-restriction-value-lambda-equation-3} are strict for $k\ne3$ by the first inequality in (d) which we have already proven. Thus the inequality in \eqref{equ-lambda-restriction-value-lambda-equation-4} is also strict for $k\ne3$. This proves that the second inequality in (d) is strict and finishes the proof of the lemma.
\end{proof}

\begin{Kor}\label{cor-potential-tail-locus-image}
We have the following inclusions:
\begin{equation*}
\phi(Y^{\tail})\subseteq\rU\subseteq\phi(Y^{\tame}).
\end{equation*}
\end{Kor}
\begin{proof}
It suffices to show these inclusions for the underlying sets of closed points in the complement of $\Gamma_0$.

Thus let $x_0\in X^{\an}\setminus\Gamma_0$ be a closed point. Because of \mbox{Theorem \ref{thm-potential-tame-locus}} it suffices to show that $x_0\in\rU$ if and only if there exists a preimage $p_0\in\phi^{-1}(x_0)$ with $p_0\in\rV$.

Let us assume first that the statement of \mbox{Lemma \ref{lem-lambda-restriction-valuative}} is true for $k_0=4$. Part (a) of \mbox{Lemma \ref{lem-lambda-restriction-valuative}} then implies that
\begin{equation*}
\frac{3\hat{b}_0(p_0)-2\hat{c}_4(p_0)}{8}=h_4(p_0)\ge h_3(p_0)=\frac{3\hat{b}_0(p_0)-2\hat{c}_3(p_0)}{6}
\end{equation*}
for any preimage $p_0\in\phi^{-1}(x_0)$. This means that $p_0\in\rV_4$. Similarly, part (b) of \mbox{Lemma \ref{lem-lambda-restriction-valuative}} implies that
\begin{IEEEeqnarray*}{cCcCc}
\frac{3\reallywidehat{\Nm b_0}(x_0)-2\reallywidehat{\Nm c_4}(x_0)}{8}&=&(\Nm h_4)(x_0)&&\\&\ge& (\Nm h_3)(x_0)&=&\frac{3\reallywidehat{\Nm b_0}(x_0)-2\reallywidehat{\Nm c_3}(x_0)}{6}.
\end{IEEEeqnarray*}
This means that $x_0\in \rU_4$. Thus the claim that $x_0\in\rU$ if and only if there exists a preimage $p_0\in\phi^{-1}(x_0)$ with $p_0\in\rV$ holds in the case that the statement of \mbox{Lemma \ref{lem-lambda-restriction-valuative}} holds for $k_0=4$.

In the case that the statement of \mbox{Lemma \ref{lem-lambda-restriction-valuative}} holds for $k_0=2$, we conclude completely analogously that $x_0\in\rU_2$ and $p_0\in\rV_2$.

Finally we consider the case that the statement of \mbox{Lemma \ref{lem-lambda-restriction-valuative}} only holds for $k_0=3$. Then part (d) of that lemma shows that we have the opposite and strict inequalities
\begin{equation*}
h_3(p_0)>h_2(p_0),h_4(p_0),\qquad(\Nm h_3)(x_0)>(\Nm h_2)(x_0),(\Nm h_4)(x_0).
\end{equation*}
This means that $x_0\not\in\rU$ and $p_0\not\in\rV$. Thus the claim is also true in this last case, completing the proof of the corollary.
\end{proof}

\section{Computing $\delta$ directly}

\label{sec-quartics-computing-delta-directly}

Let $x_0\in X^{\an}$ be a closed point different from $\infty$. In this section, we explain a method to directly compute the function $\delta_\phi$ on the interval $[x_0,\infty]\subset X^{\an}$. The approach presented here cannot replace the one developed in the previous sections, because the point $x_0$ has to stay fixed and we can only consider ``one interval $[x_0,\infty]$ at a time''. But analyzing finitely many intervals is enough to pin down the interior locus.  

Consider the coordinate $t=x-x_0$, which is a parameter for the interval $[x_0,\infty]$. Our point of departure is the minimal polynomial $G$ of the generator $y$ of $F_Y$ over $K(t)$, considered for example in \mbox{Section \ref{sec-quartics-setup}}. However, in this section we begin not by modifying $y$ using an approximation of level $m=2$, but using an approximation of level $m=1$.

Thus choose a point $p_0\in\phi^{-1}(x_0)$ and write $y_0\coloneqq y(p_0)$. Then the generator $y-y_0$ of $F_Y$ has minimal polynomial of the form
\begin{equation}\label{equ-unfinished-minimal-polynomial}
T^3+(a_2t^2+a_1t+a_0)T^2+(b_3t^3+b_2t^2+b_1t+b_0)T+c_4t^4+c_3t^3+c_2t^2+c_1t.
\end{equation}
Note that all the coefficients $a_0,a_1,a_2$, $b_0,\ldots,b_3$, and $c_0,\ldots,c_4$ are simply elements of $K$, depending on $x_0$; we do not view them as functions in $p_0$ here. (In other words, the ``generic polynomials'', which we have indicated with tildes, are absent from this section.) Finally, we consider the generator $z\coloneqq y+vt$, where $v$ is a solution to the equation
\begin{equation}\label{equ-another-cubic-equation}
v^3+a_1v^2+b_2v+c_3=0.
\end{equation}
The minimal polynomial of $z$ is then by construction of the form
\begin{equation}\label{equ-doctored-minimal-polynomial}
T^3+(a_2t^2+a_1t+a_0)T^2+(b_3t^3+b_2t^2+b_1t+b_0)T+c_4t^4+c_2t^2+c_1t.
\end{equation}
(Of course, the coefficients are not the same as the ones in \eqref{equ-unfinished-minimal-polynomial}, we just reuse the notation for convenience.) We write $A= a_2t^2+a_1t+a_0$, $B=b_3t^3+b_2t^2+b_1t+b_0$, $C=c_4t^4+c_2t^2+c_1t$. 

Let $\xi\in[x_0,\infty]$ be a point of Type II. Our goal is to compute $\delta(\xi)$. For the moment we assume that $\xi$ is a wild topological branch point so that $\phi^{-1}(\xi)=\{\eta\}$ for a Type II point $\eta\in Y^{\an}$. The completed residue fields are \emph{one-dimensional analytic fields} over $K$, meaning that they are complete valued fields over $K$ which are finite over a subfield of the form $\reallywidehat{K(x)}$, $x\not\in K$. In \cite[Section 2.4]{ctt}, a \emph{different}
\begin{equation*}
\delta^{\CTT}(\CH_2/\CH_1)
\end{equation*}
is attached to every separable extension $\CH_2/\CH_1$ of one-dimensional analytic fields over $K$. In parallel with our normalization of $\delta$ in \mbox{Section \ref{sec-greek-delta}}, we define
\begin{equation*}
\delta(\CH_2/\CH_1)\coloneqq\delta^{\add}(\CH_2/\CH_1)\coloneqq-\frac{3}{2}\log\big(\delta^{\CTT}(\CH_2/\CH_1)\big).
\end{equation*}
By the definition of $\delta$ in \cite[Section 4.1]{ctt} we then have
\begin{equation*}
\delta(\xi)=\delta(\eta)=\delta\big(\CH(\eta)/\CH(\xi)\big).
\end{equation*}
We will need the following definitions.

\begin{Def}
Let $\CH$ be a one-dimensional analytic field over $K$ with valuation $v_\CH$. An element $x\in\CH\setminus K$ such that $\CH/\reallywidehat{K(x)}$ is a finite separable extension is called a \emph{parameter} for $\CH$.
\begin{enumerate}[(a)]
\item If $x\in\CH\setminus K$ is a parameter such that $\CH/\reallywidehat{K(x)}$ is a tame extension of valued fields, then $x$ is called a \emph{tame parameter} for $\CH$. 
\item If $x\in\CH\setminus K$ is a parameter such that we have
\begin{equation*}
v_\CH\big(\sum_ia_ix^i\big)=\min_i\big(v_K(a_i)+iv_\CH(x)\big),\qquad a_i\in K,
\end{equation*}
then $x$ is called a \emph{monomial parameter} for $\CH$.
\end{enumerate}
\end{Def}

Our key tool for computing $\delta$ is the following:

\begin{Prop}\label{prop-formula-for-delta}
Let $\CH_2/\CH_1$ be a finite separable extension of one-dimensional analytic fields over $K$. Suppose that $x_1$ and $x_2$ are tame monomial parameters of $\CH_1$ and $\CH_2$ respectively. Then we have
\begin{equation*}
\frac{2}{3}\delta(\CH_2/\CH_1)=v_{\CH_2}(\tfrac{dx_1}{dx_2})+v_{\CH_2}(x_2)-v_{\CH_1}(x_1),
\end{equation*}
where $\frac{dx_1}{dx_2}$ is the unique element of $\CH_2$ for which $\frac{dx_1}{dx_2}dx_2=dx_1$ in the module of differentials $\Omega_{\CH_2/K}$.
\end{Prop}
\begin{proof}
Combine \cite[Corollary 2.4.6(ii)]{ctt} and \cite[Lemma 2.1.6(ii)]{ctt}.
\end{proof}

\begin{Lem}
\label{lem-tame-generators}
Let $\xi$ be a wild topological branch point with $\phi^{-1}(\xi)=\{\eta\}$. Then the elements $t=x-x_0\in\CH(\xi)$ and $z\in\CH(\eta)$ (constructed at the beginning of this section, with minimal polynomial \eqref{equ-doctored-minimal-polynomial}) are tame monomial parameters.
\end{Lem}
\begin{proof}
For $t$ this is clear. Indeed, $\CH(\xi)$ actually equals $\reallywidehat{K(t)}$, so $t$ is in particular a tame parameter. And since $v_\xi$ is a Gauss valuation centered at $x_0$, it follows that $t=x-x_0$ is a monomial parameter.

Let us write
\begin{equation*}
\overline{z}\coloneqq[z]_\eta=\overline{z\pi^{-s}},\qquad s=\frac{v_\xi(C)}{3}.
\end{equation*}
Exactly as in the proof of \mbox{Proposition \ref{prop-np-inseparable}} (the key assumption is that $\delta(\xi)>0$), it follows that the Newton polygon of the minimal polynomial \eqref{equ-doctored-minimal-polynomial} of $z$ is inseparable, so $\overline{z}$ satisfies the equation
\begin{equation}\label{equ-irreducible-by-construction}
T^3+\overline{C}=0,\qquad\overline{C}=\overline{C\pi^{-3s}}.
\end{equation}
However, in contrast to the situation in \mbox{Proposition \ref{prop-np-inseparable}}, we now know that $\overline{C}$ is not a third power (because by construction, $C=c_4t^4+c_2t^2+c_1t$). Thus the extension of residue fields $\kappa(\eta)/\kappa(\xi)$ is purely inseparable of degree $3$ and is generated by $\overline{z}$.

To show that $z$ is a monomial parameter, we study the extension of $v_\xi$ to $F_Y$. In the following we use MacLane's theory of inductive valuations. (See \cite{maclane} or \cite[Chapter 4]{rueth}. We give a brief introduction to it in \mbox{Section \ref{sec-applications-discoids}} below.) The Newton polygon of the minimal polynomial of $z\pi^{-s}$,
\begin{equation*}
G=T^3+A\pi^{-s}T^2+B\pi^{-2s}T+C\pi^{-3s},
\end{equation*}
shows that $z\pi^{-s}$ has valuation $0$. Thus the Gauss valuation
\begin{equation}\label{equ-first-approximant}
v_0\colon\sum_ia_iT^i\mapsto\min_i\big(v_\xi(a_i)\big),\qquad a_i\in F_X,
\end{equation}
is a first approximant of $v_\eta$. Since the polynomial \eqref{equ-irreducible-by-construction} is irreducible, it follows from \cite[Lemma 4.8]{rueth} that $G$ is a key polynomial over the first approximant \eqref{equ-first-approximant}. Thus $v_\eta$ corresponds to the pseudo-valuation $[v_0,v(G)=\infty]$. This shows that $v_\eta$ is simply the Gauss valuation with respect to $z\pi^{-s}$,
\begin{equation*}
\sum_ia_i(z\pi^{-s})^i\mapsto\min_i\big(v_K(a_i)\big).
\end{equation*}
Thus $z$ is a monomial parameter as was to be shown.

To show that $z$ is a tame parameter, it suffices to show that the degree of $\CH(\eta)/\reallywidehat{K(z)}$ is coprime to $3$. Note that the residue field $\kappa(\eta)$ is generated over the residue field $\kappa$ of $K$ by $\overline{z}=[z]_\eta$ and $\overline{t}=[t]_\eta$. The relation
\begin{equation*}
\overline{z}^3=\overline{c}_4\overline{t}^4+\overline{c}_2\overline{t}^2+\overline{c}_1\overline{t}
\end{equation*}
shows that the field extension $\kappa(\eta)/k(\overline{z})$ has degree dividing $4$. By \mbox{Remark \ref{rem-geometric-ramification-index-type-ii}}, this degree equals the degree $\CH(\eta)/\reallywidehat{K(z)}$. Thus $z$ is a tame parameter as well.
\end{proof}

\begin{Rem}\label{epp-remark}
The key step in the proof of \mbox{Lemma \ref{lem-tame-generators}} is that by construction, the element $\overline{z}$ is a generator of the field extension $\kappa(\eta)/\kappa(\xi)$. The importance of ensuring this was already foreshadowed in \mbox{Remark \ref{rem-epp-foreshadowing}}. We now provide some additional context for this result.

Suppose that $(L,v_L)/(F,v_F)$ is a finite extension of valued fields of residue characteristic $p>0$, say of ramification index $e\ge1$. The extension of valuations $v_L|v_F$ is called \emph{weakly unramified} if $e=1$. Thus if $v_L|v_F$ is weakly unramified and has separable extension of residue fields, then it is unramified. We say that a finite extension $F'/F$ \emph{eliminates the ramification} of $v_L|v_F$ if each extension $v_{L'}$ of $v_L$ to the compositum $L'\coloneqq LF'$ is weakly unramified over its restriction to $F'$. The situation is summarized in the following diagram of field and valuation extensions:

\begin{equation*}
\begin{tikzcd}
                                           &  & L' \arrow[dd, no head] &  &                                              &  & v_{L'} \arrow[dd, "e=1", no head]      \\
L \arrow[rru, no head] \arrow[dd, no head] &  &                        &  & v_L \arrow[dd, no head] \arrow[rru, no head] &  &                                        \\
                                           &  & F'                     &  &                                              &  & \rest{v_{L'}}{F'} \arrow[lld, no head] \\
F \arrow[rru, no head]                     &  &                        &  & v_F                                          &  &                                       
\end{tikzcd}
\end{equation*}
If the ramification index $e$ of $v_L|v_F$ is coprime to $p$, then it follows from \mbox{Abhyankar's} Lemma that there exists a finite extension $F'/F$ such that each extension $v_{L'}$ of $v_L$ to $L'$ is unramified over its restriction to $F'$. This result is false if $p$ divides $e$. However, under a mild assumption on the associated extension of residue fields, a theorem of Epp (\cite{epp}, \cite[Tag 09F9]{stacks}) shows that there still exists a finite extension $F'/F$ eliminating the ramification of $v_L|v_F$.

Epp's result explains why we found that the extension of residue fields $\kappa(\eta)/\kappa(\xi)$ was purely inseparable of degree $3$. Indeed, according to Epp's result, the extension of valuations $v_\eta|v_K$ becomes weakly unramified after a finite extension of $K$. Since $K$ is algebraically closed, it already is weakly unramified. This is only possible if $v_\eta|v_\xi$ has ramification index $1$, that is, if $\kappa(\eta)/\kappa(\xi)$ has degree $3$. Note the following points however.

\begin{itemize}
\item We get by using only one generator $z$ to compute $\delta$ at each point in the interval $[x_0,\infty]$; we only appropriately scale $z$ before reducing.
\item If the curve $Y/K$ is defined over some non-algebraically closed subfield $K_0\subset K$, then our generator is defined over a finite extension $K'/K_0$. In fact, if $x_0$ is a $K_0$-rational point, then $K'$ is obtained from $K_0$ by solving two cubic equations: the equation $T^3+a_0T^2+b_0T+c_0=0$ to eliminate the constant term of $C$, and the equation \eqref{equ-another-cubic-equation} to eliminate the degree-$3$ term of $C$.
\end{itemize}

In the next chapter, we will use the fact that we can compute $K'$ explicitly to compute the semistable reduction of plane quartics defined over $K_0$. See \cite[Section 6.2.1]{rueth} for a discussion of why Epp's proof is not constructive enough to determine $K'$ in general.
\end{Rem}

\begin{Kor}\label{cor-concrete-delta-formula}
We have
\begin{equation}\label{equ-concrete-delta-formula}
\delta(\xi)=\min\big(\frac{3}{2},\frac{3v_\xi(A)}{2}-\frac{v_\xi(C)}{2},\frac{3v_\xi(B)}{2}-v_\xi(C)\big),
\end{equation}
where $A,B,C$ denote the coefficients of the minimal polynomial \eqref{equ-doctored-minimal-polynomial} of the parameter $z$.
\end{Kor}
\begin{proof}
This is now simply a calculation. We first collect all the ingredients we need.

\begin{enumerate}[(i)]
\item It follows from \mbox{Proposition \ref{prop-formula-for-delta}} and \mbox{Lemma \ref{lem-tame-generators}} that
\begin{equation*}
\frac{2}{3}\delta(\xi)=v_\eta(\tfrac{dt}{dz})+v_\eta(z)-v_\xi(t).
\end{equation*}
\item Differentiating the equation $0=F(z)$ shows 
\begin{equation*}
0=dF=F_zdz+F_tdt=(3z^2+2Az+B)dz+(A'z^2+B'z+C')dt,
\end{equation*}
so
\begin{equation*}
\frac{dt}{dz}=-\frac{3z^2+2Az+B}{A'z^2+B'z+C'}.
\end{equation*}
\item Since $t$ is a monomial parameter we have
\begin{equation*}
v_\xi(A't)\ge v_\xi(A),\qquad v_\xi(B't)\ge v_\xi(B),\qquad v_\xi(C't)=v_\xi(C).
\end{equation*}
To illustrate the reason that the first two inequalities are not equalities in general, suppose that $v_\xi(B)=v_\xi(b_it^i)$ for $i=0$ or $i=3$,
while $v_\xi(B)<v_\xi(b_jt^j)$ for $j\ne i$. Then 
\begin{equation*}
B't=3b_3t^3+2b_2t^2+b_1t
\end{equation*}
has strictly larger valuation than $B$. The same reasoning applies to $A$, but not to $C=c_4t^4+c_2t^2+c_1t$, whence the equality $v_\xi(C
t)=v_\xi(C)$.
\item Since $v_\eta$ is an infinite inductive valuation obtained from a Gauss valuation (see the proof of \mbox{Lemma \ref{lem-tame-generators}}), we have
\begin{equation*}
v_\eta\big(\sum_iA_iz^i\big)=\min_i\big(v_\xi(A_i)+iv_\eta(z)\big),\qquad A_i\in K(t).
\end{equation*}
\item Because the Newton polygon of the minimal polynomial of $z$ is inseparable, we have
\begin{equation*}
v_\eta(z)=\frac{v_\xi(C)}{3}<\frac{v_\xi(B)}{2},\qquad v_\eta(z)=\frac{v_\xi(C)}{3}<v_\xi(A).
\end{equation*}
\item Combining (iii) and (v) shows
\begin{equation*}
\frac{v_\xi(C't)}{3}=\frac{v_\xi(C)}{3}<\frac{v_\xi(B)}{2}\le\frac{v_\xi(B't)}{2},
\end{equation*}
and similarly $v_\xi(C't)/3=v_\xi(C)/3<v_\xi(A't)$.
\end{enumerate}

Putting all this together, we can compute:
\begin{IEEEeqnarray*}{Clr}
&\frac{2}{3}\delta(\xi)\\=&v_\eta(\tfrac{dt}{dz})+v_\eta(z)-v_\eta(t)&\textrm{(by (i))}\\
=&v_\eta(3z^2+2Az+B)-v_\eta(A'z^2+B'z+C')+v_\eta(z)-v_\xi(t)\qquad&\textrm{(by (ii))}\\
=&v_\eta(3z^2+2Az+B)-v_\eta(A'tz^2+B'tz+C't)+v_\eta(z)&\\
=&\min(v_\eta(3z^2),v_\eta(2Az),v_\eta(B))\\
&-\min(v_\eta(A'tz^2),v_\eta(B'tz),v_\eta(C't))+v_\eta(z)\qquad&\textrm{(by (iv))}\\
=&\min\big(1+\frac{2v_\xi(C)}{3},v_\xi(A)+\frac{v_\xi(C)}{3},v_\xi(B)\big)&\textrm{(by (v))}\\&-\min\big(v_\xi(A't)+\frac{2v_\xi(C)}{3},v_\xi(B't)+\frac{v_\xi(C)}{3},v_\xi(C't)\big)+\frac{v_\xi(C)}{3}\\
=&\min\big(1+\frac{2v_\xi(C)}{3},v_\xi(A)+\frac{v_\xi(C)}{3},v_\xi(B)\big)-v_\xi(C)+\frac{v_\xi(C)}{3}&\textrm{(by (vi))}\\
=&\min\big(1, v_\xi(A)-\frac{v_\xi(C)}{3},v_\xi(B)-\frac{2v_\xi(C)}{3}\big)&
\end{IEEEeqnarray*}
\end{proof}

\begin{Rem}\label{rem-final-delta-formula}
We have derived \mbox{Corollary \ref{cor-concrete-delta-formula}} under the assumption that $\xi\in X^{\an}$ is a wild topological branch point. In particular, formula \eqref{equ-concrete-delta-formula} is true for all $\xi$ with $\delta(\xi)>0$.

If $\xi$ is not a wild topological branch point, then the right-hand side of \eqref{equ-concrete-delta-formula} is $\le0$. Indeed, the Newton polygon of $F$ with respect to $v_\xi$ cannot be inseparable, since this would imply that $\xi$ is a wild topological branch point (compare the proof of \mbox{Lemma \ref{lem-tame-generators}}). Thus we have
\begin{equation*}
v_\xi(A)\le\frac{v_\xi(C)}{3}\qquad\textrm{or}\qquad\frac{v_\xi(B)}{2}\le\frac{v_\xi(C)}{3},
\end{equation*}
which implies that the right-hand side of \eqref{equ-concrete-delta-formula} is $\le0$. All in all we thus see that for \emph{any} $\xi\in[x_0,\infty]$ we have
\begin{equation*}
\delta(\xi)=\max\Big(0,\min\big(\frac{3}{2},\frac{3v_\xi(A)}{2}-\frac{v_\xi(C)}{2},\frac{3v_\xi(B)}{2}-v_\xi(C)\big)\Big).
\end{equation*}
\end{Rem}

\begin{Ex}\label{ex-computing-delta-directly}
Let us consider the example of the smooth plane quartic curve over $\BC_3$ with equation
\begin{equation*}
Y\colon\quad F(y)=y^3 - 3y^2 + (-3x^2 - 2x)y + 3x^4 - 3x - 1=0.
\end{equation*}
The discriminant $\Delta_F$ of $F$ is a polynomial of degree $8$; by \mbox{Lemma \ref{lem-quartic-branch-locus}}, the point $\infty$ is a branch point of the morphism $\phi\colon Y\to X=\BP_K^1$ of ramification index $3$. The other branch points, each of ramification index $2$, are the eight zeros of $\Delta_F$ in $K$. Using the technique outlined in \cite[Section 4.8]{rueth}, or by applying the MCLF-function \texttt{BerkovichTree.adapt\_to\_function} to $\Delta_F$, it is not hard too see that the skeleton of $X$ spanned by the branch points of $\phi$ has the shape indicated in \mbox{Figure \ref{fig-branch-locus-tree}}.

\begin{figure}[htb]\centering\includegraphics[scale=0.38]{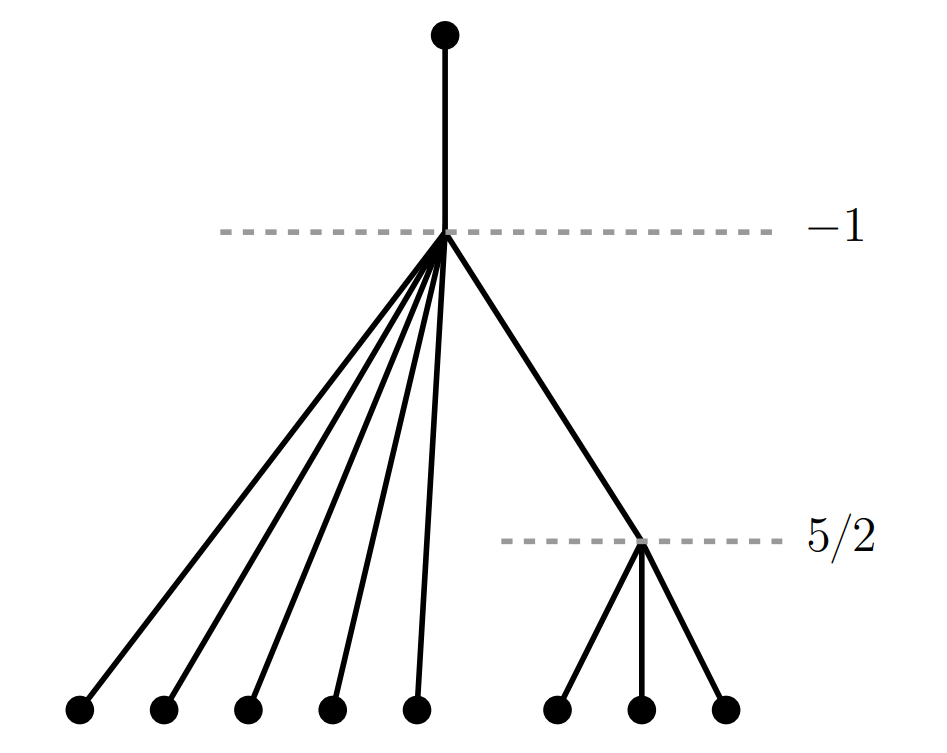}
\caption{The configuration of the branch locus of the morphism $\phi$ considered in \mbox{Example \ref{ex-computing-delta-directly}}. Recall our conventions for visualizing the tree spanned by a finite set of points (\mbox{Section \ref{sec-preliminaries-skeletons}} and particularly \mbox{Example \ref{ex-basic-tree-example}}). The point $\infty$ is at the top, the finite branch points are at the bottom. The dashed lines represent clusters of roots of the indicated radius}\label{fig-branch-locus-tree}
\end{figure}

We will describe the function $\delta$ on the interval $[x_0,\infty]$, where $x_0$ is one of the three branch points contained in the cluster of radius $5/2$. As usual, we write $\xi_r$ for the boundary point of the disk $\rD(r)$ of radius $r\in\BQ$ centered at $x_0$. The first step is to rewrite $F$ in terms of the parameter $t=x-x_0$ and to obtain a new generator $z$ as explained in the beginning of this section. Then $z$ has minimal polynomial of the form
\begin{equation*}
T^3+(a_2t^2+a_1t+a_0)T^2+(b_3t^3+b_2t^2+b_1t+b_0)T+c_4t^4+c_2t^2+c_1t.
\end{equation*}
All the coefficients lie in a finite extension of $\BQ_3$, and it is not hard to compute their valuation (see \mbox{Code Listing \ref{list-discriminant-example}}). We have
\begin{equation*}
v_K(a_0)=\frac{5}{3},\quad v_K(a_1)=2,\quad a_2=0,
\end{equation*}
\begin{equation*}
v_K(b_0)=\frac{10}{3},\quad v_K(b_1)=0,\quad v_K(b_2)=1,\quad b_3=0,
\end{equation*}
\begin{equation*}
v_K(c_1)=0,\quad v_K(c_2)=\frac{5}{3},\quad v_K(c_4)=1.
\end{equation*}
It follows from \mbox{Remark \ref{rem-final-delta-formula}} that 
\begin{equation*}
\delta(\xi_r)=\max\Big(0,\min\big(\frac{3}{2},\frac{3A(r)}{2}-\frac{C(r)}{2},\frac{3B(r)}{2}-C(r)\big)\Big).
\end{equation*}
where $A,B,C$ are the piecewise affine functions
\begin{equation*}
A(r)=\min\big(\frac{5}{3},r+1\big),\qquad B(r)=\min\big(\frac{10}{3},r,2r+1\big),
\end{equation*}
\begin{equation*}
C(r)=\min\big(r,2r+\frac{5}{3},4r+1\big).
\end{equation*}
The graph of the function $r\mapsto \delta(\xi_r)$ is plotted in \mbox{Figure \ref{fig-example-delta-graph-desmos}} --- for $r\le-1$ and $r\ge5$, the function is constant.

Recall that in \mbox{Example \ref{ex-genus-formula}} we have related the slope of $\delta$ at a point $\xi\in[x_0,\infty]$ to the number and type of branch points contained in $\rD(r)$ and the genus of the inverse image $\rC\coloneqq\phi^{-1}(\rD(r))$.

The kinks at $r=-1$ and $r=5/2$ are because of the branch points. Indeed, the slope of $\delta$ changes by $-5/2$ at the radius $r=-1$, since $\rD(-1)$ contains five fewer branch points than $\rD(r)$ for $r<-1$. And the slope changes by $-1$ at the radius $r=5/2$, since $\rD(5/2)$ contains two fewer branch points than $\rD(r)$ for $-1<r<5/2$.

Let us highlight in particular the radius $r=-2/5$. Following \mbox{Remark \ref{rem-minpoly-reduction-quartics}}, we compute the minimal polynomial of $\overline{z}\coloneqq\overline{z\pi^{-s}}$, where $s=v_r(C)/3$. It is the polynomial
\begin{equation*}
T^3+\overline{b}tT+\overline{c}t^4\in k[T],
\end{equation*}
where $\overline{b},\overline{c}\in \kappa^\times$ are certain non-zero elements we don't need to compute here. This shows that the reduction curve at the point $\eta$ above $\xi_r$ is of genus $2$, totally ramified above $\infty$ and with ordinary ramification above $t=0$.

This explains the other kinks in the graph of $\delta$. The slope at radii $0<r<5/2$ is $1/2$, an increase of $3=6/2$ over the slope at radii $-1<r<-2/5$. This is because the genus of $\rC=\phi^{-1}(\rD(r))$ is $0$ for $r>0$, while it is $3$ for $r<-2/5$. Indeed, $\rC$ contains a point of genus $2$ and a loop (lying above the interval $[-2/5,0]$) if $r<-2/5$. We will also re-recover these insights in Example \ref{ex-worked-example-2}, where we return to this example.
\end{Ex}

\begin{figure}[htb]\centering\includegraphics[scale=0.5]{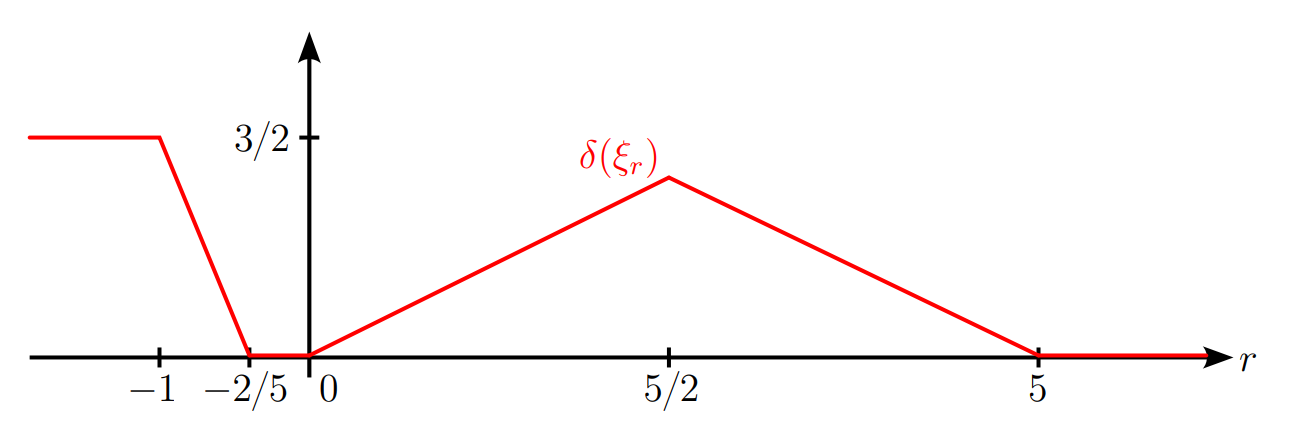}
\caption{The graph of $\delta(\xi_r)$ as a function of $r$}\label{fig-example-delta-graph-desmos}
\end{figure}

Now that we have explained how to compute $\delta$ on intervals of the form $[x_0,\infty]$, the following algorithm to compute the image of the tame locus presents itself.

\begin{Algo}\label{alg-interior-locus}Input: A smooth plane quartic curve $Y$ over $K$, equipped with a degree-$3$ morphism $\phi\colon Y\to X=\BP_K^1$, given in the normal form discussed in \mbox{Section \ref{sec-quartic-normal-form}}.
\begin{enumerate}[(1)] 
\item For each closed point $x_0\in\Gamma_0$ different from $\infty$, determine the restriction of $\delta$ to $[x_0,\infty]$ using the formula in \mbox{Remark \ref{rem-final-delta-formula}}.
\item For each $x_0$, define $\rU_{x_0}$ to be the inverse image under the canonical retraction map $\ret_{\Gamma_0}\colon X^{\an}\to\Gamma_0$ of $\{\xi\in[x_0,\infty]\mid \delta(\xi)=0\}$. It is a disjoint union of closed anuli and closed disks.
\item Compute the affinoid subdomain $\rU=\rU_2\cup\rU_4$ defined in \mbox{Section \ref{sec-quartics-norm-trick}}.
\end{enumerate}  
Output: The image $\phi(Y^{\interior})$ of the interior locus is given by $\bigcup_{x_0}\rU_{x_0}$ and the image of the tail locus $\phi(Y^{\tail})$ by $\rU\setminus\phi(Y^{\interior})$.
\end{Algo}

\begin{proof}[Proof of correctness of the algorithm]
By \mbox{Corollary \ref{cor-potential-tail-locus-image}}, it suffices to show that $\phi(Y^{\interior})=\bigcup_{x_0}\rU_{x_0}$. Let $\eta\in Y^{\an}$ be a point with retraction $\eta_1$ onto $\Sigma$, the $\phi$-minimal skeleton of $Y^{\an}$ with respect to $\Gamma_0$. To say that $\eta\in Y^{\interior}$ is to say that $\eta_1$ is not a leaf of $\Sigma$ and $\delta(\eta_1)=0$. Equivalently, $\xi_1\coloneqq\phi(\eta_1)$ lies in $\Gamma_0$ and $\delta(\xi_1)=0$. If this is the case, $\xi_1$ lies on one of the paths $[x_0,\infty]$; then $\xi\in\rU_{x_0}$. Conversely, if $\delta(\xi_1)>0$ or $\xi_1\not\in\Gamma_0$, then the retraction of $\xi$ onto $\Gamma_0$ does not lie in the locus $\{\xi\in\Gamma_0\mid\delta(\xi)=0\}$, so $\xi$ does not lie in any of the $\rU_{x_0}$.
\end{proof}

\section{Examples}\label{sec-quartics-examples}

In this section, we compute the tame locus in two examples. For illustration, we do this in relative detail, even though there is nothing here that our implementation (see \mbox{Section \ref{sec-appendix-implementation} in the appendix} and the examples in \mbox{Section \ref{sec-applications-examples}}) cannot do on its own.

\begin{Ex}\label{ex-tame-locus-1-revisited}
We revisit \mbox{Example \ref{ex-tame-locus-1}} and compute the tame locus associated to the plane quartic over $K=\BC_3$
\begin{equation*}
Y:\quad F(y)=y^3 + 3xy^2 - 3y - 2x^4 - x^2 - 1=0.
\end{equation*}
The code used for several computations in this example --- the computation of the degree-$27$ numerator of $\Nm(c_3)$, the clustering behavior of the roots of $\Nm(c_3)$, the radii $\lambda(x_0)$, and the reductions of $\tilde{H}$ --- may be found in \mbox{Code Listing \ref{list-tame-locus-1-revisited}}.

Following \mbox{Section \ref{sec-quartics-setup}}, we first compute the polynomial $\tilde{F}(T)$ defined in \eqref{equ-f-tilde}. It is
\begin{IEEEeqnarray*}{rCl}
\tilde{F}(T)&=&T^3+(3t+3x)T^2-3T\\&&-2t^4-8xt^3+(-12x^2-1)t^2+(-8x^3-2x)t-2x^4-x^2-1.
\end{IEEEeqnarray*}
Note that its coefficients are just Taylor expansions of the coefficients of $F$. Next, the approximation of order $2$ as defined in \eqref{equ-u-tilde} is
\begin{equation*}
\tilde{u}_2=y-\frac{3y^2-8x^3-2x}{3y^2+6xy-3}t.
\end{equation*}
This leads to $\tilde{H}(T)=\tilde{F}(T+\tilde{u}_2)$, a polynomial whose coefficients we denote
\begin{equation*}
\tilde{H}(T)=T^3+\tilde{A}T^2+\tilde{B}T+\tilde{C}.
\end{equation*}
Moreover, we use the usual notations $\tilde{C}=\sum_ic_it^i$ and $\tilde{B}=\sum_ib_it^i$ for the coefficients of $\tilde{C}$ and $\tilde{B}$.

To find the three points of positive genus indicated in \mbox{Figure \ref{fig-tame-locus-example-1}}, we make an educated guess. \mbox{Corollary \ref{cor-potential-tail-locus-image}} shows that on every tail disk the function $\Nm(c_3)$ has high valuation compared to $\Nm(c_2)$ or $\Nm(c_4)$. For this reason we consider the numerator of the rational function $\Nm(c_3)$,  suspecting that the tail disks are centered at its roots. It is (up to multiplication with a constant in $K$) the degree-$27$ polynomial
\begin{align*}
x^{27} &- \tfrac{27}{4}x^{26} + \tfrac{410067}{21296}x^{25} - \tfrac{4362147}{85184}x^{24} + \tfrac{5803011}{42592}x^{23} - \tfrac{43814169}{170368}x^{22}\\ &+ \tfrac{159975987}{340736}x^{21} - \tfrac{67087575}{85184}x^{20} + \tfrac{760899573}{681472}x^{19} - \tfrac{1082838429}{681472}x^{18}\\ &+ \tfrac{2331229365}{1362944}x^{17} - \tfrac{5848321365}{2725888}x^{16} + \tfrac{2354253261}{1362944}x^{15} - \tfrac{4942346193}{2725888}x^{14}\\ &+ \tfrac{5475421305}{5451776}x^{13} - \tfrac{5607571491}{5451776}x^{12} + \tfrac{1571650461}{10903552}x^{11} - \tfrac{3708417465}{10903552}x^{10}\\ &- \tfrac{2896673413}{21807104}x^9 - \tfrac{2860325271}{43614208}x^8 - \tfrac{14807853}{123904}x^7 - \tfrac{360410715}{10903552}x^6 - \tfrac{10916775}{5451776}x^5\\ &+ \tfrac{351645543}{21807104}x^4 + \tfrac{180770859}{10903552}x^3 + \tfrac{1673055}{10903552}x^2 - \tfrac{28048275}{21807104}x + \tfrac{373977}{43614208}.
\end{align*}
The denominator of $\Nm(c_3)$ serves a similar role to the \emph{monodromy polynomial} associated to a superelliptic curve defined in \cite{lehr-matignon}, which is also a polynomial at whose roots the tail disks are centered. See \cite[Section 4]{icerm} for more background on the link between monodromy polynomials and tail disks. The zeros of $\Nm(c_3)$ cluster as visualized in \mbox{Figure \ref{fig-thick-cluster-picture}}. That is, there exist three closed disks $\rD_1,\rD_2,\rD_3$, each of radius $1$, each containing nine roots of $\Nm(c_3)$. We do not claim however that these disks cut out the only proper clusters among the roots, or that the clusters cut out by $\rD_1,\rD_2,\rD_3$ have depth $1$.

\begin{figure}[htb]\centering\includegraphics[scale=0.3]{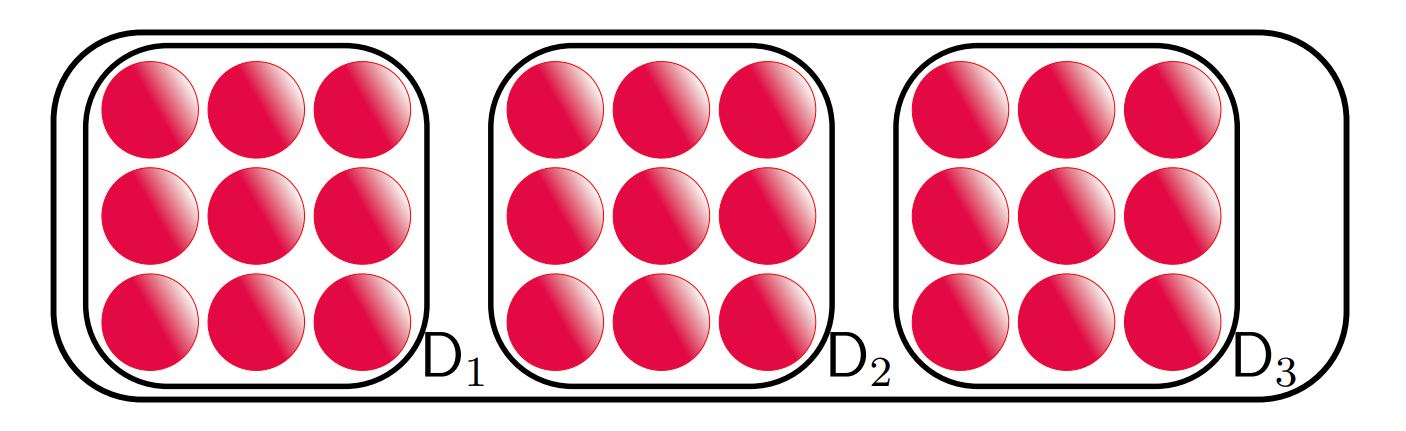}
\caption{The clustering behavior of the roots of $\Nm(c_3)$}\label{fig-thick-cluster-picture}
\end{figure}  

If the disks $\rD_1,\rD_2,\rD_3$ are concentric with three tail disks, then by \mbox{Lemma \ref{lem-lambda-restriction-valuative}(c)} we have for each Type I center $x_0$ of one of them the equality
\begin{equation*}
\lambda(x_0)=\frac{3\reallywidehat{\Nm b_0}(x_0)-2\reallywidehat{\Nm c_2}}{3\cdot4}.
\end{equation*}
For two of the disks $\rD_1,\rD_2,\rD_3$, the resulting radius is $\lambda(x_0)=3/4$, and for the third it is $\lambda(x_0)=1$. Now using \mbox{Remark \ref{rem-minpoly-reduction-quartics}} it is easy to compute the minimal polynomials of the generators $\overline{w}$ of the extension of residue fields $\kappa(\eta_i)/\kappa(\xi_i)$, where $\xi_i$ is the boundary point of $\rD_i$ and $\phi^{-1}(\xi_i)=\{\eta_i\}$. In each case, it is of the form
\begin{equation}\label{equ-artin-schreier-i-promise}
T^3+\overline{b}T+\overline{c}\overline{t}^2,\qquad\overline{b},\overline{c}\in\kappa.
\end{equation}
These are Artin-Schreier polynomials defining function fields of genus $1$ (\mbox{Lemma \ref{lem-artin-schreier-genus}}). As mentioned above, we refer to \mbox{Code Listing \ref{list-tame-locus-1-revisited}} for details on the computation of $\lambda(x_0)$ and \eqref{equ-artin-schreier-i-promise}.

This essentially shows that the tame locus has the structure claimed in Example \ref{ex-tame-locus-1}, except for the position of the ramification points. We omit this point for now, but we will once more return to this curve in \mbox{Example \ref{ex-worked-example-2}}.

\end{Ex}

\begin{Ex}\label{ex-tame-locus-2-revisited}Now we revisit \mbox{Example \ref{ex-tame-locus-2}}, computing the tame locus associated to the plane quartic over $K=\BC_3$
\begin{equation*}
Y:\quad F(y)=y^3-y^2+(3x^3+1)y+3x^4=0.
\end{equation*}
The discriminant of $F$,
\begin{equation*}
\Delta_F=-108x^9 - 243x^8 - 162x^7 - 99x^6 - 42x^4 - 30x^3 - 3,
\end{equation*}
has four zeros of valuation $0$ and five zeros of valuation $-2/5$. In particular, $\phi\colon Y\to\BP_K^1$ has ordinary ramification, one branch point being $\infty$. Among the zeros of $\Delta_F$, there is only one proper cluster, which is cut out by a closed disk of radius $0$ and contains the four zeros of $\Delta_F$ of valuation $0$. As in \mbox{Example \ref{ex-computing-delta-directly}}, we compute the function $\delta$ on the interval $[x_0,\infty]$, where $x_0$ is one of the four zeros of valuation $0$. The resulting graph is shown in \mbox{Figure \ref{fig-revisited-example-2-discriminant-graph}}; see \mbox{Code Listing \ref{list-tame-locus-2-revisited}} for this computation.

\begin{figure}[htb]\centering\includegraphics[scale=0.4]{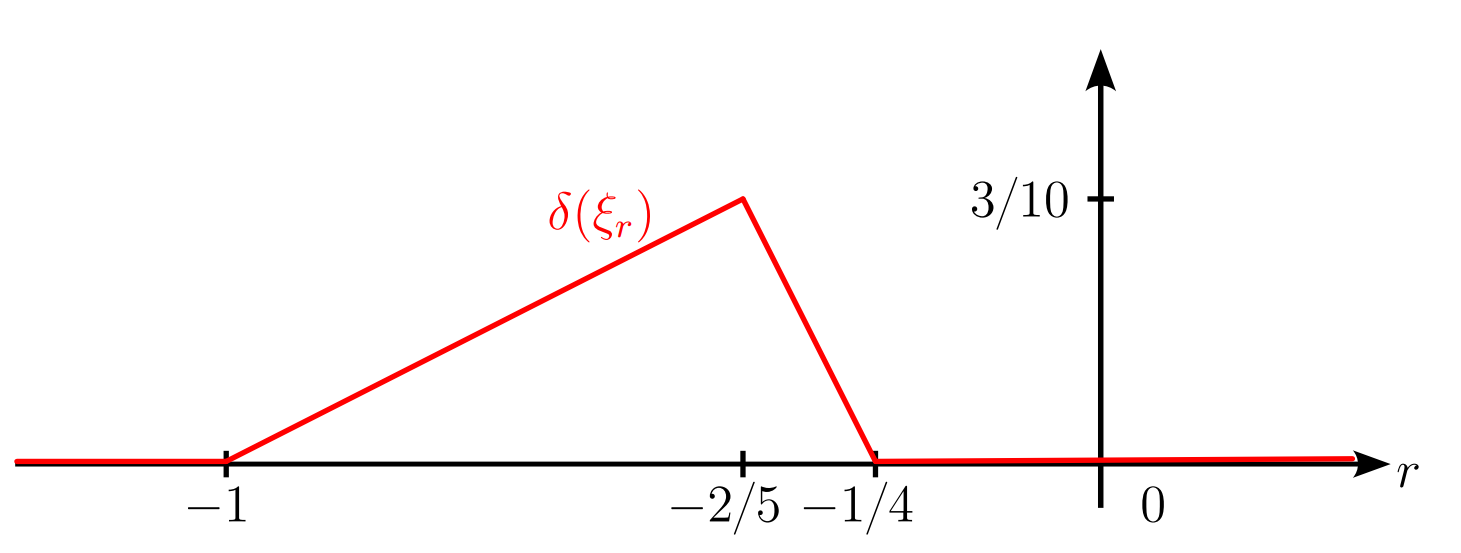}
\caption{The graph of $\delta(\xi_r)$ as a function of $r$}\label{fig-revisited-example-2-discriminant-graph}
\end{figure}  

In \mbox{Code Listing \ref{list-tame-locus-2-revisited}}, we furthermore compute the inverse images of the two points $\xi_0$ and $\xi_{-1/4}$. As can be seen from \mbox{Figure \ref{fig-revisited-example-2-discriminant-graph}}, the latter is a wild topological branch point; its unique preimage has residue field given by
\begin{equation*}
\overline{z}^3-\overline{z}^2-\overline{z}+1/\overline{t}^4=0.
\end{equation*}
This defines a curve of genus $1$, ramified only above $\infty$ (corresponding to the fact that $\delta>0$ for $-1<r<-1/4$). Two points lie above $\xi_0$, one not a topological ramification point, the other a tame topological ramification point with residue field given by
\begin{equation*}
\overline{z}^2=\overline{t}^4+\overline{t}^3+1.
\end{equation*}
This defines a curve of genus $1$, with cover to $\BP^1_\kappa$ branched at four finite points (corresponding to the cluster of four branch points of $\phi$). 

The above shows that the tame locus indeed looks as claimed in \mbox{Example \ref{ex-tame-locus-2}}. Indeed, we have recovered the two points of positive genus. They must be connected by a loop as in \mbox{Figure \ref{fig-tame-locus-example-2}}, since two branches lie above the branch pointing from $\xi_0$ towards $\infty$.

\end{Ex}

\chapter{Examples and implementation}\enlargethispage*{\baselineskip}

\label{cha-applications}

Throughout most of \mbox{Chapters \ref{cha-greek}}, \mbox{\ref{cha-analytic}, and \ref{cha-quartics}} we have assumed that our base field $K$ is algebraically closed. However, we are ultimately interested in studying curves over discretely valued fields. In this chapter, we will explain how to adapt our main results to this goal. We will continue to assume that $K$ is algebraically closed, but make the following assumption (essentially the \mbox{option (b)} in the section on \hyperref[sec-notations]{notations and conventions}):

\begin{Ass}
There exists a discretely valued subfield $K_0\subset K$ such that $K$ is the completion of an algebraic closure of $K_0$.
\end{Ass}
In \mbox{Chapter \ref{cha-quartics}}, we related the tail locus associated to a plane quartic $Y$ over $K$ to admissible functions on $Y^{\an}$, the main result being \mbox{Theorem \ref{thm-potential-tame-locus}}. In the following \mbox{Section \ref{sec-applications-discrete}}, we will see that in the case that $Y$ is defined over $K_0$ and equipped with a morphism $Y\to\BP^1_{K_0}$ over $K_0$, the construction of these valuative functions and the statement of \mbox{Theorem \ref{thm-potential-tame-locus}} remains valid while working over $K_0$.

In \mbox{Sections \ref{sec-applications-discoids}} \mbox{and \ref{sec-applications-line}} we will describe the structure of the projective line $(\BP_{K_0}^1)^{\an}$ over $K_0$. It mirrors the structure of $(\BP_K^1)^{\an}$ explained in \mbox{Section \ref{sec-preliminaries-line}}, with \emph{discoids} taking the place of closed disks. This enables us to adapt \mbox{Algorithm \ref{alg-interior-locus}} to the case of discretely valued ground field. The resulting Algorithm \ref{alg-interior-discoid-locus} is implemented in Sage. We illustrate it with various examples in \mbox{Section \ref{sec-applications-examples}}. Its output is a set of Type II valuations on the rational function field $K_0(x)$, or equivalently (using \mbox{Proposition \ref{prop-julian-theorem}}) a normal model $\CX_0$ of $\BP_{K_0}^1$. By construction, the normalization of $\CX_0$ in the function field $F_{Y_K}$ is a semistable model of $Y_K$. In \mbox{Section \ref{sec-applications-field-extension}}, we discuss the problem of finding a finite field extension $L/K_0$ such that the normalization of $\CX_0$ in $F_{Y_L}$ has semistable reduction.

\section{Discretely valued base field}

\label{sec-applications-discrete}

Suppose that we are in the situation considered in \mbox{Section \ref{sec-analytic-preliminary}}, but over the discretely valued field $K_0$. That is, we have a degree-$p$ cover
\begin{equation*}
\phi\colon Y\to X=\BP_{K_0}^1
\end{equation*}
of smooth projective and geometrically connected curves over $K_0$. Let us moreover suppose that we have an \'etale affine chart $(U,y)$ for $\phi$. This is defined just as it was in \mbox{Definition \ref{def-etale-chart}}, with $K$ replaced by $K_0$. Equivalently, $U\subseteq\BP_{K_0}^1\setminus\{\infty\}$ is a non-empty open subscheme and $y$ a generator of the field extension of function fields $F_Y/F_X$ such that the pair
\begin{equation*}
(U_K,y\otimes1)
\end{equation*}
obtained from $(U,y)$ by applying the base change functor $\cdot\otimes_{K_0}K$ is an \'etale affine chart in the previous sense.

We will now repeat the constructions of \mbox{Sections \ref{sec-analytic-preliminary}} \mbox{and \ref{sec-analytic-good-equation}}, but working with objects over $K_0$ rather than over the algebraically closed field $K$. We will see that we can describe the function $\lambda$ and the tame locus working only over $K_0$. Indeed, this is necessary for working out concrete examples and implementing our algorithms.

As before, write $A\coloneqq\CO_X(U)$ and $B\coloneqq\CO_Y(V)$, where $V=\phi^{-1}(U)$. Then $A$ and $B$ are finite type $K_0$-algebras. Furthermore write
\begin{equation*}
F(T)=T^p+f_{p-1}T^{p-1}+\ldots+f_1T+f_0
\end{equation*}
for the minimal polynomial of $y$ over $F_X$. The coefficients $f_0,\ldots,f_{p-1}$ lie in $A$ (compare \mbox{Lemma \ref{lem-equivalent-etale-chart}}), hence we have
\begin{equation*}
\tilde{f}_i(t)\coloneqq\sum_{l=0}^\infty a_{i,l}t^l\in A\llbracket t\rrbracket,\qquad\textrm{where}\quad a_{i,l}=\frac{f_i^{(l)}}{l!},\quad 0\le i\le p-1.
\end{equation*}
The $\tilde{f}_i$ are the coefficients of the ``generic'' minimal polynomial of $y$,
\begin{equation*}
\tilde{F}(T)=T^p+\tilde{f}_{p-1}T^{p-1}+\ldots+\tilde{f}_1T+\tilde{f}_0\in A\llbracket t\rrbracket[T].
\end{equation*}
The minimal polynomial
\begin{equation*}
G(T)=T^p+g_{p-1}T^{p-1}+\ldots+g_1T+g_0
\end{equation*}
of $y$ over the subfield $K_0(t)$, where $t=x-x_0$ is the parameter associated to a fixed $K_0$-rational point $x_0$, is obtained from $\tilde{F}$ by evaluation at $x_0$; we have $g_i=\sum_la_{i,l}(x_0)t^l$.

Next, we repeat the application of \mbox{Lemma \ref{lem-implicit-function-theorem}} to the polynomial $\tilde{F}$ and the zero $y\in B$ of $\tilde{F}$, resulting in a formal solution
\begin{equation*}
\tilde{y}=y+\tilde{b}_1t+\tilde{b}_2t^2+\tilde{b}_3t^3+\ldots\in B\llbracket t\rrbracket.
\end{equation*}
Consider the truncation $\tilde{u}_m=\tilde{b}_0+\tilde{b}_1t+\ldots+\tilde{b}_{m-1}t^{m-1}$ of order $m$, where $m\ge1$ is some integer. Then the coefficients $\tilde{h}_0,\ldots,\tilde{h}_p$ of the polynomial
\begin{equation*}
\tilde{H}(T)=\tilde{F}(T+\tilde{u}_m)=\sum_{i=0}^p\tilde{h}_iT^i
\end{equation*}
also lie in $B\llbracket t\rrbracket$. As in \mbox{Section \ref{sec-analytic-good-equation}}, we write
\begin{equation*}
\tilde{h}_i=\sum_{l=0}^\infty c_{i,l}t^l.
\end{equation*}
The $c_{i,l}$, $0\le i<p$, $l\ge0$, lie in $F_Y$, so define valuative functions
\begin{equation*}
\hat{c}_{i,l}\colon Y^{\an}\to\BR\cup\{\pm\infty\},\qquad\xi\mapsto v_\xi(c_{i,l}).
\end{equation*}
While we assumed that the base field was algebraically closed in \mbox{Chapter \ref{cha-greek}} (because we treated skeletons, and by extension admissible functions, only in the context of algebraically closed fields), this is not necessary for the above definition of $\hat{c}_{i,l}$. In the present context, $Y^{\an}$ is an analytic curve over the discretely valued field $K_0$, and its underlying set may be identified with pseudovaluations on the function field $F_Y$.

Proceeding as in \mbox{Lemma \ref{lem-construct-admissible-functions}} and \mbox{Remark \ref{rem-making-new-admissible-functions}}, we may then define functions
\begin{equation*}
\tilde{\lambda}=\max_{m\le k<mp}\min_{0<i<p}\min_{0\le pl<(p-i)k}\frac{p\hat{c}_{i,l}-(p-i)\hat{c}_{0,k}}{(p-i)k-pl},
\end{equation*}
\begin{equation}\label{equ-simple-lambda-tilde}
\max_{m\le k<mp}\frac{p\hat{c}_{1,0}-(p-1)\hat{c}_{0,k}}{(p-1)k}
\end{equation}
on $Y^{\an}$. The subtlety of defining these at the closed points where the constituent valuative functions have zeros or poles may be dealt with in the same manner as in the proof of \mbox{Lemma \ref{lem-construct-admissible-functions}}, namely by studying the limiting behavior of as one approaches such closed points. We do not discuss this in detail, since ultimately we need only know the values of these functions at Type II points.

Consider the morphism
\begin{equation*}
\pi\colon Y_K^{\an}\to Y^{\an}
\end{equation*}
induced by the base change morphism $Y_K\to Y$. The points lying over a point $\eta_0\in Y^{\an}$ correspond to the pseudovaluations on $F_{Y_K}$ extending the pseudovaluation $v_{\eta_0}$ on $F_Y$. Thus for every $\eta\in Y^{\an}_K$ and every valuative function $h\colon Y^{\an}\to\BR\cup\{\pm\infty\}$, say $h=\hat{f}$ for $f\in F_Y$, we have
\begin{equation}\label{equ-basechanged-valuative-function}
(\reallywidehat{f\otimes1})(\eta)=\hat{f}(\pi(\eta)).
\end{equation}
There is a similar morphism $X_K^{\an}\to X^{\an}$, which we also denote by $\pi$.

We now specialize to the case that $Y$ is a plane quartic. Using the same notation as in \mbox{Section \ref{sec-quartics-setup}}, we then have
\begin{equation*}
\tilde{H}(T)=T^3+\tilde{A}T^2+\tilde{B}T+\tilde{C},
\end{equation*}
\begin{equation*}
\tilde{A}=\sum_{l=0}^2a_lt^l,\quad\tilde{B}=\sum_{l=0}^3b_lt^l,\quad\tilde{C}=\sum_{l=0}^2c_lt^l.
\end{equation*}
We define subsets $\rV_{K_0}\subseteq Y^{\an}$ and $\rU_{K_0}\subseteq X^{\an}$ by the same formulas as in \mbox{Sections \ref{sec-quartics-tame-locus}} \mbox{and \ref{sec-quartics-norm-trick}} respectively. That is, we have
\begin{equation*}
\rV_{K_0}=\big\{\eta\in Y^{\an}\mid 3\hat{b}_0(\eta)+4\hat{c}_3(\eta)\ge6\hat{c}_2(\eta)\textrm{ or }8\hat{c}_3(\eta)\ge3\hat{b}_0(\eta)+6\hat{c}_4(\eta)\big\},
\end{equation*}
\begin{IEEEeqnarray*}{rCl}
\rU_{K_0}&=&\big\{\xi\in X^{\an}\mid 3\reallywidehat{\Nm b_0}(\xi)+4\reallywidehat{\Nm c_3}(\xi)\ge6\reallywidehat{\Nm c_2}(\xi)\\
&&\textrm{or}\quad8\reallywidehat{\Nm c_3}(\xi)\ge3\reallywidehat{\Nm b_0}(\xi)+6\reallywidehat{\Nm c_4}(\xi)\big\}.
\end{IEEEeqnarray*}

Let us write
\begin{equation*}
\rU_K\coloneqq\pi^{-1}(\rU_{K_0})\subseteq X_K^{\an},\qquad\rV_K\coloneqq\pi^{-1}(\rV_{K_0})\subseteq Y_K^{\an}.
\end{equation*}

\begin{Kor}\label{cor-potential-tail-locus-descent}
$\rU_K$ and $\rV_K$ are the same as the affinoid subdomains $\rU\subset X_K^{\an}$ (defined in \mbox{Section \ref{sec-quartics-norm-trick}}) and $\rV\subset Y_K^{\an}$ (defined in \mbox{Section \ref{sec-quartics-tame-locus}}) associated to the base change $Y_K$.
\end{Kor}
\begin{proof}
This follows immediately from the discussion in this section. Namely, a point $\eta_0\in Y^{\an}$ lies in $\rV_{K_0}$ if and only if $3\hat{b}_0(\eta_0)+4\hat{c}_3(\eta_0)\ge6\hat{c}_2(\eta_0)$ or $8\hat{c}_3(\eta_0)\ge3\hat{b}_0(\eta_0)+6\hat{c}_4(\eta_0)$. The affinoid subdomain $\rV$ of \mbox{Section \ref{sec-quartics-tame-locus}} is defined by the same equalities, except using the functions $b_0\otimes1$, $c_2\otimes1$, $c_3\otimes1$, and $c_4\otimes1$ instead of $b_0,c_2,c_3,c_4$. By \eqref{equ-basechanged-valuative-function}, a point $\eta\in\pi^{-1}(\eta_0)$ lies in $\rV$ if and only if it lies in $\rV_K$. The equality of $\rU_K$ and $\rU$ is proved in the same way.
\end{proof}

\section{Discoids}\label{sec-applications-discoids}

We have studied the structure of the Berkovich line $(\BP_K^1)^{\an}$ (where $K$ is algebraically closed) in \mbox{Section \ref{sec-preliminaries-line}}. To derive a similar description for $(\BP_{K_0}^1)^{\an}$, we need the notion of \emph{discoids}. Building on work of MacLane (\cite{maclane}, \cite{maclane2}), their connection to models and valuations has been worked out in \cite[Section 4.4]{rueth}. However, in \cite{rueth}, discoids are by definition subsets of $K$; essentially, points that are not Type I are disregarded. Following \cite[Chapter 2]{micu}, the following definition takes all points into account.

\begin{Def}
A \emph{discoid} is a subset of $(\BP_{K_0}^1)^{\an}$ of the form
\begin{equation*}
\rD[\psi,\rho]\coloneqq\{\xi\in (\BP_{K_0}^1)^{\an}\mid v_\xi(\psi)\ge\rho\},
\end{equation*}
where $\rho\in\BR$ and $\psi\in K_0[x]$ is a monic irreducible polynomial. The number $\rho$ is called the \emph{radius} of $\rD[\psi,\rho]$. The \emph{degree} of $\rD[\psi,\rho]$ is simply the degree of $\psi$. 
\end{Def}

It will be convenient to extend the notation to the case $\rho=\infty$. Then $\rD[\psi,\infty]$ is just the closed point associated to $\psi$. Note that if $\psi$ is a constant polynomial, say $\psi=x-x_0$ for some $x_0\in K_0$, then $\rD[\psi,\rho]$ is simply the closed disk of radius $\rho$ centered at $x_0$.

For a discoid $\rD\subset (\BP_{K_0}^1)^{\an}$, there is an associated Type II valuation $v_\rD$ on $K_0(x)$ (compare \cite[Lemma 4.50]{rueth}), defined by
\begin{equation*}
v_\rD(f)=\inf\{v_\xi(f)\mid\xi\in\rD\},\qquad f\in K[x].
\end{equation*}
In fact, $v_\rD$ is the valuation corresponding to the unique boundary point of $\rD$ in $(\BP_{K_0}^1)^{\an}$. We denote this boundary point by $\xi_{\rD}$ or, if $\rD=\rD[\psi,\rho]$, by $\xi_{\psi,\rho}$. If $\rho=\infty$, then we denote by $\xi_{\psi,\rho}$ the unique closed point in $\rD[\psi,\infty]$.

\begin{Prop}
The map $\rD\mapsto v_\rD$ is a bijection from the set of discoids $\rD\subset (\BP_{K_0}^1)^{\an}$ of rational radius to the set of Type II valuations on $K_0(x)$.
\end{Prop}
\begin{proof}
This follows from \cite[Theorem 4.56]{rueth}.
\end{proof}

In \mbox{Section \ref{sec-preliminaries-line}}, we saw how to describe all the points of $(\BP_K^1)^{\an}$, where $K$ was algebraically closed, in terms of nested families of disks (\mbox{Definition \ref{def-nested-family}}). The analogue of this in the present context of discoids is given by MacLane's theory of inductive valuations. We only introduce these briefly here and refer the reader to \cite{rueth} and \cite{micu} for more information.

One begins with a Gauss valuation of a certain radius $\lambda_0\in\BR$,
\begin{equation*}
v_0\colon K(x)\to\BR\cup\{\infty\},\qquad\sum_{i}a_ix^i\mapsto\min_i\{v_{K_0}(a_i)+i\lambda_0\}.
\end{equation*}
Then given so-called \emph{key polynomials} $\psi_1,\ldots,\psi_n\in K[x]$ and $\lambda_1<\lambda_2<\ldots<\lambda_n$, where $\lambda_1,\ldots,\lambda_n\in\BR\cup\{\infty\}$, one inductively defines pseudovaluations
\begin{equation*}
v_k=[v_0(x)=\lambda_0,v_1(\psi_1)=\lambda_1,\ldots,v_k(\psi_k)=\lambda_k],\qquad k=1,\ldots,n,
\end{equation*}
on $K(x)$. As the notation suggests, the pseudovaluation $v_k$ satisfies $v_k(\psi_k)=\lambda_k$; moreover, $v_k(\chi)=v_{k-1}(\chi)$ for $\chi\in K[x]$ of degree less than $\deg(\psi_k)$. The key fact (\cite[Proposition 1.110]{micu}) is then that every pseudovaluation $v$ on $K(x)$ --- with the exception of the pseudovaluation associated to the closed point $\infty\in\BP^1_{K_0}$ --- can be approximated by such inductive pseudovaluations. To be precise, we have the following (\cite[Corollary 1.117]{micu}):

\begin{itemize}
\item If $\lambda_n=\infty$, then $v_n=[v_0(x)=\lambda_0,\ldots,v_n(\psi_n)=\lambda_n]$ is a pseudo-valuation that is not a valuation, so is associated to a Type I point of $X^{\an}$.
\item If $\lambda_n\in\BQ$, then $v_n$ is a Type II valuation. In fact, it is the valuation $v_\rD$, where $\rD$ is the discoid $\rD[\psi_n,\lambda_n]$.
\item If $\lambda_n\in\BR\setminus\BQ$, then $v_n$ is associated to a Type III point of $X^{\an}$.
\item Given infinite sequences $(\psi_n)_{n\ge1}$ and $(\lambda_n)_{n\ge0}$, one may define a \emph{limit valuation}
\begin{equation*}
v_\infty\coloneqq\lim_{n\to\infty}[v_0(x)=\lambda_0,\ldots,v_n(\psi_n)=\lambda_n].
\end{equation*}
All Type IV points arise in this way. Note however that there is a definitional difference regarding Type IV points between \cite{micu} and our definition in \mbox{Section \ref{sec-preliminaries-analytic}}: In \cite{micu}, only closed points are considered Type I, while for us, a point $\xi\in X^{\an}$ such that $\CH(\xi)$ may be imbedded in $K$ is considered Type I as well. Thus the statement \cite[\mbox{Corollary 1.117(4)}]{micu} is not quite correct if one follows our definition.
\end{itemize}

It is important to understand how discoids behave under extensions of the base field. What is the preimage of a discoid $\rD[\psi,\rho]$ under the natural morphism
\begin{equation*}
\pi\colon Y^{\an}_K\to Y^{\an}?
\end{equation*}
Phrased differently, which valuations on $F_{Y_K}=F_Y\otimes_{K_0}K$ extend the valuations $\{v_\eta\mid \eta\in \rD[\psi,\rho]\}$ on $F_Y$? These questions are treated in \cite[Section 2.3]{micu}.

Let $P\in\BP_{K_0}^1\setminus\{\infty\}$ be a closed point. It corresponds to a maximal ideal in $K_0[x]$, say the one generated by the monic irreducible polynomial $\psi\in K_0[x]$. Then the points in $\BP_K^1$ above $P$ correspond to the $\deg(\psi)$ zeros of $\psi$ in $K$. We write
\begin{equation}\label{equ-numbering-of-roots}
\psi=\prod_{i=1}^s(x-\alpha_i),\qquad\alpha_1,\ldots,\alpha_s\in K,
\end{equation}
and for $r\in\BR\cup\{\infty\}$ define the sets
\begin{equation*}
I_r\coloneqq\{i\in\{1,\ldots,s\}\mid v_K(\alpha_1-\alpha_i)\ge r\},
\end{equation*}
\begin{equation*}
J_r\coloneqq\{i\in\{1,\ldots,s\}\mid v_K(\alpha_1-\alpha_i)<r\}.
\end{equation*}
Now following \cite[Definition 2.25]{micu}, we define the function
\begin{equation*}
\theta_\psi\colon\BR\cup\{\infty\}\to\BR\cup\{\infty\},\qquad r\mapsto r\abs{I_r}+\sum_{i\in J_r}v_K(\alpha_1-\alpha_i).
\end{equation*}
The function $\theta_\psi$ is independent of the numbering of the roots in \eqref{equ-numbering-of-roots} (\cite[Lemma 2.27]{micu}) and is strictly monotonously increasing (\cite[Lemma 2.28]{micu}).

\begin{Prop}\label{prop-discoid-splitting}
For any $\rho\in\BR$ we have
\begin{equation*}
\pi^{-1}(\rD[\psi,\theta_\psi(\rho)])=\bigcup_{i=1}^r\rD[\alpha_i,\rho].
\end{equation*}
\end{Prop}
\begin{proof}
This is a special case of \cite[Theorem 2.29]{micu}.
\end{proof}

Note that the union in \mbox{Proposition \ref{prop-discoid-splitting}} need not be disjoint; if the radius $\rho$ is chosen small enough, the disks $\rD[\alpha_i,\rho]$ will coalesce.

\section{The structure of $(\BP_{K_0}^1)^{\an}$}\label{sec-applications-line}

Recall that in \mbox{Section \ref{sec-preliminaries-line}} we have explained how $(\BP_K^1)^{\an}$ may be thought of as an infinitely branching tree. Chapter 2 of \cite{micu} leverages discoids to obtain a similar description of $(\BP_{K_0}^1)^{\an}$. 

Crucially, it is still true  that $(\BP^1_{K_0})^{\an}$ is a simply connected special quasipolyhedron (\cite[Theorem 4.2.1]{berkovich}). Thus there exists a unique path between any pair of points in $(\BP^1_{K_0})^{\an}$. These paths are explicitly described in \cite[Section 2.7]{micu}. For example, given two monic irreducible polynomials $\psi_1,\psi_2\in K_0[x]$, the unique path connecting the closed points $\xi_{\psi_1,\infty}$ and $\xi_{\psi_2,\infty}$ is
\begin{equation*}
\{\xi_{\psi_1,r'}\mid r'\ge r\}\cup\{\xi_{\psi_2,r'}\mid r'\ge r\},
\end{equation*}
where $r$ denotes the smallest radius for which $\xi_{\psi_2,\infty}$ is still contained in $\rD[\psi_1,r]$. The path may be visualized in the same way as in the case of algebraically closed ground field, cf.\ \mbox{Figure \ref{fig-hairy-tree-with-path}}.

The description of the affinoid subdomains of $(\BP_K^1)^{\an}$ in terms of standard  affinoid subdomains that we gave in \mbox{Section \ref{sec-quartics-norm-trick}} remains valid over $K_0$. To be more precise, we define an \emph{open discoid} to be a subset of $(\BP^1_{K_0})^{\an}$ of the form
\begin{equation*}
\rD(\psi,\rho)\coloneqq\{\xi\in(\BP^1_{K_0})^{\an}\mid v_\xi(\psi)>\rho\},
\end{equation*}
where $\psi\in K_0[x]$ is monic and irreducible and $\rho\in\BR$, or of the form
\begin{equation*}
\rD(1/x,\rho)\coloneqq\{\xi\in(\BP^1_{K_0})^{\an}\mid v_\xi(x)<-\rho\},
\end{equation*}
where $\rho\in\BR$. The latter type is just open disks centered at $\infty$. Then a \emph{standard affinoid subdomain} of $(\BP^1_{K_0})^{\an}$ is a closed discoid, or a closed disk centered at $\infty$, with a finite number of open discoids removed.

Let $S\subset(\BP_{K_0}^1)^{\an}$ be a finite set of Type I and Type II points, satisfying $\abs{S}\ge3$ or containing at least one Type II point. Using the same procedure as in the proof of \mbox{Proposition \ref{prop-tree-spanned-by-points}}, we may construct the \emph{tree spanned by $S$}. It is the union of all the paths connecting points $\xi,\xi'\in S$. (Really, this tree is a skeleton of $(\BP^1_{K_0})^{\an}$, but we have only defined skeletons in the case that the ground field is algebraically closed.)

Let $T\subset(\BP_{K_0}^1)^{\an}$ be the tree spanned by a finite set $S$. Since $(\BP_{K_0}^1)^{\an}$ is uniquely path-connected, we may in the same way as in \mbox{Definition \ref{def-retraction}} define a canonical retraction
\begin{equation*}
\ret_T\colon(\BP_{K_0}^1)^{\an}\to T.
\end{equation*}
Now let $U\subset T$ be a closed subset. Then the subset
\begin{equation*}
\ret_T^{-1}(U)\subseteq(\BP^1_{K_0})^{\an}
\end{equation*}
is an affinoid subdomain, since it is a disjoint union of standard affinoid subdomains. This way of describing affinoid subdomains by means of the retraction to a tree underlies the implementation of affinoid subdomains of $(\BP^1_{K_0})^{\an}$ in the MCLF Sage package. For details, we refer the reader to the documentation of the MCLF classes \texttt{AffinoidTree} and \texttt{AffinoidDomainOnBerkovichLine}. We illustrate the concept in \mbox{Example \ref{ex-discoid-retraction}} below.

\begin{Rem}\label{rem-mclf-affinoid-computation}\sloppy
The MCLF package also is able to compute subsets like $\rU_{K_0}\subseteq(\BP^1_{K_0})^{\an}$ defined in \mbox{Section \ref{sec-applications-discrete}} as part of the functionality of the class \texttt{PiecewiseAffineFunction}. In fact, $\rU_{K_0}$ is an affinoid subdomain. This is because the valuative functions $\reallywidehat{\Nm b_0},\reallywidehat{\Nm c_2},\reallywidehat{\Nm c_3},\reallywidehat{\Nm c_4}$ used to define $\rU_{K_0}$ are constant outside the tree spanned by the zeros and poles of $\Nm b_0,\Nm c_2,\Nm c_3,\Nm c_4$. In particular, it is clear how to represent $\rU_{K_0}$ (and other affinoid subdomains like it) using retraction to a tree.
\end{Rem}

\begin{Ex}\label{ex-discoid-retraction}
Suppose that $T$ is a tree contained in $(\BP^1_{K_0})^{\an}$ which looks like the one drawn in \mbox{Figure \ref{fig-discoid-tree-illustration}}. Let $U$ be its closed subset that is colored black. (The black-red color scheme is so chosen because in the next section, closed subsets such as $U$ will arise as the retractions of tame loci.)

\begin{figure}[htb]\centering\includegraphics[scale=0.39]{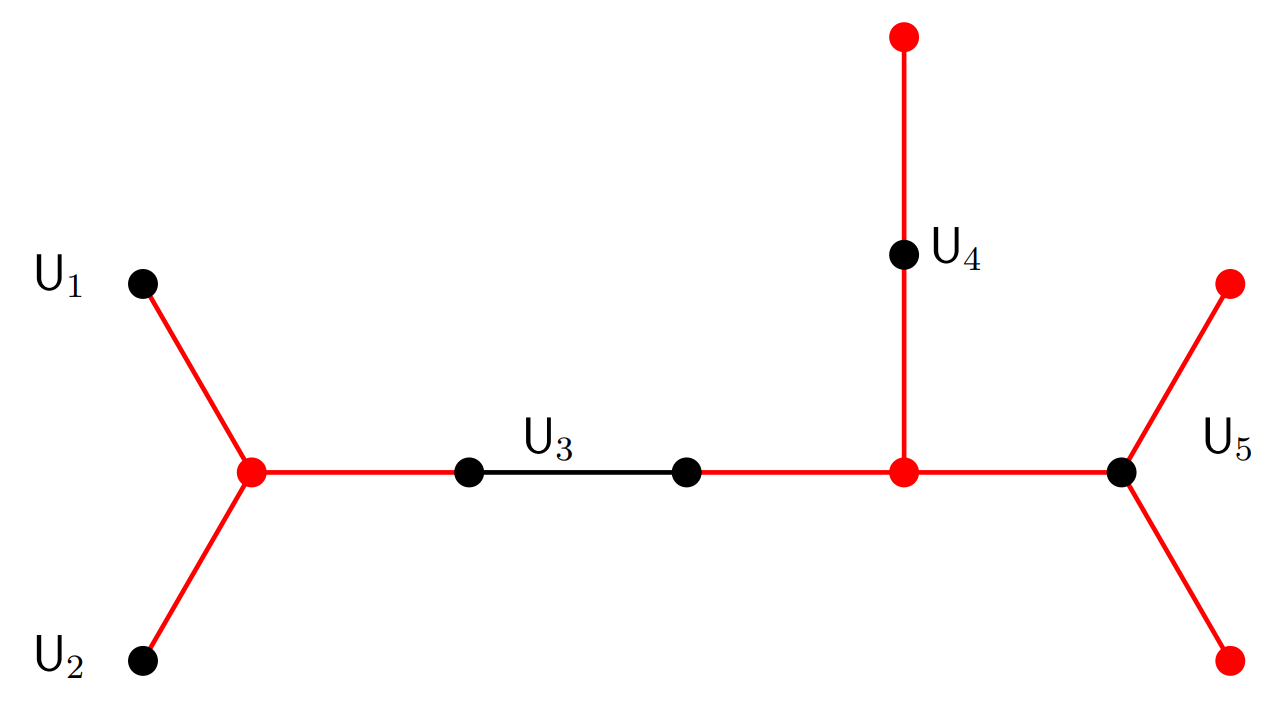}
\caption{A tree within $(\BP^1_{K_0})^{\an}$}\label{fig-discoid-tree-illustration}
\end{figure}  

The inverse image of $U$ under the retraction $\ret_T\colon(\BP^1_{K_0})^{\an}\to T$ consists of five components (or standard affinoid subdomains), numbered $\rU_1,\ldots,\rU_5$. The components $\rU_1,\rU_2$ are discoids; they are the complements of open discoids. The component $\rU_3$ is the complement of two open discoids. We call such an affinoid subdomain an \emph{anuloid} (cf.\ \cite[Definition 2.12]{micu}). The component $\rU_4$ is also an anuloid, but of radius $0$. Finally, $\rU_5$ is the complement of three discoids. One may visualize it as a closed discoid punctured by two open discoids of equal radius.
\end{Ex}

\begin{Rem}\label{rem-permanent-completion}
Let $T\subset(\BP_{K_0}^1)^{\an}$ be a tree. It may happen that the inverse image of the vertex set of $T$ under the map $\pi\colon Y^{\an}_K\to Y^{\an}$ is not the vertex set of a tree.

For example, consider the tree $T\subset(\BP^1_{\BQ_3})^{\an}$ spanned by the closed point $\xi$ with equation $x^2+3$ and the Gauss point (the boundary point of the closed disk of radius $0$ centered at $\xi$). The inverse image of this vertex set under $\pi$ consists of the Gauss point and the two closed points $\pm\sqrt{-3}$. The tree spanned by these three points also includes the Type II point that is the boundary point of the disk of radius $1/2$ centered at the closed points, see \mbox{Figure \ref{fig-permanent-completion}}.

\begin{figure}[htb]\centering\includegraphics[scale=0.3]{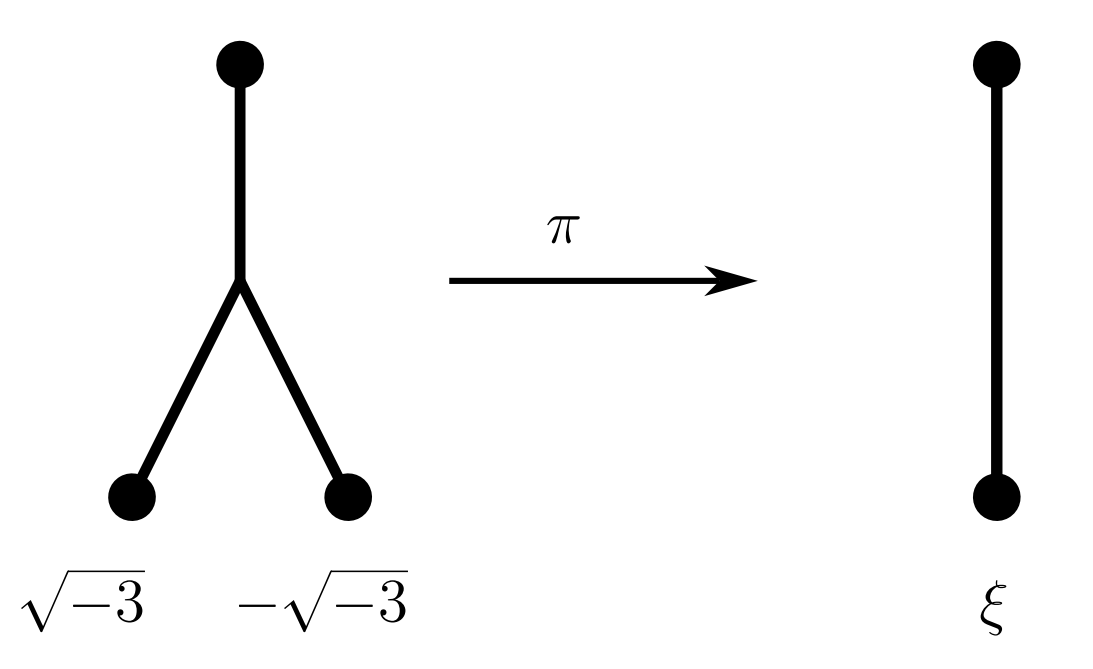}
\caption{The trees mentioned in \mbox{Remark \ref{rem-permanent-completion}}, illustrating that one needs to add another vertex to the inverse image of $T$}\label{fig-permanent-completion}
\end{figure}  

In practice, this problem will not concern us much, since MCLF has an inbuilt function \texttt{permanent\_completion} that refines a tree $T\subset(\BP^1_{K_0})^{\an}$ to a tree $T'$ such that the inverse image of the vertex set of $T'$ under $\pi$ is again a tree. We refer to the documentation of this function for details.
\end{Rem}

As an application of our discussion of discoids, we can now explain how to modify \mbox{Algorithm \ref{alg-interior-locus}} to the case of the discretely valued ground field $K_0$.

\begin{Algo}\label{alg-interior-discoid-locus}Input: A smooth plane quartic curve $Y$ over $K_0$, equipped with a degree-$3$ morphism $\phi\colon Y\to X=\BP^1_{K_0}$, given in the normal form discussed in \mbox{Section \ref{sec-quartic-normal-form}}.
\begin{enumerate}[(1)]
\item For each finite branch point $P\in X$ of $\phi$, do the following steps to compute an affinoid subdomain $\rU_P$ associated to $P$:
\begin{itemize}
\item[(1a)] Choose a root $x_0\in K$ of $\psi$, where $\psi\in K_0[x]$ is the monic irreducible polynomial generating the maximal ideal in $K_0[x]$ corresponding to $P$.
\item[(1b)] Write $t\coloneqq x-x_0$. As in \mbox{Section \ref{sec-quartics-computing-delta-directly}}, find a generator of $F_{Y_K}$ whose minimal polynomial over $K(t)$, say
\begin{equation*}
T^3+AT^2+BT+C,\qquad A,B,C\in K[t],
\end{equation*}
has coefficients of the form $A=a_2t^2+a_1t+a_0$, $B=b_3t^3+b_2t^2+b_1t+b_0$, and $C=c_4t^4+c_2t^2+c_1t$. Note that we may assume all coefficients of $A,B,C$ to lie in a finite intermediate extension $K/L/K_0$.
\item[(1c)] Using the family of Gauss valuations
\begin{equation*}
v_r\colon K(t)\to\BR\cup\{\infty\},\qquad\sum_id_it^i\mapsto\min_i\{v_K(d_i)+ir\},\qquad r\in\BR\cup\{\infty\},
\end{equation*}
define a function $\delta\colon[-\infty,\infty]\to[0,\infty]$ by
\begin{equation*}
r\mapsto\max\Big(0,\min\big(\frac{3}{2},\frac{3A(r)}{2}-\frac{C(r)}{2},\frac{3B(r)}{2}-C(r)\big)\Big).
\end{equation*}
\item[(1d)] Define
\begin{equation*}
\rU_P\coloneqq \red_{\Gamma_0}^{-1}\big(\{\xi_{\psi,\theta(r)}\mid r\in[0,\infty]:\delta(r)=0\}\big),
\end{equation*}
where $\Gamma_0\subset X^{\an}$ is the tree spanned by the branch points of $\phi$ and $\infty$.
\end{itemize}
\item Compute the affinoid subdomain $\rU_{K_0}$ as explained in Remark \ref{rem-mclf-affinoid-computation}.
\item Put $\rU^{\interior}_{K_0}\coloneqq\bigcup_P\rU_P$ and $\rU_{K_0}^{\tail}\coloneqq\rU_{K_0}\setminus\rU_{K_0}^{\interior}$.
\end{enumerate}
Result: For each point $\xi\in X_K$, the point $\xi$ lies in $\phi(Y^{\interior}_K)$ if and only if $\pi(\xi)\in\rU_{K_0}^{\interior}$ and lies in $\phi(Y^{\tail}_K)$ if and only if $\pi(\xi)\in\rU_{K_0}^{\tail}$.
\end{Algo}
\begin{proof}[Proof of correctness of the algorithm]
This follows immediately from the correctness of \mbox{Algorithm \ref{alg-interior-locus}}. Indeed, by \mbox{Proposition \ref{prop-discoid-splitting}}, the inverse image under $\pi$ of $\rU_{K_0}^{\interior}$ equals the union of the $\rU_{x_0}$ considered in \mbox{Algorithm \ref{alg-interior-locus}}, which equals the image of the interior locus $\phi(Y_K^{\interior})$. And by \mbox{Corollary \ref{cor-potential-tail-locus-descent}}, the inverse image under $\pi$ of $\rU_{K_0}\setminus\rU_{K_0}^{\interior}$ equals $\rU_K\setminus\phi(Y_K^{\interior})$, which is the same as the image of the tail locus $\phi(Y_K^{\tail})$.
\end{proof}

In summary, we have now found an affinoid subdomain
\begin{equation*}
\rU_{K_0}^{\tail}\cup\rU_{K_0}^{\interior}\subseteq(\BP_{K_0}^1)^{\an}
\end{equation*}
whose inverse image in $Y_K^{\an}$ is the tame locus associated to $\phi_K\colon Y_K\to X_K$. As explained in Remark \ref{rem-mclf-affinoid-computation}, we may describe it using retraction to a tree $T\subseteq(\BP_{K_0}^1)^{\an}$. Namely, we may take for $T$ any tree containing all boundary points of $\rU_{K_0}^{\tail}\cup\rU_{K_0}^{\interior}$. Now let $T'$ be the permanent completion (cf.\ \mbox{Remark \ref{rem-permanent-completion}}) of the tree spanned by $T$, by the branch points of $\phi$, and by $\infty$. Using the language of models, our main result is then the following. 

\begin{Kor}\label{cor-potentially-semistable-model}
Let $\CX_0$ be the model of $\BP^1_{K_0}$ that corresponds via \mbox{Proposition \ref{prop-julian-theorem}} to the set of Type II vertices of the tree $T'$. Let $\CY_0$ denote the normalization of $\CX_0$ in $F_Y$. Then the normalized base change $\CY\coloneqq(\CY_0)_K$ is a semistable model of $Y_K$.
\end{Kor}
\begin{proof}
The inverse image $\Gamma\coloneqq\pi^{-1}(T')\subset(\BP^1_K)^{\an}$ is a tree containing the tree $\Gamma_0$ spanned by the branch points of $\phi_K$ and the point $\infty$. Thus $\Sigma\coloneqq\phi_K^{-1}(\Gamma)$ is a graph containing $\phi^{-1}(\Gamma_0)$ and all boundary points of the tame locus $Y_K^{\tame}$. It follows from \mbox{Proposition \ref{prop-tame-locus-knows-all}} that $\Sigma$ is a skeleton of $Y_K$. Its set of Type II vertices corresponds via \mbox{Proposition \ref{prop-julian-theorem}} to the model $\CY$, which is therefore semistable.
\end{proof}

In the terminology of \mbox{Section \ref{sec-preliminaries-valuations}}, the model $\CY$ is \emph{potentially semistable}.

\section{Finding a field extension over which $Y$ has semistable reduction}

\label{sec-applications-field-extension}

In \mbox{Corollary \ref{cor-potentially-semistable-model}}, we have described how to find an $\CO_{K_0}$-model $\CY_0$ of $Y$ such that the normalized base change $\CY_K$ is a semistable model of $Y_K$. In this section, we will explain how to find a finite field extension $L/K_0$ such that $\CY_L$ is semistable. For this, we need not restrict ourselves to the case of plane quartic curves. Our discussion will apply to any class of curves for which we have a way to find potentially semistable $\CO_{K_0}$-models.

Let us write $G_{K_0}$ for the absolute Galois group of $K_0$ and $I_{K_0}$ for its inertia subgroup. Similarly we will write $G_L$ and $I_L$ for the absolute Galois group and inertia subgroup of a finite extension $L/K_0$.

The group $G_{K_0}$ acts naturally on $Y_K=Y\otimes_{K_0}K$ via the second factor. This action extends to certain $\CO_K$-models $\CY$ of $Y_K$. For example, if $g_Y\ge2$, there exists a \emph{stable model} $\CY^{\stab}$ of $Y_K$. It is the minimal semistable model of $Y_K$ and is characterized by the fact that all rational components of $\CY^{\stab}_s$ intersect the rest of $\CY^{\stab}_s$ in at least three points (with self-intersection points counting as two points). Uniqueness of the stable model shows that the action of $G_{K_0}$ on $Y_K$ extends uniquely to an action on $\CY^{\stab}$; see \cite[\mbox{Corollary 3.37}]{liu} for details. A similar argument shows that given any $\CO_{K_0}$-model $\CY_0$ of $Y$, the action of $G_{K_0}$ on $Y_K$ extends to the normalized base change $\CY\coloneqq(\CY_0)_K$.

\begin{Lem}\label{lem-deschamps}
Let $\CY_0$ be a potentially semistable $\CO_{K_0}$-model of a $K_0$-curve $Y$. If the inertia subgroup $I_{K_0}$ acts trivially on $\CY_s$, where $\CY=(\CY_0)_K$, then $Y$ has semistable reduction.
\end{Lem}

\begin{proof}
We use the results collected in \cite[Theorem 4.44]{liu}: Denoting by $K_0^{nr}$ the maximal unramified extension of $K_0$, there exists a minimal field extension $L/K_0^{nr}$ such that $Y_L$ has semistable reduction. Moreover, $L/K_0^{nr}$ is Galois and the action of $\Gal(L/K_0^{nr})$ on $\CY_s^{\stab}$, and hence on $\CY_s$, is faithful.

Now if the action of $I_{K_0}$ on $\CY_s$ is trivial, then the action of $\Gal(L/K_0^{nr})$ (which is a finite quotient of $I_{K_0}$) can only be faithful if $L=K_0^{nr}$. Thus $Y_{K_0}^{nr}$ has semistable reduction. In fact, because \'etale base change preserves normality (\cite[Tag 025P]{stacks}), the curve $Y$ itself must have semistable reduction: If $\CY_1$ is an $\CO_{K_0}$-model of $Y$ whose base change (with no normalization required!) $(\CY_1)_{K_0^{nr}}$ has semistable reduction, the model $\CY_1$ itself has semistable reduction.
\end{proof}

For the next proposition, we consider the familiar situation of a degree-$p$ morphism of $K_0$-curves $\phi\colon Y\to X$. Given an $\CO_K$-model $\CX$ of $X_K$, each component $Z\subseteq\CX_s$ corresponds to a Type II valuation $v_Z$ on $F_{X_K}$. We call the component $Z$ \emph{inseparable} if $v_Z$ is wildly ramified in $F_{Y_K}/F_{X_K}$. Otherwise, we call $Z$ \emph{separable}.

\begin{Prop}\label{prop-field-extension}
Let $\CX_0$ and $\CY_0$ be potentially semistable $\CO_{K_0}$-models of $X$ and $Y$ respectively such that $\CY_0$ is the normalization of $\CX_0$ in $F_{Y}$. Denote by $\CX\coloneqq(\CX_0)_K$ and $\CY\coloneqq(\CY_0)_K$ the normalized base changes to $K$. 

Suppose furthermore that $L/K_0$ is an extension with the following properties:
\begin{enumerate}[(a)]
\item The normalized base change $(\CX_0)_L$ is semistable
\item For each separable irreducible component $Z\subseteq \CX_s$ there exists an $L$-rational point $x_0$ on $X_L$ specializing to a smooth non-branch point of $\CY_s\to\CX_s$ such that the fiber $\phi^{-1}(x_0)$ consists of $L$-rational points on $Y_L$
\end{enumerate}
Then $Y_L$ has semistable reduction.
\end{Prop}
\begin{proof}
Since $(\CX_0)_L$ is semistable, the special fiber of the normalized base change $\CX$ is given simply by the base change $((\CX_0)_L)_s\otimes_{\lambda}\kappa$ (here $\lambda$ denotes the residue field of $\CO_L$). The action of $I_L$ on this special fiber is via the factor $\kappa$, so is trivial. Our goal is to show that $I_L$ also acts trivially on $\CY_s$, which by \mbox{Lemma \ref{lem-deschamps}} implies that $Y_L$ has semistable reduction.

Let $Z\subseteq\CX_s$ be an irreducible component. If $Z$ is inseparable, then clearly the unique component $W\subseteq\CY_s$ above $Z$ is fixed by the action of $I_L$. Otherwise, (b) shows that the components of $\CY_s$ above $Z$ are not permuted by the action of $I_L$. Let $W\subseteq\CY_s$ be one such component. By (b), the action of $I_L$ on $W$ fixes a non-ramification point $\overline{y}$ of the map $W\to Z$.

Assume that there exists some $\sigma\in I_L$ acting non-trivially on $W$. Then the map $W\to Z$ factors through $W/\langle\sigma\rangle$, and $\overline{y}$ is a ramification point of $W\to W/\langle\sigma\rangle$. This contradicts $\overline{y}$ being a non-ramification point, proving the claim that $I_L$ acts trivially on $\CY_s$ and finishing the proof.
\end{proof}

A word of caution: In the situation of Proposition \ref{prop-field-extension}, it is not necessarily the case that the model $(\CY_0)_L$ is semistable. This is because $\CY_s$ might contain more components than $\CY_s^{\stab}$, not all of which correspond to reduced components of the special fiber of $(\CY_0)_L$.

\begin{Ex}\label{ex-field-extension}
Let us consider once again the plane quartic from \mbox{Examples \ref{ex-lambda-as-radius-of-convergence}}, \mbox{\ref{ex-tame-locus-1}, and \ref{ex-tame-locus-1-revisited}},
\begin{equation*}
Y/\BQ_3\colon\quad y^3 + 3xy^2 - 3y - 2x^4 - x^2 - 1=0.
\end{equation*}
The $\CO_{K_0}$-model $\CX_0$ of Corollary \ref{cor-potentially-semistable-model} for which the normalization $\CY_0$ in $F_Y$ is potentially semistable corresponds to the boundary points of the three discoids
\begin{equation*}
\rD[x,1],\quad\rD[x,0],\quad\rD[x^2+1,3/4]
\end{equation*}
(see also \mbox{Example \ref{ex-worked-example-3}} below). Let us first define an extension $L_0$ of $K_0=\BQ_3$ by adjoining an element of valuation $1/4$, say $\sqrt[4]{3}$, and a square root $i$ of $-1$. Over $L_0$, the discoid $\rD[x^2+1,3/4]$ splits into the two disks $\rD[x-i,3/4]$ and $\rD[x+i,3/4]$, so $(\CX_0)_L$ is the model corresponding to the boundary points of the four disks $\rD[x,1]$, $\rD[x,0]$, $\rD[x-i,3/4]$, and $\rD[x+i,3/4]$. Since $1/4\in\Gamma_{L_0}$, the model $(\CX_0)_L$ is permanent (cf.\ the proof of \mbox{Lemma \ref{lem-ramification-index-and-multiplicity}}), and hence semistable. Thus we have dealt with Condition (a) of \mbox{Proposition \ref{prop-field-extension}}.

Next, we choose points $x_0$ as in Condition (b) of \mbox{Proposition \ref{prop-field-extension}}. Recall from \mbox{Example \ref{ex-tame-locus-1-revisited}} that above each component $Z$ corresponding to a boundary point of $\rD[x,1]$, $\rD[x-i,3/4]$, and $\rD[x+i,3/4]$ lies a single component $W$, given by an Artin-Schreier equation. In particular, $W\to Z$ is only ramified at $\infty$, and we may take for the $x_0$ points contained in the open disks $\rD(x,1)$, $\rD(x-i,3/4)$, and $\rD(x-i,3/4)$. The simplest choice is the centers $x=0$, $x=i$, and $x=-i$.

A field extension $L$ over which $Y$ acquires semistable reduction is then obtained by adjoining to $L_0$ the coordinates of the fibers $\phi^{-1}(0)$, $\phi^{-1}(i)$, and $\phi^{-1}(-i)$. Note that there is no condition for the disk $\rD[x,0]$, since the corresponding component is inseparable.
\end{Ex}

\begin{Rem}\label{rem-purity}
To find a point $x_0\in X_L(L)$ specializing to a smooth non-branch point on a given irreducible component $Z\subseteq\CX_s$ as demanded by \mbox{Proposition \ref{prop-field-extension}(b)}, one may use the theorem on \emph{purity of the branch locus} \cite[\mbox{Tag 0BMB}]{stacks}.

Indeed, any point $x_0\in X_L(L)$ whose reduction to the special fiber of $(\CX_0)_L$ is a smooth point on $Z$ and not the same as the reduction of a branch point of $\phi$ will work. To see this, denote by $\overline{x}$ the reduction of $x_0$ and let $\overline{y}$ be a preimage of $\overline{x}$ under the map $(\CY_0)_L\to(\CX_0)_L$. The local rings $\CO_{(\CX_0)_L,\overline{x}}$ and $\CO_{(\CY_0)_L,\overline{y}}$ are isomorphic to $\CO_L\llbracket T\rrbracket$ (see for example \cite[\mbox{Proposition 3.4}]{arzdorf-wewers}); in particular, they are regular local rings of dimension $2$. The prime ideals of $\CO_L\llbracket T\rrbracket$ are classified in \cite[Lemma 3.1]{helminck}. There are two types of codimension one prime ideals in $\CO_{(\CX_0)_L,\overline{x}}$: one corresponding to the generic point of $Z$, the others corresponding to closed points on $X_L$ specializing to $\overline{x}$. Analogous considerations apply to $\CO_{(\CY_0)_L,\overline{y}}$.

Since $\overline{x}$ is not the specialization of a branch point of $\phi$ and since the component $Z$ is separable, the induced map $\CO_{(\CY_0)_L,\overline{y}}\to\CO_{(\CX_0)_L,\overline{x}}$ is \'etale at codimension $1$ primes. It follows from purity of the branch locus that it is outright \'etale.
\end{Rem}

The MCLF package is able to perform the steps outlined in this section automatically. We will make extensive use of this in the examples of the following section.

\section{Examples}\label{sec-applications-examples}

In this section we present a number of worked examples. In each, we explain how to compute the semistable reduction of a plane quartic $Y$ over $\BQ_3$. All the curves and the corresponding covers $\phi\colon Y\to\BP^1_{K_0}$ are in the normal form of \mbox{Section \ref{sec-quartic-normal-form}}. A variety of different reduction types and combinatorial configurations of the tame locus is included.

At last, we let the computer do as much of the work as possible. In practice this means that we use our implementation of \mbox{Algorithm \ref{alg-interior-discoid-locus}} to compute the image of the tame locus in $(\BP_{K_0}^1)^{\an}$. The boundary points of this affinoid subdomain then determine a potentially semistable model of $Y$. Afterwards, we use the existing functionality of MCLF to compute equations for the components corresponding to the valuations we have found. We refer to \mbox{Section \ref{sec-appendix-implementation}} of the appendix for a brief overview of our implementation.

\begin{Ex}\label{ex-worked-example-1}The plane quartic
\begin{equation*}
Y/\BQ_3\colon\quad F(y)=y^3 + (2x^3 + 3x^2)y - 3x^4 - 2x^2 - 1=0.
\end{equation*}
We ask Sage to compute the image in $(\BP^1_{K_0})^{\an}$ of the tame locus associated to $Y$, using the following commands:
\begin{minted}[tabsize=2,breaklines,fontsize=\small]{text}
sage: from mclf import *
sage: R.<x,y> = QQ[]
sage: v = QQ.valuation(3)
sage: F = y^3 + (2*x^3 + 3*x^2)*y - 3*x^4 - 2*x^2 - 1
sage: Y = SmoothProjectiveCurve(F)
sage: M = Quartic3Model(Y, v)
sage: M.tame_locus()
Affinoid with 3 components:
Elementary affinoid defined by
v(1/x) >= 0
Elementary affinoid defined by
v(x^9 - 9*x^8 + 27*x^7 + 54*x^6 + 81*x^5 + 54*x^4 - 27*x^2 + 27/5) >= 6
Elementary affinoid defined by
v(x) >= 3/4
\end{minted}

Recall \mbox{Section \ref{sec-applications-line}} and \mbox{Example \ref{ex-discoid-retraction}} in particular for how we represent affinoid subdomains of $(\BP_{K_0}^1)^{\an}$ using retraction to a tree. In this case, the output is an affinoid subdomain consisting of three discoids, which may be described by retraction to the following tree:

\begin{figure}[h]\centering\includegraphics[scale=0.38]{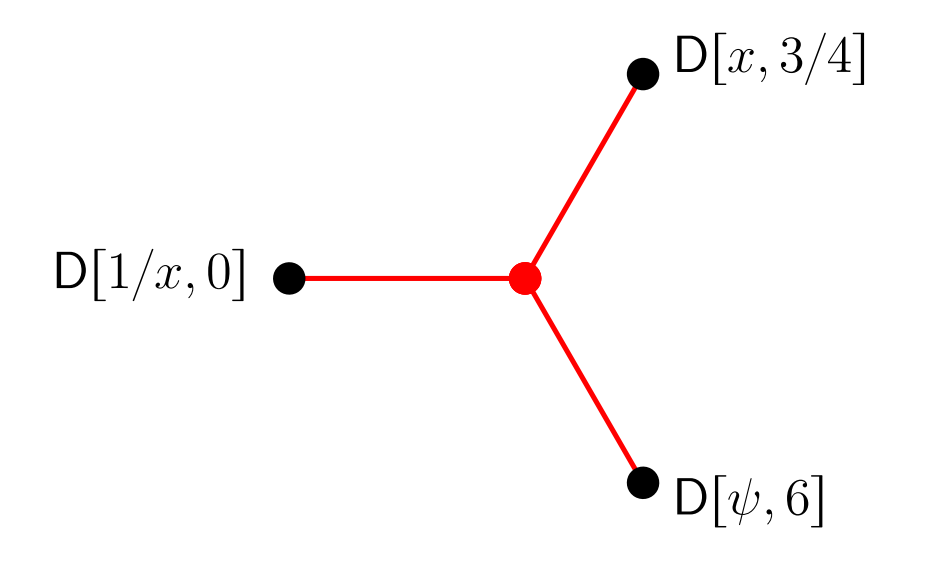}
\end{figure}  
where
\begin{equation*}
\psi=x^9 - 9x^8 + 27x^7 + 54x^6 + 81x^5 + 54x^4 - 27x^2 + \frac{27}{5}.
\end{equation*}
The polynomial $\psi$ is not equal to the discriminant $\Delta_F$ of $F$, but we have $\rD[\psi,6]=\rD[\Delta_F,6]$, and this discoid contains the branch locus except for the point $\infty$. The image of the tail locus is given by the disk of radius $3/4$:
\begin{minted}[tabsize=2,breaklines,fontsize=\small]{text}
sage: M.tail_locus()
Elementary affinoid defined by
v(x) >= 3/4
\end{minted}

It follows easily from \mbox{Proposition \ref{prop-discoid-splitting}} that $\rD[\psi,6]$ splits into $9$ disjoint disks over $K$. Thus the boundary points of these disks play no role in finding a semistable reduction of $Y$. We use MCLF to verify that the model corresponding to the boundary points of the two disks $\rD[x,0]$ and $\rD[x,3/4]$ indeed is potentially semistable. Essentially MCLF works by computing extensions of the involved valuations, that is, using \mbox{Proposition \ref{prop-julian-normalization}}. Continuing the above code fragment:
\begin{minted}[tabsize=2,breaklines,fontsize=\small]{text}
sage: X = M._X
sage: T = BerkovichTree(X)
sage: T.add_point(X.point_from_discoid(x, 0))
sage: T.add_point(X.point_from_discoid(x, 3/4))
sage: R = ReductionTree(Y, v, T)
sage: R.is_semistable()
True
\end{minted}
MCLF also provides equations for the reduction curves:
\begin{minted}[tabsize=2,breaklines,fontsize=\small]{text}
sage: R.inertial_components()[0].upper_components()[0].component()
the smooth projective curve with Function field in u1 defined by u1^3 + 2*x^3*u1 + x^2 + 2
sage: R.inertial_components()[1].upper_components()[0].component()
the smooth projective curve with Function field in u2 defined by u2^3 + u2 + x^2
\end{minted}
The reduction curve of the point above the boundary point of $\rD[x,3/4]$ is the Artin-Schreier curve of genus $1$ with equation
\begin{equation*}
y^3+y+x^2=0.
\end{equation*}
The reduction curve above the boundary point of $\rD[x,0]$ is of genus $2$ and is given by 
\begin{equation*}
y^3-x^3y+x^2-1=0.
\end{equation*}
We note that it has ordinary ramification over $x=\infty$ (the specialization of the ordinary branch point $\infty$ of $\phi$) and that it is totally ramified over $x=0$ (reflecting the fact that $\phi$ is inseparable on the branch pointing towards $x=0$).

In summary, we have found one component of genus $1$ and one component of genus $2$ on the special fiber of the model determined by the boundary points of the tame locus. To obtain the stable reduction, we can contract all other components. This corresponds to discarding the central red vertex and the boundary point of $\rD[\psi,6]$ in the tree above, keeping only the boundary points of $\rD[1/x]$ and $\rD[x,3/4]$.

\begin{Ex}\label{ex-worked-example-2}
The plane quartic
\begin{equation*}
Y/\BQ_3\colon\quad y^3 - 3y^2 + (-3x^2 - 2x)y + 3x^4 - 3x - 1=0,
\end{equation*}
previously considered in \mbox{Example \ref{ex-computing-delta-directly}}. Using the same commands as in the previous example, we find that the image in $(\BP^1_{K_0})^{\an}$ of the tame locus may be described by retraction to the following tree:
\begin{figure}[H]\centering\includegraphics[scale=0.38]{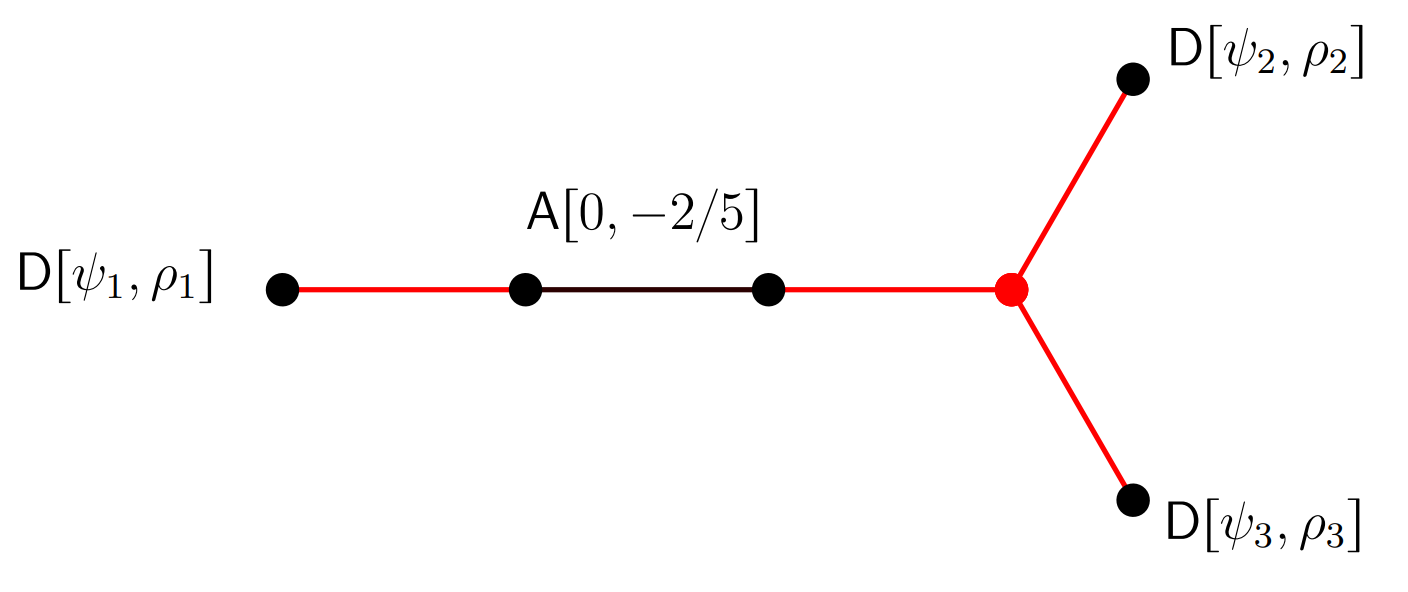}
\end{figure}  
Note that during instantiation of the \texttt{Quartic3Model} associated to $Y$, a coordinate change is performed replacing $x$ with $3x$. However, the image of the tame locus with respect to the original normalization can be obtained using the method \texttt{rescaled\_tame\_locus}:
\begin{minted}[tabsize=2,breaklines,fontsize=\small]{text}
sage: from mclf import *
sage: R.<x,y> = QQ[]
sage: v = QQ.valuation(3)
sage: F = y^3 - 3*y^2 + (-3*x^2 - 2*x)*y + 3*x^4 - 3*x - 1
sage: Y = SmoothProjectiveCurve(F)
sage: M = Quartic3Model(Y, v)
We do a base change so that the branch locus is integral and now consider the smooth projective curve with Function field in y defined by y^3 - 27*y^2 + (-27*x^2 - 54*x)*y + 27*x^4 - 729*x - 729
sage: M.tame_locus()
Affinoid with 4 components:
Elementary affinoid defined by
v(x) >= 3/5
v(1/x) >= -1
Elementary affinoid defined by
v(x^2 - x - 4) >= 3
Elementary affinoid defined by
v(x^3 + x^2 + 10*x + 8) >= 3
Elementary affinoid defined by
v(x^3 + 27/8*x^2 + 2673/2*x + 3645/22) >= 13
sage: M.rescaled_tame_locus()
Affinoid with 4 components:
Elementary affinoid defined by
v(1/x) >= 0
v(x) >= -2/5
Elementary affinoid defined by
v(x^3 - 90*x^2 - 54/7*x - 1485) >= 10
Elementary affinoid defined by
v((-63*x^2 + 3/4*x + 1)/x^2) >= 5
Elementary affinoid defined by
v((27/8*x^3 + 153*x^2 + 24*x + 1)/x^3) >= 6
\end{minted}
There are no tail discoids. The discoids $\rD[\psi_i,\rho_i]$, $i=1,2,3$, each contain branch points of $\phi$ (reflecting the structure of the branch locus of $\phi$ as explained in \mbox{Example \ref{ex-computing-delta-directly}}). The anulus
\begin{equation*}
\rA[0,-2/5]=\big\{\xi\in(\BP_{K_0}^1)^{\an}\mid 0\ge v_\xi(x)\ge-\frac{2}{5}\big\}
\end{equation*}
coincides with the interval in \mbox{Figure \ref{fig-example-delta-graph-desmos}} on which $\delta=0$. It turns out that its boundary points are all that we need:
\begin{minted}[tabsize=2,breaklines,fontsize=\small]{text}
sage: FX.<x> = FunctionField(QQ)
sage: X = BerkovichLine(FX, v)
sage: T = BerkovichTree(X)
sage: T.add_point(X.point_from_discoid(x, 0))
sage: T.add_point(X.point_from_discoid(x, -2/5))
sage: R = ReductionTree(Y, v, T)
sage: R.inertial_components()[0].upper_components()[0].component()
the smooth projective curve with Function field in u1 defined by u1^3 + x*u1 + 2
sage: R.inertial_components()[1].upper_components()[0].component()
the smooth projective curve with Function field in u2 defined by u2^3 + 1/x*u2 + 1/x^4
sage: R.is_semistable()
True
\end{minted}

Let us denote the boundary points of $\rA[0,-2/5]$ by $\xi_0$ and $\xi_{-2/5}$. The reduction curve of the point above $\xi_0$, given by
\begin{equation*}
y^3+xy+2=0,
\end{equation*}
is totally ramified at $x=0$ and has ordinary ramification at $\infty$. Conversely, the reduction curve of the point above $\xi_{-2/5}$, given by
\begin{equation*}
y^3+xy+x^4=0,
\end{equation*}
is totally ramified at $\infty$ and has ordinary ramification at $x=0$. This reflects the fact that the open interval $(\xi_0,\xi_{-2/5})$ in $\rA[0,-2/5]$ consists of tame topological branch points, while the adjacent intervals consist of wild topological branch points. The reduction curve of the point above $\xi_{-2/5}$ has genus $2$, while the reduction curve of the point above $\xi_0$ has genus $0$. A loop lies above the interval $[\xi_0,\xi_{-2/5}]$. To double check this, we let MCLF verify that there are two points above $\xi_{-1/5}$:
\begin{minted}[tabsize=2,breaklines,fontsize=\small]{text}
sage: T = BerkovichTree(X)
sage: T.add_point(X.point_from_discoid(x, -1/5))
sage: R = ReductionTree(Y, v, T)
sage: for U in R.inertial_components()[1].upper_components():
....:   U.component()
the smooth projective curve with Rational function field in x over Finite Field of size 3
the smooth projective curve with Function field in u2 defined by u2^2 + 1/x
\end{minted}

It follows that the stable reduction of $Y$ has only one component, which has geometric genus $2$ and one self-intersection point. Note that as announced in \mbox{Remark \ref{rem-mistake-in-ctt}}, the point $\xi_{-2/5}$ is not a vertex of the tree spanned by the branch points of $\phi$ even though we need its preimage to obtain a skeleton of $Y$.
\end{Ex}

\end{Ex}

\begin{Ex}\label{ex-worked-example-3}
The plane quartic
\begin{equation*}
Y/\BQ_3\colon\quad y^3 + 3xy^2 - 3y - 2x^4 - x^2 - 1=0,
\end{equation*}
previously considered in \mbox{Examples \ref{ex-lambda-as-radius-of-convergence}}, \mbox{\ref{ex-tame-locus-1}, \ref{ex-tame-locus-1-revisited}, and \ref{ex-field-extension}}. Again using the same commands as in \mbox{Example \ref{ex-worked-example-1}}, we find that the image in $(\BP^1_{K_0})^{\an}$ of the tame locus may be described by retraction to the following tree:

\begin{figure}[H]\centering\includegraphics[scale=0.38]{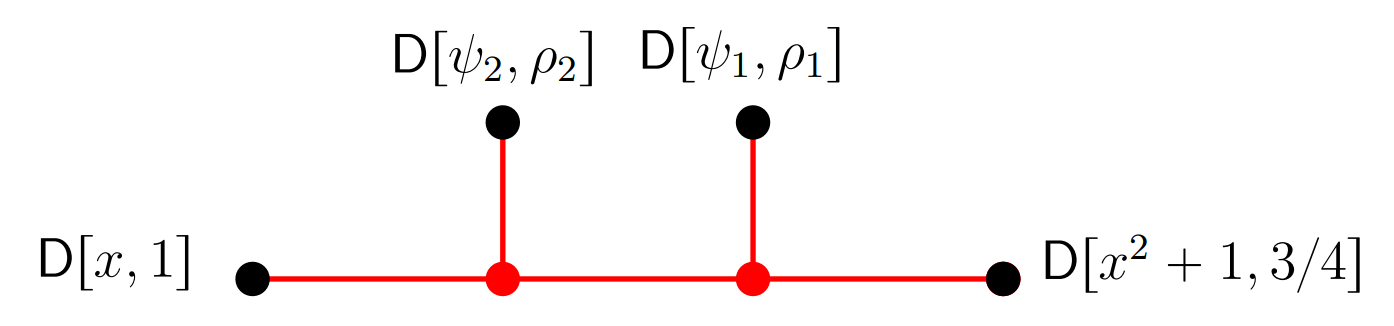}
\end{figure}  

The discoids on the left and right are tail discoids. The discoid $\rD[x^2+1,3/4]$ splits into two disks over $K$. From this the structure of the stable reduction of $Y$ is immediately clear. Since each tail component has a boundary point of genus $\ge1$, there can be no further points of positive genus nor any loops. Thus the stable reduction is a ``comb'' consisting of one inseparable rational component intersecting with three genus $1$ curves, as follows:

\begin{figure}[H]\centering\includegraphics[scale=0.3]{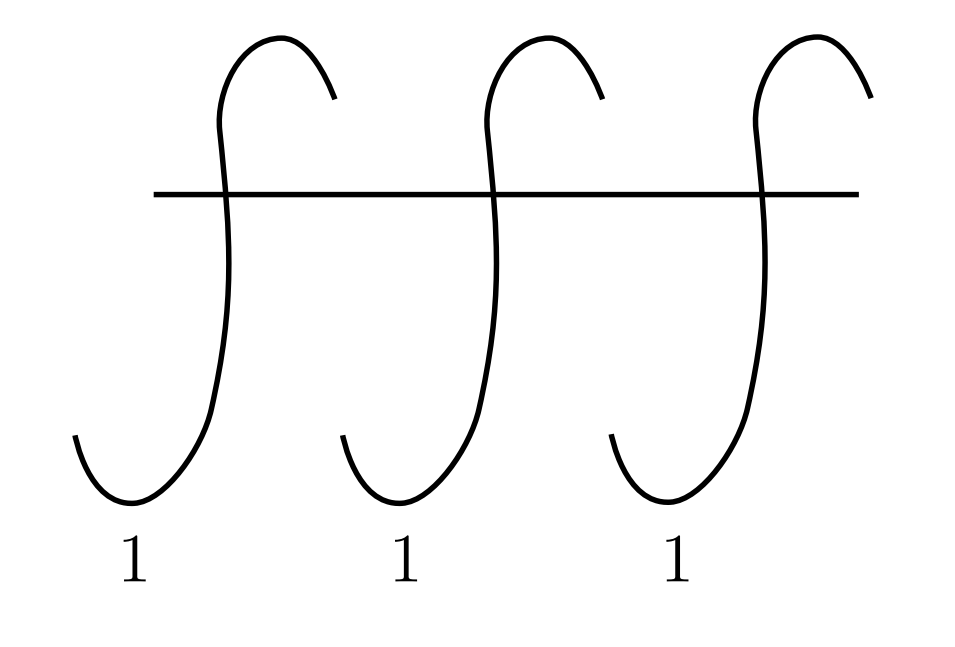}
\end{figure}  

The discoids $\rD[\psi_1,\rho_1]$ and $\rD[\psi_2,\rho_2]$ contain branch points. In fact, we already determined the structure of the branch locus of $\phi$ in \mbox{Example \ref{ex-lambda-as-radius-of-convergence}}. The discoid $\rD[\psi_1,\rho_1]$ has degree $6$, while $\rD[\psi_2,\rho_2]$ has degree $2$.

We can now also check our experimental hypothesis from \mbox{Example \ref{ex-lambda-as-radius-of-convergence}} that $\lambda(x_0)=1$, where $x_0=0$. Indeed, $x_0$ is contained in $\rD[x,1]$ and it follows from the structure of the tree above and from \mbox{Proposition \ref{prop-simple-lambda-consequence}} that $\lambda(x_0)=1$.
\end{Ex}

\begin{Ex}
\label{ex-worked-example-5}
The plane quartic
\begin{equation*}
Y/\BQ_3\colon\quad y^3 + x^3y + x^4 + 1=0,
\end{equation*}
previously considered in \mbox{Example \ref{ex-admissible-no-good}}. The image in $(\BP^1_{K_0})^{\an}$ of the tame locus may be described by retraction to the following tree:
\begin{figure}[H]\centering\includegraphics[scale=0.38]{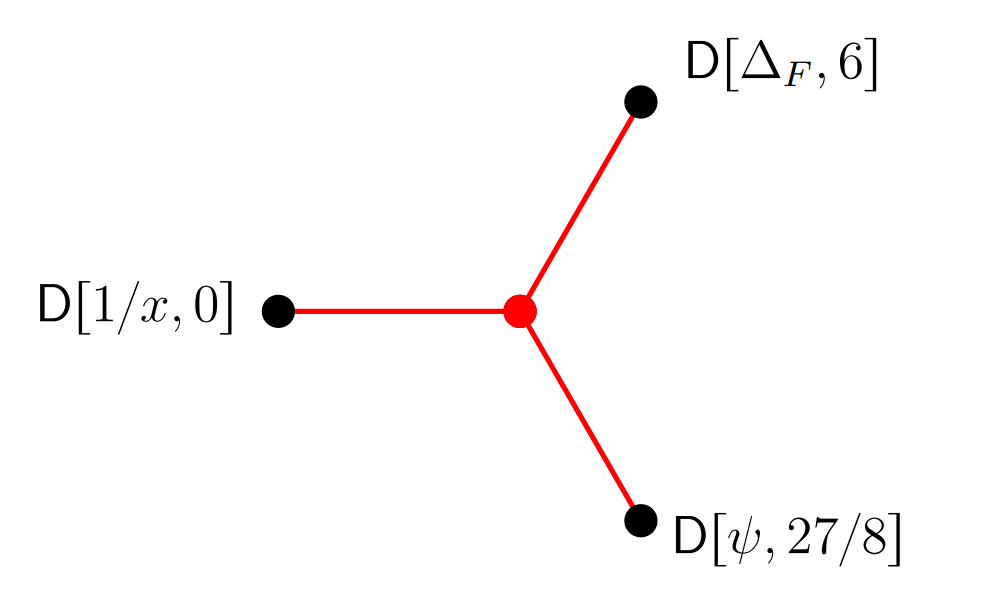}
\end{figure}  
The discoid $\rD[\Delta_F,6]$ contains the single finite branch point of $\phi$, which is of degree $9$. In fact, the (image in $(\BP^1_{K_0})^{\an}$ of the) Type II point we used in \mbox{Example \ref{ex-admissible-no-good}} lies above the path connecting $\rD[\Delta_F,6]$ to the central red vertex.

The only important component in this example is the tail discoid $\rD[\psi,27/8]$. We have
\begin{equation*}
\psi=x^9 + 27x + 54.
\end{equation*}
As expected, the reduction curve at the point above the boundary point of $\rD[\psi,27/8]$ is given by an Artin-Schreier equation
\begin{equation*}
y^3+y=x^4.
\end{equation*}
It has genus $3$ and is only ramified at $\infty$. Thus $Y$ has potentially good reduction.
\end{Ex}

\begin{Ex}\label{ex-worked-example-6}
The plane quartic
\begin{equation*}
Y/\BQ_3\colon\quad F(y)=y^3 + y^2 + (x + 1)y + 3x^4 + 2x^3 + 2x=0.
\end{equation*}
The map $\phi\colon Y\to\BP_{K_0}^1$ is totally ramified above $\infty$. There are eight finite branch points corresponding to the eight zeros of the discriminant $\Delta_F$. Two of these have valuation $1/2$, while the others have negative valuation.

The image in $(\BP^1_{K_0})^{\an}$ of the tame locus may be described by retraction to the following tree:

\begin{figure}[H]\centering\includegraphics[scale=0.38]{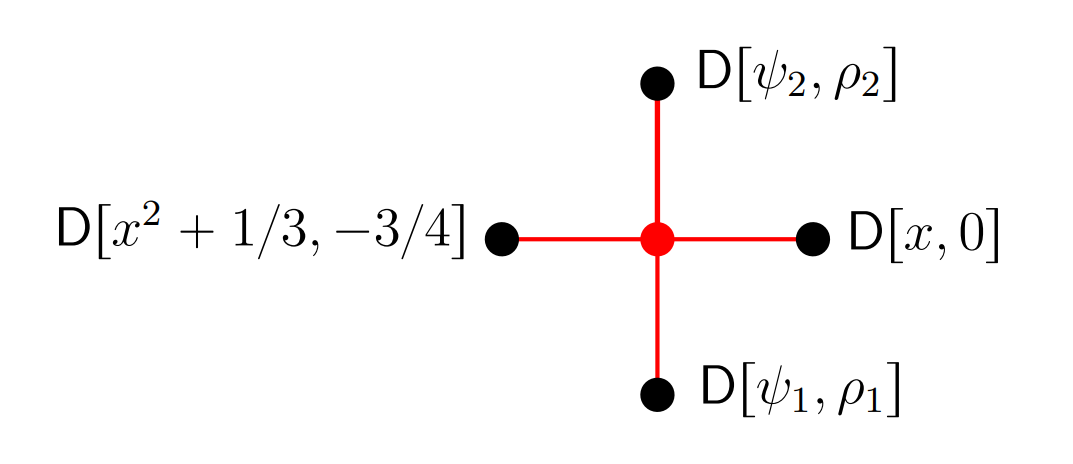}
\end{figure}  

The discoid on the left is a tail discoid. Over $K$, it splits into two disks centered at $\pm1/\sqrt{-3}$. The reduction curves at the points above the boundary points of these tail disks are the Artin-Schreier curves of genus $1$ given by
\begin{equation*}
y^2+y+x^2=0.
\end{equation*}
The discoids $\rD[\psi_1,\rho_1]$ and $\rD[\psi_2,\rho_2]$ split into $4$ and $2$ disks over $K$ respectively, so they cannot contribute to the stable reduction of $Y$. We need to keep the central red vertex however, since it is contained in the permanent completion (cf.\ \mbox{Remark \ref{rem-permanent-completion}}) of the tree spanned by $\rD[x^2+1/3,-3/4]$ and $\rD[x,0]$. It is the boundary point of the disk $\rD[x,-1/2]$. 

We furthermore add the boundary point of $\rD[x,1/2]$, which separates the two branch points of valuation $1/2$. In general it is necessary to use a tree separating the branch points of $\phi$ --- compare for example the hypotheses in \mbox{Proposition \ref{prop-tame-locus-knows-all}}, which demand that $\Sigma_1$ contain the boundary points of the tame locus \emph{and the inverse image of the tree spanned by the branch points of $\phi$}.

In summary, we take the model $\CX$ corresponding to the boundary points of the four discoids $\rD[x^2+1/3,-3/4]$, $\rD[x,-1/2]$, $\rD[x,0]$, and $\rD[x,1/2]$. We check that these determine a potentially semistable model. Note that we add the boundary point of $\rD[x,-1/2]$ using MCLF's function \texttt{permanent\_completion}.
\begin{minted}[tabsize=2,breaklines,fontsize=\small]{text}
sage: from mclf import *
sage: R.<x,y> = QQ[]
sage: v = QQ.valuation(3)
sage: F = y^3 + y^2 + (x + 1)*y + 3*x^4 + 2*x^3 + 2*x
sage: Y = SmoothProjectiveCurve(F)
sage: FX.<x> = FunctionField(QQ)
sage: X = BerkovichLine(FX, v)
sage: T = BerkovichTree(X)
sage: T.add_point(X.point_from_discoid(x, 0))
sage: T.add_point(X.point_from_discoid(x^2 + 1/3, -3/4))
sage: T.add_point(X.point_from_discoid(x, 1/2))
sage: T.permanent_completion()
sage: T.vertices()
[Point of type II on Berkovich line, corresponding to v(x) >= 0,
 Point of type II on Berkovich line, corresponding to v(1/x) >= 1/2,
 Point of type II on Berkovich line, corresponding to v(3*x^2 + 1) >= 1/4,
 Point of type II on Berkovich line, corresponding to v(x) >= 1/2]
sage: R = ReductionTree(Y, v, T)
sage: R.is_semistable()
 True
\end{minted}

The reduction curves at the points above the boundaries of each of $\rD[x,-1/2]$, $\rD[x,0]$, and $\rD[x,1/2]$ are all rational. The ones corresponding to the latter two disks intersect in two points. Thus the stable reduction of $Y$ looks as follows:
\begin{figure}[H]\centering\includegraphics[scale=0.3]{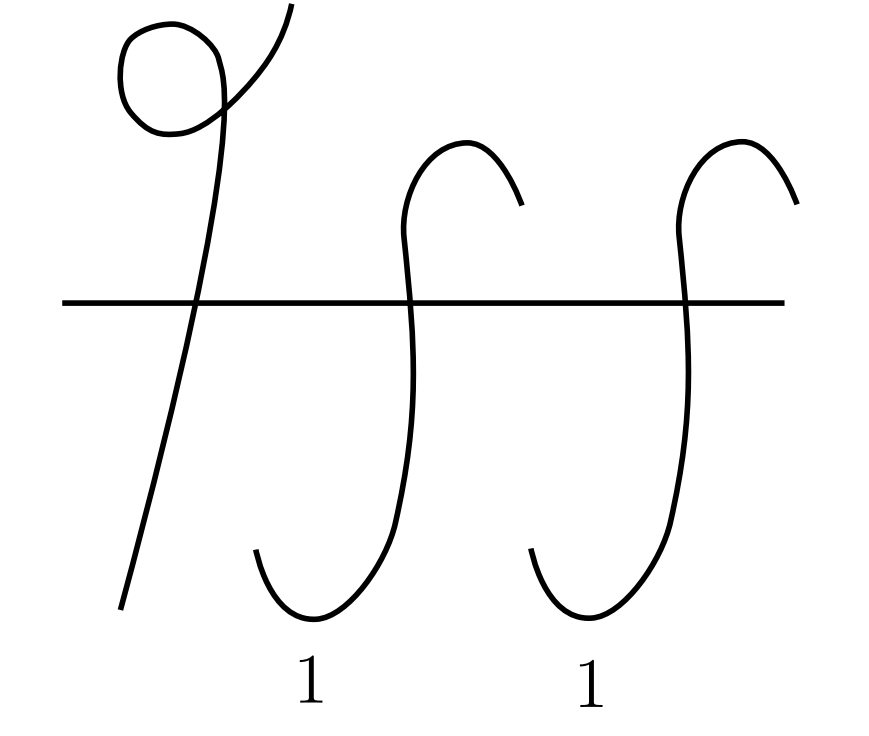}
\end{figure}

\end{Ex}

\section{An example over $\BQ_3(\zeta_3)$}\label{sec-modular-curve}

In this section we apply our results to a slightly more involved example. Namely, we compute the semistable reduction of a particular plane quartic curve that appears in the attempts of Rouse, Sutherland, and Zureick-Brown to compute the rational points on the non-split Cartan modular curve $X_{\ns}^+(27)$, see \cite{rouse-sutherland-zureickbrown}. The semistable reduction was given in \cite{ossen}, but the method of finding it (essentially \mbox{Algorithm \ref{alg-interior-discoid-locus}}) was left out.

Let $K_0=\BQ_3(\zeta_3)$ be the field obtained by adjoining to $\BQ_3$ a primitive third root of unity and let $K$ be the completion of an algebraic closure of $K_0$. As always, we normalize the valuation $v_K\colon K\to\BQ\cup\{\infty\}$ so that $v_K(3)=1$, that is to say $v_K(\zeta_3-1)=1/2$. The curve we are interested in is the plane quartic $K_0$-curve cut out of $\BP_{K_0}^2$ by the equation
\begin{equation*}
\begin{aligned}
x^4 &+ (\zeta_3 - 1)x^3y + (3\zeta_3 + 2)x^3z - 3x^2z^2 + (2\zeta_3 + 2)xy^3\\ &- 3\zeta_3 xy^2z
+ 3\zeta_3 xyz^2 - 2\zeta_3 xz^3 - \zeta_3 y^3z\\ &+ 3\zeta_3 y^2z^2 + (-\zeta_3 + 1)yz^3 + (\zeta_3 + 1)z^4=0.
\end{aligned}
\end{equation*}
Note that $[0:1:0]$ is a rational point on $Y$. We achieve the normal form discussed in \mbox{Section \ref{sec-quartic-normal-form}} by replacing $z$ with $z+x(2\zeta_3+2)/\zeta_3$. Dehomogenizing yields the affine equation
\begin{equation*}
Y\colon\quad F(y)=y^3+Ay^2+By+C=0,
\end{equation*}
where
\begin{equation*}
A=(6\zeta_3 + 12)x^2+ (36\zeta_3 + 9)x- 27,
\end{equation*}
\begin{equation*}
B=(9\zeta_3 - 18)x^3+ (-108\zeta_3 - 108)x^2+ (-162\zeta_3 + 81)x+ (81\zeta_3 + 162),
\end{equation*}
\begin{IEEEeqnarray*}{rCl}
C&=&(27\zeta_3 - 243)x^4+ (-1458\zeta_3 - 999)x^3\\&+& (-1215\zeta_3 + 1701)x^2+ (1944\zeta_3 + 2430)x+ 729\zeta_3.
\end{IEEEeqnarray*}
We use the usual projection map $\phi\colon Y\to\BP^1_{K_0}$. The discriminant $\Delta_F$ has degree $10$, so by \mbox{Lemma \ref{lem-quartic-branch-locus}} the point $\infty$ is not a branch point of $\phi$. The branch locus of $\phi$ consists of one point of degree $1$ and one point of \mbox{degree $9$}.

For convenience, the curve $Y$ is defined in \mbox{Code Listing \ref{list-modular}} and the image in $(\BP^1_{K_0})^{\an}$ of its tame locus is computed. It may be described in terms of reduction to the following tree:
\begin{figure}[H]\centering\includegraphics[scale=0.35]{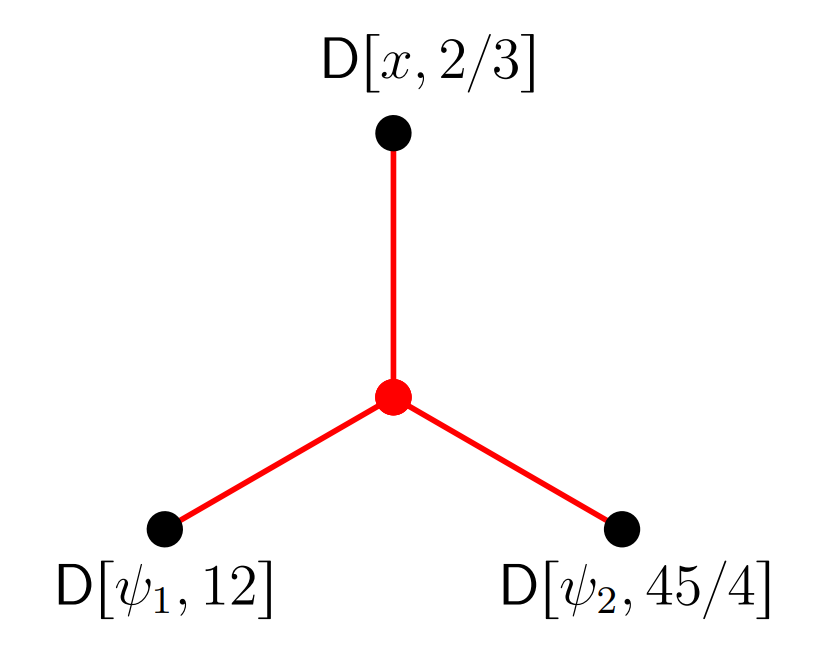}
\end{figure}  
The disk $\rD[x,2/3]$ contains $\infty$ and the branch point of degree $1$. The polynomial
\begin{equation*}
\begin{aligned}
\psi_1&=x^9 + (9\zeta_3 + 18)x^8 + (54\zeta_3 - 27/2)x^7 + (54\zeta_3 + 108)x^6 \\&+ 486\zeta_3x^5+ (-729\zeta_3 + 3645)x^4 + (729\zeta_3 + 13851)x^3 \\&+ (37179\zeta_3 + 15309)x^2 + (-6561\zeta_3/2 + 72171)x + 2187\zeta_3/13 + 15309/2
\end{aligned}
\end{equation*}
is the factor of $\Delta_F$ of degree $9$, so $\rD[\psi_1,12]$ contains the other branch point of $\phi$. The last component, the tail discoid $\rD[\psi_2,45/4]$, is the most important; we have
\begin{equation*}
\begin{aligned}
\psi_2&=x^9 + (9\zeta_3 - 9)x^8 + (54\zeta_3 + 27)x^7 + (54\zeta_3 - 27/2)x^6\\ &+ (243\zeta_3 + 972)x^5 + 729\zeta_3 x^4+ (2916\zeta_3 - 1458)x^3 + (37179\zeta_3 \\ &+ 41553)x^2 + (6561\zeta_3 + 6561/8)x - 63423\zeta_3 + 155277.
\end{aligned}
\end{equation*}
We also note that the red vertex in the middle of the tree above is the boundary point of the discoid $\rD[\psi_1,11]=\rD[\psi_2,11]$. The splitting behavior of all these discoids is particularly interesting. While the discoids $\rD[\psi_2,45/4]$ and $\rD[\psi_2,11]$ split into three disks over $K$ each containing three of the nine roots of $\psi_2$, the discoid $\rD[\psi_1,12]$ splits into nine disks, each containing one root of $\psi_1$. 

\begin{figure}[htb]\centering\includegraphics[scale=0.55]{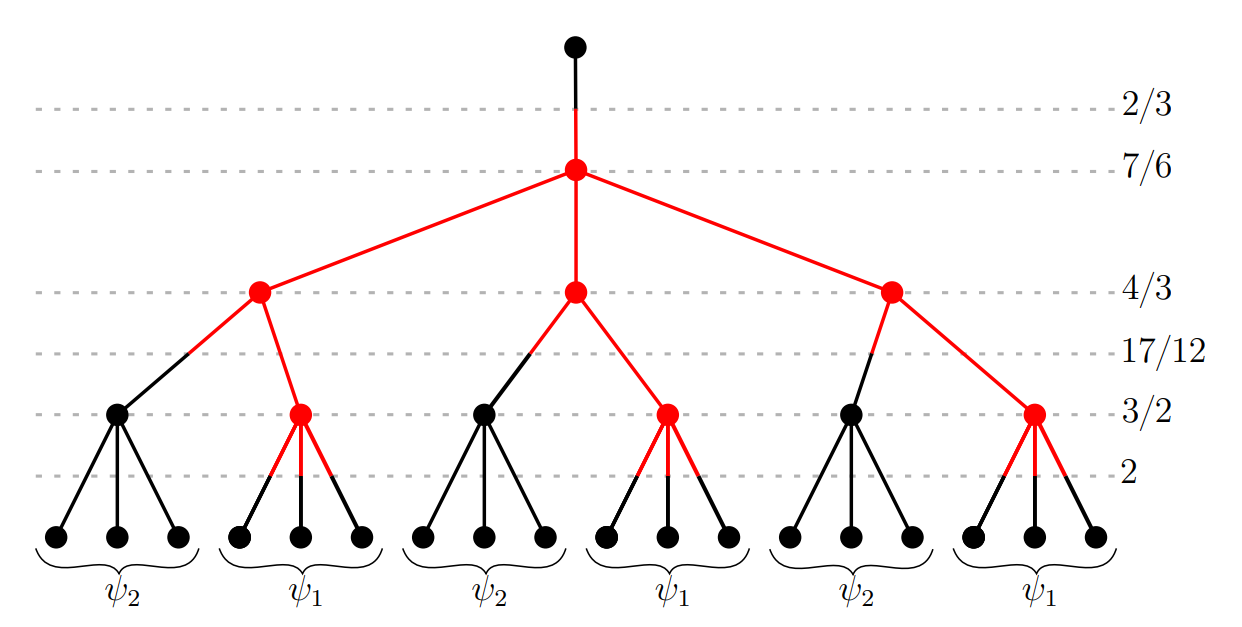}
\caption{The tree in $(\BP_K^1)^{\an}$ spanned by the zeros of $\Delta_F$ and $\psi_2$}\label{fig-big-modular-tree}
\end{figure}

The tree spanned by all the roots of $\Delta_F$ and $\psi_2$ is depicted in \mbox{Figure \ref{fig-big-modular-tree}}. The tree spanned by the zeros of $\psi_2$ already appears in \cite[Figure 3]{ossen}. As usual, the locus where $\delta>0$ is depicted in red. The dashed gray lines represent disks of the indicated radius. For example, the smallest disk containing all roots of both $\psi_1$ and $\psi_2$ has radius $7/6$. We see that the nine disks into which $\rD[\psi_1,12]$ splits have radius $2$, while the three disks into which $\rD[\psi_2,45/4]$ splits have radius $17/12$. Given a root of either $\psi_1$ or $\psi_2$, the closest other root has distance $3/2$.

Since there are three tail disks, the structure of the stable reduction of $Y$ is clear: It must be a ``comb'' (the same as in \mbox{Example \ref{ex-worked-example-3}}) consisting of one rational component intersecting three components of genus $1$, which are Artin-Schreier curves given by
\begin{equation*}
y^3-y=x^2.
\end{equation*}

\begin{Rem}
\begin{enumerate}[(a)]
\item Following the strategy outlined in \mbox{Section \ref{sec-applications-field-extension}}, a field extension of $K_0$ over which $Y$ has semistable reduction is determined in \cite{ossen}. It is of ramification index $54$ over $K_0$.
\item In \cite[Section 3]{ossen}, a coordinate transformation is used to ensure the curve $Y$ has an inflection point at infinity. It turns out that in these new coordinates, the strategy of separating branch points (\mbox{Section \ref{sec-preliminaries-admissible}}) works to compute the semistable reduction of $Y$! We have no answer to the interesting question of whether such a coordinate transformation may be found in other examples as well; in this example, this transformation was found by accident.
\end{enumerate}
\end{Rem}

\appendix

\chapter{Computer programs}

In \mbox{Section \ref{sec-appendix-implementation}}, we give a short overview of our implementation of \mbox{Algorithm \ref{alg-interior-discoid-locus}}. It is included in a branch of the MCLF Sage package available at
\begin{center}
\href{https://github.com/oossen/mclf/tree/plane-quartics}{https://github.com/oossen/mclf/tree/plane-quartics}
\end{center}
All Sage code provided in this thesis was tested on this version of MCLF with SageMath version 9.8.

\mbox{Section \ref{sec-appendix-fragments}} provides various short programs that were used for calculations in examples throughout the thesis.

\section{Implementation of Algorithm \ref{alg-interior-discoid-locus}}\label{sec-appendix-implementation}

The central part of our implementation of \mbox{Algorithm \ref{alg-interior-discoid-locus}} is the class \texttt{Quartic3Model}. It extends the MCLF class \texttt{SemistableModel}, like its siblings \texttt{AdmissibleModel} and \texttt{SuperpModel}, and provides book-keeping objects for the semistable reduction of plane quartics. For instantiation, use
\begin{minted}[tabsize=2,breaklines,fontsize=\small]{text}
sage: M = Quartic3Model(Y, vK)
\end{minted}
where \texttt{Y} is a curve over a number field $K$ and where \texttt{vK} is a valuation on $K$ with residue characteristic $3$. The two methods of \texttt{Quartic3Model} with the greatest theoretical importance are \texttt{interior\_locus} and \texttt{potential\_tail\_locus}. The former computes the image of the tame locus $Y_K^{\an}$ in $(\BP_{K_0}^1)^{\an}$, that is, the affinoid subdomain $\rU_{K_0}^{\interior}$ from \mbox{Algorithm \ref{alg-interior-discoid-locus}}. The second computes the affinoid subdomain $\rU_{K_0}$ discussed in \mbox{Section \ref{sec-applications-discrete}}. The method \texttt{tail\_locus} is available to compute the difference $\rU_{K_0}\setminus\rU_{K_0}^{\interior}$, the image of the tail locus.

The class \texttt{Quartic3Model} implements the abstract method \texttt{reduction\_tree}. The object returned essentially represents the pair of models $\CX_0$ and $\CY_0$ discussed in \mbox{Corollary \ref{cor-potentially-semistable-model}}. However, to save time in the computation of the normalization $(\CY_0)_K$, it is checked if certain vertices in the tree corresponding to $\CX_0$ can be left out. At first, it is checked if the tail components alone account for the stable reduction of $Y$. If not, the remaining boundary points of the image of the tame locus are used, and only in the end the branch points of $\phi$ are separated if necessary.

If not all branch points of $\phi$ are integral, then during the instantiation of the model a change of coordinates is performed to ensure this. For convenience, the method \texttt{rescaled\_tame\_locus} gives the tame locus with respect to the original normalization.

Our implementation is based on elementary classes \texttt{AffineMap} and \texttt{PiecewiseAffineMap} that represent affine and piecewise affine maps $\BR\to\BR$ as defined at the beginning of \mbox{Section \ref{sec-greek-definition}}.

\section{Code fragments}\label{sec-appendix-fragments}

In each of the following code listings, we note for what example the code was used; the code listing is also always mentioned in the corresponding example.

\begin{Listing}\label{list-admissible}
The following code computes the valuations of the minimal polynomial of the generator $z$ from \mbox{Example \ref{ex-admissible-no-good}}.

\begin{minted}[tabsize=2,breaklines,fontsize=\footnotesize]{text}
R.<x> = QQ[]
K.<x0> = QQ.extension(x^9 + 27/4*x^8 + 27/2*x^4 + 27/4)

R.<x> = K[]
R.<y> = R[]
F = y^3 + x^3*y + x^4 + 1
y0 = F(x=x0).roots()[1][0]
G = F(x=x+x0, y=y+y0)

v = K.valuation(3)
# G is the minimal polynomial of (y-y0)
# to compute the minimal polynomial of z = (y-y0)^{1/2}, we subtract suitable values
print([v(ai) - 1/2 for ai in G[2]])
print([v(bi) - 1 for bi in G[1]])
print([v(ci) - 3/2 for ci in G[0]])
\end{minted}
\end{Listing}

\begin{Listing}\label{list-power-series-valuations}
The following code is used in \mbox{Example \ref{ex-lambda-as-radius-of-convergence}} to compute a power series expansion and the valuations of its coefficients.
\end{Listing}

\begin{minted}[tabsize=2,breaklines,fontsize=\footnotesize]{text}
R.<a> = QQ[]
K.<a> = QQ.extension(a^3 - 3*a - 1)
v = K.valuation(3)
R.<x> = K[]
R.<y> = R[]
F = y^3 + 3*x*y^2 - 3*y - 2*x^4 - x^2 - 1
G = F(y=y+a)

approx = a

order = 100
for n in range(1, order + 1):
	next_term = -x^n*G[0][n]/G[1][0]
	G = G(y=y+next_term)
	approx += next_term

v = K.valuation(3)
for n in range(1, order + 1):
	print(n, v(approx[n]), v(approx[n])/n)
\end{minted}

\begin{Listing}\label{list-discriminant-example}
The following code computes the valuations of the coefficients of the minimal polynomial of the generator $z$ from \mbox{Example \ref{ex-computing-delta-directly}}. Note that we generate the field extension containing the point $x_0$ using the minimal polynomial of $3x_0$ instead, because Sage's functionality for extending valuations requires that the minimal polynomial be integral.

\begin{minted}[tabsize=2,breaklines,fontsize=\footnotesize]{text}
R.<x> = QQ[]
K.<a> = QQ.extension(x^8 - 22*x^6 - 114*x^5 - 237*x^4 + 346*x^3 + 2079*x^2 + 5346*x + 3645)
x0 = a/3

R.<x> = K[]
R.<y> = R[]
F = y^3 - 3*y^2 + (-3*x^2 - 2*x)*y + 3*x^4 - 3*x - 1
y0 = F(x=x0).roots()[0][0]
G = F(x=x+x0, y=y+y0)
R.<u> = K[]
f = u^3 + G[2][1]*u^2 + G[1][2]*u + G[0][3]
L.<u> = K.extension(f)

R.<x> = L[]
R.<y> = R[]
H = F(x=x+x0, y=y+y0+x*u)
v = QQ.valuation(3)
w = v.extensions(L)[0]
for h in H:
	print([w(hi) for hi in h])
\end{minted}
\end{Listing}

\begin{Listing}\label{list-tame-locus-1-revisited}
The following code computes the following, all of which are used in \mbox{Example \ref{ex-tame-locus-1-revisited}}:
\begin{itemize}
\item The ``monodromy polynomial'', which is the numerator of $\Nm(c_3)$
\item The clustering behavior of the roots of the monodromy polynomial. For this, we use the strategy explained in \cite[Section 4.8]{rueth}. We compute three Newton polygons, each showing nine roots of valuation $>1$ and eighteen roots of valuation $0$
\item For a center $x_0$ of each of the three disks $\rD_1,\rD_2,\rD_3$, the radius $\lambda(x_0)$
\item Following \mbox{Remark \ref{rem-minpoly-reduction-quartics}}, the valuations of each coefficient of the equation in whose reduction we are interested. The output shows that the constant coefficient of $\overline{B}$ and the degree-$2$ coefficient of $\overline{C}$ are the only ones not to vanish, yielding the reduction claimed in \mbox{Example \ref{ex-tame-locus-1-revisited}}
\end{itemize}

\begin{minted}[tabsize=2,breaklines,fontsize=\footnotesize]{text}
from mclf import *

R.<x,y> = QQ[]
Y = SmoothProjectiveCurve(y^3 + 3*x*y^2 - 3*y - 2*x^4 - x^2 - 1) 
FY = Y.function_field()
FX = Y.rational_function_field()
y = FY.gen()
x = FX.gen()

R.<t> = FY[]
R.<T> = R[]
A = 3*t + 3*x
B = -3
C = -2*t^4 - 8*x*t^3 + (-12*x^2 - 1)*t^2 + (-8*x^3 - 2*x)*t - 2*x^4 - x^2 - 1
F = T^3 + A*T^2 + B*T + C
w = y - (3*y^2 - 8*x^3 - 2*x)/(3*y^2 + 6*x*y - 3)*t
H = F(T+w)

A = H[2]
B = H[1]
C = H[0]
monodromy_poly = C[3].norm().numerator().monic()
print(monodromy_poly)

K.<a> = QQ.extension(monodromy_poly)
R.<y> = K[]
b = (y^3 + 3*a*y^2 - 3*y - 2*a^4 - a^2 - 1).roots()[0][0]
w = QQ.valuation(3).extensions(K)
w.append(w[1])
centers = [a, a, -a]
for i in range(3):
	print(npolygon(monodromy_poly, centers[i], w[i]))
print(w[0](a), w[1](a))

lambdas = [(3*w[i](B[0].norm().numerator()(centers[i])) - 2*w[i](C[2].norm().numerator()(centers[i])) + 2*w[i](C[2].norm().denominator()(centers[i])))/12 for i in range(3)]
print(lambdas)

h = FY.hom([b, a])
s = [2*lambdas[i]/3 for i in range(2)]
for i in range(2):
	A_val = [w[i](h(A[l])) - s[i] + l*lambdas[i] for l in range(A.degree()+1)]
	B_val = [w[i](h(B[l])) - 2*s[i] + l*lambdas[i] for l in range(B.degree()+1)]
	C_val = [w[i](h(C[l])) - 3*s[i] + l*lambdas[i] for l in range(C.degree()+1)]
	print(A_val, B_val, C_val)
\end{minted}
\end{Listing}

\begin{Listing}\label{list-tame-locus-2-revisited}
This listing accompanies \mbox{Example \ref{ex-tame-locus-2-revisited}}. We first construct an equation for $Y$ as explained in \mbox{Section \ref{sec-quartics-computing-delta-directly}}. We then compute the valuations of its coefficients, from which the graph in \mbox{Figure \ref{fig-revisited-example-2-discriminant-graph}} results by the usual formula of \mbox{Remark \ref{rem-final-delta-formula}}. We then use the functionality of MCLF to compute the points above $\xi_0$ and $\xi_{-1/4}$. Essentially this just amounts to computing extensions of the accompanying valuations to $F_Y$. But using the computer, we cannot work over the algebraically closed field $\BC_3$. Conveniently, MCLF automatically computes the necessary finite extension of $\BQ_3$.

\begin{minted}[tabsize=2,breaklines,fontsize=\footnotesize]{text}
R.<x> = QQ[]
K.<a> = QQ.extension(x^9 + 27/4*x^8 + 27/2*x^7 + 99/4*x^6 + 189/2*x^4 + 405/2*x^3 + 2187/4)
x0 = a/3

R.<x> = K[]
R.<y> = R[]
F = y^3 - y^2 + (3*x^3 + 1)*y + 3*x^4
y0 = F(x=x0).roots()[0][0]
G = F(x=x+x0, y=y+y0)
R.<u> = K[]
f = u^3 + G[2][1]*u^2 + G[1][2]*u + G[0][3]
L.<u> = K.extension(f)

R.<x> = L[]
R.<y> = R[]
H = F(x=x+x0, y=y+y0+x*u)
v = QQ.valuation(3)
w = v.extensions(L)[0]
for h in H:
	print([w(hi) for hi in h])

from mclf import *

FX.<x> = FunctionField(QQ)
R.<y> = FX[]
F = y^3 - y^2 + (3*x^3 + 1)*y + 3*x^4
FY = FX.extension(F)
Y = SmoothProjectiveCurve(FY)
X = BerkovichLine(FX, v)
T = BerkovichTree(X)
T.add_point(X.gauss_point())
T.add_point(X.point_from_discoid(x, -1/4))
R = ReductionTree(Y, v, T)
print(R.inertial_components()[0].upper_components()[0].component())
print(R.inertial_components()[1].upper_components()[0].component())
\end{minted}

\end{Listing}

\begin{Listing}\label{list-modular}
The following simply defines the curve studied in \mbox{Section \ref{sec-modular-curve}} and computes its tame locus.

\begin{minted}[tabsize=2,breaklines,fontsize=\footnotesize]{text}
from mclf import *

R.<x> = QQ[]
K.<zeta> = QQ.extension(x^2 + x + 1)
R.<x,y> = K[]
A = (6*zeta + 12)*x^2 + (36*zeta + 9)*x - 27
B = (9*zeta - 18)*x^3 + (-108*zeta - 108)*x^2 + (-162*zeta + 81)*x + (81*zeta + 162)
C = (27*zeta - 243)*x^4 + (-1458*zeta - 999)*x^3 + (-1215*zeta + 1701)*x^2 + (1944*zeta + 2430)*x + 729*zeta
F = y^3 + A*y^2 + B*y + C
Y = SmoothProjectiveCurve(F)
v = K.valuation(3)
M = Quartic3Model(Y, v)
V = M.rescaled_tame_locus()
print(V)
\end{minted}
    
\end{Listing}

\addcontentsline{toc}{chapter}{Bibliography}
\renewcommand*{\bibfont}{\small}

\printbibliography

\newpage

\begin{otherlanguage*}{german}

\section*{Zusammenfassung in deutscher Sprache}\addcontentsline{toc}{chapter}{Zusammenfassung in deutscher Sprache} \enlargethispage*{\baselineskip}

Es sei $K$ ein Körper, der vollständig ist bezüglich einer nicht-trivialen nicht-archimedischen Bewertung mit Restklassencharakteristik $p>0$. Diese Arbeit befasst sich mit der Berechnung der semistabilen Reduktion gewisser algebraischer Kurven über $K$. Es ist wohlbekannt, wie die semistabile Reduktion einer $K$-Kurve $Y$ berechnet werden kann, wenn ein Morphismus $Y\to\BP_K^1$ vom Grad $n$ kleiner als $p$ existiert. Wir behandeln aber nun den schwierigeren ``wilden'' Fall, in dem die Existenz eines Morphismus $Y\to\BP_K^1$ vom Grad genau $p$ vorausgesetzt wird.

Die analytische Geometrie über dem Körper $K$ im Sinne von Berkovich (\cite{berkovich}) ist eine nützliche Sprache für das Studium der Reduktion von Kurven über $K$. Zu jeder algebraischen $K$-Kurve $Y$ mag man ihre \emph{Analytifizierung} $Y^{\an}$ assoziieren, die die Reduktionen sämtlicher semistabiler Modelle von $Y$ in sich vereint. Die Schnittgraphen solcher Reduktionen lassen sich nämlich mit Teilräumen von $Y^{\an}$ identifizieren, den sogenannten \emph{Skeletten}. Somit verwandelt sich die ursprüngliche Problemstellung über die Bestimmung einer semistabilen Reduktion von $Y$ in die Frage, wie man ein Skelett von $Y$ finden kann.

Eine entscheidende Zutat ist die Arbeit \cite{ctt} von Cohen, Temkin, und Trushin. Ihr Gegenstand ist die zu einer Überlagerung $\phi\colon Y\to X$ assoziierte \emph{Differentenfunktion} $\delta_\phi$ auf $Y^{\an}$. Ist $\phi\colon Y\to\BP_K^1$ nun eine Überlagerung vom Grad $p$, so definieren wir ausgehend von $\delta_\phi$ eine Funktion $\lambda$ auf $(\BP_K^1)^{\an}$. Die Funktion $\lambda$ ist außerhalb eines gewissen Skeletts von $\BP_K^1$, dessen Urbild in $Y^{\an}$ ebenfalls ein Skelett ist, lokal konstant. Somit wird durch $\lambda$ ein semistabiles Modell von $Y$ bestimmt. Zu den Hauptresultaten dieser Arbeit zählen Formeln für $\lambda$ in Termen von Koeffizienten einer geeigneten Gleichung von $Y$.

Besondere Aufmerksamkeit wird dem Fall ebener quartischer Kurven über einem bezüglich einer Bewertung mit Restklassencharakteristik $3$ vollständigen Körper $K$ geschenkt. Die ebenen Quartiken bilden die einfachste Klasse von Kurven, für die noch keine allgemeine Methode zur Berechnung ihrer semistablen Reduktion zu beliebiger Restklassencharakteristik bekannt ist. Sie sind trigonal, können also mit einem Morphismus $\phi\colon Y\to\BP_K^1$ vom Grad $3$ ausgestattet werden, sodass unsere Ergebnisse auf sie angewandt werden können.

Um die Berechnung der semistabilen Reduktion ebener Quartiken in konkreten Beispielen fassbar zu machen, führen wir den Begriff des zu einer Überlagerung $\phi\colon Y\to\BP_K^1$ vom Grad $3$ assoziierten \emph{zahmen Ortes} ein. Seine Existenz entspringt den Eigenschaften der Funktion $\lambda$ und er beherrscht in gleichem Maße die semistable Reduktion der Überlagerung $\phi$.

Ein Beitrag dieser Arbeit ist auch eine Implementierung der Bestimmung des zahmen Ortes im Computeralgebrasystem SageMath, basierend auf dem MCLF-Paket (\cite{mclf}). Unter Zuhilfenahme bereits existierender MCLF-Methoden lässt sich so die semistabile Reduktion ebener Quartiken, die in einer geeigneten Normalform vorliegen, automatisch durchführen. Wir illustrieren dies in vielen Beispielen.

\end{otherlanguage*}

\section*{Brevi\=arium Lat\={\i}n\=e redditum}\addcontentsline{toc}{chapter}{Brevi\=arium Lat\={\i}n\=e redditum} 

Sit $K$ corpus expl\=etum secundum aestim\=ati\=onem n\=on trivi\=alem nec archim\=ed\=eam charact\=eristicae residu\=alis $p>0$. Pr\=opositum praesentis dissert\=ati\=onis est reducti\=onem semistabilem qu\=arundam $K$-curv\=arum comput\=are. Quae r\=es qu\=omodo effic\={\i} potest, satis n\=otum est, s\={\i} morphismus inveni\=atur $Y\to\BP_K^1$ grad\=us $n$ ips\=o $p$ min\=oris. N\=os autem tract\=abimus quaesti\=onem difficili\=orem et ``feram'' d\=e curv\=a, cui suppositus morphismus $Y\to\BP_K^1$ grad\=us ad amussim ips\=o $p$ aequ\=alis. 

Ge\=ometria analytica secundum illum Berkovich (\cite{berkovich}) super corpus $K$ instr\=umentum est \=utilissimum ad $K$-curv\=arum reducti\=on\={\i} studendum. Cuilibet $K$-curvae $Y$ adiungere licet eius \emph{analytific\=ati\=onem} $Y^{\an}$, quae omnium reducti\=on\=es exempl\=arium semistabilium ips\={\i}us $Y$ complectitur. R\=etia intersecti\=on\=alia enim eius mod\={\i} reducti\=onum subspatia ips\={\i}us $Y^{\an}$ fing\={\i} possunt, quae \emph{ossa} vocantur. Ita probl\=ema d\=e comperiend\=a reducti\=one semistabil\={\i} in quaestionem vertitur, qu\=omodo ossa ips\={\i}us $Y$ d\=etermin\=ar\={\i} possint.

Magn\={\i} m\=oment\={\i} est opus \cite{ctt} ill\=orum Cohen, Temkin, Trushin, qu\={\i} morphism\={\i}s curv\=arum $\phi\colon Y\to X$ f\={\i}n\={\i}t\={\i}s \emph{f\=uncti\=one different\={\i}} $\delta_\phi$ \=utent\=es student. Quods\={\i} $\phi\colon Y\to\BP^1_K$ est grad\=us $p$, differentem $\delta_\phi$ adhib\=emus ad f\=uncti\=onem $\lambda$ super $(\BP_K^1)^{\an}$ fabricandam, quae f\=uncti\=o constans est extr\=a quaedam ossa ips\={\i}us $\BP_K^1$, qu\=orum im\=ag\=o reciproca in $Y^{\an}$ quoque ossa suppeditat. Ita exemplar curvae $Y$ per f\=uncti\=onem $\lambda$ constituitur, quam ut coefficientibus id\=oneae aequ\=ati\=onis ips\={\i}us $Y$ d\=escr\={\i}b\=amus, in pr\={\i}m\={\i}s successibus hab\=emus.  

Animum advertimus praeciup\=e curv\={\i}s qu\=artic\={\i}s pl\=an\={\i}s super corpus expl\=etum secundum aestim\=ati\=onem charact\=eristicae residu\=alis $3$. Qu\=articae pl\=anae simplicissimae curvae sunt, qu\=arum reducti\=o semistabilis ad quamlibet charact\=eristicam residu\=alem n\=ondum semper comput\=ar\={\i} potest. Cum qu\=articae trig\=on\=al\=es sint, id est morphismum $\phi\colon Y\to\BP_K^1$ admittant, quae repperimus adhib\=er\={\i} possunt.

Ut singul\=arum qu\=artic\=arum reducti\=onem r\=ev\=er\=a attrect\=emus, \=ut\=emur \emph{loc\=o m\={\i}t\={\i}} morphism\={\i} f\={\i}n\={\i}t\={\i} $\phi\colon Y\to\BP_K^1$ grad\=us $3$. Cuius f\=ons propriet\=at\=es f\=uncti\=onis $\lambda$ sunt; reducti\=onem semistabilem ips\={\i}us $\phi$ e\=odem mod\=o regit ac ipsa $\lambda$.

Nec parva pars huius operis est illum sapientem math\=ese\=os comput\=at\=orium ``SageMath'' doc\=ere locum m\={\i}tem per s\=e s\=olum d\=etermin\=are, id quod fascicul\=o MCLF (\cite{mclf}) inn\={\i}tent\=es efficimus. Auxili\=aribus f\=uncti\=onibus veteribus fascicul\={\i} MCLF comput\=ati\=o reducti\=onis semistabilis pl\=an\=arum qu\=artic\=arum su\=a sponte fit, ut mult\={\i}s exempl\={\i}s ill\=ustr\=amus.

\end{document}